\documentclass[12pt,a4paper,oneside,reqno]{amsbook}
\usepackage[cp1251]{inputenc}
\usepackage[T2A]{fontenc}
\usepackage[english,russian]{babel}
\usepackage{amsmath,amsfonts,amssymb,amsthm,amscd,latexsym}
\newif\ifdisser
\disserfalse
\makeatletter
\def\@starttoc#1#2{%
  \begingroup
  \setTrue{#1}%
  \let\@secnumber\@empty 
  \ifx\contentsname#2%
  \else \@tocwrite{chapter}{#2}\fi
  \typeout{#2}\@xp\chaptermark\@xp{#2}%
  \@makeschapterhead{#2}\@afterheading
  \parskip\z@skip
  \makeatletter
  \@input{\jobname.#1}%
  \if@filesw
    \@xp\newwrite\csname tf@#1\endcsname
    \immediate\@xp\openout\csname tf@#1\endcsname \jobname.#1\relax
  \fi
  \global\@nobreakfalse \endgroup
  \newpage
}
\def\partrunhead#1#2#3{%
  \@ifnotempty{#2}{\uppercasenonmath{\ignorespaces#1 #2\unskip}\@ifnotempty{#3}{. }}%
  \def\@tempa{#3}%
  \ifx\@empty\@tempa\else
    \begingroup \def\\{ \ignorespaces}
    \uppercasenonmath\@tempa\@tempa
    \endgroup
  \fi
}
\def\partrunhead#1#2#3{%
 \@ifnotempty{#2}{{\ignorespaces#1 #2\unskip}\@ifnotempty{#3}{. }}%
 \def\@tempa{#3}%
 \ifx\@empty\@tempa\else
   \begingroup \def\\{ \ignorespaces}
 \sc\@tempa
   \endgroup
 \fi
}
\makeatother

\usepackage[unicode]{hyperref}
\usepackage{hypbmsec}
\ifx\undefined \pdfgentounicode
\else
\input{glyphtounicode} \pdfgentounicode=1
\fi
\ifx\pdfoutput\undefined\else
\advance\evensidemargin by -16.5mm
\advance\oddsidemargin by -16.5mm
\advance\topmargin by -8.5mm
\hypersetup{%
pdfauthor=\expandafter{Вадим В. Зудилин},%
pdftitle=\expandafter{Теорема Апери и задачи для значений
дзета-функции Римана и их \$q\$-аналогов},
colorlinks=true,
linkcolor=blue,
citecolor=blue,
urlcolor=blue,
filecolor=blue,
bookmarksnumbered=true,
pdfstartview=FitH,
pdfhighlight=/N
}
\fi
\let\le\leqslant
\let\ge\geqslant
\let\emptyset\varnothing
\def\S{\text{$\mathchar"278$}}
\renewcommand{\.}{.~\ignorespaces }
\renewcommand{\;}{\ }
\renewcommand{\:}{\colon}
\renewcommand{\d}{\mathrm{d}}
\newcommand{\wt}{\widetilde}
\newcommand{\wh}{\widehat}
\newcommand{\ol}{\overline}
\newcommand{\<}{\langle}
\renewcommand{\>}{\rangle}
\renewcommand{\[}{\lfloor}
\renewcommand{\]}{\rfloor}
\let\lf\lfloor
\let\rf\rfloor

\DeclareMathOperator{\ord}{ord}
\DeclareMathOperator{\Real}{Re}
\renewcommand{\Re}{\Real}
\DeclareMathOperator{\Imag}{Im}
\renewcommand{\Im}{\Imag}
\DeclareMathOperator{\Res}{Res}
\DeclareMathOperator{\Li}{Li}
\DeclareMathOperator{\id}{id}

\newcommand{\rom}[1]{{\rm#1}}
\newcommand{\doublesb}[2]{_{\genfrac{}{}{0pt}{1}{#1}{#2}}}
\newcommand{\ba}{{\boldsymbol a}}
\newcommand{\bb}{{\boldsymbol b}}
\newcommand{\bc}{{\boldsymbol c}}

\newcommand{\bh}{{\boldsymbol h}}

\newcommand{\balpha}{{\boldsymbol\alpha}}
\newcommand{\bbeta}{{\boldsymbol\beta}}
\newcommand{\Beta}{{\boldsymbol\eta}}
\newcommand{\fa}{{\mathfrak a}}
\newcommand{\fb}{{\mathfrak b}}
\newcommand{\fc}{{\mathfrak c}}
\newcommand{\fg}{{\mathfrak g}}
\newcommand{\fh}{{\mathfrak h}}
\newcommand{\fS}{{\mathfrak S}}
\newcommand{\fG}{{\mathfrak G}}

\newcommand{\qbinom}[2]{\genfrac[]{0pt}{}{#1}{#2}_q}

\newcommand{\xbinom}[2]{\genfrac[]{0pt}{}{#1}{#2}_x}
\begin{document}
\newtheorem{theorem}{Теорема}
\newtheorem*{theorem-apery}{Теорема Апери}
\newtheorem*{theorem-rivoal}{Теорема Ривоаля}
\newtheorem*{problem-schmidt}{Задача Шмидта}
\newtheorem{lemma}{Лемма}[chapter]
\newtheorem{proposition}{Предложение}[chapter]
\theoremstyle{remark}
\newtheorem{remark}{Замечание}[chapter]
\renewcommand{\proofname}{Доказательство}
\renewcommand{\chaptername}{Глава}
\renewcommand{\contentsname}{Содержание}
\numberwithin{equation}{chapter}
\numberwithin{section}{chapter}




\renewcommand{\cline}[1]{\hbox to\hsize{\hss#1\hss}}
\newcommand{\rline}[1]{\hbox to\hsize{\hss#1}}

\thispagestyle{empty}

\ifdisser
\cline{\large\bf Московский государственный университет}
\vskip1.5mm
\cline{\large\bf имени М.\,В.~Ломоносова}
\vskip1mm
\hbox to\hsize{\hrulefill}
\vskip2.5mm
\cline{\large\bf Механико-математический факультет}
\else
\vbox to 15mm{\vss}
\fi

\vskip45mm

\cline{\large\sl ЗУДИЛИН Вадим Валентинович}
\vskip15mm

\cline{\Large\bf Теорема Апери и задачи для значений}
\vskip2mm
\cline{\Large\bf дзета-функции Римана и их $q$-аналогов}
\vskip15mm

\ifdisser
\cline{Специальность 01.01.06 --- математическая логика,}
\smallskip
\cline{\phantom{Специальность 01.01.06 --- }алгебра и теория чисел}
\vskip15mm

\cline{Д\,И\,С\,С\,Е\,Р\,Т\,А\,Ц\,И\,Я}
\smallskip
\cline{на соискание ученой степени}
\smallskip
\cline{доктора физико-математических наук}
\vskip9mm

\rline{Научный консультант: доктор физико-математических наук,}
\rline{профессор, член-корр.\ РАН \sl НЕСТЕРЕНКО Юрий Валентинович}
\fi
\vfill

\cline{\ifdisser Москва\else Ньюкасл, Новый Южный Уэльс, Австралия\fi}
\smallskip
\cline{\ifdisser 2014\else 2013\fi}

\eject
\thispagestyle{empty}


\pdfbookmark[0]{Содержание}{Contents}
\vbox to-3mm{\vss}
\tableofcontents

\setcounter{chapter}{-1}

\begingroup

\def\chaptername{}
\def\thechapter{}

\chapter{Введение}
\label{chap:0}

\def\thechapter{0}

Изучение сумм вида
\begin{equation}
\sum_{n=1}^\infty\frac1{n^s}
\label{eq:0.1}
\end{equation}
при целых положительных значениях параметра $s$ восходит
к Л.~Эйлеру \cite{Eu1}, \cite{Eu2}.
Он, в частности, доказал расходимость ряда
в~\eqref{eq:0.1} при $s=1$ и сходимость при $s>1$,
а также знаменитые соотношения
\begin{equation}
2\sum_{n=1}^\infty\frac1{n^{2k}}
=-\frac{(2\pi i)^{2k}B_{2k}}{(2k)!}
\qquad\text{для}\quad k=1,2,3,\dots,
\label{eq:0.2}
\end{equation}
связывающие значения ряда при четных положительных~$s$
с архимедовой постоянной $\pi=3.14159265\dots$
(см.\ \cite[\S\,1.4]{Fin})
и числами Бернулли $B_s\in\nobreak\mathbb Q$;
последние могут быть определены с помощью производящей функции
\begin{equation*}
\frac z{e^z-1}
=1-\frac z2+\sum_{s=2}^\infty B_s\frac{z^s}{s!}
=1-\frac z2+\sum_{k=1}^\infty B_{2k}\frac{z^{2k}}{(2k)!}.
\end{equation*}
В 1882 году Ф.~Линдеман~\cite{Li} доказал трансцендентность числа~$\pi$
и, тем самым, трансцендентность~$\zeta(s)$ для четных~$s$.

Лишь веком спустя после Эйлера
Б.~Риман~\cite{Rie} рассмотрел ряд в~\eqref{eq:0.1}
как функцию комплексного переменного~$s$. Этот ряд представляет
в области $\Re s>1$ аналитическую функцию, которая
может быть продолжена на всю комплексную плоскость
до мероморфной функции $\zeta(s)$. Именно это аналитическое
продолжение и ряд важных свойств функции $\zeta(s)$
были открыты Риманом в его мемуаре о простых числах.
Дзета-функция Римана и ее обобщения играют неоценимую роль
в аналитической теории чисел~\cite{VK}, \cite{HW},
но тематикой настоящей \ifdisser диссертации \else монографии \fi является изучение
арифметических и аналитических свойств значений
эйлеровых сумм $\zeta(s)$ в~\eqref{eq:0.1}
при положительных $s>1$ и обобщений этих чисел.
Для краткости мы будем называть величины
\begin{equation*}
\zeta(s)=\sum_{n=1}^\infty\frac1{n^s}
\end{equation*}
при целых положительных $s$ {\it дзета-значениями},
а также {\it четными\/} и {\it нечетными дзета-значениями\/}
в зависимости от четности~$s$.

\section{Теорема Апери}
\label{sec:0.1}

Как было отмечено выше, трансцендентность (а значит, и иррациональность)
четных дзета-значений следуют из классических результатов
Эйлера и Линдемана. Формулы, подобные \eqref{eq:0.2},
для нечетных дзета-значений неизвестны,
и предположительно $\zeta(2k+1)/\pi^{2k+1}$
не является рациональным числом ни для какого целого $k\ge1$.
Арифметическая природа нечетных дзета-значений
казалась неприступной вплоть до 1978~года,
когда Р.~Апери \cite{Ap} предъявил последовательность
рациональных приближений, доказывающих иррациональность числа~$\zeta(3)$.

\begin{theorem-apery}
Число $\zeta(3)$ иррационально.
\end{theorem-apery}

История этого открытия так же,
как и строгое математическое обоснование наблюдений
Апери, изложены в~\cite{Po}. Число $\zeta(3)$ также
известно в наши дни как {\it постоянная Апери\/}
(см., например, \cite[\S\,1.6]{Fin}). В качестве
рациональных приближений к~$\zeta(3)$ Апери выбирает
последовательность $v_n/u_n\in\mathbb Q$, $n=0,1,2,\dots$,
где знаменатели $\{u_n\}=\{u_n\}_{n=0,1,\dots}$
и числители $\{v_n\}=\{v_n\}_{n=0,1,\dots}$ удовлетворяют
одной и той же полиномиальной рекурсии
\begin{equation}
(n+1)^3u_{n+1}-(2n+1)(17n^2+17n+5)u_n+n^3u_{n-1}=0
\label{eq:0.3}
\end{equation}
с начальными данными
\begin{equation}
u_0=1, \quad u_1=5, \qquad v_0=0, \quad v_1=6.
\label{eq:0.4}
\end{equation}
Тогда
\begin{equation}
\lim_{n\to\infty}\frac{v_n}{u_n}=\zeta(3),
\label{eq:0.5}
\end{equation}
но не менее важным обстоятельством являются
неожиданные (с точки зрения рекурсии~\eqref{eq:0.3})
включения
\begin{equation}
u_n=\sum_{k=0}^n\binom nk^2\binom{n+k}k^2\in\mathbb Z,
\quad D_n^3v_n\in\mathbb Z,
\qquad n=0,1,2,\dots,
\label{eq:0.6}
\end{equation}
где через $D_n$ обозначено наименьшее общее кратное
чисел $1,2,\dots,n$ (и $D_0=1$ для полноты).
Применение теоремы Пуанкаре (см., например, \cite{Ge})
к разностному уравнению~\eqref{eq:0.3} приводит к предельным
соотношениям
\begin{gather}
\lim_{n\to\infty}|u_n\zeta(3)-v_n|^{1/n}
=(\sqrt2-1)^4,
\label{eq:0.7}
\\
\lim_{n\to\infty}|u_n|^{1/n}
=\lim_{n\to\infty}|v_n|^{1/n}
=(\sqrt2+1)^4
\label{eq:0.8}
\end{gather}
согласно~\eqref{eq:0.5}, где числа $(\sqrt2-1)^4$ и $(\sqrt2+1)^4$
являются корнями характеристического многочлена $\lambda^2-34\lambda+1$
рекурсии~\eqref{eq:0.3}. Собранная информация о свойствах
последовательностей $\{u_n\}$ и $\{v_n\}$ доказывает, что
число $\zeta(3)$ не может быть рациональным. Действительно,
в предположении $\zeta(3)=a/b$, где $a,b\in\mathbb Z$,
линейные формы $r_n=bD_n^3(u_n\zeta(3)-v_n)$
являются целыми числами, ненулевыми ввиду~\eqref{eq:0.7}.
С другой стороны, $D_n^{1/n}\to e$ при $n\to\infty$ согласно
асимптотическому закону распределения простых чисел (см., например,
\cite[гл.~II, \S\,3]{VK}); следовательно,
$$
\lim_{n\to\infty}|r_n|^{1/n}
=e^3(\sqrt2-1)^4=0.59126300\ldots<1,
$$
что при достаточно большом~$n$ вступает
в противоречие с оценкой $|r_n|\ge1$ для целых ненулевых~$r_n$.
Более того, дополнительные предельные соотношения~\eqref{eq:0.8}
и стандартные аргументы (см., например, \cite[лемма~3.1]{Ha1})
позволяют измерить иррациональность постоянной Апери количественно:
\begin{equation*}
\mu(\zeta(3))
\le1+\frac{4\log(\sqrt2+1)+3}{4\log(\sqrt2+1)-3}
=13.41782023\dotsc.
\end{equation*}

Здесь и далее {\it показателем иррациональности\/} $\mu(\alpha)$
вещественного иррационального числа $\alpha$
называется величина
\begin{align*}
\mu=\mu(\alpha)
&=\inf\{c\in\mathbb R:\text{неравенство $|\alpha-a/b|\le|b|^{-c}$ имеет}
\\ &\qquad\qquad
\text{конечное число решений в $a,b\in\mathbb Z$}\};
\end{align*}
в случае $\mu(\alpha)<+\infty$ говорят, что $\alpha$~---
{\it нелиувиллево число}.

\medskip
Оригинальные рассуждения Апери (именно, соотношения
\eqref{eq:0.3}--\eqref{eq:0.8})
были настолько загадочны, что интерес к теореме
Апери не ослабевает и в наши дни.
Феномен последовательности рациональных приближений Апери
неоднократно переосмысливался
с точки зрения различных методов
(см.\ \cite{Be1}, \cite{Be2}, \cite{Be3}, \cite{Fis}, \cite{Gu},
\cite{Ha1}, \cite{Ne2}, \cite{Ne4}, \cite{Pr}, \cite{RV3},
\cite{So1}, \cite{So3}, \cite{So4}, \cite{Vi}, \cite{Ze},
\cite{Z1}, \cite{Z13}).
Новые подходы позволили усилить результат Апери количественно
--- получить лучшую оценку для показателя иррациональности числа~$\zeta(3)$
(последние этапы соревнования в этом направлении
--- работы \cite{Ha3}, \cite{RV3}).
Мы прежде всего укажем явные формулы для последовательности
$u_n\zeta(3)-v_n$, которые играют важную роль в дальнейшем изложении:
представление Бэйкерса~\cite{Be1}
\begin{equation}
u_n\zeta(3)-v_n=\iiint\limits_{[0,1]^3}
\frac{x^n(1-x)^ny^n(1-y)^nz^n(1-z)^n}
{(1-(1-xy)z)^{n+1}}\,\d x\,\d y\,\d z
\label{eq:0.9}
\end{equation}
в виде кратного вещественного интеграла,
а также ряд Гутника--Нес\-те\-ренко \cite{Gu}, \cite{Ne2}
\begin{equation}
u_n\zeta(3)-v_n=-\frac12\sum_{\nu=1}^\infty\frac{\d}{\d t}
\biggl(\frac{(t-1)(t-2)\dotsb(t-n)}
{t(t+1)(t+2)\dotsb(t+n)}\biggr)^2\bigg|_{t=\nu}
\label{eq:0.10}
\end{equation}
и ряд Болла~\cite{BR}
\begin{equation}
u_n\zeta(3)-v_n=n!^2\sum_{\nu=1}^\infty\Bigl(t+\frac n2\Bigr)
\frac{(t-1)\dotsb(t-n)\cdot(t+n+1)\dotsb(t+2n)}
{t^4(t+1)^4\dotsb(t+n)^4}\bigg|_{t=\nu}.
\label{eq:0.11}
\end{equation}

Отметим, что с помощью своего метода ``ускорения сходимости''
Апе\-ри~\cite{Ap}, \cite{Po} установил
также иррациональность числа $\zeta(2)$
без явного применения формулы $\zeta(2)=\pi^2/6$.
На этот раз, знаменатели $\{u_n'\}$ и числители $\{v_n'\}$
линейных приближающих форм $u_n'\zeta(2)-v_n'$, $n=0,1,2,\dots$,
удовлетворяют рекурсии
\begin{equation}
(n+1)^2u_{n+1}-(11n^2+11n+3)u_n-n^2u_{n-1}=0
\label{eq:0.12}
\end{equation}
с начальными данными
\begin{equation}
u_0'=1, \quad u_1'=3, \qquad v_0'=0, \quad v_1'=5;
\label{eq:0.13}
\end{equation}
при этом
\begin{equation}
u_n'=\sum_{k=0}^n\binom nk^2\binom{n+k}k\in\mathbb Z,
\quad D_n^2v_n'\in\mathbb Z,
\qquad n=0,1,2,\dots,
\label{eq:0.14}
\end{equation}
и
\begin{gather}
\lim_{n\to\infty}|u_n'\zeta(2)-v_n'|^{1/n}
=\biggl(\frac{\sqrt5-1}2\biggr)^5<e^{-2},
\label{eq:0.15}
\\
\lim_{n\to\infty}|u_n|^{1/n}
=\lim_{n\to\infty}|v_n|^{1/n}
=\biggl(\frac{\sqrt5+1}2\biggr)^5.
\label{eq:0.16}
\end{gather}
Данная последовательность приближений приводит также
к оценке
\begin{equation*}
\mu(\zeta(2))=\mu(\pi^2)
\le1+\frac{5\log((\sqrt5+1)/2)+2}{5\log((\sqrt5+1)/2)-2}
=11.85078219\dots
\end{equation*}
для показателя иррациональности числа~$\pi^2$.
Приближения Апери к~$\zeta(2)$ могут быть представлены
в виде двухкратного вещественного интеграла~\cite{Be1}
\begin{equation}
u_n'\zeta(2)-v_n'=(-1)^n\iint\limits_{[0,1]^2}
\frac{x^n(1-x)^ny^n(1-y)^n}{(1-xy)^{n+1}}\,\d x\,\d y,
\label{eq:0.17}
\end{equation}
а также в виде гипергеометрического ряда
\begin{equation}
u_n'\zeta(2)-v_n'=(-1)^n\sum_{\nu=1}^\infty
\frac{n!\cdot(t-1)(t-2)\dotsb(t-n)}
{t^2(t+1)^2(t+2)^2\dotsb(t+n)^2}\bigg|_{t=\nu}.
\label{eq:0.18}
\end{equation}

\section{Иррациональность нечетных дзета-значений}
\label{sec:0.2}

Таким образом, теорема Апери является первым существенным
продвижением в решении следующей задачи (которую по праву
можно назвать фольклорной; печатное упоминание см., например,
в~\cite{Sh}, заключительные замечания):
{\em доказать иррациональность чисел $\zeta(2k+1)$
для $k=1,2,3,\dots$\,.}

К сожалению, естественные обобщения конструкции Апери
приводят к линейным формам, содержащим значения дзета-функции
как в нечетных, так и в четных точках; это обстоятельство
не позволяло получить результаты об иррациональности~$\zeta(s)$
для нечетных $s\ge5$.
Лишь в 2000 году Т.~Ривоаль~\cite{Ri1}, используя
обобщение представления Болла \eqref{eq:0.11},
построил линейные формы, содержащие только нечетные дзета-значения
и позволяющие доказать следующий результат.

\begin{theorem-rivoal}
Среди чисел
$$
\zeta(3), \; \zeta(5), \; \zeta(7), \; \zeta(9), \; \zeta(11), \; \dots
$$
имеется бесконечно много иррациональных.
Более точно, для размерности $\delta(s)$
пространств, порожденных над~$\mathbb Q$ числами
$1,\zeta(3),\zeta(5),\dots,\linebreak[2]\zeta(s-\nobreak2),\linebreak[2]\zeta(s)$,
где $s$~нечетно, справедлива оценка
\begin{equation*}
\delta(s)
\ge\frac{\log s}{1+\log2}(1+o(1))
\qquad \text{при $s\to\infty$}.
\end{equation*}
\end{theorem-rivoal}

Линейные приближающие формы Ривоаля в~\cite{Ri1} записываются в виде
\begin{equation}
\begin{gathered}
F_n=F_{s,r,n}
=n!^{s+1-2r}\sum_{\nu=1}^\infty\Bigl(t+\frac n2\Bigr)
\frac{\prod_{j=1}^{rn}(t-j)\cdot\prod_{j=1}^{rn}(t+n+j)}
{\prod_{j=0}^n(t+j)^{s+1}}\bigg|_{t=\nu},
\\
\text{$s$ нечетно},
\end{gathered}
\label{eq:0.19}
\end{equation}
где вспомогательный параметр $r<s/2$ имеет порядок
$r\sim s/\log^2s$; в частности, ряд $F_{3,1,n}$ совпадает
с представлением \eqref{eq:0.11} для последовательности Апери.
Раскладывая рациональную функцию
параметра~$t$ под зн\'аком суммирования в сумму простейших
дробей и используя идеи работ \cite{Ni} и~\cite{Ne2},
можно показать, что справедливы включения
\begin{equation*}
2D_n^{s+1}F_n\in\mathbb Z\zeta(s)+\mathbb Z\zeta(s-2)+\dots
+\mathbb Z\zeta(5)+\mathbb Z\zeta(3)+\mathbb Z.
\end{equation*}
Кроме того, явные формулы \eqref{eq:0.19} для линейных форм
от нечетных дзета-значений позволяют вычислить асимптотическое
поведение этих форм и их коэффициентов при $n\to\infty$.
Заключительный этап доказательства теоремы Ривоаля ---
применение критерия линейной независимости Нестеренко~\cite{Ne1}.

Тот факт, что величины \eqref{eq:0.19} являются
$\mathbb Q$-линейными формами от $1$ и дзета-значений одной
четности, связан со специальной симметрией рациональной
функции параметра~$t$, стоящей в~\eqref{eq:0.19} под зн\'аком суммы.
Возможность использования менее экзотической рациональной
функции обсуждается в работах \cite{Gu}, \cite{HP1}, \cite{HP2},
\cite{Ri2}: в итоге получаются результаты о
размерности пространств, порожденных над~$\mathbb Q$
значениями полилогарифмов
\begin{equation*}
\Li_s(z)
=\sum_{n=1}^\infty\frac{z^n}{n^s}
\end{equation*}
в рациональной точке $z$, $0<|z|\le1$.

Несмотря на то, что доказательство теоремы Ривоаля использует
некоторое обобщение конструкции для доказательства теоремы Апери,
ее результат дает лишь частичное решение задачи об иррациональности дзета-значений.
Для следующего за $\zeta(3)$ иррационального нечетного $\zeta(s)$
теорема Ривоаля устанавливает лишь диапазон~\cite{BR}:
$5\le s\le169$. Дифференцирование рациональной функции под
зн\'аком суммирования (подобно представлению~\eqref{eq:0.10})
дает возможность строить $\mathbb Q$-линейные формы от нечетных
дзета-значений, не содержащие $\zeta(3)$. Это позволяет
доказать, что ``по крайней мере одно из девяти нечетных
дзета-значений $\zeta(5),\zeta(7),\dots,\zeta(21)$ иррационально''
(автор \ifdisser диссертации \else монографии \fi \cite{Z8} и Ривоаль~\cite{Ri3}
независимо установили этот результат, используя
различные обобщения конструкции из \cite{Ri1}).
Мы доказываем следующую теорему, дающую новое продвижение
в решении задачи об иррациональности нечетных дзета-значений.

\begin{theorem}
\label{th:1}
Одно из чисел
$$
\zeta(5), \; \zeta(7), \; \zeta(9), \; \zeta(11)
$$
иррационально.
\end{theorem}

Доказательству этой теоремы посвящена гл.~\ref{chap:1}
настоящей \ifdisser диссертации\else монографии\fi. Мы используем наиболее общую
форму конструкции, предложенной в работах Ривоаля, а также
арифметический метод (см., например,
\cite{Ch2}, \cite{Ru}, \cite{Ha1}),
традиционно применяемый для улучшения
оценок меры иррациональности чисел.
Отметим, что использованная техника успешно работает
и в других арифметических задачах: в~\cite{RZ} аналоги
теоремы Ривоаля и теоремы~\ref{th:1} установлены
для значений бета-функции Дирихле
$$
\beta(s)=\sum_{n=0}^\infty\frac{(-1)^n}{(2n+1)^s}
$$
в четных точках $s\ge2$. Уточнение критерия линейной независимости Нестеренко \cite{Ne1}
в работе \cite{FZ} позволило также улучшить ряд других оценок из \cite{BR} и~\cite{Z8}.

\section{Гипергеометрические ряды и кратные интегралы}
\label{sec:0.3}

Доказательство Бэйкерса~\cite{Be1} иррациональности $\zeta(2)$ и $\zeta(3)$,
использующее интегральные представления \eqref{eq:0.17} и \eqref{eq:0.9},
простое и короткое. Именно это послужило серьезным основанием
для дальнейшего применения кратных интегралов с целью количественно
усилить и обобщить результаты Апери (см.\ \cite{DV}, \cite{Ha1},
\cite{Ha2}, \cite{Ha3}, \cite{RV1}, \cite{RV2}, \cite{RV3}, \cite{Vas},
\cite{Va1}, \cite{Va2}, \cite{Vi}). О.~Василенко в~\cite{Vas}
предложил рассматривать следующее семейство $s$-кратных
интегралов, обобщающих интегралы Бэйкерса:
\begin{equation}
J_{s,n}
=\idotsint\limits_{[0,1]^s}
\frac{\prod_{j=1}^sx_j^n(1-x_j)^n}
{Q_s(x_1,\dots,x_s)^{n+1}}
\,\d x_1\dotsb\d x_s,
\label{eq:0.20}
\end{equation}
где
\begin{equation}
Q_s(x_1,\dots,x_s)
=1-x_1(1-x_2(1-\dotsb(1-x_{s-1}(1-x_s))\dotsb)).
\label{eq:0.21}
\end{equation}
Позднее Д.~Васильев \cite{Va2} изучил интегралы
$J_{4,n}$, $J_{5,n}$ и доказал, что
\begin{equation}
4D_n^4J_{4,n}\in\mathbb{Z}\zeta(4)+\mathbb{Z}\zeta(2)+\mathbb{Z},
\qquad
D_n^5J_{5,n}\in\mathbb{Z}\zeta(5)+\mathbb{Z}\zeta(3)+\mathbb{Z},
\label{eq:0.22}
\end{equation}
а также, что линейные формы в~\eqref{eq:0.20}
стремятся (достаточно быстро) к нулю при $n\to\infty$
(к сожалению, не на столько быстро, чтобы получить
новые результаты об иррациональности дзета-значений).
Включения $D_n^2J_{2,n}\in\nobreak\mathbb{Z}\zeta(2)+\nobreak\mathbb{Z}$,
$D_n^3J_{3,n}\in\mathbb{Z}\zeta(3)+\mathbb{Z}$,
доказанные Бэйкерсом в~\cite{Be1}, и \eqref{eq:0.22}
дали Васильеву основание высказать следующее предположение:
{\em имеют место включения
\begin{equation}
\begin{gathered}
2^{s-2}D_n^sJ_{s,n}
\in\mathbb{Z}\zeta(s)+\mathbb{Z}\zeta(s-2)+\dots
+\mathbb{Z}\zeta(4)+\mathbb{Z}\zeta(2)+\mathbb{Z}
\quad \text{для $s$ четного}, \\
D_n^sJ_{s,n}
\in\mathbb{Z}\zeta(s)+\mathbb{Z}\zeta(s-2)+\dots
+\mathbb{Z}\zeta(5)+\mathbb{Z}\zeta(3)+\mathbb{Z}
\quad \text{для $s$ нечетного}.
\end{gathered}
\label{eq:0.23}
\end{equation}
}

Несмотря на положительное решение этой задачи в случае $s=2,3,4,5$,
уверенность автора \cite{Va2} в справедливости \eqref{eq:0.23}
разделяли немногие. Причиной последнего
стала другая {\it неверная\/} гипотеза, именно
$2^{s-2}\*D_n^sJ_{s,n}\in\mathbb{Z}\zeta(s)+\mathbb{Z}$
для $s$ четного и
$D_n^sJ_{s,n}\in\mathbb{Z}\zeta(s)+\mathbb{Z}$
для $s$ нечетного, высказанная Васильевым в своей предыдущей
работе~\cite{Va1}.

Частичное продвижение в задаче Васильева (с точностью до умножения
на дополнительный множитель~$2D_n$) дается в следующем утверждении.

\begin{theorem}
\label{th:2}
Для каждого целого $s\ge2$ и $n=0,1,2,\dots$ справедливо тождество
\begin{equation}
J_{s,n}=F_{s,n},
\label{eq:0.24}
\end{equation}
где
\begin{equation}
F_{s,n}
=2n!^{s-1}\sum_{\nu=1}^\infty(-1)^{(s+1)(t+n+1)}\Bigl(t+\frac n2\Bigr)
\frac{\prod_{j=1}^n(t-j)\cdot\prod_{j=1}^n(t+n+j)}
{\prod_{j=0}^n(t+j)^{s+1}}\bigg|_{t=\nu}.
\label{eq:0.25}
\end{equation}
Как следствие, имеют место включения
\begin{equation}
\begin{gathered}
2^{s-2}D_n^{s+1}J_{s,n}
\in\mathbb{Z}\zeta(s)+\mathbb{Z}\zeta(s-2)+\dots
+\mathbb{Z}\zeta(4)+\mathbb{Z}\zeta(2)+\mathbb{Z}
\quad \text{для $s$ четного}, \\
D_n^{s+1}J_{s,n}
\in\mathbb{Z}\zeta(s)+\mathbb{Z}\zeta(s-2)+\dots
+\mathbb{Z}\zeta(5)+\mathbb{Z}\zeta(3)+\mathbb{Z}
\quad \text{для $s$ нечетного}.
\end{gathered}
\label{eq:0.26}
\end{equation}
\end{theorem}

Отметим, что ряд~\eqref{eq:0.25} в точности совпадает с
рядом~\eqref{eq:0.19} при $s$~нечетном и $r=1$, так что
тождество \eqref{eq:0.24} означает совпадение интегральной
конструкции $\mathbb Q$-линейных форм от дзета-значений
с конструкцией из~\cite{Ri1}.

Ряды Болла~\eqref{eq:0.11} и Ривоаля~\eqref{eq:0.19}
хорошо известны в теории обобщенных гипергеометрических
функций \cite{AAR}, \cite{Bai}, \cite{Sl}. Формально,
гипергеометрическая функция определяется рядом
\begin{equation}
{}_{m+1}F_m\biggl(\begin{matrix}
a_0, & a_1, & \dots, & a_m \\[1pt]
& b_1, & \dots, & b_m
\end{matrix}\biggm|z\biggr)
=\sum_{\nu=0}^\infty
\frac{(a_0)_\nu(a_1)_\nu\dotsb(a_m)_\nu}
{\nu!\,(b_1)_\nu\dotsb(b_m)_\nu}z^\nu,
\label{eq:0.27}
\end{equation}
где $(a)_\nu=\Gamma(a+\nu)/\Gamma(a)$ обозначает символ Похгаммера;
условие
\begin{equation}
\Re(a_0+a_1+\dots+a_m)
<\Re(b_1+\dots+b_m)
\label{eq:0.28}
\end{equation}
обеспечивает сходимость ряда~\eqref{eq:0.27}
в области $|z|\le1$ (см., например, \cite[\S\,2.1]{Bai}).
Важную роль в анализе гипергеометрических рядов
играют формулы суммирования и преобразования.
В качестве примеров укажем формулу суммирования Пфаффа--Заальшютца
\begin{equation}
{}_3F_2\biggl(\begin{matrix}-n, \, a, \, b \\
c, \, 1+a+b-c-n \end{matrix}\biggm| 1 \biggr)
=\frac{(c-a)_n(c-b)_n}{(c)_n(c-a-b)_n}
\label{eq4}
\end{equation}
($n$~--- неотрицательное целое; см., например, \cite[с.~49, формула (2.3.1.3)]{Sl}),
предельный случай теоремы Дугалла
\begin{align}
&
{}_5F_4\biggl(\begin{matrix}
a, & 1+\tfrac12a, & b, & c, & d \\[1pt]
& \tfrac12a, & 1+a-b, & 1+a-c, & 1+a-d
\end{matrix}\biggm|1\biggr)
\nonumber\\ &\quad
=\frac{\Gamma(1+a-b)\,\Gamma(1+a-c)\,\Gamma(1+a-d)\,\Gamma(1+a-b-c-d)}
{\Gamma(1+a)\,\Gamma(1+a-b-c)\,\Gamma(1+a-b-d)\,\Gamma(1+a-c-d)}
\label{eq:0.29}
\end{align}
(см.\ \cite[\S\,4.4]{Bai}),
а также преобразование Уиппла
\begin{align}
&
{}_6F_5\biggl(\begin{matrix}
a, & 1+\tfrac12a, & b, & c, & d, & e \\[1pt]
& \tfrac12a, & 1+a-b, & 1+a-c, & 1+a-d, & 1+a-e
\end{matrix}\biggm|-1\biggr)
\nonumber\\ &\quad
=\frac{\Gamma(1+a-d)\,\Gamma(1+a-e)}{\Gamma(1+a)\,\Gamma(1+a-d-e)}
\cdot{}_3F_2\biggl(\begin{matrix}
1+a-b-c, \ d, \ e \\[1pt]
1+a-b, \ 1+a-c
\end{matrix}\biggm|1\biggr)
\label{eq:0.30}
\end{align}
(см.\ \cite{Wh1} и \cite[\S\,4.4]{Bai}); ряд других примеров
приводится в главах~\ref{chap:z2} и~\ref{chap:4} (соотношение~\eqref{whipple} и предложение~\ref{prop:4.1}).
Кроме того, для гипергеометрических функций известно
множество интегральных представлений~\cite{Bai}, \cite{Sl},
из которых мы отметим классический интеграл Эйлера--Похгаммера для гауссовой функции
($m=1$ в~\eqref{eq:0.27})
\begin{equation}
{}_2F_1\biggl(\begin{matrix} a, \, b \\ c \end{matrix}
\biggm|z\biggr)
=\frac{\Gamma(c)}{\Gamma(b)\Gamma(c-b)}
\int_0^1t^{b-1}(1-t)^{c-b-1}(1-zt)^{-a}\,\d t
\label{eq5}
\end{equation}
в случае $\Re c>\Re b>0$ (см., например,
\cite[с.~20, формула (1.6.6)]{Sl}). Формула~\eqref{eq5} справедлива при
$|z|<1$, а также при любом $z\in\mathbb C$, если $a$~является неположительным целым.

В работе~\cite{Wh2} Ф.~Уиппл назвал {\it вполне уравновешенными\/}
гипергеометрические ряды~\eqref{eq:0.27}, удовлетворяющие условию
\begin{equation*}
a_0+1=a_1+b_1=\dots=a_m+b_m;
\end{equation*}
известные преобразования (например, \eqref{eq:0.29} и \eqref{eq:0.30})
относятся, как правило, именно к таким рядам.
Особую роль среди вполне уравновешенных
гипергеометрических рядов играют ряды {\it совершенно уравновешенные},
для которых выполнено дополнительное условие
\begin{equation*}
a_1=\tfrac12a_0+1, \quad b_1=\tfrac12a_0;
\end{equation*}
обзор истории и приложений вполне и совершенно уравновешенных рядов
в различных областях математики приводится в~\cite{An2}.
Ряд~\eqref{eq:0.25} (как и ряд~\eqref{eq:0.19}) является
совершенно уравновешенным:
\begin{equation}
F_{s,n}
=\frac{n!^{2s+1}(3n+2)!}{(2n+1)!^{s+2}}
\cdot{}_{s+4}F_{s+3}\biggl(\begin{matrix}
3n+2, & \tfrac32n+2, & n+1, & \dots, & n+1 \\[1pt]
& \tfrac32n+1, & 2n+2, & \dots, & 2n+2
\end{matrix}\biggm|(-1)^{s+1}\biggr).
\label{eq:0.31}
\end{equation}
Теорема~\ref{th:2} является следствием общего результата
о представлении совершенно уравновешенного ряда в виде
кратного интеграла. В гл.~\ref{chap:2} приводятся не только
подробные доказательства этого результата и теоремы~\ref{th:2},
но и ряд других важных следствий.

Именно с помощью теоремы~\ref{th:2} в последовавшей работе~\cite{KR}
задача Васильева была полностью решена в случае нечетного $s\ge3$.
Методика~\cite{KR} основана на представлении сумм \eqref{eq:0.31}
в виде кратных гипергеометрических рядов и идейно опирается
на работы \cite{Va2} и \cite{Z15}, однако техническая реализация
этих идей потребовала от авторов~\cite{KR} большой вычислительной
работы. Для рядов~\eqref{eq:0.31} возможны и другие представления
в виде кратных интегралов сорокинского типа (в~\cite{So2}, \cite{So3}
содержатся теоретико-числовые приложения подобных интегралов);
соответствующие теоремы о преобразовании кратных интегралов
установлены С.~Злобиным~\cite{Zl1}, \cite{Zl2}. Позднее Злобин \cite{Zl3}
и независимо В.~Салихов, А.~Фроловичев \cite{SF} получили другое
решение задачи Васильева.

\section({\$q\$-Аналоги дзета-значений})%
{$q$-Аналоги дзета-значений}
\label{sec:0.4}

Как обычно, величины, зависящие от числа~$q$ и превращающиеся
в классические объекты в пределе $q\to1$ (по крайней мере формально),
называются $q$-{\it аналогами\/} или $q$-{\it расширениями}.
Возможный способ $q$-рас\-ширить значения дзета-функции Римана
выглядит следующим образом (здесь $q\in\mathbb C$, $|q|<1$):
\begin{equation}
\zeta_q(s)
=\sum_{n=1}^\infty\sigma_{s-1}(n)q^n
=\sum_{\nu=1}^\infty\frac{\nu^{s-1}q^\nu}{1-q^\nu}
=\sum_{\nu=1}^\infty\frac{q^\nu\rho_s(q^\nu)}{(1-q^\nu)^s},
\qquad s=1,2,\dots,
\label{eq:0.32}
\end{equation}
где $\sigma_{s-1}(n)=\sum_{d\mid n}d^{s-1}$ обозначает сумму
степеней делителей, а многочлены
$\rho_s(x)\in\mathbb Z[x]$
могут быть определены рекурсивно с помощью формул
\begin{equation}
\rho_1=1 \qquad\text{и}\qquad
\rho_{s+1}=(1+(s-1)x)\rho_s+x(1-x)\rho_s'
\quad\text{при $s=1,2,\dots$}\,.
\label{eq:0.33}
\end{equation}
Тогда имеют место предельные соотношения
\begin{equation*}
\lim\doublesb{q\to1}{|q|<1}(1-q)^s\zeta_q(s)
=\rho_s(1)\cdot\zeta(s)
=(s-1)!\cdot\zeta(s),
\qquad s=2,3,\dots;
\end{equation*}
равенство $\rho_s(1)=(s-1)!$ следует из~\eqref{eq:0.33}.
Определенные таким образом $q$-дзета-значения~\eqref{eq:0.32}
приводят к ряду новых интересных задач в теории диофантовых
приближений и трансцендентных чисел~\cite{Z12}, которые
являются расширениями соответствующих задач
для обычных дзета-значений. Несложно показать~\cite{Z14},
что $\zeta_q(s)$ трансцендентны как функции параметра~$q$.

Для четных $s\ge2$ ряды
$E_s(q)=1-2s\zeta_q(s)/B_s$, где $B_s\in\mathbb Q$
--- числа Бернулли, известны как {\it ряды Эйзенштейна}.
Поэтому модулярное происхождение (относительно параметра
$\tau=\frac{\log q}{2\pi i}$) функций $E_4,E_6,E_8,\dots$
приводит к алгебраической независимости
$\zeta_q(2),\zeta_q(4),\zeta_q(6)$ над $\mathbb Q[q]$,
в то время как остальные четные $q$-дзета-значения являются
многочленами от~$\zeta_q(4)$ и~$\zeta_q(6)$.
В такой интерпретации следствие из теоремы Нестеренко~\cite{Ne3}
``числа $\zeta_q(2),\zeta_q(4),\zeta_q(6)$ алгебраически
независимы над~$\mathbb Q$ для алгебраического~$q$, $0<|q|<1$''
является полным $q$-расширением следствия из теоремы
Линдемана~\cite{Li} ``$\zeta(2)=\pi^2/6$ трансцендентно''.
Об арифметической природе нечетных $q$-дзета-значений
(например, $q$-аналог задачи об иррациональности дзета-значений) известно немного.
П.~Эрдёш~\cite{Er} доказал иррациональность числа $\zeta_q(1)$
($q$-{\it гармонического ряда\/})
в случае $q=p^{-1}$, где $p\in\mathbb Z\setminus\{0,\pm1\}$;
другие доказательства имеются в~\cite{Bez} и~\cite{Bo},
а в~\cite{BV} и~\cite{As2} получена оценка
\begin{equation}
\mu(\zeta_q(1))
\le\frac{2\pi^2}{\pi^2-2}=2.50828476\dots
\label{eq:0.34}
\end{equation}
для показателя иррациональности $\zeta_q(1)$ при тех же условиях
на параметр~$q$. Конструкция линейных приближающих форм для~$\zeta_q(1)$
в~\cite{BV} и~\cite{As2} непрерывно зависит от~$q$, однако
в пределе $q\to1$, $|q|<1$, получаются расходящиеся величины.
В связи с этим обстоятельством В.~Фан Ассе сформулировал в~\cite{As2}
задачу о построении линейных приближающих форм для~$\zeta_q(2)$ и $\zeta_q(3)$,
переходящих при $q\to1$ в последовательности Апери $u_n'\zeta(2)-v_n'$
и $u_n\zeta(3)-v_n$ соответственно (см.~\S\,\ref{sec:0.1}).

Методика изучения арифметических свойств
чисел $\zeta(s)$, $s=2,3,\dots$, успешно переносится
на случай $q$-дзета-значений. Именно, мы имеем в виду
гипергеометрическую конструкцию линейных форм
и арифметический метод, дополненный групповым подходом
Дж.~Рина и К.~Виолы \cite{RV2}, \cite{RV3}, \cite{Vi}.
Для каждой из этих составляющих мы указываем
необходимые $q$-расширения. Гипергеометрическим рядам при этом соответствуют
базисные гипергеометрические ряды \cite{GR}, для которых также имеют место
преобразования; кроме того, $q$-арифметический метод~\cite{Z9}, \cite{Z10}
позволяет улучшать известные меры иррациональности $q$-дзета-значений
и других $q$-постоянных (см.~\cite{BZ}, \cite{Z7}, \cite{Z9}, \cite{Z10}, \cite{Z16}, \cite{Z20a}).
Так, например, используя базисный гипергеометрический ряд, классическое преобразование
Гейне~\cite{He} и $q$-арифметический метод
мы улучшаем оценку~\eqref{eq:0.34} для показателя
иррациональности $q$-гармонического ряда (см.~\cite{Z16}):
\begin{equation*}
\mu(\zeta_q(1))\le2.46497868\dots\,.
\end{equation*}

Используя $q$-аналог гипергеометрического ${}_3F_2(1)$-ряда
и преобразование Холла \cite{GR}, мы не только решаем
упомянутую задачу Фан Ассе для~$\zeta_q(2)$, но и оптимизируем
оценку для показателя иррациональности этого числа.

\begin{theorem}
\label{th:3}
Для каждого $q=1/p$, $p\in\mathbb Z\setminus\{0,\pm1\}$,
число~$\zeta_q(2)$ является иррациональным с показателем
иррациональности, удовлетворяющим неравенству
\begin{equation}
\mu(\zeta_q(2))\le3.51887508\dots\,.
\label{eq:0.35}
\end{equation}
\end{theorem}

Количественные оценки типа~\eqref{eq:0.35} для~$\zeta_q(2)$
(доказывающие нелиувиллевость этой постоянной
в случае $q^{-1}\in\mathbb Z\setminus\{0,\pm1\}$)
до публикаций автора \cite{Z7}, \cite{Z10} не были известны, хотя, как отмечалось выше,
иррациональность~\cite{Duv} и даже трансцендентность числа~$\zeta_q(2)$
для любого алгебраического~$q$ с условием $0<|q|<1$ следует
из теоремы Нестеренко~\cite{Ne3}. В результате более аккуратного
вычисления вспомогательных параметров конструкции в работе~\cite{SVA} была
установлена оценка
\begin{equation*}
\mu(\zeta_q(2))
\le\frac{10\pi^2}{5\pi^2-24}=3.89363887\dots
\end{equation*}
лучшая, чем в \cite{Z10}. Теорема~\ref{th:3} значительно улучшает этот результат.

Доказательство теоремы~\ref{th:3} приводится в гл.~\ref{chap:3};
там же мы устанавливаем, что частный случай
\begin{equation*}
U_n(q)\zeta_q(2)-V_n(q)=(-1)^n\sum_{\nu=1}^\infty
\frac{\prod_{j=1}^n(1-q^j)\cdot\prod_{j=1}^n(1-q^jT)}
{\prod_{j=0}^n(1-q^{n+1+j}T)^2}T^{n+1}\bigg|_{T=q^\nu}
\end{equation*}
нашей гипергеометрической конструкции
также доказывает иррациональность числа~$\zeta_q(2)$
в случае $q^{-1}\in\mathbb Z\setminus\{0,\pm1\}$,
а в пределе $q\to1$ получаются рациональные
приближения Апери~\eqref{eq:0.18} к~$\zeta(2)$.
В совместной работе~\cite{KRZ} мы предъявляем
$q$-аналог последовательности
рациональных приближений~\eqref{eq:0.10}, \eqref{eq:0.11},
однако это не приводит к иррациональности величины $\zeta_q(3)$.
Вопрос об иррациональности $\zeta_q(3)$ остается открытым;
известны только частичные результаты \cite{KRZ}, \cite{JM} в духе теоремы~\ref{th:1}.

Использование $q$-арифметического метода и
гипергеометрической конструкции позволяет
получить и другие качественные и количественные результаты
для $q$-дзета-значений. Так, работа~\cite{KRZ}
содержит результат о бесконечности иррациональных
чисел среди нечетных $q$-дзета-значений
($q$-аналог теоремы Ривоаля) в случае
$q^{-1}\in\mathbb Z\setminus\{0,\pm1\}$.
Также отметим, что для $q$-дзета-значений значительный интерес представляют и вопросы
функциональной независимости \cite{Z12}, \cite{Z14}. Здесь следует упомянуть работу
Ю.~Пупырёва \cite{Pu1}, доказавшего линейную независимость $q$-дзета-значений, а также
получившего частичные результаты об их алгебраической независимости.

\section({Рациональные приближения к \$\003\266(2)\$})%
{Рациональные приближения к $\zeta(2)$}
\label{sec:0.z2}

Обзор \cite{Z19a} о гипергеометрической интерпретации
выдержанных временем мер иррациональности $\log2$, $\pi$ и $\log3$
был опубликован в~2004\,г. Всплеск активности \cite{Mar}, \cite{Sa1}, \cite{Sa2}
в последовавшее пятилетие привел не только к количественным изменениям
рассмотренных в \cite{Z19a} мер, но и к возникновению принципиально
новых конструкторских идей для их получения. Вместе с тем рекордные меры
иррациональности
$$
\mu(\zeta(2))\le5.44124250\dots
\quad\text{и}\quad
\mu(\zeta(3))\le5.51389062\dots,
$$
полученные Рином и Виолой в 1996\,г.\ и в 2001\,г.\ соответственно,
оставались незыблимыми. В основе метода Рина--Виолы --- группа
бирациональных преобразований двойных и тройных интегралов бэйкерсова
типа \eqref{eq:0.17}, \eqref{eq:0.9} (см.\ обсуждение в конце \S\,\ref{sec:2.1}
и гипергеометрическую интерпретацию групп Рина--Виолы в~\cite{Z17}),
дополненную арифметическим методом. По существу, наше доказательство
оценки меры иррациональности $q$-аналога $\zeta(2)$ в гл.~\ref{chap:3} есть
не что иное, как $q$-версия метода Рина--Виолы, переложенного на гипергеометрический язык
(см.~\cite{Z15},~\cite{Z17}). Необходимость использовать
$q$-гипергеометрические ряды вместо $q$-аналогов интегралов
(а понятие $q$-интеграла действительно существует ---
см., например, \cite{Ask}, \cite{Ex}) продиктована
отсутствием понятия замены переменных для последних.

Глава~\ref{chap:z2} настоящей \ifdisser диссертации \else монографии \fi
посвящена еще одной демонстрации успешности арифметико-гипергеометрического метода ---
доказательству новой оценки меры иррациональности $\zeta(2)=\pi^2/6$.

\begin{theorem}
\label{th:z2}
Показатель иррациональности числа $\zeta(2)=\pi^2/6$ удовлетворяет неравенству
$\mu(\zeta(2))\le5.09541178\dots$\,.
\end{theorem}

Одним из следствий теоремы~\ref{th:z2} является общая оценка
$$
\mu(\pi\sqrt{d})\le10.19082357\dots,
$$
справедливая для любого ненулевого рационального~$d$.
В то же время для некоторых частных значений $d\in\mathbb Q$ известны лучшие оценки:
результаты
$$
\mu(\pi)\le7.606308\dots,
\quad
\mu(\pi\sqrt3)\le4.601057\dots,
\quad
\mu(\pi\sqrt{10005})\le10.021363\dots
$$
получены в работах В.~Салихова~\cite{Sa2}, В.~Андросенко и Салихова \cite{AS}
и автора~\cite{Zu} соответственно.

Доказательство теоремы~\ref{th:z2} использует две гипергеометрические конструкции
построения рациональных приближений к $\zeta(2)$. Совпадение этих приближений ---
нетривиальное аналитическое тождество, не найденное (хотя и упомянутое, см.~\S\,\ref{sbailey}) в литературе.
Кроме того, гипергеометрические конструкции позволяют строить совместные рациональные
приближения к $\zeta(2)$ и $\zeta(3)$, к сожалению, недостаточно хорошие для доказательства
линейной независимости этих дзета-значений.

\section({Нижняя оценка для \$\000\134|(3/2)\000\136k\000\134|\$ и проблема Варинга})%
{Нижняя оценка для $\|(3/2)^k\|$ и проблема Варинга}
\label{sec:0.6}

Гипергеометрическая конструкция и арифметический метод, применяющиеся
для доказательства теорем~\ref{th:1}, \ref{th:3} и~\ref{th:z2} находят приложения
во многих других задачах на стыке диофантовой и аналитической теорий чисел.
Яркий пример подобного симбиоза является основным сюжетом гл.~\ref{chap:5}
настоящей \ifdisser диссертации\else монографии\fi.

Пусть $\lf\,\cdot\,\rf$ и $\{\,\cdot\,\}$ обозначают
целую и дробную части числа соответственно. Как известно~\cite{Va},
неравенство $\{(3/2)^k\}\le1-(3/4)^k$
при $k\ge6$ дает точную формулу $g(k)=2^k+\lf(3/2)^k\rf-2$
для наименьшего целого $g=g(k)$ такого, что
каждое натуральное число представимо в виде
суммы не более $g$ положительных $k$-х степеней
(проблема Варинга). К.~Малер~\cite{Ma} использовал обобщение Риду
известной теоремы Рота, чтобы показать, что неравенство
$\|(3/2)^k\|\le C^k$,
где $\|x\|=\min(\{x\},1-\{x\})$ --- расстояние от $x\in\mathbb R$
до ближайшего целого, имеет лишь конечное число решений в целых~$k$
для любого $C<1$. В частном случае $C=3/4$ получается приведенное
выше значение $g(k)$ для всех $k\ge K$, где $K$~--- некоторая абсолютная,
но неэффективная постоянная. В связи с этим возникает следующая задача:
{\em получить нетривиальную \textup(т.е.\ $C>1/2$\textup) и эффективную
\textup(в терминах~$K$\textup) оценку вида
\begin{equation}
\biggl\|\biggl(\frac32\biggr)^k\biggr\|>C^k
\qquad\text{для всех}\quad k\ge K.
\label{eq:0.40}
\end{equation}
}%
Первое продвижение в этом направлении
принадлежит А.~Бейкеру и Дж.\ Коэтсу~\cite{BC};
применив эффективные оценки линейных форм от логарифмов в
$p$-адическом случае, они показали справедливость~\eqref{eq:0.40}
с $C=2^{-(1-10^{-64})}$. Ф.~Бэйкерс~\cite{Beu} существенно
улучшил этот результат, доказав, что неравенство~\eqref{eq:0.40}
выполняется с $C=2^{-0.9}=0.5358\dots$ при $k\ge K=5000$
(хотя его доказательство давало и лучший выбор
$C=0.5637\dots$, если не заботиться о явном вычислении
эффективной границы для~$K$). Доказательство Бэйкерса
основано на приближениях Паде к остатку биномиального ряда
$(1-z)^m=\sum_{n=0}^m\binom mn(-z)^n$; позднее
А.~Дубицкас~\cite{Du} и Л.~Хабсигер~\cite{Ha}
использовали конструкцию Бэйкерса для получения оценки~\eqref{eq:0.40}
с $C=0.5769$ и $0.5770$ соответственно. Последняя работа
также содержит оценку $\|(3/2)^k\|>0.57434^k$ при $k\ge5$
на основе вычислений из~\cite{DD} и~\cite{KW}.

Модифицируя конструкцию Бэйкерса~\cite{Beu}, именно,
рассматривая приближения Паде к остатку ряда
\begin{equation}
\frac1{(1-z)^{m+1}}
=\sum_{n=0}^\infty\binom{m+n}mz^n,
\label{eq:0.41}
\end{equation}
и получая точные оценки $p$-адических порядков
возникающих биномиальных коэффициентов, мы доказываем
в гл.~\ref{chap:5} следующий результат.

\begin{theorem}
\label{th:5}
\label{TH:5}
Имеет место оценка
$$
\biggl\|\biggl(\frac32\biggr)^k\biggr\|>0.5803^k
=2^{-k\cdot0.78512916\dots}
\qquad\text{для всех}\quad k\ge K,
$$
где $K$ --- некоторая эффективная постоянная.
\end{theorem}

Конструкция гл.~\ref{chap:5} позволяет также
доказать оценки
\begin{equation}
\begin{aligned}
\biggl\|\biggl(\frac43\biggr)^k\biggr\|>0.4914^k
=3^{-k\cdot0.64672207\dots}
&\qquad\text{при}\quad k\ge K_1,
\\
\biggl\|\biggl(\frac54\biggr)^k\biggr\|>0.5152^k
=4^{-k\cdot0.47839775\dots}
&\qquad\text{при}\quad k\ge K_2,
\end{aligned}
\label{eq:0.42}
\end{equation}
где $K_1,K_2$ --- эффективные постоянные.
Наилучший результат для последовательностей
$\|(1+1/N)^k\|$ принадлежит М.~Беннетту~\cite{Ben1}:
$\|(1+1/N)^k\|>3^{-k}$ при $4\le N\le k3^k$.
Наша оценка снизу для $\|(4/3)^k\|$ дополняет результат
Беннетта~\cite{Ben2} о порядке аддитивного базиса
$\{1,N^k,(N+1)^k,(N+2)^k,\dots\}$ в случае $N=3$
(случай $N=2$ отвечает классической проблеме Варинга);
решение соответствующей задачи требует оценки
$\|(4/3)^k\|>(4/9)^k$ при $k\ge6$. Таким образом,
остается ее проверить в диапазоне $6\le k\le K_1$.
Отметим, что Пупырёв дает явное значение постоянной~$K_1$ в работе~\cite{Pu2}.

\section{Числа Апери и числа Франеля}
\label{sec:0.5}

А.~Шмидт в~\cite{Sc1} обратил внимание на тот факт,
что с последовательностью чисел Апери $\{u_n\}_{n=0,1,\dots}$
из~\eqref{eq:0.6} связано удивительное обстоятельство.
Именно, если определить числа $\{c_k\}_{k=0,1,\dots}$
последовательно с помощью равенств
\begin{equation*}
u_n=\sum_{k=0}^n\binom nk\binom{n+k}kc_k,
\qquad n=0,1,2,\dots,
\end{equation*}
то эти числа являются целыми. (Явные формулы
$$
c_n=\binom{2n}n^{-1}\sum_{k=0}^n(-1)^{n-k}
\frac{2k+1}{n+k+1}\binom{2n}{n-k}u_k,
\qquad n=0,1,2,\dots,
$$
показывают, что ожидаемыми являются включения $D_nc_n\in\mathbb Z$.)
Позднее сам Шмидт~\cite{Sc2} и независимо Ф.~Штрель~\cite{St}
получили следующее явное выражение:
\begin{equation}
c_n
=\sum_{j=0}^n\binom nj^3
=\sum_j\binom nj^2\binom{2j}n,
\qquad n=0,1,2,\dots,
\label{eq:0.36}
\end{equation}
экспериментально обнаруженное В.~Дойбером, В.~Тумзером
и Б.~Войтом. На самом деле, Штрель в~\cite{St} использовал
соответствующее тождество
\begin{equation*}
\sum_{k=0}^n\binom nk^2\binom{n+k}k^2
=\sum_{k=0}^n\binom nk\binom{n+k}k\sum_{j=0}^k\binom kj^3
\end{equation*}
в качестве модели для иллюстрации различной техники
доказательства биномиальных тождеств. Удивительным является
тот факт, что последовательность \eqref{eq:0.36}
изучалась Дж.~Франелем~\cite{Fr} еще в конце 19-го века: он
доказал, что она удовлетворяет полиномиальной рекурсии
\begin{equation*}
(n+1)^2c_{n+1}-(7n^2+7n+2)c_n-8n^2c_{n-1}=0.
\end{equation*}

Шмидт в~\cite{Sc1} отметил, что по всей видимости
феномен целочисленности, связанный с последовательностями
чисел Апери и Франеля, выполняется в следующей общей ситуации.

\begin{problem-schmidt}
Пусть для любого целого $r\ge2$ последовательность
чисел $\{c_k^{(r)}\}_{k=0,1,\dots}$, не зависимых от параметра~$n$,
определяется равенством
\begin{equation}
\sum_{k=0}^n\binom nk^r\binom{n+k}k^r
=\sum_{k=0}^n\binom nk\binom{n+k}kc_k^{(r)},
\qquad n=0,1,2,\dotsc.
\label{eq:0.37}
\end{equation}
Требуется показать, что все числа
$c_k^{(r)}$ являются целыми.
\end{problem-schmidt}

Используя алгоритм Госпера--Цайльбергера созидательного
телескопирования \cite{PWZ}, Штрель доказал в~\cite{St}
целочисленность $c_k^{(r)}$ в случае $r=3$. Задача Шмидта
была позднее сформулирована в книге~\cite{GKP} (упражнение~114
на с.~256) с указанием, что Г.~Уилф доказал включения
$c_n^{(r)}\in\mathbb Z$ для любого~$r$, но только для $n\le9$.

Глава~\ref{chap:4} посвящена решению задачи Шмидта.
Именно, мы доказываем следующие явные формулы для~$c_n^{(r)}$.

\begin{theorem}
\label{th:4}
Числа $c_n^{(r)}$ в формулировке задачи Шмидта
действительно являются целыми. Более точно, имеют место формулы
\begin{align}
c_n^{(4)}
&=\sum_{j=0}^n\binom{2j}j^3\binom nj
\sum_{k\ge j}\binom{k+j}{k-j}\binom j{n-k}\binom kj\binom{2j}{k-j},
\label{eq:0.38}
\\
c_n^{(5)}
&=\sum_{j=0}^n\binom{2j}j^4\binom nj^2
\sum_{k\ge j}\binom{k+j}{k-j}^2\binom{2j}{n-k}\binom{2j}{k-j}
\label{eq:0.39}
\end{align}
и для произвольного $s=2,3,\dots$
\begin{align*}
c_n^{(2s)}
&=\sum_{j=0}^n\binom{2j}j^{2s-1}\binom nj
\sum_{k_1\ge\dots\ge k_{s-1}\ge j}
\binom j{n-k_1}\binom{k_1}j\binom{k_1+j}{k_1-j}
\\ &\qquad\times
\binom{2j}{k_{s-1}-j}
\prod_{l=2}^{s-1}
\binom{2j}{k_{l-1}-k_l}\binom{k_l+j}{k_l-j}^2,
\displaybreak[2]\\
c_n^{(2s+1)}
&=\sum_{j=0}^n\binom{2j}j^{2s}\binom nj^2
\sum_{k_1\ge\dots\ge k_{s-1}\ge j}
\binom{2j}{n-k_1}\binom{k_1+j}{k_1-j}^2
\\ &\qquad\times
\binom{2j}{k_{s-1}-j}
\prod_{l=2}^{s-1}
\binom{2j}{k_{l-1}-k_l}\binom{k_l+j}{k_l-j}^2,
\end{align*}
где $n=0,1,2,\dots$\,.
\textup{(Биномиальные коэффициенты $\binom nk$ считаются равными
нулю в случае $k<0$ или $k>n$.)}
\end{theorem}

\section{Периоды}
\label{sec:0.8}

В качестве заключительного аккорда мы приводим результаты,
связанные с принадлежностью дзета-значений и их обобщений к так называемому классу
периодов, введенному в фундаментальной работе М.~Концевича и Д.~Загира~\cite{KZ}.
\emph{Период} --- это комплексное число,
вещественная и мнимая части которого являются значениями
абсолютно сходящихся интегралов от рациональных функций с рациональными
коэффициентами, интегрируемых по областям в $\mathbb R^n$, заданных
полиномиальными неравенствами с рациональными многочленами.
Без видимого ущерба прилагательное ``рациональные'' все три раза
может быть заменено на ``алгебраические''~\cite{KZ}. Множество периодов $\mathcal P$
счетно и имеет естественную структуру кольца. Это множество содержит
все ``важные'' математические постоянные, включая, разумеется, $\pi$ \cite{BBB}, дзета-значения~\eqref{eq:0.1}
и их многочисленные обобщения (например, \emph{кратные дзета-значения} \cite{OZ}, \cite{Z13a});
кроме того, множество $\mathcal P$ содержит все алгебраические числа.
Вместе с тем, число $1/\pi$ предположительно не принадлежит кольцу $\mathcal P$,
и многие другие примеры --- такие, как значения гипергеометрических функций в алгебраических точках
и специальные $L$-значения (еще одно обобщение дзета-значений) ---
естественным образом обитают в расширенном кольце периодов $\widehat{\mathcal P}=\mathcal P[1/\pi]=\mathcal P[(2\pi i)^{-1}]$.
В обзоре~\cite{BoCr} обсуждается важность представлений чисел в (канонической) ``замкнутой форме''
с приложениями как в чистой, так и прикладной математике.

Под $L$-функцией обычно понимается производящий ряд Дирихле
$$
L(s)=\sum_{n=1}^\infty\frac{a_n}{n^s}
$$
последовательности $a_n$, $n=1,2,\dots$, имеющей ``арифметическую значимость''; например,
связанную со счетом точек по фиксированному модулю на алгебраическом многообразии~$\mathcal M$
(обозначение для этого случая $L(s)=L(\mathcal M,s)$),
либо отвечающую (параболической) модулярной форме $f(q)=\sum_{n=1}^\infty a_nq^n$
(тогда мы пишем $L(s)=L(f,s)$). Рассматривая в качестве $\mathcal M$ эллиптическую кривую
$E:y^2=x^3-x$ (кондуктора~32), обозначим через $N_p$ количество различных точек
$(x\bmod p,y\bmod p)\in(\mathbb Z_p)^2$ на ней для каждого нечетного простого~$p$ и определим числа $a_p=p-N_p$.
Тогда
$$
L(E,s)=\prod_{p>2}\biggl(1-\frac{a_p}{p^s}+\frac1{p^{2s-1}}\biggr)^{-1}
=\sum_{n=1}^\infty\frac{a_n}{n^s}.
$$
Эта же самая $L$-функция может быть реализована как $L(f,s)$ для модулярной формы
$$
f(\tau)=\sum_{n=1}^\infty a_nq^n
=q\prod_{m=1}^\infty(1-q^{4m})^2(1-q^{8m})^2,
\quad\text{где}\;\; q=\exp(2\pi i\tau),
$$
веса~2. Подобное соответствие $L(E,s)=L(f,s)$ выполняется для всех эллиптических кривых над~$\mathbb Q$
благодаря знаменитой теореме о модулярности Э.~Уайлса~\cite{Wi}, известной ранее под именем гипотезы Таниямы--Шимуры.

Общие теоремы \cite{KZ} Бейлинсона и Денингера--Шолла устанавливают
принадлежность (некритических) значений $L$-функции, отвечающей параболической модулярной форме
$f(\tau)$ веса $k$, в точках $m\ge k$ кольцу~$\widehat{\mathcal P}$.
Несмотря на алгоритмический характер доказательства теоремы,
реализация $L$-значений явно в виде периодов является трудной задачей даже в абсолютно конкретных ситуациях.
Большинство подобных вычислений мотивировано (гипотетическими) выражениями
для логарифмической меры Малера многочленов от многих переменных, которые, в свою очередь,
эквивалентны частным случаям общих гипотез Бейлинсона--Блоха.

С целью доказательства ряда гипотез для меры Малера многочленов от двух переменных
в совместных работах \cite{RZ1}, \cite{RZ2} с М.~Роджерсом мы разработали принципиально новую методику
представления $L$-значения $L(f,2)$, отвечающего параболической модулярной форме $f(\tau)$ веса~2, в качестве периода.
В~гл.~\ref{chap:7} мы приводим обзор этой методики на конкретном примере вычисления $L(E,2)$
для эллиптической кривой $E:y^2=x^3-x$, а затем описываем общий алгоритм вычисления некритических
значений $L(f,k)$ в виде периодов и иллюстрируем его на примере $L(E,3)$.

\begin{theorem}
\label{th:7}
Для эллиптической кривой $E:y^2=x^3-x$ кондуктора~$32$ справедливы следующие
интегральные и гипергеометрические представления:
\begin{align}
L(E,2)
&=\frac{\pi}{16}\int_0^1\frac{1+\sqrt{1-x^2}}{(1-x^2)^{1/4}}\,\d x
\int_0^1\frac{\d y}{1-x^2(1-y^2)}
\label{ILE2}
\\
&=\frac{\pi^{1/2}\Gamma(\frac14)^2}{96\sqrt2}\,
\,{}_3F_2\biggl(\begin{matrix} 1, \, 1, \, \frac12 \\ \frac74, \, \frac32 \end{matrix}\biggm| 1 \biggr)
+\frac{\pi^{1/2}\Gamma(\frac34)^2}{8\sqrt2}\,
\,{}_3F_2\biggl(\begin{matrix} 1, \, 1, \, \frac12 \\ \frac54, \, \frac32 \end{matrix}\biggm| 1 \biggr),
\label{HLE2}
\\[2mm]
L(E,3)
&=\frac{\pi^2}{128}\int_0^1\frac{(1+\sqrt{1-x^2})^2}{(1-x^2)^{3/4}}\,\d x
\int_0^1\!\!\!\int_0^1\frac{\d y\,\d w}{1-x^2(1-y^2)(1-w^2)}
\label{ILE3}
\\
&=\frac{\pi^{3/2}\Gamma(\frac14)^2}{768\sqrt2}\,
{}_4F_3\biggl(\begin{matrix} 1, \, 1, \, 1, \, \frac12 \\ \frac74, \, \frac32, \, \frac32 \end{matrix}\biggm| 1 \biggr)
+\frac{\pi^{3/2}\Gamma(\frac34)^2}{32\sqrt2}\,
{}_4F_3\biggl(\begin{matrix} 1, \, 1, \, 1, \, \frac12 \\ \frac54, \, \frac32, \, \frac32 \end{matrix}\biggm| 1 \biggr)
\nonumber\\ &\qquad
+\frac{\pi^{3/2}\Gamma(\frac14)^2}{256\sqrt2}\,
{}_4F_3\biggl(\begin{matrix} 1, \, 1, \, 1, \, \frac12 \\ \frac34, \, \frac32, \, \frac32 \end{matrix}\biggm| 1 \biggr).
\label{HLE3}
\end{align}
\end{theorem}

Отметим, что литература не содержит ни одного явного интегрального представления $L(E,3)$ для
какой-либо эллиптической кривой. Полученные выражения увеличивают схожесть данных $L$-значений
с $\zeta(2)$ и $\zeta(3)$ --- ср.\ с многочисленными формулами для приближений Апери в~\S\,\ref{sec:0.1}.

\endgroup

\chapter({Одно из чисел \$\003\266(5),\003\266(7),\003\266(9),\003\266(11)\$ иррационально})%
{Одно из чисел $\zeta(5),\zeta(7),\zeta(9),\zeta(11)$ иррационально}
\label{chap:1}

Начнем с объяснения того, как обобщение ряда Болла~\eqref{eq:0.11}
и методика работы~\cite{Ne2} позволяет доказать результат
о том, что по крайней мере одно из нечетных дзета-значений
$\zeta(5),\zeta(7),\dots,\zeta(2s+1)$ для некоторого $s>\nobreak3$ иррационально.
С этой целью для каждого $n=0,1,2,\dots$ рассмотрим величину
\begin{equation}
F_n=\frac12\sum_{\nu=1}^\infty\frac{\d^2R_n(t)}{\d t^2}\bigg|_{t=\nu},
\label{eq:1.1}
\end{equation}
где
\begin{equation}
R_n(t)
=2n!^{2s-6}\Bigl(t+\frac n2\Bigr)
\frac{\prod_{j=1}^n(t-j)^3\cdot\prod_{j=1}^n(t+n+j)^3}
{\prod_{j=0}^n(t+j)^{2s}}
=\sum_{j=1}^{2s}\sum_{k=0}^n\frac{B_{jk}}{(t+k)^j}
\label{eq:1.2}
\end{equation}
--- рациональная функция аргумента~$t$. В некотором смысле
ряд~\eqref{eq:1.1} обобщает представления \eqref{eq:0.10}
и~\eqref{eq:0.11} последовательности приближений Апери.
Разложение функции $R_n(t)$ в сумму простейших
дробей влечет
\begin{align}
F_n
&=\sum_{\nu=1}^\infty
\sum_{j=1}^{2s}\frac{j(j+1)}2\sum_{k=0}^n\frac{B_{jk}}{(\nu+k)^{j+2}}
=\sum_{j=1}^{2s}\frac{j(j+1)}2\sum_{k=0}^nB_{jk}
\sum_{l=k+1}^\infty\frac1{l^{j+2}}
\nonumber\\
&=\sum_{j=1}^{2s}\frac{j(j+1)}2\sum_{k=0}^nB_{jk}
\biggl(\zeta(j+2)-\sum_{l=1}^k\frac1{l^{j+2}}\biggr)
=\sum_{j=3}^{2s+2}A_j\zeta(j)-A_0,
\label{eq:1.3}
\end{align}
где
\begin{gather}
A_j=\frac{(j-1)(j-2)}2\sum_{k=0}^nB_{j-2,k}
\label{eq:1.4}
\\
A_0=\sum_{j=1}^{2s}\frac{j(j+1)}2\sum_{k=0}^nB_{jk}
\sum_{l=1}^k\frac1{l^{j+2}}.
\label{eq:1.5}
\end{gather}
Более половины коэффициентов \eqref{eq:1.4}
равны нулю: во-первых,
\begin{equation}
A_3=\sum_{k=0}^nB_{1k}
=\sum_{k=0}^n\Res_{t=-k}R_n(t)=-\Res_{t=\infty}R_n(t)=0
\label{eq:1.6}
\end{equation}
согласно соотношению
\begin{equation}
R_n(t)=O(t^{-2})
\qquad\text{при}\quad t\to\infty
\label{eq:1.7}
\end{equation}
и теореме о вычетах;
во-вторых, рациональная функция удовлетворяет
легко проверяемому свойству
$R_n(t)=-R_n(-t-n)$,
приводящему после подстановки в разложение~\eqref{eq:1.2}
к равенству
$B_{j,k}=(-1)^{j+1}B_{j,n-k}$, $k=0,1,\dots,n$,
откуда, в частности,
\begin{align*}
A_j
&=\frac{(j-1)(j-2)}2\sum_{k=0}^nB_{j-2,k}
\\
&=(-1)^{j-1}\frac{(j-1)(j-2)}2\sum_{k=0}^nB_{j-2,n-k}
=(-1)^{j-1}A_j
\end{align*}
и, значит,
\begin{equation}
A_4=A_6=\dots=A_{2s}=A_{2s+2}=0.
\label{eq:1.8}
\end{equation}
Оставшиеся ненулевые коэффициенты $A_0,A_5,A_7,\dots,A_{2s+1}$
согласно \eqref{eq:1.4}, \eqref{eq:1.5} и формуле
\begin{equation}
B_{jk}=\frac1{(2s-j)!}\,\frac{\d^{2s-j}}{\d t^{2s-j}}
\bigl(R_n(t)(t+k)^{2s}\bigr)\big|_{t=-k}
\label{eq:1.9}
\end{equation}
являются рациональными числами, знаменатели которых можно
найти, представив исходную функцию~$R_n(t)$
в виде произведения рациональных функций $2t+n$,
\begin{equation}
\begin{gathered}
\frac{n!}{\prod_{j=0}^n(t+j)}
=\sum_{k=0}^n\frac{(-1)^k\binom nk}{t+k},
\qquad
\frac{\prod_{j=1}^n(t-j)}{\prod_{j=0}^n(t+j)}
=\sum_{k=0}^n\frac{(-1)^{n-k}\binom nk\binom{n+k}n}{t+k},
\\
\frac{\prod_{j=1}^n(t+n+j)}{\prod_{j=0}^n(t+j)}
=\sum_{k=0}^n\frac{(-1)^k\binom nk\binom{2n-k}n}{t+k}.
\end{gathered}
\label{eq:1.10}
\end{equation}
Для каждой такой функции $R(t)$ имеют место включения
\begin{equation}
D_n^j\cdot\frac1{j!}\,\frac{\d^j}{\d t^j}
\bigl(R(t)(t+k)\bigr)\big|_{t=-k}\in\mathbb Z,
\qquad k=0,1,\dots,n, \quad j=0,1,2,\dots,
\label{eq:1.11}
\end{equation}
где $D_n$ обозначает наименьшее общее кратное
чисел $1,2,\dots,n$. Действительно, если
$$
R(t)=\sum_{l=0}^n\frac{B_l}{t+l},
$$
то
$$
R(t)(t+k)=B_k+\sum\doublesb{l=0}{l\ne k}^nB_l\biggl(1+\frac{k-l}{t+l}\biggr),
$$
так что
\begin{equation*}
\frac1{j!}\,\frac{\d^j}{\d t^j}\bigl(R(t)(t+k)\bigr)\big|_{t=-k}
=-\sum\doublesb{l=0}{l\ne k}^n\frac{B_l}{(k-l)^j},
\qquad k=0,1,\dots,n, \quad j=1,2,\dots,
\end{equation*}
а в случае $j=0$ правая часть просто равна $B_k$.
Представляя в~\eqref{eq:1.9} функцию $R_n(t)(t+k)^{2s}$ в виде
произведения $2t+n$ и $2s$ сомножителей вида $R(t)(t+k)$,
где $R(t)$~--- одна из функций в~\eqref{eq:1.10},
применяя правило дифференцирования Лейбница и включения
\eqref{eq:1.11}, мы получаем
$D_n^{2s-j}\cdot B_{jk}\in\mathbb Z$ при
$k=0,1,\dots,n$ и $j=1,\dots,2s$. Окончательно,
\begin{equation}
D_n^{2s+2}F_n
\in\mathbb Z\zeta(5)+\mathbb Z\zeta(7)
+\dots+\mathbb Z\zeta(2s+1)+\mathbb Z
\label{eq:1.12}
\end{equation}
согласно \eqref{eq:1.3}, \eqref{eq:1.4}, \eqref{eq:1.5},
\eqref{eq:1.6}, \eqref{eq:1.8}.

С другой стороны, в соответствии с~\eqref{eq:1.7} и
\begin{align}
&
\Res_{t=\nu}\biggl(R_n(t)\cdot\frac{\pi^3\cos\pi t}{\sin^3\pi t}\biggr)
\nonumber\\ &\quad
=\Res_{t=\nu}
\Bigl(R_n(\nu)+R_n'(\nu)(t-\nu)
+\frac12R_n''(\nu)(t-\nu)^2
+O((t-\nu)^3)\Bigr)
\nonumber\\ &\quad\qquad\times
\bigl((t-\nu)^{-3}+O(1)\bigr)
\nonumber\\ &\quad
=\frac12R_n''(\nu),
\qquad \nu\in\mathbb Z,
\label{eq:1.13}
\end{align}
а также с учетом $R_n(\nu)=R_n'(\nu)=R_n''(\nu)=0$ для $\nu=1,2,\dots,n$
заключаем, что величина~\eqref{eq:1.1} допускает представление
в виде барнсовского интеграла~\cite{Bar}
\begin{align*}
F_n
&=-\frac1{2\pi i}\int_{t_0-i\infty}^{t_0+i\infty}
R_n(t)\biggl(\frac\pi{\sin\pi t}\biggr)^3\cos\pi t\,\d t
\\
&=\Re\biggl(\frac{(-1)^nin!^{2s-6}}{2\pi}
\int_{t_0-i\infty}^{t_0+i\infty}(2t+n)
\\ &\qquad\times
\frac{\Gamma(t)^{2s+3}\Gamma(n+1-t)^3\Gamma(2n+1+t)^3}
{\Gamma(n+1+t)^{2s+3}}
e^{-\pi it}\,\d t\biggr),
\end{align*}
где $0<t_0<n+1$. Используя асимптотику гамма-функции
\begin{equation}
\log\Gamma(z)=\Bigl(z-\frac12\Bigr)\log z-z+\frac12\log2\pi
+O\bigl(|\Re z|^{-1}\bigr),
\label{eq:1.14}
\end{equation}
где постоянная в $O(\,\cdot\,)$ абсолютная, и осуществляя
замену $t=n\tau$, находим
\begin{equation}
F_n
=\Re\biggl(\frac{(-1)^n(2\pi)^{s-1}i}{n^{s+1}}
\int_{\tau_0-i\infty}^{\tau_0+i\infty}
e^{(n+1)f(\tau)+g(\tau)}\bigl(1+O(n^{-1})\bigr)\,\d\tau\biggr),
\label{eq:1.15}
\end{equation}
где $\tau_0$~--- произвольная постоянная из интервала $0<\tau_0<1$,
\begin{align*}
f(\tau)
&=(2s+3)\tau\log\tau+3(1-\tau)\log(1-\tau)+3(2+\tau)\log(2+\tau)
\\ &\qquad
-(2s+3)(1+\tau)\log(1+\tau)-\pi i\tau,
\displaybreak[0]\\
g(\tau)
&=\log(1+2\tau)+\frac12(2s+3)\log(1+\tau)
-\frac12(2s+3)\log\tau
\\ &\qquad
-\frac32\log(1-\tau)-\frac32\log(2+\tau).
\end{align*}
Можно показать \cite[лемма~2.7]{Z8}, что
уравнение $f'(\tau)=0$ имеет два решения;
они расположены в области $\Im\tau>0$ и симметричны
относительно прямой $\Re\tau=-\frac12$.
Выберем то решение $\tau=\tau_1$ уравнения $f'(\tau)=0$,
для которого $\Re\tau>-\frac12$, и положим $\tau_0=\Re\tau_1$.
Тогда $0<\tau_0<1$. Применяя к интегралу в~\eqref{eq:1.15}
с $\tau_0=\Re\tau_1$ метод перевала, получаем
\begin{equation}
\limsup_{n\to\infty}\frac{\log|F_n|}n
=\Re f(\tau_1)
\label{eq:1.16}
\end{equation}
при условии, что $\Im f(\tau_1)\notin\pi\mathbb Z$.
Таким образом, согласно \eqref{eq:1.12} и \eqref{eq:1.16}
в случае $2s+2+\Re f(\tau_1)<0$ заключаем, что
хотя бы одно из чисел
$\zeta(5),\zeta(7),\dots,\zeta(2s+1)$
иррационально. Наконец, выбирая $s=10$, находим
\begin{gather*}
\tau_1=0.99223412\hdots+i0.01200539\hdots,
\\
f(\tau_1)=-22.02001639\hdots-i3.10440862\hdots;
\end{gather*}
следовательно, по крайней мере одно из девяти чисел
$\zeta(5),\zeta(7),\dots,\zeta(21)$ иррационально.

Дальнейшее изложение следует работам \cite{Z3}, \cite{Z4},
\cite{Z5}, \cite{Z6}, \cite{Z8}, \cite{Z17}.

\section{Арифметика простейших рациональных функций}
\label{sec:1.1}

Перед описанием общей конструкции линейных форм от $1$ и
$\zeta(5),\linebreak[3]\zeta(7),\dots$ мы приведем свойства
простейших рациональных функций
\begin{equation}
R(t)=R(a,b;t)=\begin{cases}
\dfrac{(t+b)(t+b+1)\dotsb(t+a-1)}{(a-b)!}
& \text{при $a\ge b$}, \\
\dfrac{(b-a-1)!}{(t+a)(t+a+1)\dotsb(t+b-1)}
& \text{при $a<b$},
\end{cases}
\label{eq:1.17}
\end{equation}
играющих роль строительных {\it кирпичиков\/} в
дальнейшем описании.

Следующее утверждение содержит хорошо известные свойства
целозначных многочленов.

\begin{lemma}[{см\. \cite[лемма~1]{Bak} или \cite[лемма~7]{Z2}}]
\label{lem:1.1}
Пусть $a\ge b$. Тогда для любого неотрицательного целого~$j$
выполняются включения
$$
D_{a-b}^j\cdot\frac1{j!}R^{(j)}(-k)\in\mathbb Z, \qquad k\in\mathbb Z.
$$
\end{lemma}

Следующее утверждение по существу повторяет доказательство
включений \eqref{eq:1.11} для первой из функций в~\eqref{eq:1.10}.

\begin{lemma}
\label{lem:1.2}
Пусть целые числа $a,b,a_0,b_0$ удовлетворяют неравенствам
$a_0\le a<b\le b_0$. Тогда для любого неотрицательного целого~$j$
выполняются включения
$$
D_{b_0-a_0-1}^j\cdot\frac1{j!}\bigl(R(t)(t+k)\bigr)^{(j)}\big|_{t=-k}
\in\mathbb Z,
\qquad k=a_0,a_0+1,\dots,b_0-1.
$$
\end{lemma}

Леммы \ref{lem:1.1} и \ref{lem:1.2}
дают частичную (но крайне важную) информацию
о $p$-адической оценке рациональных чисел
\begin{equation}
R^{(j)}(-k)
\quad\text{и}\quad
\bigl(R(t)(t+k)\bigr)^{(j)}\big|_{t=-k}
\label{eq:1.18}
\end{equation}
соответственно: для каждого простого~$p$ из интервала
$\sqrt N<p\le N$ имеем $\ord_pD_N=1$.
Следующие два утверждения посвящены уточнению
полученных $p$-адических оценок величин~\eqref{eq:1.18}.

\begin{lemma}
\label{lem:1.3}
Пусть $a,b,a_0,b_0$~целые, $b_0\le b<a\le a_0$, и пусть
рациональная функция $R(t)=R(a,b;t)$ определена в~\eqref{eq:1.17}.
Тогда для любого целого~$k$, $b_0\le k<a_0$,
простого $p>\sqrt{a_0-b_0-1}$ и
неотрицательного целого~$j$ имеют место оценки
\begin{align}
\ord_pR^{(j)}(-k)
&\ge-j+\biggl\[\frac{a-1-k}p\biggr\]
-\biggl\[\frac{b-1-k}p\biggr\]-\biggl\[\frac{a-b}p\biggr\]
\nonumber\\
&=-j+\biggl\[\frac{k-b}p\biggr\]
-\biggl\[\frac{k-a}p\biggr\]-\biggl\[\frac{a-b}p\biggr\].
\label{eq:1.19}
\end{align}
\end{lemma}

\begin{proof}
Выберем произвольное $p>\sqrt{a_0-b_0-1}$.
Прежде всего отметим, что согласно определению целой части числа
выполнено
$$
\[-x\]=-\[x\]-\delta_x,
\qquad\text{где}\quad
\delta_x=\begin{cases}
0, & \text{если $x\in\mathbb Z$}, \\
1, & \text{если $x\notin\mathbb Z$},
\end{cases}
$$
откуда
$$
\biggl\[-\frac sp\biggr\]=-\biggl\[\frac{s-1}p\biggr\]-1
\qquad\text{для}\quad s\in\mathbb Z.
$$
Следовательно,
\begin{equation}
\biggl\[\frac{k-b}p\biggr\]
=-\biggl\[\frac{b-1-k}p\biggr\]-1,
\qquad
\biggl\[\frac{a-1-k}p\biggr\]
=-\biggl\[\frac{k-a}p\biggr\]-1
\label{eq:1.20}
\end{equation}
для любого целого~$k$.

Непосредственные вычисления показывают, что
$$
R(-k)=\begin{cases}
\dfrac{(a-1-k)!}{(b-1-k)!\,(a-b)!}
& \text{при $k<b$}, \\
0
& \text{при $b\le k<a$}, \\
(-1)^{a-b}\dfrac{(k-b)!}{(k-a)!\,(a-b)!}
& \text{при $k\ge a$};
\end{cases}
$$
следовательно,
$$
\begin{alignedat}{2}
\ord_pR(-k)
&\ge\biggl\[\frac{a-1-k}p\biggr\]
-\biggl\[\frac{b-1-k}p\biggr\]-\biggl\[\frac{a-b}p\biggr\]
\qquad&& \text{при $k<a$},
\\
\ord_pR(-k)
&\ge\biggl\[\frac{k-b}p\biggr\]
-\biggl\[\frac{k-a}p\biggr\]-\biggl\[\frac{a-b}p\biggr\]
\qquad&& \text{при $k\ge b$},
\end{alignedat}
$$
что отвечает оценкам \eqref{eq:1.19} при $j=0$ согласно
равенству~\eqref{eq:1.20}.

Если $k<b$ или $k\ge a$, рассмотрим функцию
$$
r(t)=\frac{R'(t)}{R(t)}
=\sum_{l=b}^{a-1}\frac1{t+l},
$$
для которой при каждом целом $j\ge1$ имеют место включения
$$
r^{(j-1)}(-k)\cdot D_{\max\{a-b_0-1,a_0-b-1\}}^{j-1}\in\mathbb Z.
$$
Индукция по~$j$ и тождество
\begin{equation}
R^{(j)}(t)
=\bigl(R(t)r(t)\bigr)^{(j-1)}
=\sum_{m=0}^{j-1}\binom{j-1}mR^{(m)}(t)r^{(j-1-m)}(t)
\label{eq:1.21}
\end{equation}
с последующей подстановкой $t=-k$
приводит к требуемым оценкам~\eqref{eq:1.19}.

В случае $b\le k<a$ рассмотрим функции
$$
R_k(t)=\frac{R(t)}{t+k},
\qquad
r_k(t)=\frac{R_k'(t)}{R_k(t)}
=\sum\doublesb{l=b}{l\ne k}^{a-1}\frac1{t+l};
$$
очевидно, для каждого целого $j\ge1$ справедливы включения
$$
r_k^{(j-1)}(-k)\cdot D_{a-b-1}^{j-1}\in\mathbb Z.
$$
Тогда
$$
R^{(j)}(-k)=jR_k^{(j-1)}(-k),
$$
поскольку
$$
R_k(-k)=(-1)^{k-b}\frac{(k-b)!\,(a-1-k)!}{(a-b)!},
$$
так что индукция по~$j$ вместе с равенством~\eqref{eq:1.21}
(где мы подставляем $R_k(t),\linebreak[2]r_k(t)$
вместо $R(t),r(t)$ соответственно) показывают, что
\begin{align*}
\ord_pR^{(j)}(-k)
&\ge\ord_pR_k^{(j-1)}(-k)
\\
&\ge-(j-1)+\biggl\[\frac{k-b}p\biggr\]
+\biggl\[\frac{a-1-k}p\biggr\]-\biggl\[\frac{a-b}p\biggr\]
\end{align*}
для каждого $j\ge1$. Таким образом, применяя~\eqref{eq:1.20},
мы вновь получаем требуемые оценки~\eqref{eq:1.19}.
Лемма доказана.
\end{proof}

\begin{lemma}
\label{lem:1.4}
Пусть $a,b,a_0,b_0$~целые, $a_0\le a<b\le b_0$, и пусть
рациональная функция $R(t)=R(a,b;t)$ определена в~\eqref{eq:1.17}.
Тогда для любого целого~$k$, $a_0\le k<b_0$,
простого $p>\sqrt{b_0-a_0-1}$ и
неотрицательного целого~$j$ имеют место оценки
\begin{equation}
\ord_p\bigl(R(t)(t+k)\bigr)^{(j)}\big|_{t=-k}
\ge-j+\biggl\[\frac{b-a-1}p\biggr\]
-\biggl\[\frac{k-a}p\biggr\]-\biggl\[\frac{b-1-k}p\biggr\].
\label{eq:1.22}
\end{equation}
\end{lemma}

\begin{proof}
Зафиксируем произвольное $p>\sqrt{b_0-a_0-1}$. Имеем
$$
\bigl(R(t)(t+k)\bigr)\big|_{t=-k}=\begin{cases}
(-1)^{k-a}\dfrac{(b-a-1)!}{(k-a)!\,(b-1-k)!}
& \text{при $a\le k<b$}, \\
0 & \text{при $k<a$ или $k\ge b$},
\end{cases}
$$
что означает справедливость оценок~\eqref{eq:1.22} при $j=0$.

Рассматривая в случае $a\le k<b$ функции
$$
R_k(t)=R(t)(t+k),
\qquad
r_k(t)=\frac{R_k'(t)}{R_k(t)}
=\sum\doublesb{l=a}{l\ne k}^{b-1}\frac1{t+l},
$$
и проводя индукцию по $j\ge0$,
с помощью равенства~\eqref{eq:1.21}
(где снова функции $R(t),r(t)$ заменены на $R_k(t),r_k(t)$)
мы приходим к оценкам~\eqref{eq:1.22}.

Если же $k<a$ или $k\ge b$, заметим, что
$$
\bigl(R(t)(t+k)\bigr)^{(j)}\big|_{t=-k}
=jR^{(j-1)}(-k).
$$
Поскольку
$$
R(-k)=\begin{cases}
\dfrac{(b-a-1)!\,(a-1-k)!}{(b-1-k)!}
& \text{при $k<a$}, \\
(-1)^{b-a}\dfrac{(b-a-1)!\,(k-b)!}{(k-a)!}
& \text{при $k\ge b$},
\end{cases}
$$
индукция по~$j$ и равенства~\eqref{eq:1.20} вновь
приводят к требуемым оценкам \eqref{eq:1.22}.
Лемма доказана.
\end{proof}

Для вычисления асимптотики величин,
связанных с $p$-адически уточненными
оценками лемм~\ref{lem:1.3},~\ref{lem:1.4},
нам понадобится следующая
важная составляющая арифметического метода.

\begin{lemma}[{ср\. с~\cite[теорема~4.3 и~\S\,6]{Ch2} и \cite[лемма~3.2]{Ha1}}]
\label{lem:1.5}
Для любого $C>0$ и любого полуинтервала $[u,v)\subset(0,1)$
имеет место предельное соотношение
\begin{equation*}
\lim_{n\to\infty}\frac1n\sum\doublesb{p>\sqrt{Cn}}{\{n/p\}\in[u,v)}\log p
=\psi(v)-\psi(u)
=\int_u^v\d\psi(z),
\end{equation*}
где $\{a\}=a-\[a\]$ --- дробная часть числа и $\psi(z)$~--- логарифмическая
производная гамма-функции.
\end{lemma}

\section{Линейные формы от~$1$ и нечетных дзета-значений}
\label{sec:1.2}

Рассмотрим положительные нечетные числа $q$ и $r$,
удовлетворяющие неравенству $q\ge r+4$.
Набору положительных целочисленных параметров
$$
\bh=(h_0;h_1,\dots,h_q),
$$
связанных условием
\begin{equation}
h_1+h_2+\dots+h_q\le h_0\cdot\frac{q-r}2,
\label{eq:1.23}
\end{equation}
поставим в соответствие функцию
\begin{align}
\wt R(t)
&=\wt R(\bh;t)
\nonumber\\
&=(h_0+2t)\frac{\Gamma(h_0+t)^r\Gamma(h_1+t)\dotsb\Gamma(h_q+t)}
{\Gamma(1+t)^r\Gamma(1+h_0-h_1+t)\dotsb\Gamma(1+h_0-h_q+t)}.
\label{eq:1.24}
\end{align}
Согласно~\eqref{eq:1.23} выполнено
\begin{equation}
\wt R(t)=O(t^{-2})
\qquad\text{при}\quad t\to\infty,
\label{eq:1.25}
\end{equation}
так что ряд
\begin{equation}
\wt F(\bh)
=\frac1{(r-1)!}\sum_{t=0}^\infty\wt R^{(r-1)}(t)
\label{eq:1.26}
\end{equation}
абсолютно сходится. Если $r=1$, величина~\eqref{eq:1.26} может быть
записана как совершенно уравновешенный ряд
\begin{align*}
\wt F(\bh)
&=\frac{h_0!\,(h_1-1)!\dotsb(h_q-1)!}
{(h_0-h_1)!\dotsb(h_0-h_q)!}
\\ &\qquad\times 
{}_{q+2}F_{q+1}\biggl(\begin{matrix}
h_0, & 1+\frac12h_0, &       h_1, & \dots, &       h_q \\[1pt]
     &   \frac12h_0, & 1+h_0-h_1, & \dots, & 1+h_0-h_q
\end{matrix}\biggm|1\biggr),
\end{align*}
в то время как при $r>1$ ряд~\eqref{eq:1.26} представляется
в виде линейной комбинации $G$-функций Мейера~\cite{Lu} (интегралов
барнсовского типа), вычисляемых в точках
$e^{\pi ik}$, где $k=\pm1,\pm3,\dots,\pm(r-2)$.

Рациональная функция \eqref{eq:1.24} обладает симметрией
относительно преобразования $t\mapsto-t-h_0$:
\begin{equation}
\wt R(-t-h_0)=(-1)^{(h_0-1)(q+r)}\wt R(t)=-\wt R(t),
\label{eq:1.27}
\end{equation}
где мы пользуемся тождеством
\begin{equation*}
\Gamma(t+b)\Gamma(1-t-b)=(-1)^b\frac\pi{\sin\pi t}.
\end{equation*}
Покажем, что величина \eqref{eq:1.26} является линейной
формой от~$1$ и нечетных дзета-значений.

Для точной формулировки результата будем
предполагать набор параметров $\bh$ упорядоченным:
$$
h_1\le h_2\le\dots\le h_q<\frac12h_0;
$$
определим следующую арифметическую нормализацию
величины~\eqref{eq:1.26}:
\begin{equation}
F(\bh)
=\frac{\prod_{j=r+1}^q(h_0-2h_j)!}
{\prod_{j=1}^r(h_j-1)!^2}\cdot\wt F(\bh)
=\frac1{(r-1)!}\sum_{t=1-h_1}^\infty R^{(r-1)}(t),
\label{eq:1.28}
\end{equation}
где рациональная функция
\begin{align}
R(t)
&=(h_0+2t)\prod_{j=1}^r\frac1{(h_j-1)!}\,
\frac{\Gamma(h_j+t)}{\Gamma(1+t)}
\cdot\prod_{j=1}^r\frac1{(h_j-1)!}\,
\frac{\Gamma(h_0+t)}{\Gamma(1+h_0-h_j+t)}
\nonumber\\ &\qquad\times
\prod_{j=r+1}^q(h_0-2h_j)!\,
\frac{\Gamma(h_j+t)}{\Gamma(1+h_0-h_j+t)}
\label{eq:1.29}
\end{align}
является произведением простейших функций~\eqref{eq:1.17}.
Положим $m_0=\linebreak[4]\max\{h_r-1,h_0-2h_{r+1}\}$
и $m_j=\max\{m_0,h_0-h_1-h_{r+j}\}$
для $j=1,\dots,q-r$ и определим целочисленную величину
\begin{equation}
\Phi=\Phi(\bh)=\prod_{\sqrt{h_0}<p\le m_{q-r}}p^{\nu_p},
\label{eq:1.30}
\end{equation}
где
\begin{equation}
\nu_p
=\min_{h_{r+1}\le k\le h_0-h_{r+1}}\{\nu_{k,p}\}
\label{eq:1.31}
\end{equation}
и
\begin{align*}
\nu_{k,p}
&=\sum_{j=1}^r\biggl(\biggl\[\frac{k-1}p\biggr\]
+\biggl\[\frac{h_0-k-1}p\biggr\]
\\ &\qquad
-\biggl\[\frac{k-h_j}p\biggr\]
-\biggl\[\frac{h_0-h_j-k}p\biggr\]
-2\biggl\[\frac{h_j-1}p\biggr\]\biggr)
\\ &\phantom:\qquad
+\sum_{j=r+1}^q\biggl(\biggl\[\frac{h_0-2h_j}p\biggr\]
-\biggl\[\frac{k-h_j}p\biggr\]
-\biggl\[\frac{h_0-h_j-k}p\biggr\]\biggr).
\end{align*}
В этих обозначениях имеет место

\begin{proposition}
\label{prop:1.1}
Величина~\eqref{eq:1.28} является линейной формой от
$1,\zeta(r+\nobreak2),\linebreak[2]\zeta(r+4),\dots,\zeta(q-4),\zeta(q-2)$
с рациональными коэффициентами. Более того,
\begin{equation}
D_{m_1}^rD_{m_2}\dotsb D_{m_{q-r}}
\cdot\Phi^{-1}\cdot F(\bh)
\in\mathbb Z\zeta(q-2)+\mathbb Z\zeta(q-4)+\dots
+\mathbb Z\zeta(r+2)+\mathbb Z.
\label{eq:1.32}
\end{equation}
\end{proposition}

\begin{proof}
Применяя правило Лейбница дифференцирования произведения,
а также леммы~\ref{lem:1.1}, \ref{lem:1.2}
и леммы~\ref{lem:1.3},~\ref{lem:1.4}
к рациональной функции~\eqref{eq:1.29}, мы видим,
что числа
$$
\begin{gathered}
B_{jk}=\frac1{(q-j)!}
\cdot\bigl(R(t)(t+k)^{q-r}\bigr)^{(q-j)}\big|_{t=-k},
\\
j=r+1,\dots,q, \quad k=h_{r+1},\dots,h_0-h_{r+1},
\end{gathered}
$$
удовлетворяют соотношениям
\begin{equation}
D_{m_0}^{q-j}\cdot B_{jk}\in\mathbb Z
\label{eq:1.33}
\end{equation}
и
\begin{equation}
\ord_pB_{jk}
\ge-(q-j)+\nu_{k,p},
\label{eq:1.34}
\end{equation}
соответственно, для любого $k=h_{r+1},\dots,h_0-h_{r+1}$
и любого простого $p>\nobreak\sqrt{h_0}$.
Далее, разложение
$$
R(t)=\sum_{j=r+1}^q\sum_{k=h_j}^{h_0-h_j}\frac{B_{jk}}{(t+k)^{j-r}}
$$
приводит к ряду
\begin{align*}
F(\bh)
&=\sum_{j=r+1}^q\binom{j-2}{r-1}
\sum_{k=h_j}^{h_0-h_j}B_{jk}
\biggl(\sum_{l=1}^\infty-\sum_{l=1}^{k-h_1}\biggr)\frac1{l^{j-1}}
\\
&=\sum_{j=r+1}^qA_{j-1}\zeta(j-1)-A_0,
\end{align*}
где
\begin{align}
A_{j-1}&=\binom{j-2}{r-1}\sum_{k=h_j}^{h_0-h_j}B_{jk},
\qquad j=r+1,\dots,q,
\label{eq:1.35}
\\
A_0&=\sum_{j=r+1}^q\binom{j-2}{r-1}
\sum_{k=h_j}^{h_0-h_j}B_{jk}
\sum_{l=1}^{k-h_1}\frac1{l^{j-1}}.
\nonumber
\end{align}
Согласно~\eqref{eq:1.33} и включениям
$$
D_{m_1}^rD_{m_2}\dotsb D_{m_{j-r}}
\cdot\sum_{l=1}^{k-h_1}\frac1{l^{j-1}}\in\mathbb Z,
$$
справедливым для $k=h_j,\dots,h_0-h_j$, $j=r+1,\dots,q$,
мы получаем включения
\begin{gather*}
D_{m_0}^{q-j-1}\cdot A_j\in\mathbb Z
\qquad\text{при}\quad j=r,r+1,\dots,q-1,
\\
D_{m_1}^rD_{m_2}\dotsb D_{m_{q-r}}
\cdot A_0\in\mathbb Z,
\end{gather*}
которые несколько уточняются оценками~\eqref{eq:1.34}:
$$
\ord_pA_j\ge-(q-j-1)+\nu_p
\qquad\text{при $j=0$ и $j=r,r+1,\dots,q-1$},
$$
где показатели $\nu_p$ определены в~\eqref{eq:1.31}.
Таким образом, для завершения доказательства
мы должны показать, что
$$
A_r=0 \qquad\text{и}\qquad
A_{r+1}=A_{r+3}=\dots=A_{q-3}=A_{q-1}=0.
$$
Первое равенство следует из~\eqref{eq:1.25};
в соответствии с~\eqref{eq:1.27} получаем
$$
B_{jk}=(-1)^jB_{j,h_0-k}
\qquad\text{при $j=r+1,\dots,q$},
$$
что влечет $A_{j-1}=0$ для $j$ нечетных
согласно~\eqref{eq:1.35}.
Лемма доказана полностью.
\end{proof}

\section{Асимптотика линейных форм}
\label{sec:1.3}

Для исследования асимптотического поведения
построенных линейных форм~\eqref{eq:1.32} определим
последовательность целочисленных направлений
$\Beta=(\eta_0;\eta_1,\dots,\eta_q)$, связанных
для каждого целого $n\ge0$ с исходным набором
параметров~$\bh$ в соответствии со следующими формулами:
\begin{equation}
h_0=\eta_0n+2 \qquad\text{и}\qquad
h_j=\eta_jn+1 \quad\text{при $j=1,\dots,q$}.
\label{eq:1.36}
\end{equation}

Положим
$$
m_j=\max\{\eta_r,\eta_0-2\eta_{r+1},\eta_0-\eta_1-\eta_{r+j}\}
\qquad\text{при}\quad j=1,\dots,q-r
$$
(т.е\. старое множество параметров с этими именами
в $n$~раз больше нового). Асимптотика
величины~\eqref{eq:1.30} при $n\to\infty$
может быть вычислена с помощью леммы~\ref{lem:1.5}
в терминах целозначной функции
\begin{align*}
\varphi_0(x,y)
&=\sum_{j=1}^r\bigl(\[y\]+\[\eta_0x-y\]
-\[y-\eta_jx\]-\[(\eta_0-\eta_j)x-y\]-2\[\eta_jx\]\bigr)
\\ &\qquad
+\sum_{j=r+1}^q\bigl(\[(\eta_0-2\eta_j)x\]
-\[y-\eta_jx\]-\[(\eta_0-\eta_j)x-y\]\bigr),
\end{align*}
$1$-периодичной по каждому из аргументов $x$ и~$y$.
Согласно \eqref{eq:1.31} и~\eqref{eq:1.36} находим
$$
\nu_p=\min_{\eta_{r+1}n\le k-1\le(\eta_0-\eta_{r+1})n}
\varphi_0\biggl(\frac np,\frac{k-1}p\biggr)
\ge\varphi\biggl(\frac np\biggr),
$$
где
$$
\varphi(x)=\min_{y\in\mathbb R}\varphi_0(x,y)
=\min_{0\le y<1}\varphi_0(x,y).
$$
Кроме того,
$$
\lim_{n\to\infty}\frac{\log D_{m_jn}}n=m_j,
\qquad j=1,\dots,q-r,
$$
согласно следствию из асимптотического закона распределения
простых чисел. Поэтому имеет место следующий результат.

\begin{proposition}
\label{prop:1.2}
В приведенных выше обозначениях
\begin{align*}
&
\lim_{n\to\infty}
\frac{\log(D_{m_1n}^rD_{m_2n}\dotsb D_{m_{q-r}n}\Phi(\bh)^{-1})}n
\\ &\qquad
=rm_1+m_2+\dots+m_{q-r}
-\biggl(\int_0^1\varphi(x)\,\d\psi(x)
-\int_0^{1/m_{q-r}}\varphi(x)\,\frac{\d x}{x^2}\biggr).
\end{align*}
\end{proposition}

Перейдем теперь к асимптотической оценке величины $F(\bh)$.
Здесь мы ограничимся подробным рассмотрением случая
$r=3$, который необходим для доказательства теоремы~\ref{th:1}.
Наиболее общая ситуация исследована нами в работе~\cite{Z8}
(см\. также \cite{Z3}, \cite{Z4}, \cite{Z5}).

\begin{lemma}
\label{lem:1.6}
Максимальное значение вещественнозначной
неотрицательной функции $h(y)=|\ctg(x+iy)|$, $y\in\mathbb R$,
зависит только от $x\pmod{\pi\mathbb Z}\in\mathbb R$.
\end{lemma}

\begin{proof}
Если $x=\pi k$ для некоторого целого~$k$, то точка
$z=\nobreak x$ является (единственным) полюсом функции $\ctg z$
на прямой $\Im z=\nobreak x$. Следовательно,
в этом случае абсолютный максимум функции~$h(y)$,
равный бесконечности, достигается при $y=0$.
Поэтому считаем далее $x\notin\pi\mathbb Z$ и, значит, $\cos2x<1$.
Имеем
\begin{align}
h(y)^2
&=\biggl|\frac{\cos(x+iy)}{\sin(x+iy)}\biggr|^2
=\biggl|\frac{(e^{-y}+e^y)\cos x+i(e^{-y}-e^y)\sin x}
{(e^{-y}-e^y)\cos x+i(e^{-y}+e^y)\sin x}\biggr|^2
\nonumber\\
&=\frac{e^{-2y}+e^{2y}+2\cos2x}{e^{-2y}+e^{2y}-2\cos2x}
=1+\frac{4\cos2x}{e^{-2y}+e^{2y}-2\cos2x}
\nonumber\\
&\le1+\frac{4|\cos2x|}{2-2\cos2x},
\label{eq:1.37}
\end{align}
где мы воспользовались неравенством
$Y+1/Y\ge2$ при $Y=e^{-2y}>0$, превращающимся в равенство
только в случае $Y=1$.
Оценка в правой части~\eqref{eq:1.37} зависит только
от $x\pmod{\pi\mathbb Z}$, что и завершает доказательство леммы.
\end{proof}

\begin{lemma}
\label{lem:1.7}
Пусть $r=3$. Тогда для суммы~\eqref{eq:1.28}
справедливо интегральное представление
\begin{align}
F(\bh)
&=-\frac1{2\pi i}\int_{t_0-i\infty}^{t_0+i\infty}
R(t)\biggl(\frac\pi{\sin\pi t}\biggr)^3\cos\pi t\,\d t
\nonumber\\
&=\Re\biggl(\frac i{2\pi}\int_{t_0-i\infty}^{t_0+i\infty}
R(t)\biggl(\frac\pi{\sin\pi t}\biggr)^3e^{-\pi it}\,\d t\biggr),
\label{eq:1.38}
\end{align}
где $t_0\in\mathbb R$~--- произвольная постоянная
из интервала $1-h_1<t_0<0$.
\end{lemma}

\begin{proof}
Без ограничения общности, считаем $t_0$ нецелым.
Рассмотрим подынтегральную функцию в~\eqref{eq:1.38}
на положительно-ориен\-тируемом контуре
прямоугольника с вершинами $t_0\pm iN$, $N+1/2\pm iN$, где
{\it целое\/} число $N>0$ достаточно велико.
Согласно теореме Коши, соотношению \eqref{eq:1.13}
(примененному к~$R(t)$) и аналитичности функции~$R(t)$
внутри и на границе прямоугольника интеграл
\begin{align}
&
\frac1{2\pi i}\biggl(\int_{t_0+iN}^{t_0-iN}
+\int_{t_0-iN}^{N+1/2-iN}+\int_{N+1/2-iN}^{N+1/2+iN}
+\int_{N+1/2+iN}^{t_0+iN}\biggr)
R(t)\biggl(\frac\pi{\sin\pi t}\biggr)^3\cos\pi t\,\d t
\label{eq:1.39}
\end{align}
равен сумме вычетов подынтегральной функции в точках
$t\in\mathbb Z$, $t_0<t\le\nobreak N$, именно
\begin{equation}
\sum_{t_0<\nu\le N}\frac{R^{(2)}(\nu)}2
=\sum_{\nu=0}^N\frac{R^{(2)}(\nu)}2.
\label{eq:1.40}
\end{equation}
В последнем выражении мы также воспользовались тем,
что $R(t)$ имеет нули порядка не ниже~$3$ в точках
$t=-1,-2,\dots,1-h_1$.
Отметим теперь, что на сторонах
$[N+1/2-iN,N+1/2+iN]$, $[t_0-iN,N+1/2-iN]$
и $[N+1/2+iN,t_0+iN]$ прямоугольника выполнено
$R(t)=O(N^{-2})$ согласно~\eqref{eq:1.25}, а функция
$$
\biggl(\frac\pi{\sin\pi t}\biggr)^3\cos\pi t
=\pi^3(\ctg\pi t+\ctg^3\pi t)
$$
ограничена. Последнее утверждение
следует из ограниченности функции
$\ctg\pi t$ на рассматриваемых сторонах:
для отрезка $[N+1/2-iN,N+1/2+iN]$ мы пользуемся
леммой~\ref{lem:1.6}, а для двух других
отрезков --- тривиальной оценкой
$$
|\ctg\pi(x+iy)|
=\biggl|\frac{1+e^{\pm2\pi ix}\cdot e^{-2\pi|y|}}
{1-e^{\pm2\pi ix}\cdot e^{-2\pi|y|}}\biggr|
\le\frac{1+e^{-2\pi|y|}}{1-e^{-2\pi|y|}},
\qquad x\in\mathbb R,
$$
при $y=\pm N$. Поэтому
$$
\frac1{2\pi i}\biggl(\int_{t_0-iN}^{N+1/2-iN}
+\int_{N+1/2-iN}^{N+1/2+iN}
+\int_{N+1/2+iN}^{t_0+iN}\biggr)
R(t)\biggl(\frac\pi{\sin\pi t}\biggr)^3\cos\pi t\,\d t
=O(N^{-1})
$$
и после перехода к пределу $N\to\infty$ в~\eqref{eq:1.39}
мы получаем правую часть \eqref{eq:1.38}, в то время как
предельный переход в~\eqref{eq:1.40} дает
искомую сумму \eqref{eq:1.28}. Лемма доказана.
\end{proof}

Записывая
\begin{align*}
R(t)\biggl(\frac\pi{\sin\pi t}\biggr)^3
&=-(h_0+2t)\prod_{j=1}^3\frac{\Gamma(h_j+t)\,\Gamma(-t)}{(h_j-1)!}\,
\cdot\prod_{j=1}^3\frac1{(h_j-1)!}\,
\frac{\Gamma(h_0+t)}{\Gamma(1+h_0-h_j+t)}
\\ &\qquad\times
\prod_{j=4}^q(h_0-2h_j)!\,
\frac{\Gamma(h_j+t)}{\Gamma(1+h_0-h_j+t)},
\end{align*}
применяя на контуре в~\eqref{eq:1.38}
асимптотику гамма-функции~\eqref{eq:1.14}
и делая замену переменной $t=n(\tau-\eta_0)$,
представим интеграл из леммы~\ref{lem:1.7} в виде
\begin{equation}
F(\bh)=\Re\biggl(\frac{n^{-q/2}}{2\pi i}\int_{\tau_0-i\infty}^{\tau_0+i\infty}
e^{nf(\tau)+g(\tau)}\bigl(1+O(n^{-1})\bigr)\,\d\tau\biggr),
\label{eq:1.41}
\end{equation}
где $\tau_0$~--- произвольная вещественная постоянная из интервала
$\eta_0-\eta_1<\tau_0<\eta_0$, а функции
\begin{align*}
f(\tau)
&=3\tau\log\tau+3(\eta_0-\tau)\log(\eta_0-\tau)
\\ &\qquad
-\sum_{j=1}^q\bigl((\tau-\eta_j)\log(\tau-\eta_j)
-(\tau-\eta_0+\eta_j)\log(\tau-\eta_0+\eta_j)\bigr)
\\ &\qquad
-\pi i\tau-2\sum_{j=1}^3\eta_j\log\eta_j
+\sum_{j=4}^q(\eta_0-2\eta_j)\log(\eta_0-2\eta_j)
\end{align*}
и $g(\tau)$ голоморфны в $\tau$-плоскости с разрезами по лучам
$(-\infty,\eta_0-\eta_1]$ и $[\eta_0,+\infty)$.
Поэтому применение метода перевала (см., например, \cite[гл.~5]{Br})
приводит к следующему утверждению.

\begin{proposition}
\label{prop:1.3}
Пусть уравнение $f'(\tau)=0$ имеет единственный
корень $\tau=\tau_1$, расположенный в полосе
$\eta_0-\eta_1<\Re\tau<\eta_0$, причем
$\Im f(\tau_1)\notin\pi\mathbb Z$. Тогда
$\tau=\tau_1$ является точкой перевала для
интеграла в~\eqref{eq:1.41} с $\tau_0=\Re\tau_1$,
так что
$$
\limsup_{n\to\infty}\frac{\log|F(\bh)|}n
=\Re f(\tau_1).
$$
\end{proposition}

\begin{remark}
Практическое вычисление корня $\tau=\tau_1$
основано на том, что он является нулем многочлена
$$
(\tau-\eta_0)^3(\tau-\eta_1)\dotsb(\tau-\eta_q)
-\tau^3(\tau-\eta_0+\eta_1)\dotsb(\tau-\eta_0+\eta_q)
$$
в области $\Im\tau>0$ и $\eta_0-\eta_1<\Re\tau<\eta_0$.
Как мы показали в~\cite{Z8}, уравнение $f'(\tau)=0$
имеет в точности два корня, симметричных относительно
прямой $\tau=\frac12\eta_0$. Поскольку $\Re f(\tau)=\Re f(\eta_0-\ol\tau)$,
то корень уравнения $f'(\tau)=0$, отличный от $\tau=\tau_1$,
не влияет на асимптотику интеграла~\eqref{eq:1.41}.
\end{remark}

Подытожим результаты предложений~\ref{prop:1.1}--\ref{prop:1.3}.

\begin{proposition}
\label{prop:1.4}
В приведенных выше обозначениях пусть $r=3$ и
$$
\begin{gathered}
C_0=-\Re f(\tau_1),
\\
C_1=rm_1+m_2+\dots+m_{q-r}
-\biggl(\int_0^1\varphi(x)\,\d\psi(x)
-\int_0^{1/m_{q-r}}\varphi(x)\,\frac{\d x}{x^2}\biggr).
\end{gathered}
$$
Если $C_0>C_1$\rom, то по крайней мере одно из чисел
$$
\zeta(5), \; \zeta(7), \; \dots, \; \zeta(q-4) \;
\text{и\/} \; \zeta(q-2)
$$
иррационально.
\end{proposition}

Теперь все готово для доказательства теоремы~\ref{th:1}.

\begin{proof}[Доказательство теоремы~\rom{\ref{th:1}}]
Выбирая $r=3$, $q=13$,
$$
\eta_0=91, \qquad \eta_1=\eta_2=\eta_3=27,
\qquad \eta_j=25+j \quad\text{при}\; j=4,5,\dots,13,
$$
получаем $\tau_1=87.47900541\hdots+i\,3.32820690\dots$,
$$
\Im f(\tau_1)/\pi=-355.10280595\dots
$$
и
\begin{align*}
C_0&=-\Re f(\tau_1)=227.58019641\dots,
\\
C_1&=3\cdot35+34+8\cdot33
-\biggl(\int_0^1\varphi(x)\,\d\psi(x)
-\int_0^{1/33}\varphi(x)\,\frac{\d x}{x^2}\biggr)
\\
&=226.24944266\dots,
\end{align*}
поскольку в этом случае
$$
\varphi(x)=\nu \quad\text{при $x\in\Omega_\nu\setminus\Omega_{\nu+1}$},
\qquad \nu=0,1,\dots,9,
$$
для $x\in[0,1)$, где $\Omega_0=[0,1)$,
\begin{align*}
\Omega_1
&=\Omega_2
=\bigl[\tfrac2{91},\tfrac{36}{37}\bigr)
\cup\bigl[\tfrac{90}{91},1\bigr),
\displaybreak[0]\\
\Omega_3
&=\bigl[\tfrac2{91},\tfrac1{20}\bigr)
\cup\bigl[\tfrac5{91},\tfrac34\bigr)
\cup\bigl[\tfrac{28}{37},\tfrac{13}{14}\bigr)
\cup\bigl[\tfrac{14}{15},\tfrac{35}{37}\bigr)
\cup\bigl[\tfrac{18}{19},\tfrac{27}{28}\bigr)
\cup\bigl[\tfrac{88}{91},\tfrac{36}{37}\bigr)
\cup\bigl[\tfrac{90}{91},1\bigr),
\displaybreak[0]\\
\Omega_4
&=\bigl[\tfrac1{38},\tfrac1{22}\bigr)
\cup\bigl[\tfrac5{91},\tfrac3{26}\bigr)
\cup\bigl[\tfrac2{17},\tfrac18\bigr)
\cup\bigl[\tfrac4{31},\tfrac4{27}\bigr)
\cup\bigl[\tfrac5{33},\tfrac7{30}\bigr)
\cup\bigl[\tfrac4{17},\tfrac{12}{37}\bigr)
\cup\bigl[\tfrac{30}{91},\tfrac13\bigr)
\\ &\;\;
\cup\bigl[\tfrac{31}{91},\tfrac38\bigr)
\cup\bigl[\tfrac{14}{37},\tfrac{11}{28}\bigr)
\cup\bigl[\tfrac{13}{33},\tfrac9{22}\bigr)
\cup\bigl[\tfrac7{17},\tfrac{13}{28}\bigr)
\cup\bigl[\tfrac8{17},\tfrac12\bigr)
\cup\bigl[\tfrac{19}{37},\tfrac9{14}\bigr)
\cup\bigl[\tfrac{20}{31},\tfrac23\bigr)
\\ &\;\;
\cup\bigl[\tfrac{21}{31},\tfrac34\bigr)
\cup\bigl[\tfrac{25}{33},\tfrac{11}{14}\bigr)
\cup\bigl[\tfrac{26}{33},\tfrac{23}{28}\bigr)
\cup\bigl[\tfrac{14}{17},\tfrac{23}{27}\bigr)
\cup\bigl[\tfrac{31}{36},\tfrac{25}{27}\bigr)
\cup\bigl[\tfrac{85}{91},\tfrac{35}{37}\bigr)
\cup\bigl[\tfrac{20}{21},\tfrac{26}{27}\bigr)
\\ &\;\;
\cup\bigl[\tfrac{32}{33},\tfrac{34}{35}\bigr),
\displaybreak[0]\\
\Omega_5
&=\bigl[\tfrac1{37},\tfrac1{27}\bigr)
\cup\bigl[\tfrac1{25},\tfrac1{24}\bigr)
\cup\bigl[\tfrac5{91},\tfrac1{18}\bigr)
\cup\bigl[\tfrac2{35},\tfrac2{27}\bigr)
\cup\bigl[\tfrac3{38},\tfrac1{12}\bigr)
\cup\bigl[\tfrac8{91},\tfrac3{34}\bigr)
\cup\bigl[\tfrac2{21},\tfrac19\bigr)
\\ &\;\;
\cup\bigl[\tfrac4{33},\tfrac18\bigr)
\cup\bigl[\tfrac5{38},\tfrac4{27}\bigr)
\cup\bigl[\tfrac3{19},\tfrac16\bigr)
\cup\bigl[\tfrac5{29},\tfrac5{27}\bigr)
\cup\bigl[\tfrac4{21},\tfrac5{26}\bigr)
\cup\bigl[\tfrac6{29},\tfrac29\bigr)
\cup\bigl[\tfrac5{21},\tfrac7{27}\bigr)
\\ &\;\;
\cup\bigl[\tfrac4{15},\tfrac{10}{37}\bigr)
\cup\bigl[\tfrac27,\tfrac3{10}\bigr)
\cup\bigl[\tfrac7{23},\tfrac4{13}\bigr)
\cup\bigl[\tfrac6{19},\tfrac{12}{37}\bigr)
\cup\bigl[\tfrac{30}{91},\tfrac13\bigr)
\cup\bigl[\tfrac{10}{29},\tfrac7{20}\bigr)
\cup\bigl[\tfrac{13}{37},\tfrac5{14}\bigr)
\\ &\;\;
\cup\bigl[\tfrac{33}{91},\tfrac38\bigr)
\cup\bigl[\tfrac8{21},\tfrac5{13}\bigr)
\cup\bigl[\tfrac{13}{33},\tfrac{11}{27}\bigr)
\cup\bigl[\tfrac{12}{29},\tfrac5{12}\bigr)
\cup\bigl[\tfrac8{19},\tfrac{11}{26}\bigr)
\cup\bigl[\tfrac{14}{33},\tfrac{13}{30}\bigr)
\cup\bigl[\tfrac{40}{91},\tfrac49\bigr)
\\ &\;\;
\cup\bigl[\tfrac5{11},\tfrac{11}{24}\bigr)
\cup\bigl[\tfrac{17}{37},\tfrac6{13}\bigr)
\cup\bigl[\tfrac{17}{36},\tfrac{13}{27}\bigr)
\cup\bigl[\tfrac{16}{33},\tfrac12\bigr)
\cup\bigl[\tfrac{16}{31},\tfrac{14}{27}\bigr)
\cup\bigl[\tfrac8{15},\tfrac{19}{35}\bigr)
\cup\bigl[\tfrac{17}{31},\tfrac59\bigr)
\\ &\;\;
\cup\bigl[\tfrac{19}{33},\tfrac{15}{26}\bigr)
\cup\bigl[\tfrac{18}{31},\tfrac{16}{27}\bigr)
\cup\bigl[\tfrac{20}{33},\tfrac{17}{28}\bigr)
\cup\bigl[\tfrac{19}{31},\tfrac{17}{27}\bigr)
\cup\bigl[\tfrac{11}{17},\tfrac23\bigr)
\cup\bigl[\tfrac{17}{25},\tfrac{15}{22}\bigr)
\cup\bigl[\tfrac{20}{29},\tfrac{19}{27}\bigr)
\\ &\;\;
\cup\bigl[\tfrac{12}{17},\tfrac{17}{24}\bigr)
\cup\bigl[\tfrac{21}{29},\tfrac{20}{27}\bigr)
\cup\bigl[\tfrac{23}{31},\tfrac34\bigr)
\cup\bigl[\tfrac{69}{91},\tfrac79\bigr)
\cup\bigl[\tfrac{15}{19},\tfrac{19}{24}\bigr)
\cup\bigl[\tfrac45,\tfrac{22}{27}\bigr)
\cup\bigl[\tfrac{14}{17},\tfrac{23}{27}\bigr)
\\ &\;\;
\cup\bigl[\tfrac{25}{29},\tfrac{19}{22}\bigr)
\cup\bigl[\tfrac{27}{31},\tfrac78\bigr)
\cup\bigl[\tfrac{29}{33},\tfrac89\bigr)
\cup\bigl[\tfrac{26}{29},\tfrac9{10}\bigr)
\cup\bigl[\tfrac{28}{31},\tfrac{25}{27}\bigr)
\cup\bigl[\tfrac{31}{33},\tfrac{35}{37}\bigr)
\cup\bigl[\tfrac{87}{91},\tfrac{26}{27}\bigr)
\\ &\;\;
\cup\bigl[\tfrac{32}{33},\tfrac{33}{34}\bigr),
\displaybreak[0]\\
\Omega_6
&=\bigl[\tfrac1{36},\tfrac1{27}\bigr)
\cup\bigl[\tfrac1{17},\tfrac2{27}\bigr)
\cup\bigl[\tfrac9{91},\tfrac4{37}\bigr)
\cup\bigl[\tfrac{10}{91},\tfrac19\bigr)
\cup\bigl[\tfrac{12}{91},\tfrac4{27}\bigr)
\cup\bigl[\tfrac{16}{91},\tfrac5{27}\bigr)
\cup\bigl[\tfrac{19}{91},\tfrac8{37}\bigr)
\\ &\;\;
\cup\bigl[\tfrac5{23},\tfrac29\bigr)
\cup\bigl[\tfrac7{29},\tfrac9{37}\bigr)
\cup\bigl[\tfrac{23}{91},\tfrac7{27}\bigr)
\cup\bigl[\tfrac27,\tfrac8{27}\bigr)
\cup\bigl[\tfrac{29}{91},\tfrac{12}{37}\bigr)
\cup\bigl[\tfrac{30}{91},\tfrac13\bigr)
\cup\bigl[\tfrac{33}{91},\tfrac{10}{27}\bigr)
\\ &\;\;
\cup\bigl[\tfrac{15}{38},\tfrac{11}{27}\bigr)
\cup\bigl[\tfrac37,\tfrac{16}{37}\bigr)
\cup\bigl[\tfrac{40}{91},\tfrac49\bigr)
\cup\bigl[\tfrac9{19},\tfrac{13}{27}\bigr)
\cup\bigl[\tfrac{47}{91},\tfrac{14}{27}\bigr)
\cup\bigl[\tfrac7{13},\tfrac{20}{37}\bigr)
\cup\bigl[\tfrac{50}{91},\tfrac59\bigr)
\\ &\;\;
\cup\bigl[\tfrac{53}{91},\tfrac{16}{27}\bigr)
\cup\bigl[\tfrac8{13},\tfrac{23}{37}\bigr)
\cup\bigl[\tfrac{57}{91},\tfrac{17}{27}\bigr)
\cup\bigl[\tfrac{59}{91},\tfrac{24}{37}\bigr)
\cup\bigl[\tfrac{15}{23},\tfrac{17}{26}\bigr)
\cup\bigl[\tfrac{23}{35},\tfrac23\bigr)
\cup\bigl[\tfrac9{13},\tfrac{26}{37}\bigr)
\\ &\;\;
\cup\bigl[\tfrac{64}{91},\tfrac{19}{27}\bigr)
\cup\bigl[\tfrac{66}{91},\tfrac{19}{26}\bigr)
\cup\bigl[\tfrac{67}{91},\tfrac{20}{27}\bigr)
\cup\bigl[\tfrac{13}{17},\tfrac79\bigr)
\cup\bigl[\tfrac45,\tfrac{22}{27}\bigr)
\cup\bigl[\tfrac{76}{91},\tfrac{31}{37}\bigr)
\cup\bigl[\tfrac{16}{19},\tfrac{23}{27}\bigr)
\\ &\;\;
\cup\bigl[\tfrac{29}{33},\tfrac89\bigr)
\cup\bigl[\tfrac{31}{34},\tfrac{34}{37}\bigr)
\cup\bigl[\tfrac{23}{25},\tfrac{25}{27}\bigr)
\cup\bigl[\tfrac{31}{33},\tfrac{33}{35}\bigr)
\cup\bigl[\tfrac{87}{91},\tfrac{26}{27}\bigr),
\displaybreak[0]\\
\Omega_7
&=\bigl[\tfrac1{33},\tfrac1{27}\bigr)
\cup\bigl[\tfrac1{17},\tfrac2{27}\bigr)
\cup\bigl[\tfrac9{91},\tfrac4{37}\bigr)
\cup\bigl[\tfrac{10}{91},\tfrac19\bigr)
\cup\bigl[\tfrac{12}{91},\tfrac5{37}\bigr)
\cup\bigl[\tfrac17,\tfrac4{27}\bigr)
\cup\bigl[\tfrac{16}{91},\tfrac5{27}\bigr)
\\ &\;\;
\cup\bigl[\tfrac{19}{91},\tfrac8{37}\bigr)
\cup\bigl[\tfrac{20}{91},\tfrac29\bigr)
\cup\bigl[\tfrac{22}{91},\tfrac9{37}\bigr)
\cup\bigl[\tfrac9{35},\tfrac7{27}\bigr)
\cup\bigl[\tfrac27,\tfrac8{27}\bigr)
\cup\bigl[\tfrac{29}{91},\tfrac9{28}\bigr)
\cup\bigl[\tfrac{10}{31},\tfrac{11}{34}\bigr)
\\ &\;\;
\cup\bigl[\tfrac{33}{91},\tfrac{10}{27}\bigr)
\cup\bigl[\tfrac{36}{91},\tfrac{15}{37}\bigr)
\cup\bigl[\tfrac{37}{91},\tfrac{11}{27}\bigr)
\cup\bigl[\tfrac37,\tfrac{16}{37}\bigr)
\cup\bigl[\tfrac{40}{91},\tfrac49\bigr)
\cup\bigl[\tfrac{10}{21},\tfrac{13}{27}\bigr)
\cup\bigl[\tfrac{47}{91},\tfrac{14}{27}\bigr)
\\ &\;\;
\cup\bigl[\tfrac7{13},\tfrac{20}{37}\bigr)
\cup\bigl[\tfrac{50}{91},\tfrac59\bigr)
\cup\bigl[\tfrac{53}{91},\tfrac{16}{27}\bigr)
\cup\bigl[\tfrac8{13},\tfrac{23}{37}\bigr)
\cup\bigl[\tfrac{57}{91},\tfrac{17}{27}\bigr)
\cup\bigl[\tfrac{59}{91},\tfrac{24}{37}\bigr)
\cup\bigl[\tfrac9{13},\tfrac{26}{37}\bigr)
\\ &\;\;
\cup\bigl[\tfrac{64}{91},\tfrac{19}{27}\bigr)
\cup\bigl[\tfrac{66}{91},\tfrac{27}{37}\bigr)
\cup\bigl[\tfrac{67}{91},\tfrac{20}{27}\bigr)
\cup\bigl[\tfrac{10}{13},\tfrac79\bigr)
\cup\bigl[\tfrac{73}{91},\tfrac{30}{37}\bigr)
\cup\bigl[\tfrac{74}{91},\tfrac{22}{27}\bigr)
\cup\bigl[\tfrac{11}{13},\tfrac{23}{27}\bigr)
\\ &\;\;
\cup\bigl[\tfrac{80}{91},\tfrac89\bigr)
\cup\bigl[\tfrac{83}{91},\tfrac{34}{37}\bigr)
\cup\bigl[\tfrac{12}{13},\tfrac{25}{27}\bigr)
\cup\bigl[\tfrac{87}{91},\tfrac{26}{27}\bigr),
\displaybreak[0]\\
\Omega_8
&=\bigl[\tfrac1{31},\tfrac1{27}\bigr)
\cup\bigl[\tfrac6{91},\tfrac2{27}\bigr)
\cup\bigl[\tfrac9{91},\tfrac1{10}\bigr)
\cup\bigl[\tfrac3{29},\tfrac4{37}\bigr)
\cup\bigl[\tfrac{10}{91},\tfrac19\bigr)
\cup\bigl[\tfrac2{15},\tfrac5{37}\bigr)
\cup\bigl[\tfrac17,\tfrac4{27}\bigr)
\\ &\;\;
\cup\bigl[\tfrac3{17},\tfrac5{28}\bigr)
\cup\bigl[\tfrac7{38},\tfrac5{27}\bigr)
\cup\bigl[\tfrac7{33},\tfrac8{37}\bigr)
\cup\bigl[\tfrac{20}{91},\tfrac29\bigr)
\cup\bigl[\tfrac8{33},\tfrac9{37}\bigr)
\cup\bigl[\tfrac9{31},\tfrac7{24}\bigr)
\cup\bigl[\tfrac5{17},\tfrac8{27}\bigr)
\\ &\;\;
\cup\bigl[\tfrac4{11},\tfrac{10}{27}\bigr)
\cup\bigl[\tfrac{37}{91},\tfrac{11}{27}\bigr)
\cup\bigl[\tfrac{11}{23},\tfrac{13}{27}\bigr)
\cup\bigl[\tfrac7{13},\tfrac{20}{37}\bigr)
\cup\bigl[\tfrac{16}{29},\tfrac59\bigr)
\cup\bigl[\tfrac{53}{91},\tfrac7{12}\bigr)
\cup\bigl[\tfrac{17}{29},\tfrac{16}{27}\bigr)
\\ &\;\;
\cup\bigl[\tfrac{13}{21},\tfrac{23}{37}\bigr)
\cup\bigl[\tfrac{23}{33},\tfrac7{10}\bigr)
\cup\bigl[\tfrac{64}{91},\tfrac{19}{27}\bigr)
\cup\bigl[\tfrac{14}{19},\tfrac{20}{27}\bigr)
\cup\bigl[\tfrac{10}{13},\tfrac{27}{35}\bigr)
\cup\bigl[\tfrac{25}{31},\tfrac{30}{37}\bigr)
\cup\bigl[\tfrac{74}{91},\tfrac{22}{27}\bigr)
\\ &\;\;
\cup\bigl[\tfrac{11}{13},\tfrac{23}{27}\bigr)
\cup\bigl[\tfrac{80}{91},\tfrac{31}{35}\bigr)
\cup\bigl[\tfrac{83}{91},\tfrac{11}{12}\bigr)
\cup\bigl[\tfrac{12}{13},\tfrac{25}{27}\bigr)
\cup\bigl[\tfrac{22}{23},\tfrac{26}{27}\bigr),
\displaybreak[0]\\
\Omega_9
&=\bigl[\tfrac1{29},\tfrac1{28}\bigr)
\cup\bigl[\tfrac2{29},\tfrac1{14}\bigr)
\cup\bigl[\tfrac7{19},\tfrac{10}{27}\bigr)
\cup\bigl[\tfrac{12}{25},\tfrac{13}{27}\bigr)
\cup\bigl[\tfrac{17}{23},\tfrac{20}{27}\bigr)
\cup\bigl[\tfrac{15}{17},\tfrac{23}{26}\bigr)
\cup\bigl[\tfrac{24}{25},\tfrac{25}{26}\bigr),
\end{align*}
и $\Omega_{10}=\emptyset$.

Поэтому применение предложения~\ref{prop:1.4} завершает доказательство.
\end{proof}

Перебирая различные целочисленные наборы направлений
$\Beta=(\eta_0;\linebreak[2]\eta_1,\dots,\eta_q)$ с $q=11$,
удовлетворяющие условиям
$$
\eta_1\le\eta_2\le\dots\le\eta_q<\frac12\eta_0
\qquad\text{и}\qquad
\eta_0\le120,
$$
мы не обнаружили ни одного набора~$\Beta$, доказывающего
с помощью предложения~\ref{prop:1.4} иррациональность
хотя бы одного из трех чисел $\zeta(5)$, $\zeta(7)$ и~$\zeta(9)$.
Поэтому можно говорить о ``естественных'' границах,
достигнутых в результате теоремы~\ref{th:1} за счет применения
гипергеометрической конструкции и арифметического метода
в максимально общей форме.

\chapter[Интегральные конструкции линейных форм]%
{Интегральные конструкции линейных форм\\ от дзета-значений}
\label{chap:2}

В настоящей главе мы доказываем теорему~\ref{th:2} как частный случай
некоторого аналитического тождества, связывающего кратные
интегралы эйлерова типа и совершенно уравновешенные
гипергеометрические ряды. Отметим, что в настоящее время
известны алгоритмические методы доказательства соотношений,
подобных~\eqref{eq:0.24}. Они являются неотъемлемой частью
так называемой {\it WZ-теории}, разработанной Г.~Уилфом
и Д.~Цайльбергером \cite{WZ}, \cite{PWZ}, и заключаются
в поиске и доказательстве рекуррентных уравнений для
правой и левой части рассматриваемых тождеств. Так, например,
можно показать, что тройной интеграл~$J_{3,n}$
и соответствующий ряд~$F_{3,n}$ удовлетворяют одному
и тому же рекуррентному соотношению~\eqref{eq:0.3}.
Однако, поиск разностных уравнений для интегралов~$J_{4,n}$
и~$J_{5,n}$ с помощью разработанных алгоритмов
все еще остается за пределами возможностей
современных компьютеров. Кроме того, подобные
доказательства необходимо указывать для каждого~$s$.

В качестве исторического замечания отметим, что
совпадение интеграла $J_{3,n}$ и ряда~$F_{3,n}$
было впервые доказано в~\cite{Z17} с помощью
преобразования Бэйли \cite[\S\,6.3, формула~(2)]{Bai}
и интегральной теоремы Нестеренко \cite[теорема~2]{Ne4}.
Аналогичные аргументы (с заменой преобразования Бэйли
на преобразование Уиппла~\eqref{eq:0.30}) дают классическое
доказательство тождества $J_{2,n}=F_{2,n}$.
Эти наблюдения и сравнение асимптотик интегралов
$J_{s,n}$ и рядов~$F_{s,n}$ при $n\to\infty$ (в случае $s=4,5$)
позволили высказать в~\cite{Z17} равенство~\eqref{eq:0.24}
в качестве гипотезы. Подробная версия приводимого далее
доказательства опубликована в~\cite{Z15} (см\. также~\cite{Z11}).

\section{Основной результат и следствия из него}
\label{sec:2.1}

Основные объекты этой главы~---
общий совершенно уравновешенный гипергеометрический ряд
\begin{align}
F_s(\bh)
&=F_s(h_0;h_1,\dots,h_s)
=\frac{\Gamma(1+h_0)\,\prod_{j=1}^s\Gamma(h_j)}
{\prod_{j=1}^s\Gamma(1+h_0-h_j)}
\nonumber\\ &\qquad\times
{}_{s+2}F_{s+1}\biggl(\begin{matrix}
h_0, & 1+\tfrac12h_0, & h_1, & \dots, & h_s \\[1pt]
& \tfrac12h_0, & 1+h_0-h_1, & \dots, & 1+h_0-h_s
\end{matrix}\biggm|(-1)^{s+1}\biggr)
\nonumber\\ &\phantom:
=\sum_{\mu=0}^\infty(h_0+2\mu)
\frac{\prod_{j=0}^s\Gamma(h_j+\mu)}{\prod_{j=0}^s\Gamma(1+h_0-h_j+\mu)}
(-1)^{(s+1)\mu}
\label{eq:2.1}
\end{align}
и $s$-кратный интеграл
\begin{align}
J_s(\ba,\bb)
&=J_s\biggl(\begin{matrix}
a_0, & a_1, & \dots, & a_s \\[1pt]
& b_1, & \dots, & b_s
\end{matrix}\biggr)
\nonumber\\
&=\idotsint\limits_{[0,1]^s}
\frac{\prod_{j=1}^sx_j^{a_j-1}(1-x_j)^{b_j-a_j-1}}
{Q_s(x_1,x_2,\dots,x_s)^{a_0}}
\,\d x_1\,\d x_2\dotsb\d x_s,
\label{eq:2.2}
\end{align}
где $Q_0=1$ и
\begin{align}
Q_s
&=Q_s(x_1,x_2,\dots,x_s)
=1-x_1(1-x_2(\dotsb(1-x_{s-1}(1-x_s))\dotsb))
\nonumber\\
&=1-x_1Q_{s-1}(x_2,\dots,x_s)
=Q_{s-1}(x_1,\dots,x_{s-1})+(-1)^sx_1x_2\dotsb x_s
\label{eq:2.3}
\end{align}
(см.~\eqref{eq:0.21}).

Область сходимости обобщенного гипергеометрического
ряда~\eqref{eq:0.27} указана в~\eqref{eq:0.28}
(см\. также условие~\eqref{eq:1.23} при $r=1$).
Приведем некоторые достаточные условия для
абсолютной сходимости интеграла~\eqref{eq:2.2}.

\begin{lemma} \label{lem:2.1}
Если
\begin{equation}
\begin{gathered}
\Re b_j>\Re a_j>0, \qquad j=1,\dots,s,
\\
\Re b_1>\Re(a_1+a_0) \ \text{при $s=1$}
\quad\text{или}\quad
\Re b_1\ge\Re(a_1+a_0) \ \text{при $s>1$},
\end{gathered}
\label{eq:2.4}
\end{equation}
то интеграл в~\eqref{eq:2.2} сходится абсолютно.
\end{lemma}

\begin{proof}
При $s=1$ интеграл $J_1$ в~\eqref{eq:2.2}
является бета-инте\-гралом, так что условия~\eqref{eq:2.4}
гарантируют его абсолютную сходимость. Если же $s\ge2$,
то ввиду
\begin{equation*}
Q_s(x_1,\dots,x_s)
=1-x_1(1-x_2(1-\dotsb(1-x_{s-1}(1-x_s))\dotsb))
\ge1-x_1
\end{equation*}
выполнено
\begin{align}
&
\int_{\varepsilon_s}^{1-\varepsilon_s}\dotsi
\int_{\varepsilon_1}^{1-\varepsilon_1}
\biggl|\frac{\prod_{j=1}^sx_j^{a_j-1}(1-x_j)^{b_j-a_j-1}}
{Q_s(x_1,x_2,\dots,x_s)^{a_0}}\biggr|\,\d x_1\dotsb\d x_s
\nonumber\\ &\qquad
\le\int_{\varepsilon_s}^{1-\varepsilon_s}\dotsi
\int_{\varepsilon_1}^{1-\varepsilon_1}
\frac{\prod_{j=2}^sx_j^{\Re(a_j-1)}(1-x_j)^{\Re(b_j-a_j-1)}}
{1-x_1+x_1x_2X}\,\d x_1\dotsb\d x_s
\nonumber\\ &\qquad
\le\int_{\varepsilon_s}^{1-\varepsilon_s}\dotsi
\int_{\varepsilon_2}^{1-\varepsilon_2}
\prod_{j=2}^sx_j^{\Re(a_j-1)}(1-x_j)^{\Re(b_j-a_j-1)}
\nonumber\\ &\qquad\quad\times
\biggl(-\frac{\log(1-x_1+x_1x_2X)}{1-x_2X}\biggr)
\bigg|_{x_1=\varepsilon_1}^{1-\varepsilon_1}
\,\d x_2\dotsb\d x_s
\label{eq:2.5}
\end{align}
для любых достаточно малых $\varepsilon_1,\dots,\varepsilon_s>0$,
где мы ввели обозначение
$X=1-x_3(1-\dotsb(1-x_{s-1}(1-x_s))\dotsb)=Q_{s-2}(x_3,\dots,x_s)$.
Предельный переход $\varepsilon_1\to0$ в правой части~\eqref{eq:2.5}
приводит к интегралу
\begin{equation}
\int_{\varepsilon_s}^{1-\varepsilon_s}\dotsi
\int_{\varepsilon_2}^{1-\varepsilon_2}
\prod_{j=2}^sx_j^{\Re(a_j-1)}(1-x_j)^{\Re(b_j-a_j-1)}
\biggl(-\frac{\log(x_2X)}{1-x_2X}\biggr)\,\d x_2\dotsb\d x_s,
\label{eq:2.6}
\end{equation}
который благополучно сходится при $y=x_2X\to1$, так как
функция $(\log y)/(1-y)$ продолжается по непрерывности
в точке $y=1$. С другой стороны, при $y\to0$ справедливы
оценки
\begin{equation*}
|\log(x_2X)|<(x_2X)^{-\delta}\le x_2^{-\delta}(1-x_3)^{-\delta},
\end{equation*}
где мы выбираем $\delta=\frac12\min\{\Re a_2,\Re(b_3-a_3)\}$,
так что интеграл~\eqref{eq:2.6} действительно сходится,
а это обеспечивает абсолютную сходимость интеграла \eqref{eq:2.2}
в условиях~\eqref{eq:2.4} согласно оценке~\eqref{eq:2.5}.
\end{proof}

\begin{theorem}
\label{th:6}
Для любого $s\ge1$ имеет место тождество
\begin{align}
&
\frac{\prod_{j=1}^{s+1}\Gamma(1+h_0-h_j-h_{j+1})}
{\Gamma(h_1)\,\Gamma(h_{s+2})}
\cdot F_{s+2}(h_0;h_1,\dots,h_{s+2})
\nonumber\\ &\qquad
=J_s\biggl(\begin{matrix}
h_1, & h_2, & h_3, & \dots, & h_{s+1} \\[1pt]
& 1+h_0-h_3, & 1+h_0-h_4, & \dots, & 1+h_0-h_{s+2}
\end{matrix}\biggr)
\label{eq:2.7}
\end{align}
при условиях
\begin{gather}
1+\Re h_0>\frac2{s+1}\cdot\sum_{j=1}^{s+2}\Re h_j,
\label{eq:2.8}
\\
\Re(1+h_0-h_{j+1})>\Re h_j>0, \quad j=2,\dots,s+1,
\label{eq:2.9}
\\
\begin{gathered}
\Re(1+h_0-h_3)\ge\Re(h_1+h_2) \quad \text{при $s\ge2$},
\end{gathered}
\label{eq:2.10}
\\ \vspace{1.5pt}
h_1,h_{s+2}\ne0,-1,-2,\dots,
\label{eq:2.11}
\end{gather}
обеспечивающих абсолютную сходимость интеграла и ряда в~\eqref{eq:2.7}.
\end{theorem}

\begin{remark}
Пользуясь теоремой об аналитическом продолжении,
можно опустить условие~\eqref{eq:2.10}. Кроме того,
интерпретируя $\Gamma(h_j+\mu)/\Gamma(h_j)$ в случае $j=1,s+2$
как символ Похгаммера $(h_j)_\mu$, можно также отказаться
от ограничения~\eqref{eq:2.11} при соответствующем
суммировании в~\eqref{eq:2.1}.
\end{remark}

\begin{remark}
Как отмечает Дж.~Эндрюс в~\cite[\S\,16]{An2},
``целый обзор мог бы быть написан только для
интегралов, связанных с уравновешенными рядами''.
Теорема~\ref{th:6} немного расширила бы подобный обзор.
В~\cite{Z18} мы доказываем другие теоремы, связывающие
уравновешенные гипергеометрические объекты и эйлеровы
интегралы и связанные с общей конструкцией рациональных
приближений к дзета-значениям из гл.~\ref{chap:1}.

Другое семейство кратных интегралов
\begin{gather}
S(z)
=\idotsint\limits_{[0,1]^s}
\frac{\prod_{j=1}^sx_j^{a_j-1}(1-x_j)^{b_j-a_j-1}}
{\prod_{i=1}^m(1-zx_1x_2\dotsb x_{r_i})^{c_i}}
\,\d x_1\,\d x_2\dotsb\d x_s,
\label{eq:2.12}
\\
1\le r_1<r_2<\dots<r_m=s,
\nonumber
\end{gather}
применялось В.~Сорокиным \cite{So2}, \cite{So3}
для решения арифметических задач, связанных со значениями
полилогарифмов. Недавно С.~Злобин \cite{Zl1}, \cite{Zl2}
доказал (в наиболее общих предположениях),
что заменой переменных интегралы~\eqref{eq:2.2}
могут быть приведены к форме~\eqref{eq:2.12} с $z=1$.
Это означает, что теорема~\ref{th:6} позволяет получить
представление интегралов $S(1)$ в виде совершенно уравновешенных
гипергеометрических рядов~\eqref{eq:2.1} при некоторых
условиях на параметры $a_j$, $b_j$, $c_i$ и $r_i$ в~\eqref{eq:2.12}.
Кроме того, Злобин~\cite{Zl1} показал, что
при целых значениях параметров в~\eqref{eq:2.12}
интеграл $S(z)$ является $\mathbb Q[z^{-1}]$-линейной
комбинацией обобщенных полилогарифмов
$$
\sum_{n_1\ge n_2\ge\dots\ge n_l\ge1}
\frac{z^{n_1}}{n_1^{s_1}n_2^{s_2}\dotsb n_l^{s_l}}
\qquad\text{с}\quad
s_j\ge1, \; s_j\in\mathbb Z, \; j=1,\dots,l,
$$
где $0\le s_1+s_2+\dots+s_l\le s$ и $0\le l\le m$.
\end{remark}

В случае целочисленных параметров $\bh$
совершенно уравновешенный гипергеометрический ряд
\eqref{eq:2.1} является $\mathbb Q$-линейной формой
от четных или нечетных дзета-значений в зависимости от четности
$s\ge4$ (ср\. с предложением~\ref{prop:1.1} и \cite[\S\,9]{Z17}).
Поэтому если целые положительные параметры
$\ba$ и $\bb$ удовлетворяют дополнительному условию
\begin{equation}
b_1+a_2=b_2+a_3=\dots=b_{s-1}+a_s,
\label{eq:2.13}
\end{equation}
то интеграл \eqref{eq:2.2} есть $\mathbb Q$-линейная
форма от дзета-значений одинаковой четности.
Специализация $a_j=n+1$, $b_j=2n+2$ приводит к
тождеству~\eqref{eq:0.24} теоремы~\ref{th:1}, и
включения~\eqref{eq:0.26} следуют из рассуждений
в начале гл.~\ref{chap:1}.
(На самом деле, согласно \cite[леммы~4.2--4.4]{Z8}
мы получаем более точные включения: так, например, для $s\ge3$
нечетного выполнено
$$
2D_n^{s+1}\Phi_n^{-1}\cdot J_{s,n}
\in\mathbb Z\zeta(s)+\mathbb Z\zeta(s-2)+\dots+\mathbb Z\zeta(3)+\mathbb Z,
$$
где $\Phi_n$~--- произведение простых чисел~$p<n$, для которых
$\frac23\le\{n/p\}<1$.)
Выбор $a_j=rn+1$, $b_j=(r+1)n+2$ в~\eqref{eq:2.2}
(или, эквивалентно, $h_0=(2r+1)n+2$ и $h_j=rn+1$ для
$j=1,\dots,s+2$ в~\eqref{eq:2.1})
с целым $r\ge\nobreak1$, зависящим от заданного нечетного~$s$,
приводит к линейным формам Ривоаля~\eqref{eq:0.19}
от нечетных дзета-значений, использованных в~\cite{Ri1}
для доказательства результата о бесконечности множества
иррациональных чисел среди
$\zeta(3),\zeta(5),\zeta(7),\dots$\,.

Кроме того, следует отметить очевидную инвариантность
в предположении~\eqref{eq:2.13} величины
\begin{equation*}
\begin{split}
\frac{F_{s+2}(h_0;h_1,\dots,h_{s+2})}
{\prod_{j=1}^{s+2}\Gamma(h_j)}
&=\frac{J_s(\ba,\bb)}
{\prod_{j=2}^{s+1}\Gamma(h_j)\cdot\prod_{j=1}^{s+1}\Gamma(1+h_0-h_j-h_{j+1})}
\\
&=\frac{J_s(\ba,\bb)}
{\prod_{j=1}^s\Gamma(a_j)\cdot\Gamma(b_1+a_2-a_0-a_1)
\cdot\prod_{j=1}^s\Gamma(b_j-a_j)}
\end{split}
\end{equation*}
под действием ($\bh$-тривиальной) группы~$\fG_s$
порядка $(s+2)!$, состоящей из всех перестановок
параметров $h_1,\dots,h_{s+2}$. Этот результат также
имеет теоретико-числовые приложения.
В случаях $s=2$ и $s=3$ замена переменных
$(x_{s-1},x_s)\mapsto(1-x_s,1-x_{s-1})$ в~\eqref{eq:2.2}
реализует дополнительное преобразование~$\fc$ как интеграла \eqref{eq:2.2}
так и ряда~\eqref{eq:2.1}; при $s\ge4$ это преобразование недоступно,
поскольку нарушается условие~\eqref{eq:2.13}.
Группы $\<\fG_s,\fc\>$ порядков~$120$ и~$1920$
при $s=2$ и $s=3$ соответственно хорошо известны
(см\. \cite[\S\,3.6 и~\S\,7.5]{Bai});
Дж.~Рин и К.~Виола \cite{RV1}--\cite{RV3}, \cite{Vi}
применяют их, чтобы получить хорошие оценки
мер иррациональности чисел~$\zeta(2)$ и~$\zeta(3)$.
Группа $\<\fG_2,\fc\>$ (в терминах преобразований
$q$-гипергеометрического ряда) применяется нами далее
в гл.~\ref{chap:3} для доказательства теоремы~\ref{th:3}.
В случае $s\ge4$ группа~$\fG_s$ допускает естественную
интерпретацию как группа перестановок параметров
$$
e_{0l}=h_l-1, \quad 1\le l\le s+2,
\qquad e_{jl}=h_0-h_j-h_l, \quad 1\le j<l\le s+2;
$$
детали приведены в \cite[\S\,9]{Z17}.

\section{Гипергеометрическое доказательство}
\label{sec:2.2}

\begin{lemma} \label{lem:2.2}
Теорема~\rom{\ref{th:6}} верна при $s=1$.
\end{lemma}

\begin{proof}
Согласно предельному случаю теоремы Дугалла \eqref{eq:0.29} имеем
\begin{equation}
F_3(h_0;h_1,h_2,h_3)
=\frac{\Gamma(h_1)\,\Gamma(h_2)\,\Gamma(h_3)\,\Gamma(1+h_0-h_1-h_2-h_3)}
{\begin{aligned}
\Gamma(1+h_0-h_1-h_2)\,\Gamma(1+h_0-h_1-h_3)\,
\quad\\[-2mm] \times
\Gamma(1+h_0-h_2-h_3)
\end{aligned}}
\label{eq:2.14}
\end{equation}
при условии, что $1+\Re h_0>\Re(h_1+h_2+h_3)$
и $h_j$~не является неположительным
целым при $j=1,2,3$. С другой стороны,
интеграл в правой части~\eqref{eq:2.7} является
бета-интегралом, поэтому
\begin{equation*}
\begin{split}
J_1\biggl(\begin{matrix}
h_1, & h_2 \\[1pt] & 1+h_0-h_3
\end{matrix}\biggr)
&=\int_0^1\frac{x^{h_2-1}(1-x)^{h_0-h_2-h_3}}{(1-x)^{h_1}}\,\d x
\\
&=\frac{\Gamma(h_2)\,\Gamma(1+h_0-h_1-h_2-h_3)}{\Gamma(1+h_0-h_1-h_3)}
\end{split}
\end{equation*}
при условиях $1+\Re h_0>\Re(h_1+h_2+h_3)$ и $\Re h_2>0$.
Следовательно, умножение равенства~\eqref{eq:2.14} на необходимое
произведение гамма-множителей приводит к тождеству~\eqref{eq:2.7}
в случае $s=1$.
\end{proof}

\begin{remark}
Если определить $J_0(a_0)$ как~$1$,
то утверждение теоремы~\ref{th:6} остается справедливым
и при $s=0$, благодаря другому следствию из теоремы Дугалла
(см\. \cite[\S\,4.4, формула~(3)]{Bai}).
\end{remark}

\begin{lemma}[{\rm\cite[\S\,3.2]{Ne4}}]
\label{lem:2.3}
Пусть числа $a_0,a,b\in\mathbb C$ и $t_0\in\mathbb R$
удовлетворяют условиям
$$
\Re a_0>t_0>0, \quad \Re a>t_0>0
\quad\text{и}\quad \Re b>\Re a_0+\Re a.
$$
Тогда для любого ненулевого $z\in\mathbb C\setminus(1,+\infty)$
выполняется тождество
\begin{align}
&
\int_0^1\frac{x^{a-1}(1-x)^{b-a-1}}{(1-zx)^{a_0}}\,\d x
\nonumber\\ &\qquad
=\frac{\Gamma(b-a)}{\Gamma(a_0)}
\cdot\frac1{2\pi i}\int_{-t_0-i\infty}^{-t_0+i\infty}
\frac{\Gamma(a_0+t)\,\Gamma(a+t)\,\Gamma(-t)}{\Gamma(b+t)}
\,(-z)^t\,\d t,
\label{eq:2.15}
\end{align}
где $(-z)^t=|z|^te^{it\arg(-z)}$, $-\pi<\arg(-z)<\pi$
при $z\in\mathbb C\setminus[0,+\infty)$ и $\arg(-z)=\pm\pi$
при $z\in(0,1]$. При этом интеграл в правой части~\eqref{eq:2.15}
сходится абсолютно, и в случае $|z|\le1$ оба интеграла
в~\eqref{eq:2.15} отвечают абсолютно сходящемуся
гауссову гипергеометрическому ряду
$$
\frac{\Gamma(a)\,\Gamma(b-a)}{\Gamma(b)}
\cdot{}_2F_1\biggl(\begin{matrix}
a_0, & a \\[1pt] & b \end{matrix}\biggm|z\biggr)
=\frac{\Gamma(b-a)}{\Gamma(a_0)}\sum_{\nu=0}^\infty
\frac{\Gamma(a_0+\nu)\,\Gamma(a+\nu)}{\nu!\,\Gamma(b+\nu)}z^\nu.
$$
\end{lemma}

Определим $\epsilon_s=0$ при $s$~четном
и $\epsilon_s=1$ или~$-1$ при $s$~нечетном.

\begin{lemma} \label{lem:2.4}
Для каждого $s\ge2$ выполнено соотношение
\begin{equation*}
\begin{split}
&
J_s\biggl(\begin{matrix}
a_0, & a_1, & \dots, & a_{s-1}, & a_s \\[1pt]
& b_1, & \dots, & b_{s-1}, & b_s
\end{matrix}\biggr)
\\ &\quad
=\frac{\Gamma(b_s-a_s)}{\Gamma(a_0)}
\cdot\frac1{2\pi i}\int_{-t_0-i\infty}^{-t_0+i\infty}
\frac{\Gamma(a_0+t)\,\Gamma(a_s+t)\,\Gamma(-t)}{\Gamma(b_s+t)}
\\ &\quad\qquad\times
e^{\epsilon_s\pi it}\cdot J_{s-1}\biggl(\begin{matrix}
a_0+t, & a_1+t, & \dots, & a_{s-1}+t \\[1pt]
& b_1+t, & \dots, & b_{s-1}+t
\end{matrix}\biggr)\,\d t
\end{split}
\end{equation*}
при условии, что $\Re a_0>t_0>0$, $\Re a_s>t_0>0$,
$\Re b_s>\Re a_0+\Re a_s$ и интеграл в левой части сходится
абсолютно.
\end{lemma}

\begin{proof}
Прежде всего отметим, что первая рекурсия в~\eqref{eq:2.3}
и индуктивные рассуждения показывают справедливость неравенства
\begin{equation}
0<Q_s(x_1,x_2,\dots,x_s)<1 \qquad\text{для}\;
(x_1,x_2,\dots,x_s)\in(0,1)^s \;\text{и}\; s\ge1.
\label{eq:2.16}
\end{equation}
В соответствии со второй рекурсией в~\eqref{eq:2.3} имеем
$Q_s=Q_{s-1}\cdot(1-zx_s)$ для $s\ge2$, где
$$
z=\frac{(-1)^{s+1}x_1\dotsb x_{s-1}}{Q_{s-1}(x_1,\dots,x_{s-1})}.
$$
Для каждого набора $(x_1,\dots,x_{s-1})\in(0,1)^{s-1}$
число~$z$ является вещественным, причем $z<0$ для $s$~четных
и $0<z<1$ для $s$~нечетных, поскольку в последнем случае
$$
z=\frac{x_1\dotsb x_{s-1}}{Q_{s-1}(x_1,\dots,x_{s-2},x_{s-1})}
=\frac{x_1\dotsb x_{s-1}}{Q_{s-2}(x_1,\dots,x_{s-2})+x_1\dotsb x_{s-1}}
<1
$$
согласно~\eqref{eq:2.16}. Следовательно,
расщепляя интегрирование в~\eqref{eq:2.2}
по $(s-\nobreak1)$-мерному кубу
$(0,1)^{s-1}$ и интервалу $(0,1)$ и применяя лемму~\ref{lem:2.3}
к интегралу
$$
\int_0^1\frac{x_s^{a_s-1}(1-x_s)^{b_s-a_s-1}}{(1-zx_s)^{a_0}}\,\d x_s,
$$
мы получаем необходимое соотношение.
\end{proof}

\begin{proof}[Доказательство теоремы~\rom{\ref{th:6}}]
Случай $s=1$ рассмотрен в лемме~\ref{lem:2.2}.
Поэтому считаем далее, что $s\ge2$, тождество~\eqref{eq:2.7}
выполнено при~$s-1$ и, дополнительно,
\begin{equation}
1+\Re h_0>\frac 2s\cdot\sum_{j=1}^{s+1}\Re h_j,
\qquad \Re h_{s+2}<1.
\label{eq:2.17}
\end{equation}
Ограничения~\eqref{eq:2.17} без труда снимаются
в окончательном результате с помощью теоремы об аналитическом
продолжении.

В соответствии с индуктивным предположением
для $t\in\mathbb C$ с $\Re t<0$ выполнено
\begin{align}
&
J_{s-1}\biggl(\begin{matrix}
h_1+t, & h_2+t, & h_3+t, & \dots, & h_s+t \\[1pt]
& 1+h_0-h_3+t, & 1+h_0-h_4+t, & \dots, & 1+h_0-h_{s+1}+t
\end{matrix}\biggr)
\nonumber\\ &\qquad
=\frac{\prod_{j=1}^s\Gamma(1+h_0-h_j-h_{j+1})}
{\Gamma(h_1+t)\,\Gamma(h_{s+1}+t)}
\cdot F_{s+1}(h_0+2t;h_1+t,\dots,h_{s+1}+t)
\nonumber\\ &\qquad
=\frac{\prod_{j=1}^s\Gamma(1+h_0-h_j-h_{j+1})}
{\Gamma(h_1+t)\,\Gamma(h_{s+1}+t)}
\cdot\frac1{2\pi i}\int_{-q_0-i\infty}^{-q_0+i\infty}
(h_0+2t+2q)
\nonumber\\ &\qquad\quad\times
\frac{\Gamma(h_0+2t+q)\,\prod_{j=1}^{s+1}\Gamma(h_j+t+q)\,\Gamma(-q)}
{\prod_{j=1}^{s+1}\Gamma(1+h_0-h_j+t+q)}
\,e^{\epsilon_{s-1}\pi iq}\,\d q,
\label{eq:2.18}
\end{align}
где вещественное число $q_0>0$ удовлетворяет условиям
$$
\begin{gathered}
\Re(h_0+2t)>q_0, \quad \Re(1+\tfrac12h_0+t)>q_0,
\\
\Re(h_j+t)>q_0 \quad\text{при}\; j=1,\dots,s+1,
\end{gathered}
$$
а абсолютная сходимость комплексного интеграла барнсовского типа
следует из \cite[лемма~3]{Ne4}.
Сдвигая переменную интегрирования $t+q\mapsto q$ в~\eqref{eq:2.18}
(и используя при этом равенство
$e^{\epsilon_s\pi it}\cdot e^{\epsilon_{s-1}\pi iq}
=e^{\epsilon_{s-1}\pi i(t+q)}\cdot e^{\epsilon_1\pi it}$),
применяя лемму~\ref{lem:2.4} и меняя местами интегрирования
в полученном двойном комплексном интеграле
(законность гарантируется абсолютной сходимостью интегралов),
мы получаем соотношение
\begin{align}
&
J_s\biggl(\begin{matrix}
h_1, & h_2, & h_3, & \dots, & h_s, & h_{s+1} \\[1pt]
& 1+h_0-h_3, & 1+h_0-h_4, & \dots, & 1+h_0-h_{s+1}, & 1+h_0-h_{s+2}
\end{matrix}\biggr)
\nonumber\\ &\qquad
=\frac{\prod_{j=1}^{s+1}\Gamma(1+h_0-h_j-h_{j+1})}{\Gamma(h_1)}
\nonumber\\ &\qquad\quad\times
\frac1{2\pi i}\int_{-q_1-i\infty}^{-q_1+i\infty}
(h_0+2q)\frac{\prod_{j=1}^{s+1}\Gamma(h_j+q)}
{\prod_{j=1}^{s+1}\Gamma(1+h_0-h_j+q)}e^{\epsilon_{s-1}\pi iq}
\nonumber\\ &\qquad\quad\times
\frac1{2\pi i}\int_{-t_0-i\infty}^{-t_0+i\infty}
\frac{\Gamma(-q+t)\,\Gamma(h_0+q+t)\,\Gamma(-t)}
{\Gamma(1+h_0-h_{s+2}+t)}
e^{\epsilon_1\pi it}\,\d t\,\d q,
\label{eq:2.19}
\end{align}
где $q_1=q_0+t_0$.
Поскольку $\Re h_{s+2}<1$ и $h_{s+2}\ne0,-1,-2,\dots$,
внутренний барнсовский интеграл сворачивается
с помощью леммы~\ref{lem:2.3}:
\begin{align*}
&
\frac1{2\pi i}\int_{-t_0-i\infty}^{-t_0+i\infty}
\frac{\Gamma(-q+t)\,\Gamma(h_0+q+t)\,\Gamma(-t)}
{\Gamma(1+h_0-h_{s+2}+t)}
e^{\pm\pi it}\,\d t
\\ &\qquad
=\frac{\Gamma(-q)}{\Gamma(1-h_{s+2}-q)}
\int_0^1\frac{x^{h_0+q-1}(1-x)^{-h_{s+2}-q}}{(1-x)^{-q}}\,\d x
\displaybreak[0]\\ &\qquad
=\frac{\Gamma(-q)}{\Gamma(1-h_{s+2}-q)}
\cdot\frac{\Gamma(h_0+q)\,\Gamma(1-h_{s+2})}{\Gamma(1+h_0-h_{s+2}+q)}
\displaybreak[0]\\ &\qquad
=\frac{\Gamma(h_0+q)\,\Gamma(h_{s+2}+q)\,\Gamma(-q)}
{\Gamma(h_{s+2})\,\Gamma(1+h_0-h_{s+2}+q)}
\cdot\frac{\sin\pi(h_{s+2}+q)}{\sin\pi h_{s+2}}
\\ &\qquad
=\frac{\Gamma(h_0+q)\,\Gamma(h_{s+2}+q)\,\Gamma(-q)}
{\Gamma(h_{s+2})\,\Gamma(1+h_0-h_{s+2}+q)}
\\ &\qquad\quad\times
\biggl(e^{\pi iq}\cdot\frac{1-i\ctg\pi h_{s+2}}2
+e^{-\pi iq}\cdot\frac{1+i\ctg\pi h_{s+2}}2\biggr).
\end{align*}
Подставляя последнее выражение в~\eqref{eq:2.19}, находим
\begin{equation*}
\begin{split}
&
J_s\biggl(\begin{matrix}
h_1, & h_2, & h_3, & \dots, & h_s, & h_{s+1} \\[1pt]
& 1+h_0-h_3, & 1+h_0-h_4, & \dots, & 1+h_0-h_{s+1}, & 1+h_0-h_{s+2}
\end{matrix}\biggr)
\\ &\quad
=\frac{\prod_{j=1}^{s+1}\Gamma(1+h_0-h_j-h_{j+1})}
{\Gamma(h_1)\,\Gamma(h_{s+2})}
\\ &\quad\quad\times
\biggl(\frac{1-i\ctg\pi h_{s+2}}{4\pi i}
\int_{-q_1-i\infty}^{-q_1+i\infty}
(h_0+2q)
\\ &\quad\quad\;\qquad\times
\frac{\prod_{j=0}^{s+2}\Gamma(h_j+q)\,\Gamma(-q)}
{\prod_{j=1}^{s+2}\Gamma(1+h_0-h_j+q)}
e^{(\epsilon_{s-1}+1)\pi iq}\,\d q
\\ &\quad\quad\;
+\frac{1+i\ctg\pi h_{s+2}}{4\pi i}
\int_{-q_1-i\infty}^{-q_1+i\infty}
(h_0+2q)
\\ &\quad\quad\;\qquad\times
\frac{\prod_{j=0}^{s+2}\Gamma(h_j+q)\,\Gamma(-q)}
{\prod_{j=1}^{s+2}\Gamma(1+h_0-h_j+q)}
e^{(\epsilon_{s-1}-1)\pi iq}\,\d q
\biggr).
\end{split}
\end{equation*}
Если $s$~четно, полагаем $\epsilon_{s-1}=-1$ в первом интеграле
и $\epsilon_{s-1}=1$ во втором. Поэтому оба интеграла равны
$$
\int_{-q_1-i\infty}^{-q_1+i\infty}
(h_0+2q)\frac{\prod_{j=0}^{s+2}\Gamma(h_j+q)\,\Gamma(-q)}
{\prod_{j=1}^{s+2}\Gamma(1+h_0-h_j+q)}
e^{\epsilon_s\pi iq}\,\d q
=2\pi i\cdot F_{s+2}(h_0;h_1,\dots,h_{s+2}),
$$
что представляет собой требуемое равенство~\eqref{eq:2.7}.
Теорема~\ref{th:6} доказана.
\end{proof}

\chapter({Иррациональность \$q\$-дзета-значений})%
{Иррациональность $q$-дзета-значений}
\label{chap:3}

Исследование арифметических свойств $q$-дзета-значений
основано на применении $q$-гипергеометрической конструкции
и $q$-арифметике, основные компоненты которой
приводятся в следующей таблице.
Здесь, как и раньше, $\[\,\cdot\,\]$~--- целая часть числа и
сокращение ``НОК'' используется для обозначения
наименьшего общего кратного.

\smallskip
\hbox to\hsize{\hss\vbox{\offinterlineskip
\halign to135.35mm{\strut\tabskip=1mm minus 1mm
\strut\vrule#&\hbox to60mm{\hss#\hss}&%
\vrule#&\hbox to75mm{\hss#\hss}&%
\vrule#\tabskip=0pt\cr\noalign{\hrule}
height9.2pt&обычные объекты&%
&$q$-расширения, $p=1/q\in\mathbb Z\setminus\{0,\pm1\}$&\cr
\noalign{\hrule\vskip1pt\hrule}
&числа $n\in\mathbb Z$&%
&``числа''
$[n]_p=\dfrac{p^n-1\vphantom{\big|}}{p-1\vphantom{\big|}}
\in\mathbb Z[p]$&\cr
\noalign{\hrule}
&простые $l\in\{2,3,5,7,\dots\}\subset\mathbb Z$&%
&$\gathered
\text{неприводимые круговые многочлены} \vphantom{\big|} \\ \vspace{-5.5mm}
\Phi_l(p)=\prod\doublesb{k=1}{(k,l)=1}^l(p-e^{2\pi ik/l})
\in\mathbb Z[p] \\ \vspace{-1mm}
\endgathered$&\cr
\noalign{\hrule}
&гамма-функция Эйлера $\Gamma(t)$&%
&$\gathered
\text{$q$-гамма-функция Джексона} \vphantom{\big|} \\ \vspace{-2mm}
\Gamma_q(t)=\frac{\prod_{\nu=1}^\infty(1-q^\nu)}
{\prod_{\nu=1}^\infty(1-q^{t+\nu-1})\vphantom{\Big|_{0_0}}}\,(1-q)^{1-t}
\endgathered$&\cr
\noalign{\hrule}
&$\gathered
\text{факториал $n!=\Gamma(n+1)$} \vphantom{\big|} \\ \vspace{-3.5mm}
n!=\prod_{\nu=1\vphantom{{}_{\displaystyle0_a}}}^n\nu\in\mathbb Z
\endgathered$&%
&$\gathered
\text{$q$-факториал $[n]_q!=\Gamma_q(n+1)$} \vphantom{\big|^0}
\\ \vspace{-3.5mm}
[n]_p!=\prod_{\nu=1}^n\frac{p^\nu-1}{p-1}
=p^{n(n-1)/2}[n]_q!\in\mathbb Z[p]
\\ \vspace{-2mm}
\endgathered$&\cr
\noalign{\hrule}
&$\displaystyle
\ord_ln!=\biggl\[\frac nl\biggr\]
+\biggl\[\frac n{l^2\vphantom{|_{\displaystyle0_0}}}\biggr\]+\dotsb$&%
&$\displaystyle
\ord_{\Phi_l(p)}[n]_p!=\biggl\[\frac{n\vphantom{\big|}}l\biggr\]$, \
$l=2,3,4,\dots$&\cr
\noalign{\hrule}
&$\aligned
D_n&=\text{НОК}(1,\dots,n) \\ \vspace{-3pt}
&=\prod_{\text{простые $l\le n$}}l^{\[\log n/\log l\]}
\in\mathbb Z
\endaligned$&%
&$\aligned
D_n(p)&=\text{НОК}(p-1,[1]_p,\dots,[n]_p) \vphantom{\big|^0} \\ \vspace{-2.5mm}
&=\prod_{l=1}^n\Phi_l(p)\in\mathbb Z[p]
\endaligned$&\cr
\noalign{\hrule}
&$\gathered
\text{асимптотический закон} \\[-1.2mm]
\text{распределения простых чисел} \\
\endgathered$&&формула Мертенса&\cr
&$\displaystyle
\lim_{n\to\infty}\frac{\log D_n\vphantom{\big|}}n=1
$&&$\displaystyle
\lim_{n\to\infty}\frac{\log|D_n(p)|}{n^2\log|p|}
=\frac3{\pi^2\vphantom{|_{\displaystyle0_0}}}
$&\cr
\noalign{\hrule}
}}\hss}

\smallskip
Арифметическая мотивировка и обоснование второго
столбца таблицы приводятся в следующем параграфе.
Настоящая глава основана на работах \cite{Z7}, \cite{Z9},
\cite{Z10}, \cite{Z12}.

\section({\$q\$-Арифметика}){$q$-Арифметика}
\label{sec:3.1}

Напомним стандартные $q$-обозначения (см.~\cite[гл.~1]{GR}):
\begin{gather*}
(a;q)_n=\prod_{\nu=1}^n(1-aq^{\nu-1}),
\qquad
(a_1,a_2,\dots,a_m;q)_n=(a_1;q)_n(a_2;q)_n\dotsb(a_m;q)_n,
\\
\Gamma_q(t)=\frac{(q;q)_\infty}{(q^t;q)_\infty}(1-q)^{1-t},
\qquad
[n]_q!=\Gamma_q(n+1)=\frac{(q;q)_n}{(1-q)^n},
\\
\qbinom nk=\frac{[n]_q!}{[k]_q!\,[n-k]_q!}
=\frac{(q;q)_n}{(q;q)_k\cdot(q;q)_{n-k}},
\end{gather*}
где $k=0,1,\dots,n$ и $n=0,1,2,\dots$\,.

Рассмотрим {\it круговые многочлены\/} \cite{Wa}
\begin{equation}
\Phi_l(x)=\prod\doublesb{k=1}{(k,l)=1}^l(x-e^{2\pi ik/l}),
\quad
\deg\Phi_l(x)=\varphi(l),
\qquad l=1,2,\dots,
\label{eq:3.1}
\end{equation}
где $\varphi(\,\cdot\,)$~--- функция Эйлера. Хорошо
известно, что коэффициенты многочленов~\eqref{eq:3.1}
являются целыми числами
и для каждого $l=1,2,\dots$ многочлен $\Phi_l(x)$
неприводим над~$\mathbb Z$;
кроме того, имеет место формула
\begin{equation}
x^n-1=\prod_{l\mid n}\Phi_l(x).
\label{eq:3.2}
\end{equation}
Поскольку
$$
(x;x)_n=(1-x)(1-x^2)\dotsb(1-x^n)
=\pm\prod_{k=1}^n\prod_{l\mid k}\Phi_l(x),
$$
мы получаем, что разложение $(x;x)_n$ в произведение
неприводимых многочленов содержит только многочлены~\eqref{eq:3.1},
причем
\begin{equation}
\ord_{\Phi_l(x)}(x;x)_n=\biggl\[\frac nl\biggr\],
\qquad l=1,2,\dots,
\label{eq:3.3}
\end{equation}
где $\[\,\cdot\,\]$ обозначает целую часть числа.
Простым следствием формулы~\eqref{eq:3.3} является
\begin{equation}
\ord_{\Phi_l(x)}\xbinom nk
=\biggl\[\frac nl\biggr\]
-\biggl\[\frac kl\biggr\]-\biggl\[\frac{n-k}l\biggr\],
\label{eq:3.4}
\end{equation}
что делает возможным рассматривать круговые многочлены
как $q$-аналоги простых чисел. В свою очередь,
формула~\eqref{eq:3.4} влечет включения
\begin{equation}
\xbinom nk\in\mathbb Z[x],
\qquad k=0,1,\dots,n, \quad n=0,1,2,\dots,
\label{eq:3.5}
\end{equation}
которые обычно доказываются с помощью $q$-версии
треугольника Паскаля.

Из разложения~\eqref{eq:3.2} следует, что многочлен
$$
D_n(x)=\prod_{l=1}^n\Phi_l(x)\in\mathbb Z[x]
$$
является наименьшим общим кратным многочленов $x-1,x^2-1,\dots,x^n-1$,
иными словами, $D_n(x)$~есть многочлен наименьшей степени,
удовлетворяющий включениям
$$
D_n(x)\cdot\frac1{x^k-1}\in\mathbb Z[x],
\qquad k=1,2,\dots,n.
$$
Формула Мертенса~\cite{Me} (см.\ также \cite[\S\,18.5, теорема~330]{HW})
\begin{equation}
\sum_{l\le n}\varphi(l)=\frac3{\pi^2}n^2+O(n\log n),
\label{eq:3.6}
\end{equation}
являющаяся в нашем случае $q$-аналогом асимптотического
закона распределения простых чисел, приводит к следующим
утверждениям.

\begin{lemma}[{см.\ \cite[\S\,2]{BV} и \cite[лемма~2]{As2}}]
\label{lem:3.1}
Для любого $p\in\mathbb Z\setminus\{0,\pm1\}$
справедливо предельное соотношение
$$
\lim_{n\to\infty}\frac{\log|D_n(p)|}{n^2}=\frac3{\pi^2}\log|p|.
$$
\end{lemma}

\begin{lemma}[\cite{Z9}]
\label{lem:3.2}
Для любого $p\in\mathbb Z\setminus\{0,\pm1\}$ и любого
полуинтервала $[u,v)\subset(0,1)$, $u,v\in\mathbb Q$,
справедливо предельное соотношение
\begin{align*}
\lim_{n\to\infty}\frac1{n^2}
\sum_{l:\{n/l\}\in[u,v)}\log|\Phi_l(p)|
&=\frac3{\pi^2}\bigl(\psi'(u)-\psi'(v)\bigr)\log|p|
\\
&=\frac{3\log|p|}{\pi^2}\int_u^v\d(-\psi'(x)),
\end{align*}
где $\psi(x)$~--- логарифмическая производная гамма-функции.
\end{lemma}

\begin{proof}
Согласно формуле~\eqref{eq:3.6}
\begin{align*}
\frac1{n^2}\sum_{l:\{n/l\}\in[u,v)}\log|\Phi_l(p)|
&\sim\frac1{n^2}\sum_{l:\{n/l\}\in[u,v)}\log|\Phi_l(p)|
\sim\frac{\log|p|}{n^2}\sum_{l:\{n/l\}\in[u,v)}\varphi(l)
\\
&\sim\frac{\log|p|}{n^2}\cdot\sum_{N=0}^\infty
\sum_{\frac n{v+N}<l\le\frac n{u+N}}\varphi(l)
\\
&\sim\frac{3\log|p|}{\pi^2n^2}\cdot\sum_{N=0}^\infty
\biggl(\frac{n^2}{(u+N)^2}-\frac{n^2}{(v+N)^2}\biggr)
\end{align*}
при $n\to\infty$. Остается воспользоваться формулой
$$
\psi'(x)=\frac{\d^2}{\d x^2}\log\Gamma(x)
=\frac1{x^2}+\frac1{(x+1)^2}+\frac1{(x+2)^2}+\frac1{(x+3)^2}+\dotsb.
$$
\vskip-9.5mm
\end{proof}

\bigskip
Как несложно видеть, лемма~\ref{lem:3.2} является
$q$-аналогом леммы~\ref{lem:1.5}.

\section({\$q\$-Гипергеометрическая конструкция}){$q$-Гипергеометрическая конструкция}
\label{sec:3.2}

Зафиксируем целые параметры
\begin{equation}
(\ba,\bb)
=\biggl(\begin{matrix}
a_1, & a_2, & a_3 \\[1pt]
b_1, & b_2, & b_3
\end{matrix}\biggr),
\label{eq:3.7}
\end{equation}
удовлетворяющие условиям
\begin{equation}
\min\{a_2,a_3\}\ge b_1=1, \quad
a_2<b_2, \quad a_3<b_3, \quad
a_1+a_2+a_3\le b_2+b_3-1,
\label{eq:3.8}
\end{equation}
и рассмотрим $q$-базисный гипергеометрический ряд~\cite{GR}
\begin{align}
G_q(\ba,\bb)
&=\frac{\Gamma_q(a_2)\,\Gamma_q(a_3)\,\Gamma_q(b_2-a_2)\,\Gamma_q(b_3-a_3)}
{(1-q)^2\Gamma_q(b_2)\,\Gamma_q(b_3)}
\nonumber\\ &\qquad\times 
{}_3\phi_2\biggl(\begin{matrix}
q^{a_1}, & q^{a_2}, & q^{a_3} \\[1pt]
         & q^{b_2}, & q^{b_3}
\end{matrix}\biggm|q,q^{b_2+b_3-a_1-a_2-a_3}\biggr)
\nonumber\\ &\phantom:
=\frac{\Gamma_q(a_2)\,\Gamma_q(a_3)\,\Gamma_q(b_2-a_2)\,\Gamma_q(b_3-a_3)}
{(1-q)^2\Gamma_q(b_2)\,\Gamma_q(b_3)}
\nonumber\\ &\qquad\times 
\sum_{t=0}^\infty
\frac{(q^{a_1},q^{a_2},q^{a_3};q)_t}{(q^{b_1},q^{b_2},q^{b_3};q)_t}
q^{t(b_2+b_3-a_1-a_2-a_3)},
\label{eq:3.9}
\end{align}
абсолютно сходящийся в области $|q|<1$.
Очевидная симметрия величины $G_q(\ba,\bb)$ приводит
к следующему утверждению.

\begin{lemma} \label{lem:3.3}
Величина $G_q(\ba,\bb)$ инвариантна относительно преобразования
$$
\sigma\:\biggl(\begin{matrix}
a_1, & a_2, & a_3 \\[1pt]
  1, & b_2, & b_3
\end{matrix}\biggr)
\mapsto\biggl(\begin{matrix}
a_1, & a_3, & a_2 \\[1pt]
  1, & b_3, & b_2
\end{matrix}\biggr).
$$
\end{lemma}

Из тождества Холла \cite[формула~(3.2.10)]{GR}
\begin{align*}
&
{}_3\phi_2\biggl(\begin{matrix}
q^{a_1}, & q^{a_2}, & q^{a_3} \\[1pt]
         & q^{b_2}, & q^{b_3}
\end{matrix}\biggm|q,q^{b_2+b_3-a_1-a_2-a_3}\biggr)
\\ &\quad
=\frac{\Gamma_q(b_2)\,\Gamma_q(b_3)\,\Gamma_q(b_2+b_3-a_1-a_2-a_3)}
{\Gamma_q(a_2)\,\Gamma_q(b_2+b_3-a_2-a_3)\,\Gamma_q(b_2+b_3-a_1-a_2)}
\\ &\quad\qquad\times 
{}_3\phi_2\biggl(\begin{matrix}
q^{b_3-a_2}, & q^{b_2-a_2},         & q^{b_2+b_3-a_1-a_2-a_3} \\[1pt]
             & q^{b_2+b_3-a_2-a_3}, & q^{b_2+b_3-a_1-a_2}
\end{matrix}\biggm|q,q^{a_2}\biggr)
\end{align*}
вытекает также ``нетривиальное'' преобразование
величины $G_q(\ba,\bb)$.

\begin{lemma} \label{lem:3.4}
Величина $G_q(\ba,\bb)$ инвариантна относительно преобразования
$$
\tau\:\biggl(\begin{matrix}
a_1, & a_2, & a_3 \\[1pt]
  1, & b_2, & b_3
\end{matrix}\biggr)
\mapsto\biggl(\begin{matrix}
b_3-a_2, & b_2-a_2,         & b_2+b_3-a_1-a_2-a_3 \\[1pt]
      1, & b_2+b_3-a_2-a_3, & b_2+b_3-a_1-a_2
\end{matrix}\biggr).
$$
\end{lemma}

Следующее утверждение содержит рекуррентные соотношения
для величины \eqref{eq:3.9}, которые являются $q$-аналогом
тождеств, полученных в доказательстве теоремы~2.1
в~\cite[с.~31]{RV2}.

\begin{lemma}
\label{lem:3.5}
Справедливо тождество
\begin{align}
&
q^{b_2+b_3-a_1-a_2-a_3}
G_q\biggl(\begin{matrix}
a_1, & a_2, & a_3 \\[1pt]
b_1, & b_2, & b_3
\end{matrix}\biggr)
\nonumber\\ &\quad
=G_q\biggl(\begin{matrix}
a_1, & a_2-1, & a_3-1 \\[1pt]
b_1, & b_2-1, & b_3-1
\end{matrix}\biggr)
-G_q\biggl(\begin{matrix}
a_1-1, & a_2-1, & a_3-1 \\[1pt]
b_1  , & b_2-1, & b_3-1
\end{matrix}\biggr).
\label{eq:3.10}
\end{align}
\end{lemma}

\begin{proof}
Поскольку
\begin{align*}
\frac{(q^{a_1};q)_{t+1}}{(q;q)_{t+1}}
-\frac{(q^{a_1-1};q)_{t+1}}{(q;q)_{t+1}}
&=\frac{(q^{a_1};q)_t}{(q;q)_t}
\cdot\frac{(1-q^{a_1+t})-(1-q^{a_1-1})}{1-q^{t+1}}
\\
&=q^{a_1-1}\frac{(q^{a_1};q)_t}{(q;q)_t},
\end{align*}
заключаем, что
\begin{align}
&
\sum_{t=0}^\infty
\frac{(q^{a_1},q^{a_2-1},q^{a_3-1};q)_t}{(q,q^{b_2-1},q^{b_3-1};q)_t}
q^{tc}
-\sum_{t=0}^\infty
\frac{(q^{a_1-1},q^{a_2-1},q^{a_3-1};q)_t}{(q,q^{b_2-1},q^{b_3-1};q)_t}
q^{tc}
\nonumber\\ &\quad
=\sum_{t=0}^\infty\biggl(
\frac{(q^{a_1},q^{a_2-1},q^{a_3-1};q)_{t+1}}
{(q,q^{b_2-1},q^{b_3-1};q)_{t+1}}
-\frac{(q^{a_1-1},q^{a_2-1},q^{a_3-1};q)_{t+1}}
{(q,q^{b_2-1},q^{b_3-1};q)_{t+1}}
\biggr)q^{(t+1)c}
\nonumber\\ &\quad
=q^{c+a_1-1}\frac{(1-q^{a_2-1})(1-q^{a_3-1})}{(1-q^{b_2-1})(1-q^{b_3-1})}
\sum_{t=0}^\infty
\frac{(q^{a_1},q^{a_2},q^{a_3};q)_t}{(q,q^{b_2},q^{b_3};q)_t}
q^{tc},
\label{eq:3.11}
\end{align}
где $c=b_2+b_3-a_1-a_2-a_3$. С другой стороны,
\begin{align}
&
\sum_{t=0}^\infty
\frac{(q^{a_1-1},q^{a_2-1},q^{a_3-1};q)_t}{(q,q^{b_2-1},q^{b_3-1};q)_t}
q^{tc}
-\sum_{t=0}^\infty
\frac{(q^{a_1-1},q^{a_2-1},q^{a_3-1};q)_t}{(q,q^{b_2-1},q^{b_3-1};q)_t}
q^{t(c+1)}
\nonumber\\ &\quad
=\sum_{t=0}^\infty
\frac{(q^{a_1-1},q^{a_2-1},q^{a_3-1};q)_t}{(q,q^{b_2-1},q^{b_3-1};q)_t}
(1-q^t)q^{tc}
\nonumber\\ &\quad
=\sum_{t=0}^\infty
\frac{(q^{a_1-1},q^{a_2-1},q^{a_3-1};q)_{t+1}}
{(q,q^{b_2-1},q^{b_3-1};q)_{t+1}}
(1-q^{t+1})q^{(t+1)c}
\nonumber\\ &\quad
=q^c\frac{(1-q^{a_1-1})(1-q^{a_2-1})(1-q^{a_3-1})}
{(1-q^{b_2-1})(1-q^{b_3-1})}
\sum_{t=0}^\infty
\frac{(q^{a_1},q^{a_2},q^{a_3};q)_t}{(q,q^{b_2},q^{b_3};q)_t}
q^{tc}.
\label{eq:3.12}
\end{align}
Складывая левые и правые части соотношений~\eqref{eq:3.11} и \eqref{eq:3.12},
получаем
\begin{align}
&
\sum_{t=0}^\infty
\frac{(q^{a_1},q^{a_2-1},q^{a_3-1};q)_t}{(q,q^{b_2-1},q^{b_3-1};q)_t}
q^{tc}
-\sum_{t=0}^\infty
\frac{(q^{a_1-1},q^{a_2-1},q^{a_3-1};q)_t}{(q,q^{b_2-1},q^{b_3-1};q)_t}
q^{t(c+1)}
\nonumber\\ &\quad
=q^c\frac{(1-q^{a_2-1})(1-q^{a_3-1})}{(1-q^{b_2-1})(1-q^{b_3-1})}
\sum_{t=0}^\infty
\frac{(q^{a_1},q^{a_2},q^{a_3};q)_t}{(q,q^{b_2},q^{b_3};q)_t}
q^{tc}.
\label{eq:3.13}
\end{align}
После домножения обеих частей получившегося равенства~\eqref{eq:3.13} на
$$
\frac{\Gamma_q(a_2-1)\,\Gamma_q(a_3-1)\,
\Gamma_q(b_2-a_2)\,\Gamma_q(b_3-a_3)}
{(1-q)^2\Gamma_q(b_2-1)\,\Gamma_q(b_3-1)}
$$
мы приходим к требуемому соотношению~\eqref{eq:3.10}.
\end{proof}

В следующем параграфе мы показываем,
что построенные нами величины~\eqref{eq:3.9}
являются линейными формами от~$1$ и~$\zeta_q(2)$.

\section{Арифметика линейных форм}
\label{sec:3.3}

Свяжем с параметрами~\eqref{eq:3.7} другой набор из десяти
целых чисел:
\begin{equation}
\begin{aligned}
c_{00}&=(b_1+b_2+b_3)-(a_1+a_2+a_3)-2,
\\
c_{jk}&=\begin{cases}
a_j-b_k &\text{для $k=1$}, \\
b_k-a_j-1 &\text{для $k=2,3$},
\end{cases}
\qquad j,k=1,2,3.
\end{aligned}
\label{eq:3.14}
\end{equation}
Согласно \eqref{eq:3.8} набор
\begin{equation}
\{c_{00},c_{21},c_{22},c_{33},c_{31}\}
\label{eq:3.15}
\end{equation}
состоит из целых неотрицательных чисел, в то время как
целые параметры в наборе
\begin{equation}
\{c_{11},c_{23},c_{13},c_{12},c_{32}\}
\label{eq:3.16}
\end{equation}
могут быть и отрицательными.
Отметим также, что исходные параметры~\eqref{eq:3.7}
однозначно восстанавливаются по любому из наборов~\eqref{eq:3.15}
и \eqref{eq:3.16}:
\begin{gather*}
\begin{alignedat}{3}
a_1&=c_{22}+c_{33}-c_{00}+1, \quad&
a_2&=c_{21}+1, \quad&
a_3&=c_{31}+1,
\\
b_1&=1, \quad&
b_2&=c_{21}+c_{22}+2, \quad&
b_3&=c_{31}+c_{33}+2;
\end{alignedat}
\\
\begin{alignedat}{3}
a_1&=c_{11}+1, \quad&
a_2&=c_{11}+c_{13}-c_{23}+1, \quad&
a_3&=c_{11}+c_{12}-c_{32}+1,
\\
b_1&=1, \quad&
b_2&=c_{11}+c_{12}+2, \quad&
b_3&=c_{11}+c_{13}+2.
\end{alignedat}
\end{gather*}

Согласно леммам~\ref{lem:3.3},~\ref{lem:3.4} и формулам
\eqref{eq:3.14} действие преобразований $\sigma,\tau$ на параметры~$\bc$
определяется следующим образом:
\begin{align}
\sigma\:\biggl(\begin{matrix}
c_{00}, & c_{21}, & c_{22}, & c_{33}, & c_{31} \\[1pt]
c_{11}, & c_{23}, & c_{13}, & c_{12}, & c_{32}
\end{matrix}\biggr)
&\mapsto\biggl(\begin{matrix}
c_{00}, & c_{31}, & c_{33}, & c_{22}, & c_{21} \\[1pt]
c_{11}, & c_{32}, & c_{12}, & c_{13}, & c_{23}
\end{matrix}\biggr),
\label{eq:3.17}
\\
\tau\:\biggl(\begin{matrix}
c_{00}, & c_{21}, & c_{22}, & c_{33}, & c_{31} \\[1pt]
c_{11}, & c_{23}, & c_{13}, & c_{12}, & c_{32}
\end{matrix}\biggr)
&\mapsto\biggl(\begin{matrix}
c_{21}, & c_{22}, & c_{33}, & c_{31}, & c_{00}  \\[1pt]
c_{23}, & c_{13}, & c_{12}, & c_{32}, & c_{11}
\end{matrix}\biggr).
\label{eq:3.18}
\end{align}
Таким образом, преобразования $\sigma,\tau$ задают перестановки
$10$-элементного множества
\begin{equation}
\bc=\biggl(\begin{matrix}
c_{00}, & c_{21}, & c_{22}, & c_{33}, & c_{31} \\[1pt]
c_{11}, & c_{23}, & c_{13}, & c_{12}, & c_{32}
\end{matrix}\biggr)
\label{eq:3.19}
\end{equation}
и не меняют величины
\begin{equation}
H_q(\bc)=G_q(\ba,\bb).
\label{eq:3.20}
\end{equation}
Обозначим через $\fG_0\subset\fS_{10}$ группу, порожденную
перестановками \eqref{eq:3.17}, \eqref{eq:3.18}; отметим, что
порядок~$\sigma$ равен~$2$, а порядок~$\tau$ равен~$5$.

Отметим также, что формулы \eqref{eq:3.14} позволяют записать рекуррентные
соотношения леммы~\ref{lem:3.5} в виде
\begin{align}
q^{c_{00}+1}H_q(\bc)
&=H_q\biggl(\begin{matrix}
c_{00}, & c_{21}-1, & c_{22},   & c_{33},   & c_{31}-1 \\[1pt]
c_{11}, & c_{23},   & c_{13}-1, & c_{12}-1, & c_{32}
\end{matrix}\biggr)
\nonumber\\ &\qquad
-H_q\biggl(\begin{matrix}
c_{00}+1, & c_{21}-1, & c_{22}, & c_{33}, & c_{31}-1 \\[1pt]
c_{11}-1, & c_{23},   & c_{13}, & c_{12}, & c_{32}
\end{matrix}\biggr).
\label{eq:3.21}
\end{align}

Для каждого набора параметров~\eqref{eq:3.7} и соответствующего
набора~\eqref{eq:3.14} определим величину
\begin{align}
m=m(\bc)
&=c_{00}+c_{21}+c_{22}+c_{33}+c_{31}
=c_{11}+c_{23}+c_{13}+c_{12}+c_{32}
\nonumber\\
&=2(b_1+b_2+b_3)-(a_1+a_2+a_3)-7;
\label{eq:3.22}
\end{align}
через $m_1=m_1(\bc)$ и $m_2=m_2(\bc)$ обозначим
два максимальных элемента, стоящих на разных местах
в наборе~\eqref{eq:3.16};
тот факт, что $m_1\ge0$ и $m_2\ge0$, доказан
в~\cite[теорема~2.1]{RV2}. Положим
\begin{align}
N(\bc)
&=\frac{c_{00}(c_{00}-1)+c_{21}(c_{21}-1)+c_{22}(c_{22}-1)
+c_{33}(c_{33}-1)+c_{31}(c_{31}-1)}2,
\label{eq:3.43}
\\
L(\bc)
&=N(\bc)-(c_{00}c_{22}+c_{21}c_{33}+c_{22}c_{31}+c_{33}c_{00}+c_{31}c_{21})-1
\nonumber\\
&=\frac{a_1(a_1-1)}2+\frac{a_2(a_2-1)}2+\frac{a_3(a_3-1)}2
\nonumber\\ &\quad
-(b_1-1)(b_2-1)-(b_2-1)(b_3-1)-(b_3-1)(b_1-1).
\label{eq:3.23}
\end{align}
Как следует из определений, величины
$$
H_q(\bc), \quad m(\bc), \quad m_1(\bc), \quad m_2(\bc), \quad N(\bc), \quad L(\bc)
$$
инвариантны относительно действия группы~$\fG_0$.

Отметим, что $\zeta_q(2)$ не является рациональной
(и даже алгебраической) функцией над полем $\mathbb C(q)$.

\begin{proposition} \label{prop:3.1}
Справедливо включение
\begin{equation}
D_{m_1}(x)D_{m_2}(x)\cdot H_q(\bc)
\in\mathbb Z[x]\zeta_q(2)+\mathbb Z[x],
\quad\text{где}\;\; x=q^{-1}.
\label{eq:3.24}
\end{equation}
\end{proposition}

\begin{proof}
Будем проводить доказательство индукцией по величине
\begin{equation}
m(\bc)=c_{00}+c_{21}+c_{22}+c_{33}+c_{31},
\label{eq:3.25}
\end{equation}
где каждое слагаемое в сумме~\eqref{eq:3.25} неотрицательно.

В качестве базы индукции рассмотрим случаи, когда
по крайней мере три из параметров~\eqref{eq:3.15} нулевые.
Здесь возможны две ситуации: три нулевых параметра
расположены или не расположены по порядку в циклическом
наборе~\eqref{eq:3.15}
(т.е.\ параметр $c_{00}$ считается следующим за параметром~$c_{31}$).
Первая ситуация с помощью одного или нескольких применений
циклической перестановки~\eqref{eq:3.18} может быть сведена
к случаю
\begin{equation}
c_{22}=c_{33}=c_{31}=0,
\label{eq:3.26}
\end{equation}
а вторая --- к случаю
\begin{equation}
c_{00}=c_{22}=c_{33}=0.
\label{eq:3.27}
\end{equation}
Рассмотрим сначала вторую возможность.

В случае~\eqref{eq:3.27} имеем
$$
c_{11}=c_{22}+c_{33}-c_{00}=0,
$$
откуда
$$
a_1=c_{11}+1=1, \quad b_2=c_{22}+a_2+1=a_2+1, \quad b_3=c_{33}+a_3+1=a_3+1.
$$
Поэтому
\begin{align*}
H_q(\bc)
&=G_q\biggl(\begin{matrix}
1, & a_2,   & a_3   \\[1pt]
1, & a_2+1, & a_3+1
\end{matrix}\biggr)
\\
&=\frac{\Gamma_q(a_2)\,\Gamma_q(a_3)}
{(1-q)^2\Gamma_q(a_2+1)\,\Gamma_q(a_3+1)}
\cdot{}_3\phi_2\biggl(\begin{matrix}
q, & q^{a_2},   & q^{a_3}   \\[1pt]
   & q^{a_2+1}, & q^{a_3+1}
\end{matrix}\biggm|q,q\biggr)
\\
&=\frac1{(1-q^{a_2})(1-q^{a_3})}
\sum_{t=0}^\infty
\frac{(q^{a_2},q^{a_3};q)_t}{(q^{a_2+1},q^{a_3+1};q)_t}q^t
\\
&=\sum_{t=0}^\infty\frac{q^t}{(1-q^{t+a_2})(1-q^{t+a_3})}.
\end{align*}
Если $a_2=a_3$, то
\begin{align*}
H_q(\bc)
&=q^{-a_2}\sum_{t=0}^\infty\frac{q^{t+a_2}}{(1-q^{t+a_2})^2}
=q^{-a_2}\biggl(\sum_{n=1}^\infty-\sum_{n=1}^{a_2-1}\biggr)
\frac{q^n}{(1-q^n)^2}
\\
&=q^{-a_2}\zeta_q(2)
-q^{-a_2}\sum_{n=1}^{a_2-1}\frac{q^n}{(1-q^n)^2}
=x^{a_2}\zeta_q(2)
-x^{a_2}\sum_{n=1}^{a_2-1}\frac{x^n}{(x^n-1)^2},
\end{align*}
откуда вытекает включение
\begin{equation}
D_{a_2-1}(x)^2\cdot H_q(\bc)
\in\mathbb Z[x]\zeta_q(2)+\mathbb Z[x].
\label{eq:3.28}
\end{equation}
Если $a_2\ne a_3$, то (считая для определенности $a_2<a_3$)
\begin{align*}
H_q(\bc)
&=\frac1{q^{a_2}-q^{a_3}}
\sum_{t=0}^\infty\biggl(\frac1{1-q^{t+a_2}}-\frac1{1-q^{t+a_3}}\biggr)
\\
&=\frac1{q^{a_2}-q^{a_3}}\biggl(\sum_{n=a_2}^{a_3-1}\frac1{1-q^n}+a_2-a_3\biggr)
\\
&=\frac{x^{a_3}}{x^{a_3-a_2}-1}
\biggl(x^{a_2}\sum_{n=a_2}^{a_3-1}\frac{x^{n-a_2}}{x^n-1}+a_2-a_3\biggr)
\end{align*}
и, значит,
\begin{equation}
D_{a_3-a_2}(x)D_{a_3-1}(x)\cdot H_q(\bc)
\in\mathbb Z[x].
\label{eq:3.29}
\end{equation}
Оба включения \eqref{eq:3.28}, \eqref{eq:3.29} означают, что
в случае~\eqref{eq:3.27} требуемое утверждение~\eqref{eq:3.24} выполнено,
поскольку
$$
\{c_{11},c_{23},c_{13},c_{12},c_{32}\}
=\{0,a_3-a_2,a_3-1,a_2-1,a_2-a_3\}.
$$

Рассмотрим теперь случай~\eqref{eq:3.26}. Возможность $c_{00}=0$
была исследована выше, поэтому считаем далее $c_{00}>0$,
т.е.\ $a_1\le0$ и ряд в~\eqref{eq:3.9} содержит лишь
конечное число слагаемых. По аналогичной причине мы можем
считать $c_{21}>0$ (иначе после применения перестановки~$\tau$
мы вновь приходим к случаю~\eqref{eq:3.27}), так что $a_2>1$.
Имеем
$$
b_2=c_{22}+a_2+1=a_2+1, \quad a_3=c_{31}+1=1, \quad
b_3=c_{33}+a_3+1=2,
$$
откуда
\begin{align}
H_q(\bc)
&=G_q\biggl(\begin{matrix}
a_1, & a_2,   & 1   \\[1pt]
  1, & a_2+1, & 2
\end{matrix}\biggr)
\nonumber\\
&=\frac{\Gamma_q(a_2)}{(1-q)^2\Gamma_q(a_2+1)}
\cdot{}_3\phi_2\biggl(\begin{matrix}
q^{a_1}, & q^{a_2},   & q   \\[1pt]
         & q^{a_2+1}, & q^2
\end{matrix}\biggm|q,q^{-a_1+2}\biggr)
\nonumber\\
&=\frac1{(1-q)(1-q^{a_2})}
\sum_{t=0}^\infty\frac{(q^{a_1},q^{a_2};q)_t}{(q^2,q^{a_2+1};q)_t}
q^{t(-a_1+2)}.
\label{eq:3.30}
\end{align}
Полагая $n=-a_1\ge0$, перепишем получившийся ряд в виде
\begin{align}
H_q(\bc)
&=\sum_{t=0}^n\frac{(q^{-n};q)_t}{(q;q)_{t+1}}
\cdot\frac{q^{t(n+2)}}{1-q^{t+a_2}}
=\sum_{t=0}^n\frac{(q^{-n};q)_t}{(q;q)_t}
\cdot\frac{q^{t(n+2)}}{(1-q^{t+1})(1-q^{t+a_2})}
\nonumber\\
&=\sum_{t=0}^n(-1)^tx^{t(t+1)/2}\xbinom nt
\cdot\frac{x^{(t+1)+(t+a_2)-t(n+2)}}{(x^{t+1}-1)(x^{t+a_2}-1)}
\nonumber\\
&=x^{a_2+1-n(n-1)/2}\sum_{t=0}^n(-1)^t\xbinom nt
\cdot\frac{x^{(n-t)(n-t-1)/2}}{(x^{t+1}-1)(x^{t+a_2}-1)}.
\label{eq:3.31}
\end{align}
Согласно~\eqref{eq:3.5} формула~\eqref{eq:3.31} влечет включение
\begin{equation}
x^{-(a_2+1)+n(n-1)/2}\cdot D_{n+1}(x)D_{a_2+n}(x)\cdot H_q(\bc)
\in\mathbb Z[x].
\label{eq:3.32}
\end{equation}
С другой стороны, сумма в правой части~\eqref{eq:3.30}
имеет другое представление:
\begin{align*}
H_q(\bc)
&=\frac{q^{a_1-2}}{(1-q^{a_1-1})(1-q^{a_2-1})}
\sum_{t=0}^\infty
\frac{(q^{a_1-1},q^{a_2-1};q)_{t+1}}{(q,q^{a_2};q)_{t+1}}
q^{(t+1)(-a_1+2)}
\\
&=\frac{q^{a_1-2}}{(1-q^{a_1-1})(1-q^{a_2-1})}
\biggl({}_2\phi_1\biggl(\begin{matrix}
q^{a_1-1}, & q^{a_2-1} \\[1pt]
           & q^{a_2}
\end{matrix}\biggm|q,q^{-a_1+2}\biggr)-1\biggr).
\end{align*}
Теперь $q$-аналог Гейне \cite[формула~(1.5.2)]{GR}
для формулы суммирования Гаусса позволяет свернуть
последний $q$-базисный ряд:
\begin{align*}
H_q(\bc)
&=\frac{q^{a_1-2}}{(1-q^{a_1-1})(1-q^{a_2-1})}
\biggl(\frac{(q;q)_{-a_1+1}}{(q^{a_2};q)_{-a_1+1}}-1\biggr)
\\
&=-q^{-1}\frac{(q;q)_n}{(q^{a_2-1};q)_{n+2}}
-\frac{q^{-n-2}}{(1-q^{-n-1})(1-q^{a_2-1})}
\\
&=-x^{a_2(n+2)}\frac{(x;x)_n}{(x^{a_2-1};x)_{n+2}}
-\frac{x^{a_2+n+1}}{(1-x^{n+1})(x^{a_2-1}-1)},
\end{align*}
откуда следует включение
\begin{equation}
(x^{n+1}-1)(x;x)_{a_2+n}\cdot H_q(\bc)
\in\mathbb Z[x].
\label{eq:3.33}
\end{equation}
Поскольку многочлен $x$ взаимно прост с многочленами
$D_{n+1}(x)$, $D_{a_2+n}(x)$,
$x^{n+1}-1$ и $(x;x)_{a_2+n}$,
включения~\eqref{eq:3.32},~\eqref{eq:3.33} можно записать в виде
$$
D_{n+1}(x)D_{a_2+n}(x)\cdot H_q(\bc)
\in\mathbb Z[x]
$$
или после возврата к параметру $a_2=-n$:
$$
D_{-a_1+1}(x)D_{a_2-a_1}(x)\cdot H_q(\bc)
\in\mathbb Z[x].
$$
Поскольку
$$
\{c_{11},c_{23},c_{13},c_{12},c_{32}\}
=\{a_1-1,-a_2+1,-a_1+1,a_2-a_1,a_2-1\},
$$
требуемое включение~\eqref{eq:3.24} также оказывается выполненным.
Это завершает доказательство базы индукции.

Предположим теперь, что по крайней мере три параметра
в наборе \eqref{eq:3.15} ненулевые и для наборов~$\bc'$
с условием $m(\bc')<m(\bc)$ требуемое включение~\eqref{eq:3.24}
доказано. Среди трех ненулевых
параметров в циклическом наборе~\eqref{eq:3.15} всегда можно
выбрать два, не являющихся соседями; циклическая
перестановка~\eqref{eq:3.18} позволяет перейти к
$\fG_0$-эквивалентному набору, у которого $c_{21}>0$
и $c_{31}>0$. Тогда включение~\eqref{eq:3.24} вытекает
из рекуррентных соотношений~\eqref{eq:3.21} и индукционного
предположения, поскольку для наборов
\begin{equation*}
\begin{aligned}
\bc'&=\biggl(\begin{matrix}
c_{00}, & c_{21}-1, & c_{22},   & c_{33},   & c_{31}-1 \\[1pt]
c_{11}, & c_{23},   & c_{13}-1, & c_{12}-1, & c_{32}
\end{matrix}\biggr),
\\
\bc''&=\biggl(\begin{matrix}
c_{00}+1, & c_{21}-1, & c_{22}, & c_{33}, & c_{31}-1 \\[1pt]
c_{11}-1, & c_{23},   & c_{13}, & c_{12}, & c_{32}
\end{matrix}\biggr)
\end{aligned}
\end{equation*}
выполнено
\begin{equation*}
m_1(\bc')\le m_1(\bc), \quad m_2(\bc')\le m_2(\bc),
\qquad
m_1(\bc'')\le m_1(\bc), \quad m_2(\bc'')\le m_2(\bc).
\end{equation*}
Это завершает доказательство индукционного перехода
и предложения~\ref{prop:3.1}.
\end{proof}

В виду трансцендентности $q$-функции $\zeta_q(2)$ над $\mathbb C(q)$, представление \eqref{eq:3.24}
определяет полиномиальные (по $x=q^{-1}$) коэффициенты линейной формы $D_{m_1}(x)D_{m_2}(x)\cdot H_q(\bc)$ единственным образом.
Обозначим через $M=M(\bc)=M(\ba,\bb)$ минимальный порядок нуля по~$x$ этих двух многочленов; иными словами,
максимальное целое $M\ge0$, для которого
\begin{equation}
x^{-M}\cdot D_{m_1}(x)D_{m_2}(x)\cdot H_q(\bc)
\in\mathbb Z[x]\zeta_q(2)+\mathbb Z[x],
\quad\text{где}\;\; x=q^{-1}.
\label{eq:3.34}
\end{equation}
$\fG_0$-Инвариантность величины $H_q(\bc)$ означает, что новая характеристика $M(\bc)$ также $\fG_0$-инвариантна.

\section({Групповая структура для \$\003\266\_q(2)\$})%
{Групповая структура для $\zeta_q(2)$}
\label{sec:3.4}

Потребуем теперь более жестких, чем \eqref{eq:3.8},
ограничений на параметры~\eqref{eq:3.7}:
\begin{equation}
\{b_1=1\}\le\{a_1,a_2,a_3\}<\{b_2,b_3\},
\quad
a_1+a_2+a_3\le b_1+b_2+b_3-2.
\label{eq:3.38}
\end{equation}
Тогда все параметры~\eqref{eq:3.14} неотрицательны,
и помимо преобразований~$\sigma,\tau$ можно также рассматривать
всевозможные перестановки параметров
$a_1,\linebreak[2]a_2,a_3$, с одной
стороны, и перестановку параметров $b_2,b_3$, с другой стороны,
не меняющие величины
\begin{equation}
\frac{\Gamma_q(a_1)}{\Gamma_q(b_2-a_2)\,\Gamma_q(b_3-a_3)}
\cdot G_q(\ba,\bb)
=\frac{[c_{11}]_q!}{[c_{22}]_q!\,[c_{33}]_q!}\cdot H_q(\bc).
\label{eq:3.39}
\end{equation}
Таким образом, \eqref{eq:3.39} инвариантна относительно
действия ``$(\ba,\bb)$-триви\-альной''
группы~$\fG_1$, порожденной перестановками
\begin{equation}
\fa_1=(a_1 \; a_3),
\quad
\fa_2=(a_2 \; a_3),
\quad
\fb=(b_2 \; b_3).
\label{eq:3.40}
\end{equation}
В терминах параметров~$\bc$ действие перестановок~\eqref{eq:3.40}
записывается следующим образом:
\begin{equation}
\begin{gathered}
\fa_1=(c_{11} \; c_{31})(c_{12} \; c_{32})(c_{13} \; c_{33}),
\quad
\fa_2=(c_{21} \; c_{31})(c_{22} \; c_{32})(c_{23} \; c_{33}),
\\
\fb=(c_{12} \; c_{13})(c_{22} \; c_{23})(c_{32} \; c_{33}).
\end{gathered}
\label{eq:3.41}
\end{equation}
Рассматривая теперь группу
$$
\fG=\<\fG_0,\fG_1\>=\<\sigma,\tau,\fa_1,\fa_2,\fb\>
$$
как группу перестановок 10-элементного множества~$\bc$,
а также учитывая $\fG_0$-инва\-риантность величины~\eqref{eq:3.20}
и $\fG_1$-инвариантность величины~\eqref{eq:3.39},
мы приходим к следующему утверждению.

\begin{lemma}[{ср.\ с~\cite[\S\,3]{RV2} и \cite[лемма~14]{Z17}}]
\label{lem:3.6}
Величина
$$
\frac{H_q(\bc)}{\Pi_q(\bc)},
\qquad\text{где}\quad
\Pi_q(\bc)=[c_{00}]_q!\,[c_{21}]_q!\,[c_{22}]_q!\,[c_{33}]_q!\,[c_{31}]_q!,
$$
инвариантна относительно действия группы~$\fG$.
\end{lemma}

В~\cite{RV2} доказано, что группа $\fG\subset\fS_{10}$ имеет
порядок~$120$. В качестве образующих группы~$\fG$ можно выбрать
(см.\ \cite[\S\,6]{Z17}) перестановки~\eqref{eq:3.41} и
$$
\fh=(c_{00} \; c_{22})(c_{11} \; c_{33})(c_{13} \; c_{31})
$$
второго порядка; при этом $\sigma=\fa_2\,\fb$ и
$\tau=(\fa_2\,\fa_1\,\fb\,\fh)^2$.

\begin{lemma} \label{lem:3.7}
Величины $m(\bc)$, $L(\bc)$ и $M(\bc)+N(\bc)$
инвариантны относительно действия группы~$\fG$.
\end{lemma}

\begin{proof}
Тот факт, что эти величины $\fG_0$-инвариантны был показан в~\S\,\ref{sec:3.3};
$\fG_1$-инвариантность $m(\bc)$ и $L(\bc)$ следует из определений \eqref{eq:3.22}
и \eqref{eq:3.23} этих характеристик через параметры $(\ba,\bb)$.
С учетом непосредственно проверяемого равенства
$$
[n]_q!=x^{-n(n-1)/2}[n]_x!,
\quad\text{где}\;\; x=q^{-1},
$$
из леммы~\ref{lem:3.6} вытекает $\fG_1$-инвариантность величины
\begin{equation}
\frac{H_q(\bc)}{x^{-N(\bc)}\Pi_x(\bc)},
\label{eq:3.42}
\end{equation}
где $N(\bc)$ определено в~\eqref{eq:3.43}; отсюда и из \eqref{eq:3.34}
заключаем, что $M(\bc)+N(\bc)$ инвариантны относительно действия группы~$\fG_1$.
\end{proof}

Обозначая через~$m_1^*\ge m_2^*$ два максимальных элемента,
стоящих на разных местах в $10$-элементном множестве~$\bc$,
в соответствии с~\eqref{eq:3.34} получаем включение
\begin{equation}
x^{-M}\cdot D_{m_1^*}(x)D_{m_2^*}(x)\cdot H_q(\bc)
\in\mathbb Z[x]\zeta_q(2)+\mathbb Z[x],
\quad\text{где}\;\; x=q^{-1}.
\label{eq:3.44}
\end{equation}
Кроме того, величины $m_1^*,m_2^*$ (в отличие от введенных
в~\S\,\ref{sec:3.3} величин $m_1,m_2$) являются $\fG$-инвариантными.

Для заданного набора параметров~$(\ba,\bb)$, удовлетворяющего
условиям \eqref{eq:3.38}, и соответствующего набора~\eqref{eq:3.14}
определим многочлен
$$
\Omega(x)=\prod_{l=1}^{m_1^*}\Phi_l^{\nu_l}(x)\in\mathbb Z[x],
$$
где
\begin{equation}
\nu_l=\max_{\fg\in\fG}\ord_{\Phi_l(x)}\frac{\Pi_x(\bc)}{\Pi_x(\fg\bc)},
\qquad l=1,2,\dots\,.
\label{eq:3.45}
\end{equation}

\begin{proposition} \label{prop:3.2}
Справедливо включение
\begin{equation}
x^{-M}\cdot D_{m_1^*}(x)D_{m_2^*}(x)\cdot\Omega^{-1}(x)\cdot H_q(\bc)
\in\mathbb Z[x]\zeta_q(2)+\mathbb Z[x],
\quad\text{где}\;\; x=q^{-1}.
\label{eq:3.46}
\end{equation}
\end{proposition}

\begin{proof}
Ввиду $\fG$-инвариантности величин $M(\bc)+N(\bc)$, $m_1^*(\bc)$,
$m_2^*(\bc)$ и \eqref{eq:3.42}
согласно включению~\eqref{eq:3.44} заключаем, что
для любой перестановки $\fg\in\fG$ линейная форма
\begin{align*}
&
x^{-M(\bc)}\cdot D_{m_1^*(\bc)}(x)D_{m_2^*(\bc)}(x)
\cdot\frac{\Pi_x(\fg\bc)}{\Pi_x(\bc)}\cdot H_q(\bc)
\\ &\qquad
=x^{-M(\bc)-N(\bc)+N(\fg\bc)}\cdot D_{m_1^*(\bc)}(x)D_{m_2^*(\bc)}(x)
\cdot H_q(\fg\bc)
\\ &\qquad
=x^{-M(\fg\bc)}\cdot D_{m_1^*(\fg\bc)}(x)D_{m_2^*(\fg\bc)}(x)
\cdot H_q(\fg\bc)
\end{align*}
лежит в $\mathbb Z[x]\zeta_q(2)+\mathbb Z[x]$. Воспользуемся теперь тем фактом,
что $\zeta_q(2)$ как функция от~$x=q^{-1}$ иррациональна над $\mathbb Q(x)$;
кроме того, в разложении $\Pi_x(\fg\bc)$, $\fg\in\fG$, на
неприводимые множители участвуют только многочлены~\eqref{eq:3.1}.
С учетом доказанных в~\cite{RV2} неравенств
$$
\nu_l=0 \quad\text{для $l>m_1^*$}, \qquad
\nu_l\le1 \quad\text{для $m_2^*<l\le m_1^*$}
$$
мы получаем требуемые включения~\eqref{eq:3.46}.
Предложение доказано.
\end{proof}

Ввиду $\fG_0$-инвариантности величины~$H_q(\bc)$ для определения
показателей~\eqref{eq:3.45} нам достаточно рассмотреть
действие представителей левых смежных классов
факторгруппы~$\fG/\fG_0$ (порядка~$12$)
на $5$-элементное (упорядоченное) множество
$\bc'=(c_{00},c_{21},c_{22},c_{33},c_{31})$:
\begin{align}
\nu_l
&=\max_{\fg\in\fG/\fG_0}
\ord_{\Phi_l(x)}\frac{\Pi_x(\bc)}{\Pi_x(\fg\bc)}
\nonumber\\
&=\max_{\fg\in\fG/\fG_0}\biggl(
\biggl\[\frac{c_{00}}l\biggr\]
+\biggl\[\frac{c_{21}}l\biggr\]
+\biggl\[\frac{c_{22}}l\biggr\]
+\biggl\[\frac{c_{33}}l\biggr\]
+\biggl\[\frac{c_{31}}l\biggr\]
\nonumber\\ &\qquad
-\biggl\[\frac{\fg c_{00}}l\biggr\]
-\biggl\[\frac{\fg c_{21}}l\biggr\]
-\biggl\[\frac{\fg c_{22}}l\biggr\]
-\biggl\[\frac{\fg c_{33}}l\biggr\]
-\biggl\[\frac{\fg c_{31}}l\biggr\]
\biggr),
\quad l=1,2,\dots,m_1^*,
\label{eq:3.47}
\end{align}
согласно~\eqref{eq:3.3}.
В качестве таких представителей выберем
\begin{equation}
\begin{gathered}
\fg_0=\id,
\quad
\fg_1=\fa_1\,\fa_2\,\fa_1,
\quad
\fg_2=\fa_1,
\quad
\fg_3=\fa_2,
\\
\fg_4=\fa_1\,\fa_2,
\quad
\fg_5=\fa_2\,\fa_1,
\quad
\fg_6=\fh\,\fa_1\,\fa_2\,\fa_1,
\quad
\fg_7=\fh\,\fa_2,
\\
\fg_8=\fh\,\fa_1\,\fa_2,
\quad
\fg_9=\fh\,\fa_2\,\fa_1,
\quad
\fg_{10}=\fh\,\fa_1\,\fa_2\,\fa_1\,\fh\,\fa_2,
\quad
\fg_{11}=\fh\,\fa_1\,\fa_2\,\fa_1\,\fh\,\fa_2\,\fa_1.
\end{gathered}
\label{eq:3.48}
\end{equation}
Тогда
\begin{equation}
\begin{gathered}
\fg_0\bc'=(c_{00},c_{21},c_{22},c_{33},c_{31}),
\qquad
\fg_1\bc'=(c_{00},c_{11},c_{12},c_{33},c_{31}),
\\
\fg_2\bc'=(c_{00},c_{21},c_{22},c_{13},c_{11}),
\qquad
\fg_3\bc'=(c_{00},c_{31},c_{32},c_{23},c_{21}),
\\
\fg_4\bc'=(c_{00},c_{11},c_{12},c_{23},c_{21}),
\qquad
\fg_5\bc'=(c_{00},c_{31},c_{32},c_{13},c_{11}),
\\
\fg_6\bc'=(c_{22},c_{33},c_{12},c_{11},c_{13}),
\qquad
\fg_7\bc'=(c_{22},c_{13},c_{32},c_{23},c_{21}),
\\
\fg_8\bc'=(c_{22},c_{33},c_{12},c_{23},c_{21}),
\qquad
\fg_9\bc'=(c_{22},c_{13},c_{32},c_{31},c_{33}),
\\
\fg_{10}\bc'=(c_{12},c_{23},c_{32},c_{31},c_{33}),
\qquad
\fg_{11}\bc'=(c_{12},c_{23},c_{32},c_{13},c_{11}).
\end{gathered}
\label{eq:3.49}
\end{equation}

\section{Оценки линейных форм и их коэффициентов}
\label{sec:3.5}

В этом параграфе мы также будем считать, что набор целочисленных
параметров~\eqref{eq:3.7} удовлетворяет условиям~\eqref{eq:3.38}.
Используя явное выражение величины~\eqref{eq:3.20} в виде
\begin{equation}
\begin{gathered}
G_q(\ba,\bb)=A\zeta_q(2)-B,
\\
A=A_q(\ba,\bb)=A_q(\bc)\in\mathbb Q(q), \quad
B=B_q(\ba,\bb)=B_q(\bc)\in\mathbb Q(q),
\end{gathered}
\label{eq:3.50}
\end{equation}
мы получим оценки для $|G_q(\ba,\bb)|$ и~$|A|$
при $|q|\le1/2$.

Начнем с иного представления величины~\eqref{eq:3.9}.
Именно, рассмотрим функцию
\begin{align}
R_q(t)
=R_q(\ba,\bb;t)
&=\frac{\Gamma_q(b_2-a_2)\,\Gamma_q(b_3-a_3)}
{(1-q)^2\Gamma_q(a_1-b_1+1)}
\cdot q^{t(b_1+b_2+b_3-a_1-a_2-a_3-2)}
\nonumber\\ &\phantom:\qquad\times 
\frac{\Gamma_q(t+a_1)\,\Gamma_q(t+a_2)\,\Gamma_q(t+a_3)}
{\Gamma_q(t+b_1)\,\Gamma_q(t+b_2)\,\Gamma_q(t+b_3)}
\label{eq:3.51}
\end{align}
и запишем~\eqref{eq:3.9} в виде
\begin{equation}
G_q(\ba,\bb)=\sum_{t=0}^\infty R_q(t)q^t.
\label{eq:3.52}
\end{equation}

\begin{proposition} \label{prop:3.3}
Пусть $c_{00}=b_2+b_3-a_1-a_2-a_3-1\ge5$ и $|q|\le1/2$.
Тогда справедливы оценки
\begin{equation}
3^{-3(b_2+b_3)}<|G_q(\ba,\bb)|<3^{3(b_2+b_3)}.
\label{eq:3.53}
\end{equation}
\end{proposition}

\begin{proof}
Из функционального уравнения для $\Gamma_q$-функции
\begin{equation}
\Gamma_q(t+1)=\frac{1-q^t}{1-q}\Gamma_q(t)
\label{eq:3.54}
\end{equation}
получаем
$$
\frac{R_q(t+1)}{R_q(t)}
=\frac{(1-q^{t+a_1})(1-q^{t+a_2})(1-q^{t+a_3})}
{(1-q^{t+b_1})(1-q^{t+b_2})(1-q^{t+b_3})}\cdot q^{c_{00}},
$$
откуда
\begin{equation}
\frac{|R_q(t+1)q^{t+1}|}{|R_q(t)q^t|}
\le\frac{(1+|q|)^3}{(1-|q|)^3}\cdot|q|^{c_{00}+1}
\le3^3\cdot2^{-(c_{00}+1)}<\frac12.
\label{eq:3.55}
\end{equation}
Применяя теперь полученную оценку~\eqref{eq:3.55}
к ряду в~\eqref{eq:3.52} находим, что, с одной стороны,
\begin{align}
|G_q(\ba,\bb)|
&\le|R_q(0)|\cdot\biggl(1+\frac{|R_q(1)q|}{|R_q(0)|}
+\frac{|R_q(2)q^2|}{|R_q(0)|}+\frac{|R_q(3)q^3|}{|R_q(0)|}+\dotsb\biggr)
\nonumber\\
&<|R_q(0)|\cdot\biggl(1+\frac12+\frac1{2^2}+\frac1{2^3}+\dotsb\biggr)
=2|R_q(0)|
\label{eq:3.56}
\end{align}
и, с другой стороны,
\begin{equation}
|G_q(\ba,\bb)|
\ge|R_q(0)|\cdot\biggl(1-\frac{|R_q(1)q|}{|R_q(0)|}\biggr)
>\frac12|R_q(0)|.
\label{eq:3.57}
\end{equation}
Воспользуемся тривиальными неравенствами
$$
3^{-n}\le\biggl(\frac{1-|q|}{1+|q|}\biggr)^n
\le|\Gamma_q(n+1)|
\le\biggl(\frac{1+|q|}{1-|q|}\biggr)^n\le3^n,
\qquad n=0,1,2,\dots,
$$
для оценки всех $\Gamma_q$-множителей, входящих в~$R_q(0)$.
С учетом
\begin{multline*}
(b_2-a_2-1)+(b_3-a_3-1)+(a_1-b_1)
+(a_1+a_2+a_3-3)+(b_1+b_2+b_3-3)
\\
<3(b_2+b_3)-1
\end{multline*}
это приводит к оценкам
\begin{equation}
3^{-3(b_2+b_3)+1}<|R_q(0)|<3^{3(b_2+b_3)-1}.
\label{eq:3.58}
\end{equation}
Собирая теперь вместе \eqref{eq:3.56}--\eqref{eq:3.58},
получаем требуемые неравенства
\linebreak[4]
\eqref{eq:3.53}.
Предложение доказано.
\end{proof}

\begin{proposition} \label{prop:3.4}
Для коэффициента $A=A_q(\ba,\bb)\in\mathbb Q(q)$
в представлении~\eqref{eq:3.50} при $|q|\le1/2$
справедлива оценка
\begin{equation}
|A|<3^{2(b_2+b_3)}
\cdot|q|^{L},
\label{eq:3.59}
\end{equation}
где показатель $L=L(\ba,\bb)=L(\bc)$ определен в~\eqref{eq:3.23}.
Кроме того, для для характеристики $M=M(\ba,\bb)=M(\bc)$,
определенной в \eqref{eq:3.34}, имеет место неравенство
\begin{align}
M\ge M_0(\ba,\bb)
&=(a_2-1)b_2+(a_3-1)b_3+\frac{a_1^2-3a_2^2-3a_3^2}4
\nonumber\\ &\quad
-\frac{a_1a_2+a_2a_3+a_3a_1-a_1-a_2-a_3}2+1.
\label{eq:3.59a}
\end{align}
\end{proposition}

\begin{proof}
Нам будет удобно ввести упорядоченную версию $(\ba^*,\bb^*)$
набора параметров $(\ba,\bb)$, именно
$$
\begin{gathered}
b_1^*=1<a_1^*\le a_2^*\le a_3^*<b_2^*\le b_3^*,
\\{}
\{a_1^*,a_2^*,a_3^*\}=\{a_1,a_2,a_3\},
\quad
\{b_1^*,b_2^*,b_3^*\}=\{b_1,b_2,b_3\}.
\end{gathered}
$$

Согласно функциональному уравнению~\eqref{eq:3.54}
выполнено
$$
\frac{\Gamma_q(t+a_j)}{\Gamma_q(t+b_j)}=\begin{cases}
\dfrac{(1-q^{t+b_j})(1-q^{t+b_j+1})\dotsb(1-q^{t+a_j-1})}{(1-q)^{a_j-b_j}}
& \text{при $j=1$}, \\
\dfrac{(1-q)^{b_j-a_j}}{(1-q^{t+a_j})(1-q^{t+a_j+1})\dotsb(1-q^{t+b_j-1})}
& \text{при $j=2,3$};
\end{cases}
$$
поэтому $R_q(t)$ в~\eqref{eq:3.51} является рациональной
функцией параметра $T=q^t$ над полем $\mathbb Q(q)=\mathbb Q(q^{-1})$:
\begin{align}
R_q(t)
&=\frac{[b_2-a_2-1]_q!\,[b_3-a_3-1]_q!}{[a_1-b_1]_q!}
\cdot\frac{(1-q^{b_1}T)\dotsb(1-q^{a_1-1}T)}{(1-q)^{a_1-b_1}}
\nonumber\\ &\qquad\times 
\frac{(1-q)^{b_2-a_2-1}}{(1-q^{a_2}T)\dotsb(1-q^{b_2-1}T)}
\cdot\frac{(1-q)^{b_3-a_3-1}}{(1-q^{a_3}T)\dotsb(1-q^{b_3-1}T)}
\nonumber\\ &\qquad\times 
T^{b_2+b_3-a_1-a_2-a_3-1}.
\label{eq:3.60}
\end{align}
Поскольку степень числителя функции~\eqref{eq:3.60} как функции от~$T$
на двойку меньше степени знаменателя, имеем
\begin{equation}
R_q(t)=O(T^{-2}) \qquad\text{при}\quad T\to\infty,
\label{eq:3.61}
\end{equation}
и $R_q(t)$ может быть разложена в сумму простейших дробей:
$$
R_q(t)
=\sum_{k=a_3^*}^{b_2^*-1}\frac{A_k}{(1-q^kT)^2}
+\sum_{k=a_2^*}^{b_3^*-1}\frac{B_k}{1-q^kT}.
$$
Из условия \eqref{eq:3.61} следует, что
$$
\sum_{k=a_2^*}^{b_3^*-1}B_kq^{-k}
=-\sum_{k=a_2^*}^{b_3^*-1}\Res_{T=q^{-k}}R_q(t)
=\Res_{T=\infty}R_q(t)
=0;
$$
поэтому
\begin{align*}
G_q(\ba,\bb)
&=\sum_{t=0}^\infty\biggl(
\sum_{k=a_3^*}^{b_2^*-1}\frac{A_kq^t}{(1-q^{t+k})^2}
+\sum_{k=a_2^*}^{b_3^*-1}\frac{B_kq^t}{1-q^{t+k}}\biggr)
\\
&=\sum_{k=a_3^*}^{b_2^*-1}A_kq^{-k}
\sum_{t=0}^\infty\frac{q^{t+k}}{(1-q^{t+k})^2}
+\sum_{k=a_3^*}^{b_2^*-1}B_kq^{-k}
\sum_{t=0}^\infty\frac{q^{t+k}}{1-q^{t+k}}
\\
&=\sum_{k=a_3^*}^{b_2^*-1}A_kq^{-k}
\biggl(\sum_{l=1}^\infty-\sum_{l=1}^{k-1}\biggr)\frac{q^l}{(1-q^l)^2}
+\sum_{k=a_2^*}^{b_3^*-1}B_kq^{-k}
\biggl(\sum_{l=1}^\infty-\sum_{l=1}^{k-1}\biggr)\frac{q^l}{1-q^l}
\\
&=A\zeta_q(2)-B,
\end{align*}
где
\begin{gather}
A=\sum_{k=a_3^*}^{b_2^*-1}A_kq^{-k},
\label{eq:3.62}
\\
B=\sum_{k=a_3^*}^{b_2^*-1}A_kq^{-k}
\sum_{l=1}^{k-1}\frac{q^l}{(1-q^l)^2}
+\sum_{k=a_2^*}^{b_3^*-1}B_kq^{-k}
\sum_{l=1}^{k-1}\frac{q^l}{1-q^l}
\nonumber
\end{gather}
--- рациональные функции переменной~$q$.
Для коэффициентов $A_k$, $a_3^*\le k<\nobreak b_2^*$, с помощью
представления~\eqref{eq:3.60} находим явные формулы
(напомним, что $b_1=1$):
\begin{align}
A_k
&=R_q(t)(1-q^kT)^2\big|_{T=q^{-k}}
=R_q(t)(1-q^{t+k})^2\big|_{t=-k}
\nonumber\\
&=(-1)^{a_1-b_1}q^{(a_1-b_1)(a_1+b_1-2k-1)/2}\qbinom{k-b_1}{a_1-b_1}
\nonumber\\ &\qquad\times 
(-1)^{k-a_2}q^{(k-a_2)(k-a_2+1)/2}\qbinom{b_2-a_2-1}{k-a_2}
\nonumber\\ &\qquad\times 
(-1)^{k-a_3}q^{(k-a_3)(k-a_3+1)/2}\qbinom{b_3-a_3-1}{k-a_3}
\cdot
q^{-k(b_2+b_3-a_1-a_2-a_3-1)}
\nonumber\\
&=(-1)^{a_1+a_2+a_3-1}
q^{a_1(a_1-1)/2+a_2(a_2-1)/2+a_3(a_3-1)/2-k(b_2+b_3-3)+k^2}
\nonumber\\ &\qquad\times
\qbinom{k-b_1}{a_1-b_1}
\qbinom{b_2-a_2-1}{k-a_2}\qbinom{b_3-a_3-1}{k-a_3},
\qquad a_3^*\le k<b_2^*.
\label{eq:3.62a}
\end{align}
Функция $k^2-k(b_2+b_3-2)$ убывает при изменении~$k$
от~$a_3^*$ до $b_2^*-1=\min\{b_2,b_3\}-1$ и принимает
в указанном промежутке наименьшее значение
$-(b_2-1)(b_3-1)$ при $k=b_2^*-1$. Кроме того,
\begin{align*}
\biggl|\qbinom{k-b_1}{a_1-b_1}\biggr|
&=\biggl|\frac{(1-q^{k-a_1+1})\dotsb(1-q^{k-b_1})}
{(1-q)(1-q^2)\dotsb(1-q^{a_1-b_1})}\biggr|
\le\biggl(\frac{1+|q|}{1-|q|}\biggr)^{a_1-b_1},
\\
\biggl|\qbinom{b_j-a_j-1}{k-a_j}\biggr|
&=\biggl|\frac{(1-q)(1-q^2)\dotsb(1-q^{b_j-a_j-1})}
{(1-q)\dotsb(1-q^{b_j-k-1})\cdot(1-q)\dotsb(1-q^{k-a_j})}\biggr|
\\
&\le\biggl(\frac{1+|q|}{1-|q|}\biggr)^{b_j-a_j-1},
\qquad j=2,3.
\end{align*}
Следовательно,
\begin{align*}
|A_kq^{-k}|
&\le\biggl(\frac{1+|q|}{1-|q|}\biggr)^{a_1-a_2-a_3+b_2+b_3-3}
\\ &\qquad\times 
|q|^{a_1(a_1-1)/2+a_2(a_2-1)/2+a_3(a_3-1)/2-(b_2-1)(b_3-1)},
\qquad a_3^*\le k<b_2^*.
\end{align*}
С учетом неравенств $a_1<b_2$ и
$$
\frac{1+|q|}{1-|q|}\le3
$$
при $|q|\le1/2$, а также того факта, что
в суммировании~\eqref{eq:3.62} участвует не более $b_3<3^{b_3}$
слагаемых, мы окончательно приходим к требуемой оценке~\eqref{eq:3.59}.

Записывая формулу \eqref{eq:3.62a} в виде
$$
A_k=(-1)^{a_1+a_2+a_3-1}x^{\mu(k)}
\xbinom{k-b_1}{a_1-b_1}\xbinom{b_2-a_2-1}{k-a_2}\xbinom{b_3-a_3-1}{k-a_3},
$$
где
\begin{align*}
\mu(k)
&=k^2-k(a_1+a_2+a_3)+\frac{a_1(a_1-1)}2-\frac{a_2(a_2-1)}2-\frac{a_3(a_3-1)}2
\\ &\quad
+a_2(b_2-1)+a_3(b_3-1)
\\
&\ge\mu\biggl(\frac{a_1+a_2+a_3}2\biggr)
\ge\mu\biggl(\frac{a_1+a_2+a_3}2\biggr)-c_{00}
=M_0(\ba,\bb),
\end{align*}
и пользуясь \eqref{eq:3.62}, заключаем, что $x^{-M_0}A\in\mathbb Z[x]$.
Аналогичные рассуждения показывают, что $x^{-M_0}B\in\mathbb Z[x]$. Поскольку
величина $M$ определена в~\eqref{eq:3.34} как максимально возможное целое,
для которого $x^{-M_0}A\in\mathbb Z[x]$ и $x^{-M_0}B\in\mathbb Z[x]$, заключаем,
что $M\ge M_0$.

Предложение доказано полностью.
\end{proof}

Формула
\begin{align*}
M_0(\bc)+N(\bc)
&=\frac34(a_1^2+a_2^2+a_3^2)+\frac12(a_1a_2+a_2a_3+a_3a_1)
\\ &\quad
-(a_1+a_2+a_3+1)(b_2+b_3-2)
\\ &\quad
+b_2^2+b_3^2+b_2b_3-3(b_2+b_3)+4
\displaybreak[2]\\
&=\frac34(c_{00}^2+c_{21}^2+c_{22}^2+c_{33}^2+c_{31}^2)
\\ &\quad
+\frac12(c_{00}c_{21}+c_{21}c_{22}+c_{22}c_{33}+c_{33}c_{31}+c_{31}c_{00})
\\ &\quad
-\frac12(c_{00}c_{22}+c_{21}c_{33}+c_{22}c_{31}+c_{33}c_{00}+c_{31}c_{21})
\\ &\quad
-\frac12(c_{00}+c_{21}+c_{22}+c_{33}+c_{31})-\frac14
\end{align*}
показывает, что выражение $M_0(\bc)+N(\bc)$ инвариантно относительно действия группы~$\fG$
(ср.\ с леммой~\ref{lem:3.7}).

\section({Мера иррациональности \$\003\266\_q(2)\$})%
{Мера иррациональности $\zeta_q(2)$}
\label{sec:3.6}

Зафиксируем теперь набор целочисленных параметров
({\it направлений\/}) $(\balpha,\bbeta)$,
удовлетворяющих условиям
\begin{equation}
\{\beta_1=0\}
\le\{\alpha_1,\alpha_2,\alpha_3\}
\le\{\beta_2,\beta_3\},
\quad
\alpha_1+\alpha_2+\alpha_3
\le\beta_1+\beta_2+\beta_3,
\label{eq:3.63}
\end{equation}
и для каждого $n=0,1,2,\dots$, свяжем их с исходным
набором параметров~\eqref{eq:3.7} по правилу
\begin{equation}
\begin{alignedat}{3}
a_1&=\alpha_1n+1, \quad&
a_2&=\alpha_2n+1, \quad&
a_3&=\alpha_3n+1,
\\
b_1&=\beta_1n+1, \quad&
b_2&=\beta_2n+2, \quad&
b_3&=\beta_3n+2.
\end{alignedat}
\label{eq:3.64}
\end{equation}
Определяя множество дополнительных параметров~$\bc$ равенствами
\begin{equation}
\begin{aligned}
c_{00}&=(\beta_1+\beta_2+\beta_3)-(\alpha_1+\alpha_2+\alpha_3),
\\
c_{jk}&=\begin{cases}
\alpha_j-\beta_k & \text{для $k=1$}, \\
\beta_k-\alpha_j & \text{для $k=2,3$},
\end{cases}
\qquad j,k=1,2,3;
\end{aligned}
\label{eq:3.65}
\end{equation}
мы получаем, что $10$-элементное множество $\bc\cdot n$
отвечает набору параметров~\eqref{eq:3.64} в соответствии с~\eqref{eq:3.14}.
С набором~\eqref{eq:3.65} свяжем определенные ранее характеристики
$m(\bc)$, $m_1^*(\bc)$, $m_2^*(\bc)$, $M(\bc)$ и рассмотрим
функцию
\begin{align}
\varphi(x)
&=\max_{0\le k\le11}\bigl(
\[c_{00}\cdot x\]+\[c_{21}\cdot x\]+\[c_{22}\cdot x\]
+\[c_{33}\cdot x\]+\[c_{31}\cdot x\]
\nonumber\\ &\qquad
-\[\fg_kc_{00}\cdot x\]-\[\fg_kc_{21}\cdot x\]-\[\fg_kc_{22}\cdot x\]
-\[\fg_kc_{33}\cdot x\]-\[\fg_kc_{31}\cdot x\]\bigr),
\label{eq:3.66}
\end{align}
где представители $\fg_k$, $k=0,1,\dots,11$, левых смежных классов
факторгруппы $\fG/\fG_0$ и их действие
на параметры $c_{00},c_{21},c_{22},c_{33},c_{31}$
указаны в \eqref{eq:3.48}, \eqref{eq:3.49}. Отметим, что ввиду
$\fG$-инвариантности характеристики $m(\bc)$ функция \eqref{eq:3.66}
является $1$-периодической.

\begin{proposition} \label{prop:3.5}
В приведенных выше обозначениях положим
$$
C_0=\mu-\frac3{\pi^2}\biggl({m_1^*}^2+{m_2^*}^2
+\int_0^1\varphi(x)\,\d\psi'(x)\biggr),
\qquad
C_1=\beta_2\beta_3-\frac{\alpha_1^2+\alpha_2^2+\alpha_3^2}2,
$$
где
$$
\mu=\alpha_2\beta_2+\alpha_3\beta_3+\frac{\alpha_1^2-3\alpha_2^2-3\alpha_3^2}4-\frac{\alpha_1\alpha_2+\alpha_2\alpha_3+\alpha_3\alpha_1}2.
$$
Если $C_0>0$, то число $\zeta_q(2)$ иррационально
для любого $q=1/p$, $p\in\mathbb Z\setminus\{0,\pm1\}$,
и справедлива оценка
\begin{equation}
\mu(\zeta_q(2))\le\frac{C_1}{C_0}
\label{eq:3.67}
\end{equation}
для показателя иррациональности.
\end{proposition}

\begin{proof}
Пусть $q^{-1}=p\in\mathbb Z\setminus\{0,\pm1\}$.
Для заданного набора направлений $(\balpha,\bbeta)$
и соответствующего набора~\eqref{eq:3.65} рассмотрим
последовательности
$$
\begin{gathered}
H_n=H(\bc n), \qquad
L_n=p^{-M_0(\bc n)}\cdot D_{m_1^*n}(p)D_{m_2^*n}(p)
\cdot\prod_{l=1}^{m_1^*n}\Phi_l^{-\omega_0(n/l)}(p),
\\
n=0,1,2,\dots\,.
\end{gathered}
$$
Поскольку
$m_1^*(\bc n)=m_1^*n$,
$m_2^*(\bc n)=m_2^*n$ и $\nu_l=\omega_0(n/l)$ согласно~\eqref{eq:3.47},
$n=0,1,2,\dots$, предложение~\ref{prop:3.2} влечет включения
$$
\wt H_n=L_nH_n\in\mathbb Z[p]\zeta_q(2)+\mathbb Z[p]
\subset\mathbb Z\zeta_q(2)+\mathbb Z,
\qquad n=0,1,2,\dots\,.
$$
С другой стороны, записывая линейные формы~$H_n$
в виде $H_n=A_n\zeta_q(2)-B_n$, $n=0,1,2,\dots$,
и применяя предложения~\ref{prop:3.3},~\ref{prop:3.4}, а также
$M_0(\bc n)=\mu n^2+o(n)$, находим
\begin{equation}
\lim_{n\to\infty}\frac{\log|H_n|}{n^2}=0,
\qquad
\varlimsup_{n\to\infty}\frac{\log|A_n|}{n^2}\le C_1\log|p|,
\label{eq:3.68}
\end{equation}
а асимптотическое поведение последовательности~$L_n$,
$n=0,1,2,\dots$, определяется с помощью лемм~\ref{lem:3.1},~\ref{lem:3.2}:
\begin{equation}
\lim_{n\to\infty}\frac{\log|L_n|}{n^2}=-C_0\log|p|.
\label{eq:3.69}
\end{equation}
Поэтому в случае $C_0>0$ иррациональность числа~$\zeta_q(2)$
вытекает из оценок
$$
0<|\wt H_n|<|p|^{-(C_0-\varepsilon)n^2},
$$
справедливых для всех $n\ge n_0(\varepsilon)$,
где в качестве~$\varepsilon>0$
можно взять $C_0/2$. Неравенство~\eqref{eq:3.67} выводится
из предельных соотношений~\eqref{eq:3.68},~\eqref{eq:3.69}
стандартным образом (см., например, \cite[лемма~3.1]{Ha1}
или \cite[лемма~2]{Da}). Это завершает
доказательство предложения.
\end{proof}

\begin{proof}[Доказательство теоремы~\ref{th:3}]
Реализуя перебор по всем целочисленным направлениям~$(\balpha,\bbeta)$,
удовлетворяющих условиям~\eqref{eq:3.63} и $\beta_2+\beta_3\le100$,
с помощью программы для калькулятора \texttt{GP-PARI}, мы обнаружили,
что наилучшая оценка~\eqref{eq:0.35}
для показателя иррациональности $\zeta_q(2)$
достигается (с точностью до действия группы~$\fG$ и умножения
вектора направлений на положительное целое) на наборе
$$
\alpha_1=10, \quad \alpha_2=11, \quad \alpha_3=12, \qquad
\beta_2=24, \quad \beta_3=25.
$$
В этом случае $\mu=837/4$, $m_1^*=16$, $m_2^*=15$,
$$
\varphi(x)=\begin{cases}
2, &\text{если $x\in\bigl[\frac3{16},\frac15\bigr)
\cup\bigl[\frac5{13},\frac25\bigr)
\cup\bigl[\frac6{13},\frac7{15}\bigr)
\cup\bigl[\frac9{16},\frac47\bigr)
\cup\bigl[\frac7{11},\frac9{14}\bigr)$},
\\
& \kern18mm
\text{$\cup\bigl[\frac9{13},\frac7{10}\bigr)
\cup\bigl[\frac{10}{13},\frac{11}{14}\bigr)
\cup\bigl[\frac{11}{13},\frac67\bigr)
\cup\bigl[\frac{12}{13},\frac{13}{14}\bigr)$},
\\
1, &\text{если $x\in\bigl[\frac1{16},\frac1{14}\bigr)
\cup\bigl[\frac1{13},\frac1{10}\bigr)
\cup\bigl[\frac18,\frac17\bigr)
\cup\bigl[\frac2{13},\frac3{16}\bigr)
\cup\bigl[\frac15,\frac3{14}\bigr)$}
\\
& \kern18mm
\text{$\cup\bigl[\frac3{13},\frac27\bigr)
\cup\bigl[\frac4{13},\frac5{14}\bigr)
\cup\bigl[\frac38,\frac5{13}\bigr)
\cup\bigl[\frac25,\frac37\bigr)
\cup\bigl[\frac5{11},\frac6{13}\bigr)$}
\\
& \kern18mm
\text{$\cup\bigl[\frac7{15},\frac12\bigr)
\cup\bigl[\frac7{13},\frac9{16}\bigr)
\cup\bigl[\frac47,\frac35\bigr)
\cup\bigl[\frac8{13},\frac7{11}\bigr)
\cup\bigl[\frac9{14},\frac23\bigr)$}
\\
& \kern18mm
\text{$\cup\bigl[\frac{11}{16},\frac9{13}\bigr)
\cup\bigl[\frac7{10},\frac57\bigr)
\cup\bigl[\frac8{11},\frac{11}{15}\bigr)
\cup\bigl[\frac34,\frac{10}{13}\bigr)
\cup\bigl[\frac{11}{14},\frac45\bigr)$}
\\
& \kern18mm
\text{$\cup\bigl[\frac9{11},\frac{11}{13}\bigr)
\cup\bigl[\frac67,\frac{13}{15}\bigr)
\cup\bigl[\frac{10}{11},\frac{12}{13}\bigr)
\cup\bigl[\frac{13}{14},\frac{14}{15}\bigr)$},
\\
0 &\text{в остальных случаях},
\end{cases}
$$
для $x\in[0,1)$. Следовательно,
\begin{align*}
C_0&=\frac{837}4-\frac3{\pi^2}(16^2+15^2-182.92436375\dots)
=118.64587116\dots,
\\
C_1&=14\cdot15-\frac{5^2+6^2+7^2}2=\frac{835}2,
\end{align*}
и согласно предложению~\ref{prop:3.5} мы получаем
требуемую оценку~\eqref{eq:0.35}.
Теорема доказана.
\end{proof}

\section({\$q\$-Аналог последовательности Апери})%
{$q$-Аналог последовательности Апери}
\label{sec:3.7}

Выбор параметров
\begin{equation}
a_1=a_2=a_3=n+1, \quad
b_1=1, \quad b_2=b_3=2n+2,
\qquad\text{где}\quad n=0,1,2,\dots,
\label{eq:3.70}
\end{equation}
приводит к величинам $C_0=5/4-6/\pi^2>0$, $C_1=5/2$
в обозначениях предложения~\ref{prop:3.5}, а значит, к иррациональности
числа~$\zeta_q(2)$ для $q^{-1}\in\mathbb Z\setminus\{0,\pm1\}$.
Соответствующая оценка для показателя иррациональности
в этом случае имеет вид
$$
\mu(\zeta_q(2))
\le\frac{10\pi^2}{5\pi^2-24}=3.89363887\dots\,.
$$
Цель этого параграфа --- показать, что случай~\eqref{eq:3.70}
является точным $q$-аналогом оригинального доказательства
Апери~\cite{Ap} иррациональности~$\zeta(2)$.

Зафиксируем целое $n\ge0$ и
запишем разложение рациональной функции~\eqref{eq:3.51}
в сумму простейших дробей относительно параметра~$T=q^t$:
\begin{align}
R_q(t)
&=\frac{(1-qT)\dotsb(1-q^nT)}{(1-q^{n+1}T)\dotsb(1-q^{2n+1}T)}
\cdot\frac{(q;q)_nT^n}{(1-q^{n+1}T)\dotsb(1-q^{2n+1}T)}
\nonumber
\displaybreak[2]\\
&=(-1)^n\sum_{k=0}^n\qbinom{k+n}k\qbinom nk
\frac{(-1)^kq^{k(k+1)/2-kn-n(n+1)/2}}{1-q^{k+n+1}T}
\nonumber\\ &\qquad\times 
\sum_{j=0}^n\qbinom nj
\frac{(-1)^jq^{j(j+1)/2-jn-n(n+1)}}{1-q^{j+n+1}T}
\nonumber
\displaybreak[2]\\
&=(-1)^n\sum_{k=0}^n\qbinom{k+n}k\qbinom nk^2
\frac{q^{k(k+1)-2kn-3n(n+1)/2}}{(1-q^{t+k+n+1})^2}
\nonumber\\ &\quad
+(-1)^n\sum_{k=0}^n\sum\doublesb{j=0}{j\ne k}^n
\qbinom{k+n}k\qbinom nk\qbinom nj
q^{k(k+1)/2+j(j+1)/2-(k+j)n-3n(n+1)/2}
\nonumber\\ &\quad\qquad\times
\frac{(-1)^{k+j}}{q^k-q^j}
\biggl(\frac1{1-q^{t+k+n+1}}-\frac1{1-q^{t+j+n+1}}\biggr).
\label{eq:3.71}
\end{align}
Учитывая равенства
$$
R_q(t)=0 \qquad\text{для}\quad t=-1,-2,\dots,-n
$$
и
$$
\begin{aligned}
\sum_{t=-n}^\infty\frac{q^t}{(1-q^{t+k+n+1})^2}
&=q^{-(k+n+1)}\biggl(\zeta_q(2)-\sum_{l=1}^k\frac{q^l}{(1-q^l)^2}\biggr),
\\
\sum_{t=-n}^\infty\frac{q^t}{1-q^{t+k+n+1}}
&=q^{-(k+n+1)}\biggl(\zeta_q(1)-\sum_{l=1}^k\frac{q^l}{1-q^l}\biggr),
\end{aligned}
\qquad k=0,1,\dots,n,
$$
из~\eqref{eq:3.71} мы получаем линейную форму
\begin{align}
H_n(q)
&=(-1)^nq^{(3n+2)(n+1)/2}\sum_{t=0}^\infty R_q(t)q^t
=(-1)^nq^{(3n+2)(n+1)/2}\sum_{t=-n}^\infty R_q(t)q^t
\nonumber\\
&=A_n(q)\zeta_q(2)-B_n(q)
\label{eq:3.72}
\end{align}
(коэффициент при $\zeta_q(1)$ в~\eqref{eq:3.72}
равен~$0$ согласно предложению~\ref{prop:3.1}), где
\begin{align*}
A_n(q)
&=\sum_{k=0}^n\qbinom{k+n}k\qbinom nk^2q^{k^2-2kn},
\\
B_n(q)
&=\sum_{k=0}^n\qbinom{k+n}k\qbinom nk^2q^{k^2-2kn}
\sum_{l=1}^k\frac{q^l}{(1-q^l)^2}
\\ &\quad
+\sum_{k=0}^n\sum\doublesb{j=0}{j\ne k}^n
\qbinom{k+n}k\qbinom nk\qbinom nj
q^{k(k+1)/2+j(j+1)/2-(k+j)n}
\\ &\quad\qquad\times
\frac{(-1)^{k+j}}{q^k-q^j}
\biggl(q^{-k}\sum_{l=1}^k\frac{q^l}{1-q^l}
-q^{-j}\sum_{l=1}^j\frac{q^l}{1-q^l}\biggr).
\end{align*}
Осуществляя теперь предельный переход $q\to1$,
имеем
\begin{align*}
u_n'
&=\lim_{q\to1}A_n(q)
=\sum_{k=0}^n\binom{k+n}k\binom nk^2,
\\
v_n'
&=\lim_{q\to1}(1-q)^2B_n(q)
=\sum_{k=0}^n\binom{k+n}k\binom nk^2
\sum_{l=1}^k\frac1{l^2}
\\ &\phantom{=\lim_{q\to1}(1-q)^2B_n(q)}\quad
+\sum_{k=0}^n\sum\doublesb{j=0}{j\ne k}^n
\binom{k+n}k\binom nk\binom nj
\frac{(-1)^{k+j}}{j-k}
\biggl(\sum_{l=1}^k\frac1l
-\sum_{l=1}^j\frac1l\biggr).
\end{align*}
Остается заметить (см., например, \cite[\S\,4]{As1}),
что последовательность линейных форм
$$
u_n'\zeta(2)-v_n', \qquad n=0,1,2,\dots,
$$
есть в точности последовательность
Апери~\eqref{eq:0.12}--\eqref{eq:0.18}.

\chapter({Мера иррациональности \$\003\266(2)\$})%
{Мера иррациональности $\zeta(2)$}
\label{chap:z2}

Изложение доказательства теоремы~\ref{th:z2} в этой главе организовано следующим образом.
В~\S\,\ref{BSI} мы приводим некоторые вспомогательные аналитические и арифметические составляющие,
в то время как \S\,\ref{gc1} и~\S\,\ref{app} посвящены деталям первой гипергеометрической конструкции
рациональных приближений к~$\zeta(2)$.
В~\S\,\ref{sbailey} обсуждается тождество между двумя гипергеометрическими интегралами,
которое мотивирует вторую гипергеометрическую конструкцию приближений к~$\zeta(2)$; именно
эту конструкцию мы изучаем в~\S\,\ref{s1}. В заключительном параграфе, \S\,\ref{epi}, мы собираем
воедино полученную информацию и доказываем теорему~\ref{th:z2}, а также приводим комментарии
о гипергеометрических конструкциях, связанных с использованными в этой главе.
Основной результат этой главы опубликован в~\cite{Z100}.

\section{Прелюдия: вспомогательные леммы}
\label{BSI}

В этом параграфе приводятся вспомогательные результаты о разложении интегралов
типа Меллина--Барнса и о специальной арифметике целозначных многочленов.

\begin{lemma}
\label{barnes}
При $\ell=0,1,2,\dots$ выполнено
\begin{equation}
\frac1{2\pi i}\int_{1/2-i\infty}^{1/2+i\infty}\biggl(\frac\pi{\sin\pi t}\biggr)^2
\frac{(t-1)(t-2)\dotsb(t-\ell)}{\ell!}\,\d t
=\frac{(-1)^\ell}{\ell+1}.
\label{expp}
\end{equation}
\end{lemma}

\begin{proof}
Подынтегральное выражение представляется в виде
$$
\frac1{\ell!}\biggl(\frac\pi{\sin\pi t}\biggr)^2
\frac{\Gamma(t)}{\Gamma(t-\ell)}
=\frac{(-1)^\ell}{\ell!}\,\Gamma(t)^2\Gamma(1-t)\,\Gamma(1+\ell-t);
$$
вычисление \eqref{expp} соответствующего интеграла следует непосредственно из первой леммы Барнса \cite[\S\,4.2.1]{Sl}.
\end{proof}

\begin{lemma}
\label{nest}
При $k=1,2,\dots$ имеет место разложение
\begin{equation}
\frac1{2\pi i}\int_{1/2-i\infty}^{1/2+i\infty}\biggl(\frac\pi{\sin\pi t}\biggr)^2\frac{\d t}{t+k}
=\sum_{m=1}^\infty\frac1{(m+k)^2}
=\zeta(2)-\sum_{\ell=1}^k\frac1{\ell^2}.
\label{exp2}
\end{equation}
\end{lemma}

\begin{proof}
Поскольку
$$
\biggl(\frac\pi{\sin\pi t}\biggr)^2\d t
=\d(-\pi\,\cot\pi t),
$$
интегрирование по частям левой части~\eqref{exp2} преобразует интеграл в
$$
-\frac1{2\pi i}\int_{1/2-i\infty}^{1/2+i\infty}\frac{\pi\,\cot\pi t\,\d t}{(t+k)^2}.
$$
Рассматривая сначала в качестве контура интегрирования прямоугольник с вершинами
$(1/2\pm iN,N+1/2\pm iN)$, где $N>0$ целое, применяя теорему о вычетах подобно тому,
как это было сделано в доказательства леммы~\ref{lem:1.7} (см.\ также \cite[лемма~2.4]{Z8}),
и переходя к пределу $N\to\infty$, мы получаем требуемое разложение~\eqref{exp2}.
\end{proof}

\begin{remark}
Как и в случае леммы~\ref{lem:1.7}, формулировка и доказательство леммы~\ref{nest}
подсказаны леммой~2 работы Ю.~Нестеренко~\cite{Ne2}. По-существу приводимая нами лемма
является частным случаем леммы~2.4 из \cite{Z8}, в которой используется искусственное
ограничение на порядок роста правильной рациональной функции на бесконечности.
Указанное ограничение устраняется интегрированием по частям, как это сделано выше.

Нестеренко также указывает альтернативный способ доказательства леммы~\ref{nest},
без интегрирования по частям. Основную сложность (ср.\ с леммой~\ref{lem:1.7} и леммой~2.4 из \cite{Z8})
представляет оценка интеграла по вертикальному отрезку, проходящему через точку $N+1/2$.
Выделим в нем подотрезок длины $2\sqrt N$ с центром в этой точке: интеграл по этому
подотрезку стремится к нулю, так как длина существенно меньше знаменателя.
Интеграл по оставшейся части отрезка также стремится к нулю: здесь маленькой оказывается дробь с синусом.
Таким образом, к нулю стремится и весь интеграл при $N\to\infty$.
\end{remark}

Как и ранее, через $D_n$ обозначается наименьшее общее кратное чисел $1,2,\dots,n$.

\begin{lemma}
\label{iv1}
Для целых $b<a$ определим
$$
R(t)=R(a,b;t)=\frac{(t+b)(t+b+1)\dotsb(t+a-1)}{(a-b)!}.
$$
Тогда
$$
R(k)\in\mathbb Z, \quad D_{a-b}\cdot\frac{\d R(t)}{\d t}\bigg|_{t=k}\in\mathbb Z
\quad\text{и}\quad
D_{a-b}\cdot\frac{R(k)-R(\ell)}{k-\ell}\in\mathbb Z
$$
при всех $k,\ell\in\mathbb Z$, $\ell\ne k$.
\end{lemma}

\begin{proof}
Обозначим через $m=b-a$ степень многочлена $R(t)$.
Первые два семейства включений уже обсуждались нами в лемме~\ref{lem:1.1}.
Что касается последнего семейства, введем еще один многочлен
$P(t)=(R(t)-R(\ell))/(t-\ell)$ степени $m-1$. Поскольку $D_m\cdot1/(k-\ell)$
является целым при $k=\ell+1,\ell+2,\dots,\ell+m$, а также $R(k)-R(\ell)\in\mathbb Z$, мы можем заключить, что
$D_m\cdot P(k)\in\mathbb Z$ при тех значениях~$k$. Последнее означает, что многочлен $D_mP(t)$
степени $m-1$ принимает целые значения при $m$~последовательных значениях аргумента;
в соответствии с \cite[отд.~8, гл.~2, задача~87]{PS} этот многочлен целозначный.
\end{proof}

\begin{lemma}
\label{iv2}
Пусть $R(t)$ представляет собой произведение нескольких целозначных многочленов
$$
R_j(t)=R(a_j,b_j;t)=\frac{(t+b_j)(t+b_j+1)\dotsb(t+a_j-1)}{(a_j-b_j)!}, \quad\text{где}\; b_j<a_j;
$$
положим $m=\max_j\{a_j-b_j\}$. Тогда
\begin{equation}
R(k)\in\mathbb Z, \quad D_m\cdot\frac{\d R(t)}{\d t}\bigg|_{t=k}\in\mathbb Z
\quad\text{и}\quad
D_m\cdot\frac{R(k)-R(\ell)}{k-\ell}\in\mathbb Z
\label{iv-incl}
\end{equation}
при всех $k,\ell\in\mathbb Z$, $\ell\ne k$.
\end{lemma}

\begin{proof}
Достаточно доказать это утверждение для произведения двух многочленов $R(t)$ и $\wt R(t)$,
удовлетворяющих включениям~\eqref{iv-incl}, и применить принцип математической индукции по количеству множителей.
Имеем
\begin{align*}
\bigl(R(t)\wt R(t)\bigr)\big|_{t=k}
&=R(k)\wt R(k)
\in\mathbb Z,
\\
D_m\frac{\d(R(t)\wt R(t))}{\d t}\bigg|_{t=k}
&=R(k)\cdot D_m\frac{\d\wt R(t)}{\d t}\bigg|_{t=k}
+D_m\frac{\d R(t)}{\d t}\bigg|_{t=k}\cdot\wt R(k)
\in\mathbb Z
\\
D_m\frac{R(k)\wt R(k)-R(\ell)\wt R(\ell)}{k-\ell}
&=D_m\frac{R(k)-R(\ell)}{k-\ell}\cdot\wt R(k)
\\ &\qquad
+R(\ell)\cdot D_m\frac{\wt R(k)-\wt R(\ell)}{k-\ell}
\in\mathbb Z,
\end{align*}
откуда следует требуемый результат.
\end{proof}

\section{Первая гипергеометрическая конструкция}
\label{gc1}

Частный случай общей конструкции из этого параграфа рассматривался нами в~\cite[\S\,2]{Z23}.

Набору целочисленных параметров
\begin{equation*}
(\ba,\bb)=\biggl(\begin{matrix} a_1, \, a_2, \, a_3, \, a_4 \\ b_1, \, b_2, \, b_3, \, b_4 \end{matrix}\biggr),
\end{equation*}
удовлетворяющих условиям
\begin{equation}
\begin{gathered}
b_1,b_2,b_3\le a_1,a_2,a_3,a_4<b_4,
\\
d=(a_1+a_2+a_3+a_4)-(b_1+b_2+b_3+b_4)\ge0,
\end{gathered}
\label{cond1}
\end{equation}
поставим в соответствие рациональную функцию
\begin{align}
R(t)
&=R(\ba,\bb;t)
=\frac{(t+b_1)\dotsb(t+a_1-1)}{(a_1-b_1)!}
\cdot\frac{(t+b_2)\dotsb(t+a_2-1)}{(a_2-b_2)!}
\nonumber\\ &\phantom{=R(\ba,\bb;t)} \qquad\times
\frac{(t+b_3)\dotsb(t+a_3-1)}{(a_3-b_3)!}
\cdot\frac{(b_4-a_4-1)!}{(t+a_4)\dotsb(t+b_4-1)}
\label{eq:gc}
\\[-4mm]
&=\Pi(\ba,\bb)\cdot\frac{\Gamma(t+a_1)\,\Gamma(t+a_2)\,\Gamma(t+a_3)\,\Gamma(t+a_4)}
{\Gamma(t+b_1)\,\Gamma(t+b_2)\,\Gamma(t+b_3)\,\Gamma(t+b_4)},
\label{eq:P0}
\end{align}
где
$$
\Pi(\ba,\bb)
=\frac{(b_4-a_4-1)!}{(a_1-b_1)!\,(a_2-b_2)!\,(a_3-b_3)!}.
$$
Также введем упорядоченные версии $a_1^*\le a_2^*\le a_3^*\le a_4^*$ параметров $a_1,a_2,a_3,a_4$
и $b_1^*\le b_2^*\le b_3^*$ параметров $b_1,b_2,b_3$, так что наборы $\{a_1^*,a_2^*,a_3^*,a_4^*\}$ и $\{b_1^*,b_2^*,b_3^*\}$
совпадают с $\{a_1,a_2,a_3,a_4\}$ и $\{b_1,b_2,b_3\}$ соответственно.
В этих обозначениях рациональная функция $R(t)$ имеет полюсы при $t=-k$, где $k=a_4^*,a_4^*+1,\dots,b_4-1$,
нули при $t=-\ell$, где $\ell=b_1^*,b_1^*+1,\dots,a_3^*-1$, и нули кратности~2 при
$t=-\ell$, где $\ell=b_2^*,b_2^*+1,\dots,a_2^*-1$.

Разложение функции \eqref{eq:gc} в сумму простейших дробей имеет вид
\begin{equation}
R(t)=\sum_{k=a_4^*}^{b_4-1}\frac{C_k}{t+k}+P(t),
\label{eq:P1}
\end{equation}
где $P(t)$ --- многочлен степени $d$ (см.~\eqref{cond1}) и
\begin{align}
C_k
&=\bigl(R(t)(t+k)\bigr)|_{t=-k}
\nonumber\\
&=(-1)^{d+b_4+k}\binom{k-b_1}{k-a_1}\binom{k-b_2}{k-a_2}\binom{k-b_3}{k-a_3}\binom{b_4-a_4-1}{k-a_4}\in\mathbb Z
\label{eq:P2}
\end{align}
при $k=a_4^*,a_4^*+1,\dots,b_4-1$.

\begin{lemma}
\label{lem:ap}
Положим $c=\max\{a_1-b_1,a_2-b_2,a_3-b_3\}$. Тогда
$D_cP(t)$ является целозначным многочленом степени~$d$.
\end{lemma}

\begin{proof}
Запишем $R(t)=R_1(t)R_2(t)$, где
$$
R_1(t)=\frac{\prod_{j=b_1}^{a_1-1}(t+j)}{(a_1-b_1)!}
\cdot\frac{\prod_{j=b_2}^{a_2-1}(t+j)}{(a_2-b_2)!}
\cdot\frac{\prod_{j=b_3}^{a_3-1}(t+j)}{(a_3-b_3)!}
$$
--- произведение трех целозначных многочленов и
$$
R_2(t)=\frac{(b_4-a_4-1)!}{\prod_{j=a_4}^{b_4-1}(t+j)}
=\sum_{k=a_4}^{b_4-1}\frac{(-1)^{k-a_4}\binom{b_4-a_4-1}{k-a_4}}{t+k}.
$$

Из леммы~\ref{iv2} следует, что
\begin{equation}
\begin{gathered}
D_c\cdot\frac{\d R_1(t)}{\d t}\bigg|_{t=j}\in\mathbb Z
\quad\text{при}\; j\in\mathbb Z,
\\
D_c\cdot\frac{R_1(j)-R_1(m)}{j-m}\in\mathbb Z
\quad\text{при}\; j,m\in\mathbb Z, \; j\ne m.
\end{gathered}
\label{eq:P3}
\end{equation}

Кроме того, отметим явные формулы
\begin{align*}
C_k
&=R_1(-k)\cdot\bigl(R_2(t)(t+k)\bigr)\big|_{t=-k}
\\
&=R_1(-k)\cdot(-1)^{k-a_4}\binom{b_4-a_4-1}{k-a_4}
\quad\text{при}\; k\in\mathbb Z;
\end{align*}
выражения обращаются в нуль для значений $k$ за пределами диапазона $a_4^*\le k\le b_4-1$.

При $\ell\in\mathbb Z$ имеем
\begin{align*}
&
\frac{\d}{\d t}\bigl(R(t)(t+\ell)\bigr)\bigg|_{t=-\ell}
=\frac{\d}{\d t}\bigl(R_1(t)\cdot R_2(t)(t+\ell)\bigr)\bigg|_{t=-\ell}
\\ &\quad
=\frac{\d R_1(t)}{\d t}\bigg|_{t=-\ell}
\cdot\bigl(R_2(t)(t+\ell)\bigr)\big|_{t=-\ell}
+R_1(-\ell)\cdot\frac{\d}{\d t}\bigr(R_2(t)(t+\ell)\bigr)\bigg|_{t=-\ell}
\displaybreak[2]\\ &\quad
=\frac{\d R_1(t)}{\d t}\bigg|_{t=-\ell}\cdot(-1)^{\ell-a_4}\binom{b_4-a_4-1}{\ell-a_4}
\\ &\quad\qquad
+R_1(-\ell)\cdot\frac{\d}{\d t}\sum_{k=a_4}^{b_4-1}(-1)^{k-a_4}\binom{b_4-a_4-1}{k-a_4}
\biggl(1-\frac{-\ell+k}{t+k}\biggr)\bigg|_{t=-\ell}
\displaybreak[2]\\ &\quad
=\frac{\d R_1(t)}{\d t}\bigg|_{t=-\ell}\cdot(-1)^{\ell-a_4}\binom{b_4-a_4-1}{\ell-a_4}
+R_1(-\ell)\sum_{\substack{k=a_4\\k\ne\ell}}^{b_4-1}\frac{(-1)^{k-a_4}\binom{b_4-a_4-1}{k-a_4}}{-\ell+k}
\end{align*}
и
\begin{align*}
\frac{\d}{\d t}\biggl(\sum_{k=a_4^*}^{b_4-1}\frac{C_k}{t+k}\cdot(t+\ell)\biggr)\bigg|_{t=-\ell}
&=\frac{\d}{\d t}\biggl(\sum_{k=a_4}^{b_4-1}\frac{C_k}{t+k}\cdot(t+\ell)\biggr)\bigg|_{t=-\ell}
\\
&=\frac{\d}{\d t}\sum_{k=a_4}^{b_4-1}C_k\biggl(1-\frac{-\ell+k}{t+k}\biggr)\bigg|_{t=-\ell}
=\sum_{\substack{k=a_4\\k\ne\ell}}^{b_4-1}\frac{C_k}{-\ell+k}
\\
&=\sum_{\substack{k=a_4\\k\ne\ell}}^{b_4-1}\frac{R_1(-k)\cdot(-1)^{k-a_4}\binom{b_4-a_4-1}{k-a_4}}{-\ell+k}.
\end{align*}
Следовательно,
\begin{align*}
P(-\ell)
&=\frac{\d}{\d t}\bigl(P(t)(t+\ell)\bigr)\big|_{t=-\ell}
=\frac{\d}{\d t}\biggl(R(t)(t+\ell)
-\sum_{k=a_4^*}^{b_4-1}\frac{C_k}{t+k}\cdot(t+\ell)\biggr)\bigg|_{t=-\ell}
\\
&=\frac{\d R_1(t)}{\d t}\bigg|_{t=-\ell}\cdot(-1)^{\ell-a_4}\binom{b_4-a_4-1}{\ell-a_4}
\\ &\qquad
+\sum_{\substack{k=a_4\\k\ne\ell}}^{b_4-1}(-1)^{k-a_4}\binom{b_4-a_4-1}{k-a_4}\frac{R_1(-\ell)-R_1(-k)}{-\ell+k},
\end{align*}
что в соответствии с включениями \eqref{eq:P3} влечет требуемое
$D_cP(-\ell)\in\mathbb Z$ при всех $\ell\in\mathbb Z$.
\end{proof}

Наконец, определим величину
\begin{equation}
r(\ba,\bb)
=\frac{(-1)^d}{2\pi i}\int_{C-i\infty}^{C+i\infty}\biggl(\frac\pi{\sin\pi t}\biggr)^2R(\ba,\bb;t)\,\d t,
\label{eq:P4}
\end{equation}
где $C$ --- произвольное число из интервала $-a_2^*<C<1-b_2^*$. Определение не зависит от выбора $C$,
поскольку подынтегральное выражение не имеет особых точек в полосе $-a_2^*<\Re t<1-b_2^*$.

\begin{proposition}
\label{prop1}
Справедливо разложение
\begin{equation}
r(\ba,\bb)=q(\ba,\bb)\zeta(2)-p(\ba,\bb)
\qquad\text{с}\quad
q(\ba,\bb)\in\mathbb Z, \quad D_{c_1}D_{c_2}p(\ba,\bb)\in\mathbb Z,
\label{eq:P5}
\end{equation}
где
\begin{equation*}
c_1=\max\{a_1-b_1,a_2-b_2,a_3-b_3,b_4-a_2^*-1\}
\quad\text{и}\quad
c_2=\max\{d+1,b_4-a_2^*-1\}.
\end{equation*}
Кроме того, величина $r(\ba,\bb)/\Pi(\ba,\bb)$ инвариантна относительно любой перестановки
параметров $a_1,a_2,a_3,a_4$.
\end{proposition}

\begin{proof}
Полагаем $C=1/2-a_2^*$ в~\eqref{eq:P4} и записываем \eqref{eq:P1} в виде
\begin{equation*}
R(t)=\sum_{k=a_4^*}^{b_4-1}\frac{C_k}{t+k}+\sum_{\ell=0}^dA_\ell P_\ell(t+a_2^*),
\end{equation*}
где
\begin{equation*}
P_\ell(t)=\frac{(t-1)(t-2)\dotsb(t-\ell)}{\ell!}
\end{equation*}
и $D_cA_\ell\in\mathbb Z$ в соответствии с леммой~\ref{lem:ap}.
Применяя леммы~\ref{barnes} и \ref{nest}, получаем
\begin{align*}
r(\ba,\bb)
&=\frac{(-1)^d}{2\pi i}\int_{1/2-i\infty}^{1/2+i\infty}\biggl(\frac\pi{\sin\pi t}\biggr)^2R(t-a_2^*)\,\d t
\\
&=\zeta(2)\cdot(-1)^d\sum_{k=a_4^*}^{b_4-1}C_k
-(-1)^d\sum_{k=a_4^*}^{b_4-1}C_k\sum_{\ell=1}^{k-a_2^*}\frac1{\ell^2}
+\sum_{\ell=0}^d\frac{(-1)^{d+\ell}A_\ell}{\ell+1}.
\end{align*}
Это представление означает, что $r(\ba,\bb)$ действительно имеет указанную форму~\eqref{eq:P5}.
Инвариантность $r(\ba,\bb)/\Pi(\ba,\bb)$ относительно перестановок
параметров $a_1,a_2,a_3,a_4$ непосредственно вытекает из~\eqref{eq:P0} и определения \eqref{eq:P4} величины $r(\ba,\bb)$.
\end{proof}

\section{Арифметика и асимптотика линейных форм}
\label{app}

В частном случае
\begin{equation}
\begin{alignedat}{4}
a_1&=7n+1, &\quad a_2&=6n+1, &\quad a_3&=5n+1, &\quad a_4&=\phantom08n+1,
\\
b_1&=1, &\quad b_2&=\phantom1n+1, &\quad b_3&=2n+1, &\quad b_4&=14n+2
\end{alignedat}
\label{P-ex}
\end{equation}
согласно предложению~\ref{prop1} получаем
\begin{equation}
r_n=r(\ba,\bb)=q_n\zeta(2)-p_n,
\qquad\text{где}\quad
q_n, \, D_{9n}D_{8n}p_n\in\mathbb Z.
\label{P-fin}
\end{equation}

Асимптотическое поведение величин $r_n$ и $q_n$ в общем случае
\begin{equation}
\begin{alignedat}{4}
a_1&=\alpha_1n+1, &\quad a_2&=\alpha_2n+1, &\quad a_3&=\alpha_3n+1, &\quad a_4&=\alpha_4n+1,
\\
b_1&=\beta_1n+1, &\quad b_2&=\beta_2n+1, &\quad b_3&=\beta_3n+1, &\quad b_4&=\beta_4n+2,
\end{alignedat}
\label{P-gen}
\end{equation}
когда целые параметры $\alpha_j$ и $\beta_j$ удовлетворяют условиям
$$
\beta_1,\beta_2,\beta_3<\alpha_1,\alpha_2,\alpha_3,\alpha_4<\beta_4,
\quad
\alpha_1+\alpha_2+\alpha_3+\alpha_4>\beta_1+\beta_2+\beta_3+\beta_4
$$
(для обеспечения условий \eqref{cond1}),
получается стандартным образом.

\begin{lemma}
\label{lem:Pan}
Кубическое уравнение $\prod_{j=1}^4(\tau-\alpha_j)-\prod_{j=1}^4(\tau-\beta_j)=0$
имеет один вещественный корень $\tau_1$ и пару комплексно-сопряженных корней $\tau_0,\ol{\tau_0}$.
Тогда
$$
\limsup_{n\to\infty}\frac{\log|r_n|}n=\Re f_0(\tau_0)
\quad\text{и}\quad
\lim_{n\to\infty}\frac{\log|q_n|}n=\Re f_0(\tau_1),
$$
где
\begin{align*}
f_0(\tau)
&=\sum_{j=1}^4\bigl(\alpha_j\log(\tau-\alpha_j)-\beta_j\log(\tau-\beta_j)\bigr)
\\[-4.5pt] &\qquad
-\sum_{j=1}^3(\alpha_j-\beta_j)\log(\alpha_j-\beta_j)+(\beta_4-\alpha_4)\log(\beta_4-\alpha_4).
\end{align*}
\end{lemma}

В случае асимптотики $r_n$ подобное утверждение было нами показано в предложении~\ref{prop:1.3} гл.~\ref{chap:1}.
Доказательство леммы~\ref{lem:Pan} мало отличается от доказательства аналогичных результатов в~\cite{Z8}, \cite{Z15}, \cite{Z17}.
Альтернативный способ вычисления асимптотического поведения $r_n$ и $q_n$ основан
на применении теоремы Пуанкаре к явным рекуррентным соотношениям для обеих последовательностей;
мы затрагиваем вопрос получения последних для нашего частного выбора параметров~\eqref{P-ex}
в~\S\,\ref{sbailey}.

\medskip
Когда параметры выбраны в соответствии с~\eqref{P-ex}, мы получаем
\begin{equation}
\begin{aligned}
-\limsup_{n\to\infty}\frac{\log|r_n|}n=C_0&=15.88518998\dots,
\\ 
\lim_{n\to\infty}\frac{\log|q_n|}n=C_1&=23.22906071\dotsc.
\end{aligned}
\label{eq:Pan}
\end{equation}

\medskip
В общей ситуации \eqref{P-gen} величины $c_1$ и $c_2$ в предложении~\ref{prop1}
имеют вид $\gamma_1n$ и $\gamma_2n$, где целые $\gamma_1$ и $\gamma_2$ однозначно
определяются через данные $\alpha_j,\beta_j$, $j=1,\dots,4$;
для простоты упорядочим их: $\gamma_1\ge\gamma_2$. Напомним, что через
$\lf\,\cdot\,\rf$ обозначается целая часть вещественного числа.

\begin{lemma}
\label{lem:Par}
В приведенных обозначениях справедливы включения
\begin{equation}
\Phi_n^{-1}q_n,\,\Phi_n^{-1}D_{\gamma_1n}D_{\gamma_2n}p_n\in\mathbb Z,
\label{acorr}
\end{equation}
где $\Phi_n=\prod_{p\le\gamma_2n}p^{\varphi(n/p)}$ \textup(произведение по простым $p$\textup) и
\begin{align*}
\varphi(x)=\max_{\boldsymbol\alpha'=\sigma\boldsymbol\alpha:\sigma\in\mathfrak S_4}
&\biggl(\lf(\beta_4-\alpha_4)x\rf-\lf(\beta_4-\alpha_4')x\rf
\\[-2pt] &\qquad
-\sum_{j=1}^3\bigl(\lf(\alpha_j-\beta_j)x\rf-\lf(\alpha_j'-\beta_j)x\rf\bigr)\biggr);
\end{align*}
максимум здесь берется по всевозможным перестановкам $(\alpha_1',\alpha_2',\alpha_3',\alpha_4')$
набора $(\alpha_1,\alpha_2,\alpha_3,\alpha_4)$. Кроме того,
$$
\lim_{n\to\infty}\frac{\log\Phi_n}n
=\int_0^1\varphi(x)\,\d\psi(x)-\int_0^{1/\gamma_2}\varphi(x)\,\frac{\d x}{x^2},
$$
где $\psi(x)$ --- логарифмическая производная гамма-функции.
\end{lemma}

\begin{proof}
Арифметическая ``корректировка'' в~\eqref{acorr} использует теперь уже стандартное рассуждение,
основанное на инвариантности относительно перестановок из предложения~\ref{prop1};
мы подробно обсуждаем этот метод (даже в более сложной $q$-ситуации) в доказательстве предложения~\ref{prop:3.2}.
Функция $\varphi(x)$ ``считает'' максимум
$$
\varphi\biggl(\frac np\biggr)
=\max_{\sigma\in\mathfrak S_4}\ord_p\frac{\Pi(\ba,\bb)}{\Pi(\sigma\ba,\bb)}.
\qedhere
$$
\end{proof}

\begin{remark}
Существует альтернативный способ вычисления $\varphi(x)$ с помощью
\begin{align}
\varphi(x)=\min_{0\le y<1}\biggl(\sum_{j=1}^3
&
\bigl(\lf y-\beta_jx\rf-\lf y-\alpha_jx\rf-\lf(\alpha_j-\beta_j)x\rf\bigr)
\nonumber\\[-6pt] &\qquad
+\lf(\beta_4-\alpha_4)x\rf-\lf\beta_4x-y\rf-\lf y-\alpha_4x\rf\biggr),
\nonumber
\end{align}
хотя далеко не очевидно, что новое выражение представляет ту же самую функцию $\varphi(x)$,
что и в лемме~\ref{lem:Par}. Именно этой технологией мы пользовались
в доказательстве предложения~\ref{prop:1.2} (см.\ также \cite[\S\,4]{Z8}
и \cite[\S\,2]{Ne5}). Именно ей мы воспользуемся далее в~\S\,\ref{s1}.
\end{remark}

При выборе~\eqref{P-ex} мы получаем $\gamma_1=9$, $\gamma_2=8$ и
\begin{equation}
\varphi(x)=\begin{cases}
2 &\text{при $x\in\bigl[\frac16,\frac15\bigr)\cup\bigl[\frac14,\frac27\bigr)\cup\bigl[\frac12,\frac47\bigr)\cup\bigl[\frac56,\frac67\bigr)$}, \\
1 &\text{при $x\in\bigl[\frac18,\frac17\bigr)\cup\bigl[\frac15,\frac14\bigr)\cup\bigl[\frac27,\frac37\bigr)\cup\bigl[\frac47,\frac56\bigr)\cup\bigl[\frac67,\frac89\bigr)$}, \\
0 &\text{в остальных случаях},
\end{cases}
\label{eq:P-ar}
\end{equation}
так что
$$
\lim_{n\to\infty}\frac{\log\Phi_n}n=8.12793878\dotsc.
$$
Полагая
$$
C_2=\lim_{n\to\infty}\frac{\log(\Phi_n^{-1}D_{9n}D_{8n})}n=9+8-8.12793878\hdots=8.87206121\dots
$$
и пременяя лемму~2.1 из \cite{Ha1a}, мы приходим к следующей оценке показателя иррациональности $\zeta(2)$:
\begin{equation*}
\mu(\zeta(2))\le\frac{C_0+C_1}{C_0-C_2}=5.57728968\dotsc.
\end{equation*}
Разумеется, этот результат хуже оценки, установленной ранее в~\cite{RV2} Рином и Виолой.
Далее мы увидим, что включения \eqref{acorr} могут быть еще уточнены при нашем выборе
параметров~\eqref{P-ex}.

\section{Интерлюдия: гипергеометрический интеграл}
\label{sbailey}

\begin{proposition}
\label{prop:int}
Для каждого $n=0,1,2,\dots$ справедливо равенство
\begin{align}
&
\frac{(6n)!}{(7n)!\,(5n)!\,(3n)!}\,\frac1{2\pi i}\int_{-i\infty}^{i\infty}
\frac{\Gamma(7n+1+t)\,\Gamma(6n+1+t)\,\Gamma(5n+1+t)}{\Gamma(1+t)\,\Gamma(n+1+t)\,\Gamma(2n+1+t)}
\nonumber\\ & \phantom{\frac{(6n)!}{(7n)!\,(5n)!\,(3n)!}\,} \qquad\times
\frac{\Gamma(8n+1+t)}{\Gamma(14n+2+t)}\biggl(\frac\pi{\sin\pi t}\biggr)^2\d t
\nonumber\\ &\;
=\frac{(6n)!^2}{(9n)!\,(3n)!}\,
\frac1{2\pi i}\int_{-i\infty}^{i\infty}
\frac{\Gamma(11n+2+2t)\,\Gamma(3n+1+t)}{\Gamma(2n+2+2t)\,\Gamma(1+t)}
\nonumber\\ &\; \phantom{=\frac{(6n)!^2}{(9n)!\,(3n)!}\,} \qquad\times
\frac{\Gamma(4n+1+t)\,\Gamma(5n+1+t)}{\Gamma(10n+2+t)\,\Gamma(11n+2+t)}
\,\frac\pi{\sin2\pi t}\,\d t,
\label{P=T}
\end{align}
где пути интегрирования разделяют две группы полюсов подынтегральных выражений.
\end{proposition}

\begin{proof}
Применяя алгоритм Госпера--Цайльбергера созидательного телескопирования
к рациональным функциям
$$
R(t)=\frac{\prod_{j=1}^{7n}(t+j)}{(7n)!}\,\frac{\prod_{j=1}^{5n}(t+n+j)}{(5n)!}\,
\frac{\prod_{j=1}^{3n}(t+2n+j)}{(3n)!}\,\frac{(6n)!}{\prod_{j=1}^{6n+1}(t+8n+j)}
$$
и
$$
\hat R(t)
=\frac{\prod_{j=2}^{9n+1}(2t+2n+j)}{(9n)!}\,\frac{\prod_{j=1}^{3n}(t+j)}{(3n)!}\,
\frac{(6n)!}{\prod_{j=1}^{6n+1}(t+4n+j)}\,\frac{(6n)!}{\prod_{j=1}^{6n+1}(t+5n+j)}
$$
подобно тому, как это сделано в доказательстве теоремы~5.4 в~\cite{BBBC},
получаем, что интегралы
$$
r_n=\frac1{2\pi i}\int_{-i\infty}^{i\infty}R(t)\biggl(\frac\pi{\sin\pi t}\biggr)^2\d t
\quad\text{и}\quad
\hat r_n=\frac1{2\pi i}\int_{-i\infty}^{i\infty}\hat R(t)\,\frac\pi{\sin2\pi t}\,\d t
$$
удовлетворяют \emph{одному и тому же} рекуррентному уравнению
$$
s_0(n)r_{n+3}+s_1(n)r_{n+2}+s_2(n)r_{n+1}+s_3(n)r_n=0
\quad\text{при}\; n=0,1,2,\dots,
$$
где $s_0(n)$, $s_1(n)$, $s_2(n)$ и $s_3(n)$ --- многочлены от $n$ степени~64.
Проверяя равенство в~\eqref{P=T} непосредственно при $n=0,1,2$, заключаем, что оно имеет место при всех~$n$.
\end{proof}

Схожие приложения алгоритма созидательного телескопирования к доказательству
тождеств для интегралов типа Меллина--Барнса обсуждаются в \cite{Gui}, \cite{Sta}.

\begin{remark}
Отметим, что левая часть~\eqref{P=T} есть не что иное, как линейная форма из \S\,\ref{gc1},
отвечающая нашему частному выбору \eqref{P-ex} параметров.
Характеристический многочлен получающегося рекуррентного уравнения равен
\begin{align*}
&
2^2\,3^{12}\,7^{14}\,\lambda^3+3^3\,7^7\,794493690983053821271\,\lambda^2
-2^{20}\,3^4\,7^5\,2687491277\,\lambda+2^{48},
\end{align*}
и его нули определяют асимптотики \eqref{eq:Pan} последовательностей $r_n$ и $q_n$ с помощью
теоремы Пуанкаре.
\end{remark}

Для ``достаточно общего'' выбора \emph{целых} параметров ожидается верным следующее тождество:
\begin{align}
&
\frac1{2\pi i}\int_{-i\infty}^{i\infty}
\frac{\Gamma(a+t)\,\Gamma(b+t)\,\Gamma(e+t)\,\Gamma(f+t)}{\Gamma(1+t)\,\Gamma(1+a-e+t)\,\Gamma(1+a-f+t)\,\Gamma(g+t)}
\biggl(\frac\pi{\sin\pi t}\biggr)^2\d t
\nonumber\\ &\;
=(-1)^{a+b+e+f}\frac{\Gamma(e+f-a)\,\Gamma(e)\,\Gamma(f)}{\Gamma(g-b)}
\nonumber\\ &\;\quad\times
\frac1{2\pi i}\int_{-i\infty}^{i\infty}
\frac{\Gamma(a-b+g+2t)\,\Gamma(a+t)\,\Gamma(e+t)\,\Gamma(f+t)}{\Gamma(1+a+2t)\,\Gamma(1+a-b+t)\,\Gamma(e+f+t)\,\Gamma(g+t)}
\,\frac\pi{\sin2\pi t}\,\d t.
\label{bmiss}
\end{align}
Более того, предполагается верным и другое ``двойственное'' тождество, в котором
$(\pi/\sin\pi t)^2$ и $\pi/\cos2\pi t$ заменены соответственно на
$\pi^3\cos\pi t/\allowbreak(\sin\pi t)^3$ и $(\pi/\sin\pi t)^2$;
интегралы в двойственном тождестве представляют рациональные приближения к~$\zeta(3)$ (см.\ \cite{Z90}).
По всей видимости, эти равенства могут быть доказаны в подобной общности путем применения
так называемых смежных соотношений для интегралов в обеих частях;
этот метод, однако, не выглядит простым.

Предложение~\ref{prop:int} является частным случаем \eqref{bmiss}, когда
$$
a=8n+1, \quad b=5n+1, \quad e=6n+1, \quad f=7n+1, \quad g=14n+2.
$$

Тождество \eqref{bmiss} и двойственное к нему, судя по всему, являются
частными случаями еще более общего тождества между гипергеометрическими интегралами,
имеющего место для \emph{комплексных} параметров. Мы не смогли найти его в литературе,
хотя упоминание о нем содержится в конце статьи Бейли~\cite{Bai0}, опубликованной в~1932\,г.:
\begin{quote}
``Формула (1.4)\footnote{Мы приводим эту формулу далее под именем соотношения \eqref{whipple}.}
и следующая за ней чрезвычайно затруднительны для обобщений,
и окончательный результат оказался неожиданным.
Полученные формулы содержат пять гипергеометрических рядов вместо трех или четырех, как в предыдущих результатах.
В каждом из случаев два ряда являются почти уравновешенными и второго типа,
один --- почти уравновешенным и первого типа, а оставшиеся два --- сбалансированными.
При проведении этих исследований некоторые интегралы барнсова типа
вычисляются аналогично известным суммированиям гипергеометрических рядов.
Соображения места, однако, не позволяют этим результатам быть детально сформулированными.''\footnote%
{Оригинальный текст из~\cite{Bai0}:
``The formula (1.4) and its successor are rather more troublesome to generalize, and the final result was unexpected.
The formulae obtained involve five series instead of three or four as previously obtained.
In each case two of the series are nearly-poised and of the second kind, one is nearly-poised
and of the first kind, and the other two are Saalsch\"utzian in type.
In the course of these investigations some integrals of Barnes's type
are evaluated analogous to known sums of hypergeometric series. Considerations of space,
however, prevent these results being given in detail.''}
\end{quote}
Чем не в духе знаменитого ``Я нашел этому поистине чудесное доказательство, но поля книги слишком узки для него'' П.~Ферма?
Добавим еще, что в последнем абзаце гл.~6 своей книги \cite{Bai} Бейли опять упоминает
о затруднительном обобщении без каких-либо деталей. Было ли Бейли известно разыскиваемое нами общее тождество?

\section{Вторая гипергеометрическая конструкция}
\label{s1}

Материал предыдущего параграфа подсказывает другую конструкцию рациональных приближений к $\zeta(2)$.
На этот раз мы рассматриваем рациональную функцию
\begin{align}
\hat R(t)
&=\hat R(\hat\ba,\hat\bb;t)
=\frac{(2t+\hat b_0)(2t+\hat b_0+1)\dotsb(2t+\hat a_0-1)}{(\hat a_0-\hat b_0)!}
\cdot\frac{(t+\hat b_1)\dotsb(t+\hat a_1-1)}{(\hat a_1-\hat b_1)!}
\nonumber\\ &\phantom{=R(\ba,\bb;t)} \qquad\times
\frac{(\hat b_2-\hat a_2-1)!}{(t+\hat a_2)\dotsb(t+\hat b_2-1)}
\cdot\frac{(\hat b_3-\hat a_3-1)!}{(t+\hat a_3)\dotsb(t+\hat b_3-1)}
\nonumber\displaybreak[2]\\
&=\hat\Pi(\hat\ba,\hat\bb)\cdot\frac{\Gamma(2t+\hat a_0)\,\Gamma(t+\hat a_1)\,\Gamma(t+\hat a_2)\,\Gamma(t+\hat a_3)}
{\Gamma(2t+\hat b_0)\,\Gamma(t+\hat b_1)\,\Gamma(t+\hat b_2)\,\Gamma(t+\hat b_3)},
\nonumber
\end{align}
где
$$
\hat\Pi(\hat\ba,\hat\bb)
=\frac{(\hat b_2-\hat a_2-1)!\,(\hat b_3-\hat a_3-1)!}{(\hat a_0-\hat b_0)!\,(\hat a_1-\hat b_1)!}
$$
и целые параметры
\begin{equation*}
(\hat\ba,\hat\bb)
=\biggl(\begin{matrix} \hat a_0; \hat a_1, \, \hat a_2, \, \hat a_3 \\ \hat b_0; \, \hat b_1, \, \hat b_2, \, \hat b_3 \end{matrix}\biggr)
\end{equation*}
удовлетворяют требованиям
\begin{equation}
\begin{gathered}
\tfrac12\hat b_0,\hat b_1\le\tfrac12\hat a_0,\hat a_1,\hat a_2,\hat a_3<\hat b_2,\hat b_3,
\\
\hat a_0+\hat a_1+\hat a_2+\hat a_3=\hat b_0+\hat b_1+\hat b_2+\hat b_3-2.
\end{gathered}
\label{cond2}
\end{equation}
Последнее условие означает, что $\hat R(t)=O(1/t^2)$ при $t\to\infty$.
Хотя здесь будут менее важны детали арифметики линейных приближающих форм
в духе \S\,\ref{gc1}, удобно ввести упорядоченные версии
$\hat a_1^*\le\hat a_2^*\le\hat a_3^*$ и $\hat b_2^*\le\hat b_3^*$
параметров $\hat a_1,\hat a_2,\hat a_3$ и $\hat b_2,\hat b_3$ соответственно.
Такое упорядочение и условия~\eqref{cond2} позволяют заключить, что
рациональная функция $\hat R(t)$ имеет полюсы при $t=-k$, где $\hat a_2^*\le k\le\hat b_3^*-1$,
двойные полюсы при $t=-k$, где $\hat a_3^*\le k\le\hat b_2^*-1$, и нули при
$t=-\ell/2$, где $\hat b_0\le\ell\le\hat a_0^*-1$ в обозначении $\hat a_0^*=\min\{\hat a_0,2\hat a_2^*\}$.

Разложение правильной рациональной функции $\hat R(t)$ в сумму простейших дробей имеет вид
\begin{equation*}
\hat R(t)=\sum_{k=\hat a_3^*}^{\hat b_2^*-1}\frac{A_k}{(t+k)^2}+\sum_{k=\hat a_2^*}^{\hat b_3^*-1}\frac{B_k}{t+k},
\end{equation*}
где
\begin{align}
A_k
&=\bigl(\hat R(t)(t+k)^2\bigr)|_{t=-k}
\nonumber\\
&=(-1)^{\hat d}\binom{2k-\hat b_0}{2k-\hat a_0}\binom{k-\hat b_1}{k-\hat a_1}
\binom{\hat b_2-\hat a_2-1}{k-\hat a_2}\binom{\hat b_3-\hat a_3-1}{k-\hat a_3}\in\mathbb Z
\label{eq:T2}
\end{align}
с $\hat d=\hat b_2+\hat b_3$,
при $k=\hat a_3^*,\hat a_3^*+1,\dots,\hat b_2^*-1$, и аналогично
$$
B_k
=\frac{\d}{\d t}\bigl(\hat R(t)(t+k)^2\bigr)|_{t=-k}
$$
при $k=\hat a_2^*,\hat a_2^*+1,\dots,\hat b_3^*-1$.
Более того,
\begin{equation}
\sum_{k=\hat a_2^*}^{\hat b_3^*-1}B_k
=-\Res_{t=\infty}\hat R(t)=0
\label{eq:T4}
\end{equation}
согласно теореме о сумме вычетов.

Включения
\begin{equation}
D_{\max\{\hat a_0-\hat b_0,\hat a_1-\hat b_1,\hat b_3^*-\hat a_2-1,\hat b_3^*-\hat a_3-1\}}\cdot B_k\in\mathbb Z
\label{eq:T3}
\end{equation}
доказываются стандартным образом, как, например, это было сделано в доказательстве предложения~\ref{prop:1.1} гл.~\ref{chap:1}.
Более важными для нас являются оценки порядка вхождения простых~$p$:
\begin{align}
\ord_pA_k, \, 1+\ord_pB_k
&\ge\biggl\lf\frac{2k-\hat b_0}p\biggr\rf
-\biggl\lf\frac{2k-\hat a_0}p\biggr\rf
-\biggl\lf\frac{\hat a_0-\hat b_0}p\biggr\rf
\nonumber\\ &\quad
+\biggl\lf\frac{k-\hat b_1}p\biggr\rf
-\biggl\lf\frac{k-\hat a_1}p\biggr\rf
-\biggl\lf\frac{\hat a_1-\hat b_1}p\biggr\rf
\nonumber\\ &\quad
+\sum_{j=2}^3\biggl(\biggl\lf\frac{\hat b_j-\hat a_j-1}p\biggr\rf
-\biggl\lf\frac{k-\hat a_j}p\biggr\rf
-\biggl\lf\frac{\hat b_j-1-k}p\biggr\rf\biggr)
\label{eq:T3aaa}
\end{align}
при $k=\hat a_2^*,\dots,\hat b_3^*-1$, которые следуют из лемм~\ref{lem:1.3} и~\ref{lem:1.4}.

В этом параграфе мы рассматриваем величину
\begin{equation*}
\hat r(\hat\ba,\hat\bb)
=\frac{(-1)^{\hat d}}{2\pi i}\int_{C/2-i\infty}^{C/2+i\infty}\frac\pi{\sin2\pi t}\,\hat R(\hat\ba,\hat\bb;t)\,\d t,
\end{equation*}
где $C$ выбирается произвольным из интервала $-\hat a_0^*<C<1-\hat b_0$.

\begin{proposition}
\label{prop2}
Справедливо представление
\begin{equation}
\hat r(\hat\ba,\hat\bb)=\hat q(\hat\ba,\hat\bb)\zeta(2)-\hat p(\hat\ba,\hat\bb)
\qquad\text{с}\quad
\hat q(\hat\ba,\hat\bb)\in\mathbb Z, \quad D_{\hat c_1}D_{\hat c_2}\hat p(\hat\ba,\hat\bb)\in\mathbb Z,
\label{eq:T6}
\end{equation}
где
\begin{equation*}
\begin{gathered}
\hat c_1=\max\{\hat a_0-\hat b_0,\hat a_1-\hat b_1,\hat b_3^*-\hat a_2-1,\hat b_3^*-\hat a_3-1,2\hat b_2^*-\hat a_0^*-2\},
\\
\hat c_2=2\hat b_3^*-\hat a_0^*-2.
\end{gathered}
\end{equation*}
\end{proposition}

\begin{proof}
Воспользуемся
$$
\Res_{t=m/2}\frac{\pi\hat R(t)}{\sin2\pi t}
=\frac{(-1)^m}2\,\hat R(t)\big|_{t=m/2}
$$
при $m\ge1-\hat a_0^*$, чтобы записать
\begin{align*}
\hat r(\hat\ba,\hat\bb)
&=-\frac{(-1)^{\hat d}}2\sum_{m=1-\hat a_0^*}^\infty(-1)^m\hat R(t)\big|_{t=m/2}
\\
&=(-1)^{\hat d}\sum_{k=\hat a_3^*}^{\hat b_2^*-1}2A_k\sum_{m=1-\hat a_0^*}^\infty\frac{(-1)^{m-1}}{(m+2k)^2}
+(-1)^{\hat d}\sum_{k=\hat a_2^*}^{\hat b_3^*-1}B_k\sum_{m=1-\hat a_0^*}^\infty\frac{(-1)^{m-1}}{m+2k}
\displaybreak[2]\\
&=2\sum_{\ell=1}^\infty\frac{(-1)^{\ell-1}}{\ell^2}\cdot(-1)^{\hat d}\sum_{k=\hat a_3^*}^{\hat b_2^*-1}A_k
-(-1)^{\hat d}\sum_{k=\hat a_3^*}^{\hat b_2^*-1}2A_k\sum_{\ell=1}^{2k-\hat a_0^*}\frac{(-1)^{\ell-1}}{\ell^2}
\\ &\qquad
+\sum_{\ell=1}^\infty\frac{(-1)^{\ell-1}}{\ell}\cdot(-1)^{\hat d}\sum_{k=\hat a_2^*}^{\hat b_3^*-1}B_k
-(-1)^{\hat d}\sum_{k=\hat a_2^*}^{\hat b_3^*-1}B_k\sum_{\ell=1}^{2k-\hat a_0^*}\frac{(-1)^{\ell-1}}{\ell}
\displaybreak[2]\\
&=\zeta(2)\cdot(-1)^{\hat d}\sum_{k=\hat a_3^*}^{\hat b_2^*-1}A_k
-(-1)^{\hat d}\sum_{k=\hat a_3^*}^{\hat b_2^*-1}2A_k\sum_{\ell=1}^{2k-\hat a_0^*}\frac{(-1)^{\ell-1}}{\ell^2}
\\ &\qquad
-(-1)^{\hat d}\sum_{k=\hat a_2^*}^{\hat b_3^*-1}B_k\sum_{\ell=1}^{2k-\hat a_0^*}\frac{(-1)^{\ell-1}}{\ell},
\end{align*}
где мы применили равенство \eqref{eq:T4}.
В соответствии с включениями \eqref{eq:T2}, \eqref{eq:T3}
полученное представление величины $\hat r(\hat\ba,\hat\bb)$ действительно имеет заявленную форму~\eqref{eq:T6}.
\end{proof}

\begin{remark}
Биномиальные выражения \eqref{eq:P2} и \eqref{eq:T2} позволяют записать коэффициенты
$$
q(\ba,\bb)=(-1)^d\sum_{k=a_4^*}^{b_4-1}C_k
\quad\text{и}\quad
\hat q(\hat\ba,\hat\bb)=(-1)^{\hat d}\sum_{k=\hat a_3^*}^{\hat b_2^*-1}A_k
$$
в виде некоторых ${}_4F_3$- и ${}_5F_4$-гипергеометрических рядов соответственно.
Тогда классическое преобразование Уиппла \cite[с.~65, (2.4.2.3)]{Sl},
\begin{multline}
{}_4F_3\biggl(\begin{matrix} f, \, 1+f-h, \, h-a, \, -N \\
h, \, 1+f+a-h, \, g \end{matrix}\biggm|1\biggr)
=\frac{(g-f)_N}{(g)_N}
\\ \times
{}_5F_4\biggl(\begin{alignedat}{5}
a&, \, & -N&, \, & 1+f-g&, \, & \tfrac12f&, \, & \tfrac12f+\tfrac12 \\
&& h&, \, & 1+f+a-h&, \, & \tfrac12(1+f-N-g)&, \, & \tfrac12(1+f-N-g)+\tfrac12
\end{alignedat}\biggm|1\biggr),
\label{whipple}
\end{multline}
может быть сформулировано в виде следующего тождества:
\begin{equation*}
q\biggl(\begin{alignedat}{4} a_1&, &\, a_2&, &\, a_3&, &\, a_4& \\ 1&, &\, a_4-a_1&+1, &\, a_4-a_2&+1, &\, b_4& \end{alignedat}\biggr)
=\hat q\biggl(\begin{alignedat}{4}  b_4-a_3+a_4&; &\, a_2&, &\, a_1&, &\, a_4& \\ a_4+1&; &\, a_4-a_3&+1, &\, a_1\,+\,&a_2, &\, b_4& \end{alignedat}\biggr).
\end{equation*}
Отметим, что \eqref{bmiss} эквивалентно
\begin{equation*}
r\biggl(\begin{alignedat}{4} a_1&, &\, a_2&, &\, a_3&, &\, a_4& \\ 1&, &\, a_4-a_1&+1, &\, a_4-a_2&+1, &\, b_4& \end{alignedat}\biggr)
=\hat r\biggl(\begin{alignedat}{4}  b_4-a_3+a_4&; &\, a_2&, &\, a_1&, &\, a_4& \\ a_4+1&; &\, a_4-a_3&+1, &\, a_1\,+\,&a_2, &\, b_4& \end{alignedat}\biggr),
\end{equation*}
так что именно преобразование Уиппла \eqref{whipple} натолкнуло нас на мысль о совпадении
двух разных семейств линейных форм от $1$ и $\zeta(2)$.
\end{remark}

Как и в \S\,\ref{app}, мы делаем следующий выбор параметров:
\begin{equation}
\begin{alignedat}{4}
\hat a_0&=\hat\alpha_0n+2, &\quad \hat a_1&=\hat\alpha_1n+1, &\quad \hat a_2&=\hat\alpha_2n+1, &\quad \hat a_3&=\hat\alpha_3n+1,
\\
\hat b_0&=\hat\beta_0n+2, &\quad \hat b_1&=\hat\beta_1n+1, &\quad \hat b_2&=\hat\beta_2n+2, &\quad \hat b_3&=\hat\beta_3n+2,
\end{alignedat}
\label{T-gen}
\end{equation}
где целые $\hat\alpha_j$ и $\hat\beta_j$, $j=0,\dots,3$, удовлетворяют соотношениям
$$
\begin{gathered}
\tfrac12\hat\beta_0,\hat\beta_1<\tfrac12\hat\alpha_0,\hat\alpha_1,\hat\alpha_2,\hat\alpha_3<\hat\beta_2,\hat\beta_3,
\quad
\hat\alpha_0+\hat\alpha_1+\hat\alpha_2+\hat\alpha_3=\hat\beta_0+\hat\beta_1+\hat\beta_2+\hat\beta_3,
\end{gathered}
$$
обеспечивающим выполнение требований \eqref{cond2}. Для конструкции этого параграфа
несложно привести аналог леммы~\ref{lem:Pan}, однако нас прежде всего интересуют арифметические аспекты.

\begin{lemma}
\label{lem:Tar}
Предполагая параметры выбранными в соответствии с \eqref{T-gen}, для каждого простого $p$ имеем
$$
\ord_p\hat q(\hat\ba,\hat\bb)\ge\hat\varphi(n/p)
\quad\text{и}\quad
\ord_p\hat p(\hat\ba,\hat\bb)\ge-2+\hat\varphi(n/p),
$$
где \textup($1$-периодическая и целозначная\textup) функция $\hat\varphi(x)$ задается выражением
\begin{align*}
\hat\varphi(x)=\min_{0\le y<1}\biggl(
&
\lf 2y-\hat\beta_0x\rf-\lf 2y-\hat\alpha_0x\rf-\lf(\hat\alpha_0-\hat\beta_0)x\rf
\\[-4.5pt] &\;
+\lf y-\hat\beta_1x\rf-\lf y-\hat\alpha_1x\rf-\lf(\hat\alpha_1-\hat\beta_1)x\rf
\\[-3pt] &\;
+\sum_{j=2}^3\bigl(\lf(\hat\beta_j-\hat\alpha_j)x\rf-\lf\hat\beta_jx-y\rf-\lf y-\hat\alpha_jx\rf\bigr)\biggr).
\end{align*}
\end{lemma}

\begin{proof}
Оценки следуют из неравенств \eqref{eq:T3aaa}, явных выражений для $\hat q(\hat\ba,\hat\bb)$
и $\hat p(\hat\ba,\hat\bb)$, полученных в доказательстве предложения~\ref{prop2}: мы просто вводим обозначение
$y=(k-1)/p$ и минимизируем по всем возможным индексам~$k$.
\end{proof}

Укажем теперь, что в случае специального выбора параметров $(\hat\ba,\hat\bb)$,
\begin{equation*}
\begin{alignedat}{4}
\hat a_0&=11n+2, &\quad \hat a_1&=3n+1, &\quad \hat a_2&=\phantom04n+1, &\quad \hat a_3&=\phantom05n+1,
\\
\hat b_0&=\phantom02n+2, &\quad \hat b_1&=1, &\quad \hat b_2&=10n+2, &\quad \hat b_3&=11n+2,
\end{alignedat}
\end{equation*}
мы получаем линейные формы
\begin{equation*}
\hat r_n=\hat r(\hat\ba,\hat\bb)=\hat q_n\zeta(2)-\hat p_n,
\end{equation*}
совпадающие согласно предложению \ref{prop:int} с линейными формами \eqref{P-ex}, \eqref{P-fin}:
$$
r_n=q_n\zeta(2)-p_n=\hat q_n\zeta(2)-\hat p_n,
$$
откуда следует совпадение последовательностей коэффициентов, $q_n=\hat q_n$ и $p_n=\hat p_n$ для $n=0,1,2,\dots$\,.
Из этого совпадения и леммы~\ref{lem:Tar} заключаем, что
$$
\hat\Phi_n^{-1}q_n=\hat\Phi_n^{-1}\hat q_n\in\mathbb Z,
\quad
\hat\Phi_n^{-1}D_{9n}D_{8n}p_n=\hat\Phi_n^{-1}D_{9n}D_{8n}\hat p_n\in\mathbb Z,
$$
где $\hat\Phi_n=\prod_{p\le8n}p^{\hat\varphi(n/p)}$ и
\begin{align}
\hat\varphi(x)
&=\min_{0\le y<1}\bigl(\lf 2y-2x\rf-\lf 2y-11x\rf-\lf9x\rf
+\lf y\rf-\lf y-3x\rf-\lf3x\rf
\nonumber\\ &\quad
+\lf6x\rf-\lf10x-y\rf-\lf y-4x\rf
+\lf6x\rf-\lf11x-y\rf-\lf y-5x\rf\bigr)
\nonumber\\
&=\begin{cases}
2 &\text{при $x\in\bigl[\frac16,\frac29\bigr)\cup\bigl[\frac12,\frac59\bigr)\cup\bigl[\frac56,\frac78\bigr)$}, \\
1 &\text{при $x\in\bigl[\frac29,\frac49\bigr)\cup\bigl[\frac59,\frac79\bigr)\cup\bigl[\frac78,\frac89\bigr)$}, \\
0 &\text{в остальных случаях},
\end{cases}
\label{eq:T-ar}
\end{align}
так что
$$
\lim_{n\to\infty}\frac{\log\hat\Phi_n}n=7.03418177\dotsc.
$$
Сравнивая \eqref{eq:P-ar} и \eqref{eq:T-ar}, находим, что $\varphi(x)\ge\hat\varphi(x)$ для всех $x\in[0,1)$
за исключением $x\in\bigl[\frac15,\frac29\bigr)\cup\bigl[\frac37,\frac49\bigr)\cup\bigl[\frac67,\frac78\bigr)$.
Последнее означает, что с выбором $\tilde\Phi_n=\prod_{p\le8n}p^{\tilde\varphi(n/p)}$, где
\begin{align*}
\tilde\varphi(x)
&=\max\{\varphi(x),\hat\varphi(x)\}
\\
&=\begin{cases}
2 &\text{при $x\in\bigl[\frac16,\frac29\bigr)\cup\bigl[\frac14,\frac27\bigr)\cup\bigl[\frac12,\frac47\bigr)\cup\bigl[\frac56,\frac78\bigr)$}, \\
1 &\text{при $x\in\bigl[\frac18,\frac17\bigr)\cup\bigl[\frac29,\frac14\bigr)\cup\bigl[\frac27,\frac49\bigr)\cup\bigl[\frac47,\frac56\bigr)\cup\bigl[\frac78,\frac89\bigr)$}, \\
0 &\text{в остальных случаях},
\end{cases}
\end{align*}
справедливы включения
\begin{equation}
\tilde\Phi_n^{-1}q_n,\,\tilde\Phi_n^{-1}D_{9n}D_{8n}p_n\in\mathbb Z
\quad\text{for}\; n=0,1,2,\dots,
\label{eq:far}
\end{equation}
и, кроме того,
$$
\lim_{n\to\infty}\frac{\log\tilde\Phi_n}n
=8.79117698\dotsc.
$$

\section{Финал: доказательство теоремы~\ref{th:z2} и заключение}
\label{epi}

\begin{proof}[Доказательство теоремы~\textup{\ref{th:z2}}]
В этой главе мы построили линейные формы $r_n=q_n\zeta(2)-p_n$, $n=0,1,2,\dots$,
такие, что рациональные коэффициенты $q_n$ и $p_n$ удовлетворяют включениям~\eqref{eq:far},
в то время как асимптотика величин $r_n$ и $q_n$ при $n\to\infty$ задается соотношениями \eqref{eq:Pan}.
Полагая
$$
\tilde C_2=\lim_{n\to\infty}\frac{\log(\tilde\Phi_n^{-1}D_{9n}D_{8n})}n=8.20882301\dots
$$
и применяя лемму 2.1 из \cite{Ha1a}, мы приходим к оценке
$$
\mu(\zeta(2))\le\frac{C_0+C_1}{C_0-\tilde C_2}=5.09541178\dots
$$
для показателя иррациональности $\zeta(2)=\pi^2/6$.
\end{proof}

Как было отмечено ранее, последовательность $r_n=q_n\zeta(2)-p_n$ приближений,
построенных в доказательства теоремы~\ref{th:z2}, имеет двойственную последовательность
$r_n'=q_n\zeta(3)-p_n'$ рациональных приближений к~$\zeta(3)$ с тем же коэффициентом~$q_n$.
Можно показать, что $\tilde\Phi_n^{-1}D_{9n}D_{8n}^2p_n'\in\mathbb Z$ при $n=0,1,2,\dots$ и
$$
\limsup_{n\to\infty}\frac{\log|r_n'|}n=-C_0=-15.88518998\dots
$$
(ср.\ с~\eqref{eq:Pan}). Поскольку
$$
\lim_{n\to\infty}\frac{\log(\tilde\Phi_n^{-1}D_{9n}D_{8n}^2)}n=16.20882301\hdots>C_0,
$$
линейные формы $\tilde\Phi_n^{-1}D_{9n}D_{8n}^2r_n'\in\mathbb Z\zeta(3)+\mathbb Z$
неограниченно возрастают при $n\to\infty$, поэтому не могут быть использованы для доказательства иррациональности
числа $\zeta(3)$
(которая в этой ситуации означала бы также линейную независимость над $\mathbb Q$ чисел $1$, $\zeta(2)$ и $\zeta(3)$).
С помощью рекуррентного уравнения, использованного в доказательстве предложения~\ref{prop:int},
которому одновременно удовлетворяют последовательности $q_n$, $r_n=q_n\zeta(2)-p_n$ и $r_n'=q_n\zeta(3)-p_n'$,
мы вычислили первые триста значений последовательности
$$
\Lambda_n=\text{НОД}(\tilde\Phi_n^{-1}D_{9n}D_{8n}^2q_n,\tilde\Phi_n^{-1}D_{9n}D_{8n}^2p_n,\tilde\Phi_n^{-1}D_{9n}D_{8n}^2p_n'),
\quad n=0,1,2,\dotsc.
$$
Простые, участвующие в разложении $\Lambda_n$ на множители, не имеют явной характеризации в зависимости от~$n$;
для большинства значений $n$ эти простые $p$ попадают в (асимптотически невидимый) диапазон $p\le\sqrt{8n}$.
Тем не менее, абсолютные значения форм
$$
(\Lambda_n\tilde\Phi_n)^{-1}D_{9n}D_{8n}^2r_n\in\mathbb Z\zeta(2)+\mathbb Z
\quad\text{и}\quad
(\Lambda_n\tilde\Phi_n)^{-1}D_{9n}D_{8n}^2r_n'\in\mathbb Z\zeta(3)+\mathbb Z
$$
оказываются одновременно меньшими единицы для
\begin{align*}
n&=1,\dots,21,23,\dots,35,37,38,39,41,42,43,47,\dots,50,53,54,64,68,
\\ &\qquad
71,\dots,74,79,80,81,84,85,89,101,102,106,110,113,128,129,178,228,
\end{align*}
где $n\le300$.

\medskip
Нам представляется интересным исследовать арифметически другие классические гипергеометрические тождества,
приводимые в книгах Бейли \cite{Bai} и Слейтер \cite{Sl}: философия здесь заключается в том, что
за любым гипергеометрическим преобразованием кроится интересная арифметика.
Установленная ранее мера иррациональности $\zeta(2)$ в \cite{RV2}
и наилучшая известная мера иррациональности $\zeta(3)$ в \cite{RV3}
имеют глубокие гипергеометрические корни (см.~\cite{Z17}).
Еще одной иллюстрацией является гипергеометрическая конструкция рациональных приближений
к $\zeta(4)$ из работы~\cite{Z15}.

\chapter({Оценка снизу для расстояния от \$(3/2)\000\136k\$ до ближайшего целого})%
{Оценка снизу для расстояния от $(3/2)^k$\\ до ближайшего целого}
\label{chap:5}

Настоящая глава посвящена доказательству теоремы~\ref{th:5}
и оценок~\eqref{eq:0.42} (теорема~\ref{th:8} далее).
Приводимое далее доказательство опубликовано в работе~\cite{Z22}.

\section{Приближения Паде биномиального ряда}
\label{sec:5.1}

Зафиксируем два положительных целых параметра $a$ и $b$,
удовлетворяющих $2a\le b$. Из формулы~\eqref{eq:0.41} получаем
\begin{align}
\label{eq6}
\biggl(\frac32\biggr)^{3(b+1)}
&=\biggl(\frac{27}8\biggr)^{b+1}
=3^{b+1}\biggl(1-\frac19\biggr)^{-(b+1)}
\\
&=3^{b+1}\sum_{k=0}^\infty\binom{b+k}b\biggl(\frac19\biggr)^k
=3^{b-2a+1}\sum_{k=0}^\infty\binom{b+k}b3^{2(a-k)}
\nonumber\\
&=\text{целое число}+3^{b-2a+1}\sum_{k=a}^\infty\binom{b+k}b3^{2(a-k)}
\nonumber\\
&\equiv3^{b-2a+1}\sum_{\nu=0}^\infty\binom{a+b+\nu}b3^{-2\nu}\pmod{\mathbb Z}.
\nonumber
\end{align}
Это представление мотивирует (ср\. с~\cite{Beu}) построение приближений Паде
функции
\begin{equation}
F(z)=F(a,b;z)
=\sum_{\nu=0}^\infty\binom{a+b+\nu}bz^\nu
=\binom{a+b}b\sum_{\nu=0}^\infty
\frac{(a+b+1)_\nu}{(a+1)_\nu}z^\nu
\label{eq7}
\end{equation}
и последующее их применение с выбором $z=1/9$.

Выберем произвольное неотрицательное целое $n\le b$ и воспользуемся
общим рецептом из работ~\cite{Ch}, \cite{Zu}. Именно, рассмотрим многочлен
\begin{align}
\label{eq8}
Q_n(x)
&=\binom{a+b+n}{a+b}
{}_2F_1\biggl(\begin{matrix}-n, \, a+n \\ a+b+1 \end{matrix}
\biggm|x\biggr)
\\
&=\sum_{\mu=0}^n\binom{a+n-1+\mu}\mu\binom{a+b+n}{n-\mu}(-x)^\mu
=\sum_{\mu=0}^nq_\mu x^\mu
\in\mathbb Z[x]
\nonumber
\end{align}
степени~$n$. Тогда выполнено
\begin{align}
\label{eq9}
Q_n(z^{-1})F(z)
&=\sum_{\mu=0}^nq_{n-\mu}z^{\mu-n}
\cdot\sum_{\nu=0}^\infty\binom{a+b+\nu}bz^\nu
\\
&=\sum_{l=0}^\infty z^{l-n}
\sum_{\substack{\mu=0\\\mu\le l}}^nq_{n-\mu}\binom{a+b+l-\mu}b
\nonumber\\
&=\sum_{l=0}^{n-1}r_lz^{l-n}
+\sum_{l=n}^\infty r_lz^{l-n}
=P_n(z^{-1})+R_n(z);
\nonumber
\end{align}
здесь многочлен
\begin{equation}
P_n(x)=\sum_{l=0}^{n-1}r_lx^{n-l}\in\mathbb Z[x],
\qquad\text{где}\quad
r_l=\sum_{\mu=0}^lq_{n-\mu}\binom{a+b+l-\mu}b,
\label{eq10}
\end{equation}
имеет степень не выше~$n$, а коэффициенты остатка
$$
R_n(z)=\sum_{l=n}^\infty r_lz^{l-n}
$$
имеют следующее представление:
\begin{align*}
r_l
&=\sum_{\mu=0}^nq_{n-\mu}\binom{a+b+l-\mu}b
\\
&=\sum_{\mu=0}^n(-1)^{n-\mu}\binom{a+2n-1-\mu}{n-\mu}
\binom{a+b+n}\mu\binom{a+b+l-\mu}b
\\
&=(-1)^n\frac{(a+b+n)!}{(a+n-1)!n!b!}
\sum_{\mu=0}^n(-1)^\mu\binom n\mu
\frac{(a+2n-1-\mu)!(a+b+l-\mu)!}{(a+l-\mu)!(a+b+n-\mu)!}
\\
&=(-1)^n\frac{(a+b+n)!}{(a+n-1)!n!b!}
\\ &\qquad\times 
\frac{(a+2n-1)!(a+b+l)!}{(a+l)!(a+b+n)!}
\sum_{\mu=0}^n\frac{(-n)_\mu(-a-l)_\mu(-a-b-n)_\mu}
{\mu!(-a-2n+1)_\mu(-a-b-l)_\mu}
\\
&=(-1)^n\frac{(a+2n-1)!(a+b+l)!}{(a+n-1)!(a+l)!n!b!}
\cdot{}_3F_2\biggl(\begin{matrix}-n, \, -a-l, \, -a-b-n \\
-a-2n+1, \, -a-b-l \end{matrix}\biggm| 1 \biggr).
\end{align*}
Применяя формулу суммирования \eqref{eq4} Пфаффа--Заальшютца
с выбором $-a-l$, $-a-b-n$ и $-a-b-l$ в качестве $a$, $b$ и $c$ соответственно,
получаем
$$
r_l
=(-1)^n\frac{(a+2n-1)!(a+b+l)!}{(a+n-1)!(a+l)!n!b!}
\cdot\frac{(-b)_n(n-l)_n}{(-a-b-l)_n(a+n)_n}.
$$
Исходное условие $n\le b$ гарантирует, что коэффициенты
$r_l$ не обращаются в нуль тождественно
(в противном случае $(-b)_n=0$).
Кроме того, $(n-\nobreak l)_n=\nobreak0$ для $l$ из интервала
$n\le l\le2n-1$, так что $r_l=0$ для таких~$l$,
в то время как
\begin{align*}
r_l
&=\frac{(a+2n-1)!(a+b+l)!}{(a+n-1)!(a+l)!n!b!}
\\ &\qquad\times 
\frac{b!/(b-n)!\cdot(l-n)!/(l-2n)!}
{(a+b+l)!/(a+b+l-n)!\cdot(a+2n-1)!/(a+n-1)!}
\\
&=\frac{(a+b+l-n)!(l-n)!}{n!(b-n)!(a+l)!(l-2n)!}
\qquad\text{при}\quad l\ge2n.
\end{align*}
Окончательно,
\begin{align}
\label{eq11}
R_n(z)
&=\sum_{l=2n}^\infty r_lz^{l-n}
=z^n\sum_{\nu=0}^\infty r_{\nu+2n}z^\nu
\\
&=z^n\frac1{n!(b-n)!}\sum_{\nu=0}^\infty
\frac{(a+b+n+\nu)!(n+\nu)!}{\nu!(a+2n+\nu)!}z^\nu
\nonumber\\
&=z^n\binom{a+b+n}{b-n}
\cdot{}_2F_1\biggl(\begin{matrix} a+b+n+1, \, n+1  \\
a+2n+1 \end{matrix}\biggm| z \biggr).
\nonumber
\end{align}

Применяя интегральное представление~\eqref{eq5} к многочлену \eqref{eq8}
и остатку \eqref{eq11}, мы приходим к следующему утверждению.

\begin{lemma}
\label{l1}
Справедливы интегральные представления
$$
Q_n(z^{-1})
=\frac{(a+b+n)!}{(a+n-1)!n!(b-n)!}
\int_0^1t^{a+n-1}(1-t)^{b-n}(1-z^{-1}t)^n\,\d t
$$
и
$$
R_n(z)
=\frac{(a+b+n)!}{(a+n-1)!n!(b-n)!}\,z^n
\int_0^1t^n(1-t)^{a+n-1}(1-zt)^{-(a+b+n+1)}\,\d t.
$$
\end{lemma}

Нам также понадобится линейная независимость пары соседних
приближений Паде, вытекающая из следующего утверждения.

\begin{lemma}
\label{l2}
Выполнено
\begin{equation}
Q_{n+1}(x)P_n(x)-Q_n(x)P_{n+1}(x)
=(-1)^n\binom{a+2n+1}{a+n}\binom{a+b+n}{b-n}x.
\label{eq12}
\end{equation}
\end{lemma}

\begin{proof}
Очевидно, левая часть~\eqref{eq12} является многочленом,
свободный член которого равен нулю ввиду $P_n(0)=P_{n+1}(0)=0$
в соответствии с~\eqref{eq10}. С другой стороны,
\begin{align*}
&
Q_{n+1}(z^{-1})P_n(z^{-1})-Q_n(z^{-1})P_{n+1}(z^{-1})
\\ &\qquad
=Q_{n+1}(z^{-1})\bigl(Q_n(z^{-1})F(z)-R_n(z)\bigr)
\\ &\qquad\qquad
-Q_n(z^{-1})\bigl(Q_{n+1}(z^{-1})F(z)-R_{n+1}(z)\bigr)
\\ &\qquad
=Q_n(z^{-1})R_{n+1}(z)-Q_{n+1}(z^{-1})R_n(z),
\end{align*}
а из \eqref{eq8}, \eqref{eq11} заключаем, что в последней сумме
возникают только отрицательные степени переменной~$z$:
\begin{align*}
&
-Q_{n+1}(z^{-1})R_n(z)
\\ &\qquad
=(-1)^n\binom{a+2n+1}{a+n}z^{-n-1}\bigl(1+O(z)\bigr)
\cdot\binom{a+b+n}{b-n}z^n\bigl(1+O(z)\bigr)
\\ &\qquad
=(-1)^n\binom{a+2n+1}{a+n}\binom{a+b+n}{b-n}\frac1z+O(1)
\qquad\text{при}\quad z\to0.
\end{align*}
\vskip-10.2mm
\end{proof}

\section{Арифметические составляющие}
\label{sec:5.2}

Сразу отметим, что для простого $p>\sqrt{N}$ имеют место формулы
$$
\ord_pN!=\biggl\lf\frac Np\biggr\rf
\qquad\text{и}\qquad
\ord_pN=\lambda\biggl(\frac Np\biggr),
$$
где
$$
\lambda(x)
=1-\{x\}-\{-x\}
=1+\lf x\rf+\lf-x\rf
=\begin{cases}
1, & \text{если $x\in\mathbb Z$}, \\
0, & \text{если $x\notin\mathbb Z$}.
\end{cases}
$$

Для каждого простого $p>\sqrt{a+b+n}$ определим показатели
\begin{align}
\label{eq13}
e_p
&=\min_{\mu\in\mathbb Z}\biggl(-\biggl\{-\frac{a+n}p\biggr\}
+\biggl\{-\frac{a+n+\mu}p\biggr\}+\biggl\{\frac\mu p\biggr\}
\\ &\qquad
-\biggl\{\frac{a+b+n}p\biggr\}
+\biggl\{\frac{a+b+\mu}p\biggr\}
+\biggl\{\frac{n-\mu}p\biggr\}\biggr)
\nonumber\\
&=\min_{\mu\in\mathbb Z}\biggl(\biggl\lf-\frac{a+n}p\biggr\rf
-\biggl\lf-\frac{a+n+\mu}p\biggr\rf-\biggl\lf\frac\mu p\biggr\rf
\nonumber\\ &\qquad
+\biggl\lf\frac{a+b+n}p\biggr\rf
-\biggl\lf\frac{a+b+\mu}p\biggr\rf
-\biggl\lf\frac{n-\mu}p\biggr\rf\biggr)
\nonumber\\
&\le\min_{0\le\mu\le n}\ord_p\frac{a+n}{a+n+\mu}
\binom{a+n+\mu}\mu\binom{a+b+n}{n-\mu}
\nonumber\\
&=\min_{0\le\mu\le n}\ord_p\binom{a+n-1+\mu}\mu\binom{a+b+n}{n-\mu}
\nonumber
\end{align}
и
\begin{align}
\label{eq14}
e_p'
&=\min_{\mu\in\mathbb Z}\biggl(-\biggl\{\frac{a+n+\mu}p\biggr\}
+\biggl\{\frac{a+n}p\biggr\}+\biggl\{\frac\mu p\biggr\}
\\ &\qquad
-\biggl\{\frac{a+b+n}p\biggr\}
+\biggl\{\frac{a+b+\mu}p\biggr\}
+\biggl\{\frac{n-\mu}p\biggr\}\biggr)
\nonumber
\displaybreak[2]
\\
&=\min_{\mu\in\mathbb Z}\biggl(\biggl\lf\frac{a+n+\mu}p\biggr\rf
-\biggl\lf\frac{a+n}p\biggr\rf-\biggl\lf\frac\mu p\biggr\rf
\nonumber\\ &\qquad
+\biggl\lf\frac{a+b+n}p\biggr\rf
-\biggl\lf\frac{a+b+\mu}p\biggr\rf
-\biggl\lf\frac{n-\mu}p\biggr\rf\biggr)
\nonumber\\
&\le\min_{0\le\mu\le n}\ord_p\binom{a+n+\mu}\mu\binom{a+b+n}{n-\mu}.
\nonumber
\end{align}
Положим
\begin{equation}
\Phi=\Phi(a,b,n)
=\prod_{p>\sqrt{a+b+n}}p^{e_p}
\qquad\text{и}\qquad
\Phi'=\Phi'(a,b,n)
=\prod_{p>\sqrt{a+b+n}}p^{e_p'}.
\label{eq15}
\end{equation}
Из \eqref{eq8}, \eqref{eq10} и \eqref{eq13}, \eqref{eq15}
получаем следующий результат.

\begin{lemma}
\label{l3}
Справедливы включения
$$
\Phi^{-1}\cdot\binom{a+n-1+\mu}\mu\binom{a+b+n}{n-\mu}
\in\mathbb Z \qquad\text{при}\quad \mu=0,1,\dots,n,
$$
откуда
$$
\Phi^{-1}Q_n(x)\in\mathbb Z[x]
\qquad\text{и}\qquad
\Phi^{-1}P_n(x)\in\mathbb Z[x].
$$
\end{lemma}

Арифметическая информация в случае $n+1$ вместо $n$
дается в следующем утверждении.

\begin{lemma}
\label{l4}
Справедливы включения
\begin{equation}
(n+1){\Phi'}^{-1}\cdot\binom{a+n+\mu}\mu\binom{a+b+n+1}{n+1-\mu}
\in\mathbb Z \qquad\text{при}\quad \mu=0,1,\dots,n+1,
\label{eq16}
\end{equation}
откуда
$$
(n+1){\Phi'}^{-1}Q_{n+1}(x)\in\mathbb Z[x]
\qquad\text{и}\qquad
(n+1){\Phi'}^{-1}P_{n+1}(x)\in\mathbb Z[x].
$$
\end{lemma}

\begin{proof}
Запишем
$$
\binom{a+n+\mu}\mu\binom{a+b+n+1}{n+1-\mu}
=\binom{a+n+\mu}\mu\binom{a+b+n}{n-\mu}
\cdot\frac{a+b+n+1}{n+1-\mu}.
$$
Тогда для $p\nmid n+1-\mu$ выполнено
\begin{equation}
\ord_p\binom{a+n+\mu}\mu\binom{a+b+n+1}{n+1-\mu}
\ge\ord_p\binom{a+n+\mu}\mu\binom{a+b+n}{n-\mu}
\ge e_p';
\label{eq17}
\end{equation}
в противном случае $\mu\equiv n+1\pmod p$, так что
$\mu/p-(n+1)/p\in\mathbb Z$ или
\begin{align}
\label{eq18}
&
\ord_p\binom{a+n+\mu}\mu\binom{a+b+n+1}{n+1-\mu}
\\ &\qquad
=-\biggl\{\frac{a+n+\mu}p\biggr\}
+\biggl\{\frac{a+n}p\biggr\}+\biggl\{\frac\mu p\biggr\}
\nonumber\\ &\qquad\quad
-\biggl\{\frac{a+b+n+1}p\biggr\}
+\biggl\{\frac{a+b+\mu}p\biggr\}
+\biggl\{\frac{n+1-\mu}p\biggr\}
\nonumber\\ &\qquad
=-\biggl\{\frac{a+2n+1}p\biggr\}
+\biggl\{\frac{a+n}p\biggr\}+\biggl\{\frac{n+1}p\biggr\}
\nonumber\\ &\qquad
=\ord_p\binom{a+2n+1}{n+1}
=\ord_p\binom{a+2n}n+\ord_p\frac{a+2n+1}{n+1}
\nonumber\\ &\qquad
=\ord_p\binom{a+n+\mu}\mu\binom{a+b+n}{n-\mu}\bigg|_{\mu=n}
+\ord_p\frac{a+2n+1}{n+1}
\nonumber\\ &\qquad
\ge e_p'-\ord_p(n+1).
\nonumber
\end{align}
Комбинируя~\eqref{eq17} и~\eqref{eq18}, приходим к требуемым
включениям~\eqref{eq16}.
\end{proof}

\section{Доказательство теоремы~\ref{th:5}}
\label{sec:5.3}

Целочисленные параметры $a$, $b$ и $n$ будут теперь зависеть
от возрастающего параметра $m\in\mathbb N$ следующим образом:
$$
a=\alpha m, \qquad b=\beta m, \qquad
n=\gamma m \quad\text{или}\quad n=\gamma m+1,
$$
где с выбором целочисленных параметров $\alpha$,
$\beta$, $\gamma$, удовлетворяющих $2\alpha\le\nobreak\beta$ и
$\gamma<\beta$, мы определимся позднее. Тогда лемма~\ref{l1}
и метод Лапласа приводят к асимптотике
\begin{align}
\label{eq19}
C_0(z)
&=\lim_{m\to\infty}\frac{\log|R_n(z)|}m
\\
&=(\alpha+\beta+\gamma)\log(\alpha+\beta+\gamma)
-(\alpha+\gamma)\log(\alpha+\gamma)
\nonumber\\ &\;\quad
-\gamma\log\gamma
-(\beta-\gamma)\log(\beta-\gamma)
+\gamma\log|z|
\nonumber\\ &\;\quad
+\max_{0\le t\le1}\Re\bigl(\gamma\log t
+(\alpha+\gamma)\log(1-t)-(\alpha+\beta+\gamma)\log(1-zt)\bigr)
\nonumber
\end{align}
и
\begin{align}
\label{eq20}
C_1(z)
&=\lim_{m\to\infty}\frac{\log|Q_n(z^{-1})|}m
\\
&=(\alpha+\beta+\gamma)\log(\alpha+\beta+\gamma)
-(\alpha+\gamma)\log(\alpha+\gamma)
\nonumber\\ &\;\quad
-\gamma\log\gamma
-(\beta-\gamma)\log(\beta-\gamma)
\nonumber\\ &\;\quad
+\max_{0\le t\le1}\Re\bigl((\alpha+\gamma)\log t
+(\beta-\gamma)\log(1-t)+\gamma\log(1-z^{-1}t)\bigr).
\nonumber
\end{align}
Кроме того, из \eqref{eq13}--\eqref{eq15} и леммы~\ref{lem:1.5} имеем
\begin{equation}
\begin{aligned}
C_2
&=\lim_{m\to\infty}\frac{\log\Phi(\alpha m,\beta m,\gamma m)}m
=\int_0^1\varphi(x)\,\d\psi(x),
\\
C_2'
&=\lim_{m\to\infty}\frac{\log\Phi'(\alpha m,\beta m,\gamma m)}m
=\int_0^1\varphi'(x)\,\d\psi(x),
\end{aligned}
\label{eq21}
\end{equation}
где, как и ранее, $\psi(x)$~--- логарифмическая производная
гамма-функции, а $1$-периодические функции $\varphi(x)$ и $\varphi'(x)$
определяются следующим образом:
\begin{gather*}
\varphi(x)
=\min_{0\le y<1}\wh\varphi(x,y),
\qquad
\varphi'(x)
=\min_{0\le y<1}\wh\varphi'(x,y),
\\
\begin{aligned}
\wh\varphi(x,y)
&=-\{-(\alpha+\gamma)x\}
+\{-(\alpha+\gamma)x-y\}+\{y\}
\\ &\qquad
-\{(\alpha+\beta+\gamma)x\}
+\{(\alpha+\beta)x+y\}+\{\gamma x-y\},
\\
\wh\varphi'(x,y)
&=-\{(\alpha+\gamma)x+y\}
+\{(\alpha+\gamma)x\}+\{y\}
\\ &\qquad
-\{(\alpha+\beta+\gamma)x\}
+\{(\alpha+\beta)x+y\}+\{\gamma x-y\}.
\end{aligned}
\end{gather*}

\begin{lemma}
\label{l5}
Функции $\varphi(x)$ и $\varphi'(x)$
отличаются на множестве меры~$0$, так что
\begin{equation}
C_2=C_2'.
\label{eq22}
\end{equation}
\end{lemma}

\begin{proof}
Сразу отметим, что $\wh\varphi(x,0)=\wh\varphi'(x,0)\in\{0,1\}$,
откуда следует, что $\varphi(x)\in\{0,1\}$ и $\varphi'(x)\in\{0,1\}$.

В соответствии с определением функций
$$
\Delta(x,y)=\wh\varphi(x,y)-\wh\varphi'(x,y)
=\lambda((\alpha+\gamma)x)-\lambda((\alpha+\gamma)x+y).
$$
Предположим, что $(\alpha+\gamma)x\notin\mathbb Z$; в этом случае
имеем $\Delta(x,y)=-\lambda((\alpha+\gamma)x+y)\le0$,
откуда $\Delta(x,y)=0$, если $y\ne-\{(\alpha+\gamma)x\}$.
Таким образом, единственная возможность получить $\varphi(x)\ne\varphi'(x)$
следующая:
\begin{equation}
\label{eq22a}
\begin{aligned}
\wh\varphi(x,y)&\ge1 \quad\text{при $y\ne-\{(\alpha+\gamma)x\}$},
\\
\wh\varphi(x,y)&=0 \quad\text{при $y=-\{(\alpha+\gamma)x\}$}.
\end{aligned}
\end{equation}
Поскольку $\wh\varphi'(x,0)\ge1$, выполнено
$\{(\alpha+\beta)x\}+\{\gamma x\}=1$; далее,
\begin{align*}
&
\{(\alpha+\beta)x+y\}+\{\gamma x-y\}
-\{(\alpha+\beta+\gamma)x\}
\\ &\qquad
=\begin{cases}
1, & \text{если $0\le y<1-\{(\alpha+\beta)x\}$}, \\
0, & \text{если $1-\{(\alpha+\beta)x\}\le y\le\{\gamma x\}$}, \\
1, & \text{если $\{\gamma x\}<y<1$},
\end{cases}
\end{align*}
и
$$
\{-(\alpha+\gamma)x-y\}+\{y\}-\{-(\alpha+\gamma)x\}
=\begin{cases}
0, & \text{если $0\le y\le\{-(\alpha+\gamma)x\}$}, \\
1, & \text{если $\{-(\alpha+\gamma)x\}<y<1$}.
\end{cases}
$$
Условия \eqref{eq22a} на функцию $\wh\varphi(x,y)$
влекут $1-\{(\alpha+\beta)x\}=\{-(\alpha+\gamma)x\}$
или, что то же самое, $(\beta-\gamma)x\in\mathbb Z$.

Наконец, множества $\{x\in\mathbb R:(\alpha+\gamma)x\in\mathbb Z\}$
и $\{x\in\mathbb R:(\beta-\gamma)x\in\mathbb Z\}$ имеют меру~$0$,
что завершает доказательство леммы.
\end{proof}

Нашей целью является получение оценки снизу для модуля величины
$\varepsilon_k$, где
$$
\biggl(\frac32\biggr)^k=M_k+\varepsilon_k,
\qquad M_k\in\mathbb Z, \quad 0<|\varepsilon_k|<\frac12.
$$
Представим $k\ge3$ в виде $k=3(\beta m+1)+j$
с неотрицательными целыми $m$ и $j<3\beta$.
Умножая обе части~\eqref{eq9} на ${\wt\Phi}^{-1}3^{b-2a+j+1}$, где
$\wt\Phi=\Phi(\alpha m,\beta m,\gamma m)$ в случае $n=\gamma m$
и $\wt\Phi=\Phi'(\alpha m,\beta m,\gamma m)/(\gamma m+1)$
в случае $n=\gamma m+1$, и подставим $z=1/9$:
\begin{align}
\label{eq23}
&
Q_n(9){\wt\Phi}^{-1}2^j
\cdot\biggl(\frac32\biggr)^j3^{b-2a+1}F\biggl(a,b;\frac19\biggr)
\\ &\qquad
=P_n(9){\wt\Phi}^{-1}3^{b-2a+j+1}
+R_n\biggl(\frac19\biggr){\wt\Phi}^{-1}3^{b-2a+j+1}.
\nonumber
\end{align}
Из \eqref{eq6}, \eqref{eq7} мы видим, что
$$
\biggl(\frac32\biggr)^j3^{b-2a+1}F\biggl(a,b;\frac19\biggr)
\equiv\biggl(\frac32\biggr)^{3(b+1)+j}\pmod{\mathbb Z}
=\biggl(\frac32\biggr)^k,
$$
так что левая часть представляется в виде $M_k'+\varepsilon_k$ для
некоторого $M_k'\in\mathbb Z$ и равенство~\eqref{eq23}
может быть записано в виде
\begin{equation}
Q_n(9){\wt\Phi}^{-1}2^j\cdot\varepsilon_k
=M_k''+R_n\biggl(\frac19\biggr){\wt\Phi}^{-1}3^{b-2a+j+1},
\label{eq24}
\end{equation}
где
$$
M_k''=P_n(9){\wt\Phi}^{-1}3^{b-2a+j+1}
-Q_n(9){\wt\Phi}^{-1}2^jM_k'
\in\mathbb Z
$$
согласно леммам~\ref{l3} и~\ref{l4}. В свою очередь, лемма~\ref{l2}
гарантирует $M_k''\ne0$ для хотя бы одного из чисел $n=\gamma m$ или $\gamma m+1$;
мы фиксируем соответствующий выбор~$n$.
Предполагая далее, что
\begin{equation}
C_0\biggl(\frac19\biggr)-C_2+(\beta-2\alpha)\log3
<0,
\label{eq25}
\end{equation}
из \eqref{eq19} и \eqref{eq21} получаем
$$
\biggl|R_n\biggl(\frac19\biggr){\wt\Phi}^{-1}3^{b-2a+j+1}\biggr|
<\frac12
\qquad\text{для всех}\quad m\ge N_1,
$$
где $N_1>0$ является некоторой эффективной абсолютной постоянной.
Следовательно, в соответствии с~\eqref{eq24} и $|M_k''|\ge1$ имеем
$$
|Q_n(9){\wt\Phi}^{-1}2^j|\cdot|\varepsilon_k|
\ge|M_k''|
-\biggl|R_n\biggl(\frac19\biggr){\wt\Phi}^{-1}3^{b-2a+j+1}\biggr|
>\frac12,
$$
так что согласно~\eqref{eq19}, \eqref{eq20} заключаем, что
$$
|\varepsilon_k|
>\frac{\wt\Phi}{2^{j+1}|Q_n(9)|}
\ge\frac{\wt\Phi}{2^{3\beta}|Q_n(9)|}
>e^{-m(C_1(1/9)-C_2+\delta)}
$$
для любого $\delta>0$ и $m>N_2(\delta)$ при условии, что
$C_1(1/9)-C_2+\delta>0$;
здесь $N_2(\delta)$ эффективно зависит от~$\delta$.
Наконец, ввиду $k>3\beta m$ получаем оценку
\begin{equation}
|\varepsilon_k|
>e^{-k(C_1(1/9)-C_2+\delta)/(3\beta)},
\label{eq26}
\end{equation}
справедливую при всех $k\ge K_0(\delta)$, где постоянная
$K_0(\delta)$ может быть явно выражена в терминах
$\max(N_1,N_2(\delta))$.

Выбирая $\alpha=\gamma=9$ и $\beta=19$
(что является оптимальным выбором целочисленных параметров
$\alpha,\beta,\gamma$, по крайней мере при ограничении
$\beta\le\nobreak100$), находим
$$
C_0\biggl(\frac19\biggr)
=3.28973907\dots,
\qquad
C_1\biggl(\frac19\biggr)
=35.48665992\dots,
$$
и
$$
\varphi(x)=\begin{cases}
1, & \text{если}\ \{x\}\in\bigl[\frac2{37},\frac1{18}\bigr]
\cup\bigl[\frac3{37},\frac1{10}\bigr)
\cup\bigl[\frac4{37},\frac19\bigr)
\cup\bigl[\frac6{37},\frac16\bigr]
\cup\bigl[\frac7{37},\frac15\bigr)
\\ &\phantom{\text{если}\ \{x\}\in}
\cup\bigl[\frac8{37},\frac29\bigr)
\cup\bigl[\frac{10}{37},\frac5{18}\bigr]
\cup\bigl[\frac{11}{37},\frac3{10}\bigr)
\cup\bigl[\frac{12}{37},\frac13\bigr)
\cup\bigl[\frac{14}{37},\frac7{18}\bigr]
\\ &\phantom{\text{если}\ \{x\}\in}
\cup\bigl[\frac{16}{37},\frac49\bigr)
\cup\bigl[\frac{18}{37},\frac12\bigr)
\cup\bigl[\frac{20}{37},\frac59\bigr)
\cup\bigl[\frac{22}{37},\frac35\bigr)
\cup\bigl[\frac{24}{37},\frac23\bigr)
\\ &\phantom{\text{если}\ \{x\}\in}
\cup\bigl[\frac{28}{37},\frac79\bigr)
\cup\bigl[\frac{32}{37},\frac89\bigr)
\cup\bigl[\frac{36}{37},1\bigr),
\\
0 & \text{в противном случае},
\end{cases}
$$
откуда
$$
C_2=C_2'=4.46695926\dotsc.
$$
С помощью этих вычислений проверяем условие~\eqref{eq25},
$$
C_0\biggl(\frac19\biggr)-C_2+(\beta-2\alpha)\log3
=-0.07860790\dots,
$$
и с выбором $\delta=0.00027320432\dots$ имеем
$$
e^{-(C_1(1/9)-C_2+\delta)/(3\beta)}
=0.5803.
$$
Согласно~\eqref{eq26} последнее соотношение завершает
доказательство теоремы~\ref{th:5}.
\qed

\section{Смежные результаты}
\label{sec:5.4}

Конструкция этой главы позволяет нам получить схожие результаты
для последовательностей $\|(4/3)^k\|$ и $\|(5/4)^k\|$, используя
представления
\begin{align}
\label{eq27}
\biggl(\frac43\biggr)^{2(b+1)}
&=2^{b+1}\biggl(1+\frac18\biggr)^{-(b+1)}
\\
&\equiv(-1)^a2^{b-3a+1}
F\biggl(a,b;-\frac18\biggr)\pmod{\mathbb Z},
\nonumber\\ &\qquad\qquad
\text{где}\quad 3a\le b,
\nonumber
\end{align}
и
\begin{align}
\label{eq28}
\biggl(\frac54\biggr)^{7b+3}
&=2\cdot5^b\biggl(1+\frac3{125}\biggr)^{-(2b+1)}
\\
&\equiv(-1)^a2\cdot3^a\cdot5^{b-3a}
F\biggl(a,2b;-\frac3{125}\biggr)\pmod{\mathbb Z},
\nonumber\\ &\qquad\qquad
\text{где}\quad 3a\le b.
\nonumber
\end{align}
Именно, выбирая $a=5m$, $b=15m$, $n=6m(+1)$
в случае~\eqref{eq27} и $a=3m$, $b=9m$, $n=7m(+1)$
в случае~\eqref{eq28} и повторяя рассуждения \S\,\ref{sec:5.3},
мы приходим к следующему результату.

\begin{theorem}
\label{th:8}
Справедливы оценки
\begin{align*}
\biggl\|\biggl(\frac43\biggr)^k\biggr\|>0.4914^k
=3^{-k\cdot0.64672207\dots}
&\qquad\text{при}\quad k\ge K_1,
\\
\biggl\|\biggl(\frac54\biggr)^k\biggr\|>0.5152^k
=4^{-k\cdot0.47839775\dots}
&\qquad\text{при}\quad k\ge K_2,
\end{align*}
где $K_1$ и $K_2$ --- некоторые эффективные постоянные.
\end{theorem}

Общий случай последовательности $\|(1+1/N)^k\|$
с целым $N\ge5$ может быть рассмотрен так же, как и в~\cite{Beu}, \cite{Ben1},
с использованием представления
$$
\biggl(\frac{N+1}N\biggr)^{b+1}
=\biggl(1-\frac1{N+1}\biggr)^{-(b+1)}
\equiv F\biggl(0,b;\frac1{N+1}\biggr)\pmod{\mathbb Z}.
$$

\chapter{Решение задачи А.~Шмидта}
\label{chap:4}

Правая часть~\eqref{eq:0.37} определяет так называемое
{\it лежандрово преобразование\/} последовательности
$\{c_k^{(r)}\}_{k=0,1,\dots}$.
В общей ситуации если
\begin{equation*}
a_n=\sum_{k=0}^n\binom nk\binom{n+k}kc_k
=\sum_{k=0}^n\binom{2k}k\binom{n+k}{n-k}c_k,
\end{equation*}
то согласно известному соотношению между обратными
лежандровыми парами выполнено
\begin{equation*}
\binom{2n}nc_n=\sum_k(-1)^{n-k}d_{n,k}a_k,
\end{equation*}
где
\begin{equation*}
d_{n,k}=\binom{2n}{n-k}-\binom{2n}{n-k-1}
=\frac{2k+1}{n+k+1}\binom{2n}{n-k}.
\end{equation*}
Следовательно, полагая
\begin{equation}
t_{n,j}^{(r)}=\sum_{k=j}^n(-1)^{n-k}d_{n,k}\binom{k+j}{k-j}^r,
\label{eq:4.1}
\end{equation}
получаем
\begin{equation}
\binom{2n}nc_n^{(r)}
=\sum_{j=0}^n\binom{2j}j^rt_{n,j}^{(r)}.
\label{eq:4.2}
\end{equation}
Случай $r=1$ задачи~\ref{prob:1} тривиален (потому и не включен
в формулировку задачи), в то время как случаи
$r=2$ и $r=3$ реализованы в \cite{Sc2}, \cite{St} благодаря тому, что
$t_{n,j}^{(2)}$ и $t_{n,j}^{(3)}$ имеют {\it замкнутую форму}.
Именно, с помощью алгоритма Цайльбергера созидательного телескопирования
для этих последовательностей несложно показать, что
они удовлетворяют простым разностным уравнениям первого порядка
как по~$n$, так и по~$j$. К сожалению, в случае $r\ge4$
подобная аргументация уже недоступна.

Ф.~Штрель в~\cite[\S\,4.2]{St} заметил, что
требуемая в теореме~\ref{th:4} целозначность
следует из делимости произведения
$\binom{2j}j^r\cdot t_{n,j}^{(r)}$ на $\binom{2n}n$
при всех $j$, $0\le j\le n$. Он высказал гипотезу
о справедливости результата, обобщающего теорему~\ref{th:4},
которое мы и будем доказывать.

\begin{theorem}
\label{th:9}
Числа $\binom{2n}n^{-1}\binom{2j}jt_{n,j}^{(r)}$
являются целыми.
\end{theorem}

Наша стратегия доказательства теоремы~\ref{th:9}
(а, значит, и теоремы~\ref{th:4}) заключается в следующем:
переписать \eqref{eq:4.1} в гипергеометрическом виде
и применить подходящие формулы суммирования и преобразования
(предложения~\ref{prop:4.1} и~\ref{prop:4.2} ниже).
Приводимое далее доказательство опубликовано
в~\cite{Z19} и~\cite{Z20}.

\section{Совершенно уравновешенные ряды}
\label{sec:4.1}

Заменяя $n-k$ на~$l$ в~\eqref{eq:4.1}, имеем
\begin{equation*}
t_{n,j}^{(r)}=\sum_{l\ge0}(-1)^l
\frac{2n-2l+1}{2n-l+1}\binom{2n}l\binom{n-l+j}{n-l-j}^r,
\end{equation*}
где ряд справа обрывается (т.е\. содержит конечное число
членов). Будет удобно записывать подобные обрывающиеся ряды
кратко в виде $\sum_l$, что, на самом деле, является
стандартным соглашением (см., например, \cite{PWZ}).
Отношение двух последовательных членов последней суммы равно
\begin{equation*}
\frac{-(2n+1)+l}{1+l}
\cdot\frac{-\frac12(2n-1)+l}{-\frac12(2n+1)+l}
\cdot\biggl(\frac{-(n-j)+l}{-(n+j)+l}\biggr)^r,
\end{equation*}
поэтому
\begin{equation*}
t_{n,j}^{(r)}
=\binom{n+j}{n-j}^r
\cdot{}_{r+2}F_{r+1}\biggl(\begin{matrix}
-(2n+1), & -\tfrac12(2n-1), & -(n-j), & \dots, & -(n-j) \\[2.5pt]
& -\tfrac12(2n+1), & -(n+j), & \dots, & -(n+j)
\end{matrix}\biggm|1\biggr)
\end{equation*}
является совершенно уравновешенным гипергеометрическим рядом.
В да\-льнейшем будем опускать аргумент $z=1$ у гипергеометрических рядов.

Следующие два классических результата --- суммирование Дугалла
\linebreak[4]
${}_5F_4(1)$-ряда (доказанное в~1907) и преобразование Уиппла
${}_7F_6(1)$-ряда (доказанное в~1926) --- потребуются в случаях
$r=3,4,5$ для теорем~\ref{th:4} и~\ref{th:9}.

\begin{proposition}[{\cite[\S\,4.3]{Bai}}]
\label{prop:4.1}
Имеют место равенства
\begin{equation}
{}_5F_4\biggl(\begin{matrix}
a, & 1+\frac12a, & c, & d, & -m \\[2.5pt]
& \frac12a, & 1+a-c, & 1+a-d, & 1+a+m
\end{matrix}\biggr)
=\frac{(1+a)_m\,(1+a-c-d)_m}{(1+a-c)_m\,(1+a-d)_m}
\label{eq:4.3}
\end{equation}
и
\begin{align}
&
{}_7F_6\biggl(\begin{matrix}
a, & 1+\frac12a, & b, & c, & d, & e, & -m \\[2.5pt]
& \frac12a, & 1+a-b, & 1+a-c, & 1+a-d, & 1+a-e, & 1+a+m
\end{matrix}\biggr)
\nonumber\\ &\qquad
=\frac{(1+a)_m\,(1+a-d-e)_m}{(1+a-d)_m\,(1+a-e)_m}
\cdot{}_4F_3\biggl(\begin{gathered}
1+a-b-c, \, d, \, e, \, -m \\
1+a-b, \, 1+a-c, \, d+e-a-m
\end{gathered}\biggr),
\label{eq:4.4}
\end{align}
где $m$ неотрицательное целое.
\end{proposition}

Как следствие~\eqref{eq:4.3} получаем
(без всякого созидательного телескопирования)
\begin{equation*}
t_{n,j}^{(3)}
=\binom{n+j}{n-j}^3
\cdot\frac{(-2n)_{n-j}(-2n+2(n-j))_{n-j}}{(-2n+(n-j))_{n-j}^2}
=\frac{(2n)!}{(3j-n)!\,(n-j)!^3},
\end{equation*}
что в точности совпадает с выражением, найденным
в \cite[\S\,4.2]{St}. Следовательно, из~\eqref{eq:4.2}
вытекает следующее точная формула:
\begin{equation*}
c_n^{(3)}
=\binom{2n}n^{-1}\sum_j\binom{2j}j^3
\frac{(2n)!}{(3j-n)!\,(n-j)!^3}
=\sum_j\binom{2j}j^2\binom{2j}{n-j}\binom nj^2.
\end{equation*}

В случае $r=5$ применим преобразование~\eqref{eq:4.4}:
\begin{align*}
t_{n,j}^{(5)}
&=\binom{n+j}{n-j}^5
\cdot\frac{(-2n)_{n-j}(-2n+2(n-j))_{n-j}}{(-2n+(n-j))_{n-j}^2}
\\ &\qquad\times 
{}_4F_3\biggl(\begin{gathered}
-2j, \, -(n-j), \, -(n-j), \, -(n-j) \\
-(n+j), \, -(n+j), \, 3j-n+1
\end{gathered}\biggr)
\\
&=\binom{n+j}{n-j}^2\frac{(2n)!}{(3j-n)!\,(n-j)!^3}
\sum_l\frac{(-2j)_l\,(-(n-j))_l^3}{l!\,(-(n+j))_l^2(3j-n+1)_l}
\\
&=\frac{(2n)!}{(2j)!\,(n-j)!^2}
\sum_l\binom{n-l+j}{n-l-j}^2\binom{2j}l
\binom{2j}{n-l-j}
\\
&=\frac{(2n)!}{(2j)!\,(n-j)!^2}
\sum_k\binom{k+j}{k-j}^2\binom{2j}{n-k}
\binom{2j}{k-j},
\end{align*}
значит,
\begin{equation*}
\binom{2n}n^{-1}\binom{2j}jt_{n,j}^{(5)}
=\binom nj^2
\sum_k\binom{k+j}{k-j}^2\binom{2j}{n-k}
\binom{2j}{k-j}
\end{equation*}
являются целыми, и в соответствии с~\eqref{eq:4.2}
мы получаем~\eqref{eq:0.39}.

В случае $r=4$ применим формулу \eqref{eq:4.4} с $b=(1+a)/2$
(так что ряд в левой части редуцируется
к совершенно уравновешенному ${}_6F_5(1)$-ряду):
\begin{align*}
t_{n,j}^{(4)}
&=\binom{n+j}{n-j}^4
\cdot\frac{(-2n)_{n-j}(-2n+2(n-j))_{n-j}}{(-2n+(n-j))_{n-j}^2}
\\ &\qquad\times 
{}_4F_3\biggl(\begin{gathered}
-j, \, -(n-j), \, -(n-j), \, -(n-j) \\
-n, \, -(n+j), \, 3j-n+1
\end{gathered}\biggr)
\displaybreak[2]\\
&=\binom{n+j}{n-j}\frac{(2n)!}{(3j-n)!\,(n-j)!^3}
\sum_l\frac{(-j)_l\,(-(n-j))_l^3}{l!\,(-n)_l\,(-(n+j))_l(3j-n+1)_l}
\\
&=\frac{(2n)!\,j!}{n!\,(n-j)!\,(2j)!}
\sum_l\binom{n-l+j}{n-l-j}\binom jl\binom{n-l}j\binom{2j}{n-l-j}
\\
&=\frac{(2n)!\,j!}{n!\,(n-j)!\,(2j)!}
\sum_k\binom{k+j}{k-j}\binom j{n-k}\binom kj\binom{2j}{k-j},
\end{align*}
откуда вновь $\binom{2n}n^{-1}\binom{2j}jt_{n,j}^{(4)}\in\mathbb Z$,
и мы получаем формулу~\eqref{eq:0.38}.

\section{Кратное преобразование Эндрюса}
\label{sec:4.2}

``Классические'' гипергеометрические тождества справляются
только в случаях\footnote{Как заметил К.~Краттенталер,
это не совсем верно, так как приводимое далее преобразование
Эндрюса является следствием классических преобразования Уиппла~\eqref{eq:4.4}
и формулы Пфаффа--Заальшютца~\eqref{eq4}.}
$r=2,3,4,5$ в теоремах \ref{th:4} и~\ref{th:9}.
Для того чтобы доказать эти теоремы в общем случае,
нам понадобится многомерное обобщение преобразования
Уиппла~\eqref{eq:4.4}. Требуемое обобщение установлено
Дж.~Эндрюсом в~\cite[теорема~4]{An1}.
Переходя к пределу при $q\to1$ в теореме Эндрюса,
мы получаем следующий результат.

\begin{proposition}
\label{prop:4.2}
Для $s\ge1$ и $m$ неотрицательного целого
\begin{align*}
&
{}_{2s+3}F_{2s+2}\biggl(\begin{matrix}
a, & 1+\frac12a, & b_1, & c_1, & b_2, & c_2, & \dots \\[2.5pt]
& \frac12a, & 1+a-b_1, & 1+a-c_1, & 1+a-b_2, & 1+a-c_2, & \dots
\end{matrix}
\\ &\qquad\qquad\qquad\qquad\qquad\qquad\qquad
\begin{matrix}
\dots, & b_s, & c_s, & -m \\[1.5pt]
\dots, & 1+a-b_s, & 1+a-c_s, & 1+a+m
\end{matrix}\biggr)
\\ &\quad
=\frac{(1+a)_m(1+a-b_s-c_s)_m}{(1+a-b_s)_m(1+a-c_s)_m}
\sum_{l_1\ge0}\frac{(1+a-b_1-c_1)_{l_1}(b_2)_{l_1}(c_2)_{l_1}}
{l_1!\,(1+a-b_1)_{l_1}(1+a-c_1)_{l_1}}
\\ &\quad\qquad\times
\sum_{l_2\ge0}\frac{(1+a-b_2-c_2)_{l_2}(b_3)_{l_1+l_2}(c_3)_{l_1+l_2}}
{l_2!\,(1+a-b_2)_{l_1+l_2}(1+a-c_2)_{l_1+l_2}}
\dotsb
\\ &\quad\qquad\times
\sum_{l_{s-1}\ge0}
\frac{(1+a-b_{s-1}-c_{s-1})_{l_{s-1}}
(b_s)_{l_1+\dots+l_{s-1}}(c_s)_{l_1+\dots+l_{s-1}}}
{l_{s-1}!\,(1+a-b_{s-1})_{l_1+\dots+l_{s-1}}
(1+a-c_{s-1})_{l_1+\dots+l_{s-1}}}
\\ &\quad\qquad\qquad\times
\frac{(-m)_{l_1+\dots+l_{s-1}}}
{(b_s+c_s-a-m)_{l_1+\dots+l_{s-1}}}.
\end{align*}
\end{proposition}

\begin{proof}[Доказательство теоремы~{\rm\ref{th:9}}]
Как и в~\S\,\ref{sec:4.1}, будем различать случаи
в зависимости от четности~$r$.

Если $r=2s+1$, то, полагая $a=-(2n+1)$ и
$b_1=c_1=\dots=b_s=c_s=-m=-(n-j)$
в предложении~\ref{prop:4.2}, находим
\begin{align*}
t_{n,j}^{(2s+1)}
&=\binom{n+j}{n-j}^{2s-2}\frac{(2n)!}{(3j-n)!\,(n-j)!^3}
\sum_{l_1}\binom{2j}{l_1}
\biggl(\frac{(-(n-j))_{l_1}}{(-(n+j))_{l_1}}\biggr)^2
\\ &\qquad\times
\sum_{l_2}\binom{2j}{l_2}
\biggl(\frac{(-(n-j))_{l_1+l_2}}{(-(n+j))_{l_1+l_2}}\biggr)^2
\dotsb
\\ &\qquad\times
\sum_{l_{s-1}}\binom{2j}{l_{s-1}}
\biggl(\frac{(-(n-j))_{l_1+\dots+l_{s-1}}}
{(-(n+j))_{l_1+\dots+l_{s-1}}}\biggr)^2
\\ &\qquad\qquad\times
\frac{(-1)^{l_1+\dots+l_{s-1}}(-(n-j))_{l_1+\dots+l_{s-1}}}
{(3j-n+1)_{l_1+\dots+l_{s-1}}}
\displaybreak[2]\\
&=\frac{(2n)!}{(2j)!\,(n-j)!^2}
\sum_{l_1}\binom{2j}{l_1}\binom{n-l_1+j}{n-l_1-j}^2
\sum_{l_2}\binom{2j}{l_2}\binom{n-l_1-l_2+j}{n-l_1-l_2-j}^2
\dotsb
\\ &\qquad\times
\sum_{l_{s-1}}\binom{2j}{l_{s-1}}
\binom{n-l_1-\dotsb-l_{s-1}+j}{n-l_1-\dotsb-l_{s-1}-j}^2
\cdot\binom{2j}{n-l_1-\dotsb-l_{s-1}-j}.
\end{align*}

Если $r=2s$, то применим предложение~\ref{prop:4.2}
с выбором $a=-(2n+1)$, $b_1=(a+1)/2=-n$ и
$c_1=b_2=\dots=b_s=c_s=-m=-(n-j)$:
\begin{align*}
t_{n,j}^{(2s)}
&=\binom{n+j}{n-j}^{2s-3}\frac{(2n)!}{(3j-n)!\,(n-j)!^3}
\sum_{l_1}\binom j{l_1}\frac{(-(n-j))_{l_1}}{(-n)_{l_1}}
\,\frac{(-(n-j))_{l_1}}{(-(n+j))_{l_1}}
\\ &\qquad\times
\sum_{l_2}\binom{2j}{l_2}
\biggl(\frac{(-(n-j))_{l_1+l_2}}{(-(n+j))_{l_1+l_2}}\biggr)^2
\dotsb
\\ &\qquad\times
\sum_{l_{s-1}}\binom{2j}{l_{s-1}}
\biggl(\frac{(-(n-j))_{l_1+\dots+l_{s-1}}}
{(-(n+j))_{l_1+\dots+l_{s-1}}}\biggr)^2
\\ &\qquad\qquad\times
\frac{(-1)^{l_1+\dots+l_{s-1}}(-(n-j))_{l_1+\dots+l_{s-1}}}
{(3j-n+1)_{l_1+\dots+l_{s-1}}}
\displaybreak[2]\\
&=\frac{(2n)!\,j!}{n!\,(n-j)!\,(2j)!}
\sum_{l_1}\binom j{l_1}\binom{n-l_1}j\binom{n-l_1+j}{n-l_1-j}
\\ &\qquad\times
\sum_{l_2}\binom{2j}{l_2}\binom{n-l_1-l_2+j}{n-l_1-l_2-j}^2
\dotsb
\\ &\qquad\times
\sum_{l_{s-1}}\binom{2j}{l_{s-1}}
\binom{n-l_1-\dotsb-l_{s-1}+j}{n-l_1-\dotsb-l_{s-1}-j}^2
\cdot\binom{2j}{n-l_1-\dotsb-l_{s-1}-j}.
\end{align*}

В обоих случаях требуемая целочисленность
\begin{equation*}
\binom{2n}n^{-1}\binom{2j}jt_{n,j}^{(r)}\in\mathbb Z,
\qquad j=0,1,\dots,n,
\end{equation*}
следует из полученных формул, что и доказывает
теорему~\ref{th:9}.
\end{proof}

Теорема~\ref{th:4} по сути была доказана во время доказательства
теоремы~\ref{th:9} с явными формулами для $c_n^{(4)}$, $c_n^{(5)}$
и в общем для $c_n^{(r)}$, $r\ge2$.

\chapter({Интегральные представления \$L\$-рядов эллиптических кривых})%
{Интегральные представления $L$-рядов эллиптических кривых}
\label{chap:7}

Содержание этой главы устроено следующим образом.
Мы демонстрируем метод из \cite{RZ1}, \cite{RZ2} на примере вычисления
$L$-значения $L(E,2)$ в \S\,\ref{ch7:s1}, а затем в \S\,\ref{ch7:s2} описываем
общий алгоритм, на котором основан метод. В~\S\,\ref{ch7:s3}
мы иллюстрируем этот алгоритм на примере вычисления $L(E,3)$ как периода;
наконец, в \S\,\ref{ch7:s4} мы строим гипергеометрические представления
для найденных ранее интегральных, которые, в свою очередь, приводят
к иным интегральным представлениям --- всё это составляет предмет теоремы~\ref{th:7}.
В примерах \S\,\ref{ch7:s1}, \S\,\ref{ch7:s3} (и~\S\,\ref{ch7:s4}),
$E$~обозначает эллиптическую кривую кондуктора~32. Выбор этой кривой обусловлен
двумя причинами. Во-первых, она не обсуждается в наших совместных работах \cite{RZ1}, \cite{RZ2}
с М.~Роджерсом; во-вторых, возникающие модулярные параметризации наиболее
простые и классические. Основные результаты главы опубликованы в работе автора~\cite{Z96}.

\medskip
Всюду в этой главе мы используем обозначение $q=e^{2\pi i\tau}$ для $\tau$
из верхней полуплоскости $\operatorname{Im}\tau>0$, так что $|q|<1$. Нашим
основным конструктором для построения модулярных форм служит
эта-функция Дедекинда
\begin{equation*}
\eta(\tau)=q^{1/24}\prod_{m=1}^{\infty}(1-q^m)
=\sum_{n=-\infty}^{\infty}(-1)^nq^{(6n+1)^2/24}
\end{equation*}
с ее модулярной инволюцией
\begin{equation}
\eta(-1/\tau)=\sqrt{-i\tau}\eta(\tau).
\label{mod_inv}
\end{equation}
Для краткости мы используем обозначение $\eta_k=\eta(k\tau)$. Отметим, что
бесконечное произведение для эта-функции в более общей форме уже появлялось
в~гл.~\ref{chap:3} под именем $q$-гамма функции, а логарифмическая производная
эта-функции (с точностью до нормировки) есть не что иное, как $\zeta_q(2)$,
также обсуждавшееся в~гл.~\ref{chap:3}.

Для функций параметра $\tau$ или $q=e^{2\pi i\tau}$ будем использовать дифференциальный оператор
$$
\delta=\frac1{2\pi i}\,\frac{\d}{\d\tau}=q\,\frac{\d}{\d q},
$$
а также обозначим через $\delta^{-1}$ обратное интегральное преобразование,
нормализованное нулем в $\tau=i\infty$ (или в $q=0$):
$$
\delta^{-1}f=\int_0^qf\,\frac{\d q}q.
$$
В частности, для (параболической) модулярной формы $f(\tau)=\sum_{n=1}^\infty a_nq^n$ справедливо представление
\begin{equation}
L(f,m)=\frac1{(m-1)!}\int_0^1f\log^{m-1}q\,\frac{\d q}q
=\sum_{n=1}^\infty\frac{a_n}{n^m}=(\delta^{-m}f)|_{q=1}
\label{Lfm}
\end{equation}
в случае, когда указанная сумма имеет смысл.

\section({\$L(E,2)\$})%
{$L(E,2)$}
\label{ch7:s1}

Как уже было отмечено в~\S\,\ref{sec:0.8},
$L$-ряд эллиптической кривой кондуктора 32 совпадает с $L$-рядом
параболической формы $f(\tau)=\eta_4^2\eta_8^2$
(см.\ наш пример в~\S\,\ref{sec:0.8}).

Нам понадобится следующее разложение в ряд Ламберта
(являющийся частным случаем разложения модулярной формы веса~1 в ряд Эйзенштейна):
\begin{gather*}
\frac{\eta_8^4}{\eta_4^2}
=\sum_{m\ge1}\biggl(\frac{-4}m\biggr)\frac{q^m}{1-q^{2m}}
=\sum_{\substack{m,n\ge1\\\text{$n$ неч.}}}\biggl(\frac{-4}m\biggr)q^{mn}
=\sum_{m,n\ge1}a(m)b(n)q^{mn},
\\
\text{где}\quad
a(m)=\biggl(\frac{-4}m\biggr), \quad b(n)=n\bmod2,
\end{gather*}
а через $\bigl(\frac{-4}m\bigr)$ обозначен квадратичный характер по модулю~4.

Тогда
\begin{align}
f(it)
&=\frac{\eta_8^4}{\eta_4^2}\,\frac{\eta_4^4}{\eta_8^2}\bigg|_{\tau=it}
=\frac{\eta_8^4}{\eta_4^2}\bigg|_{\tau=it}
\cdot\frac1{2t}\,\frac{\eta_8^4}{\eta_4^2}\bigg|_{\tau=i/(32t)}
\nonumber\\
&=\frac1{2t}\sum_{m_1,n_1\ge1}a(m_1)b(n_1)e^{-2\pi m_1n_1t}
\sum_{m_2,n_2\ge1}b(m_2)a(n_2)e^{-2\pi m_2n_2/(32t)},
\label{f32}
\end{align}
где $t>0$ и мы воспользовались инволюцией~\eqref{mod_inv}.

Далее
\begin{align*}
L(E,2)
&=L(f,2)
=\int_0^1f\,\log q\,\frac{\d q}q
=-4\pi^2\int_0^\infty f(it)t\,\d t
\displaybreak[2]\\
&=-2\pi^2\int_0^\infty\sum_{m_1,n_1,m_2,n_2\ge1}a(m_1)b(n_1)b(m_2)a(n_2)
\\ &\qquad\qquad\times
\exp\biggl(-2\pi\biggl(m_1n_1t+\frac{m_2n_2}{32t}\biggr)\biggr)\d t
\displaybreak[2]\\
&=-2\pi^2\sum_{m_1,n_1,m_2,n_2\ge1}a(m_1)b(n_1)b(m_2)a(n_2)
\\ &\qquad\times
\int_0^\infty\exp\biggl(-2\pi\biggl(m_1n_1t+\frac{m_2n_2}{32t}\biggr)\biggr)\d t.
\end{align*}
Новизна нашего метода --- в элементарном аналитическом преобразовании, которое представляет
собой замену переменной $t=n_2u/n_1$ во внутреннем интеграле. Такая замена не влияет
на подынтегральное выражение, но изменяет дифференциальную форму:
\begin{align*}
L(E,2)
&=-2\pi^2\sum_{m_1,n_1,m_2,n_2\ge1}\frac{a(m_1)b(n_1)b(m_2)a(n_2)n_2}{n_1}
\\ &\qquad\times
\int_0^\infty\exp\biggl(-2\pi\biggl(m_1n_2u+\frac{m_2n_1}{32u}\biggr)\biggr)\d u
\displaybreak[2]\\
&=-2\pi^2\int_0^\infty
\sum_{m_1,n_2\ge1} a(m_1)a(n_2)n_2 e^{-2\pi m_1n_2u}
\\ &\qquad\qquad\times
\sum_{m_2,n_1\ge1} \frac{b(m_2)b(n_1)}{n_1} e^{-2\pi m_2n_1/(32u)}\d u.
\end{align*}
Первый двойной ряд под интегралом отвечает
\begin{equation*}
\sum_{m,n\ge1}a(m)a(n)n\,q^{mn}
=\sum_{m,n\ge1}\biggl(\frac{-4}{mn}\biggr)n\,q^{mn}
=\sum_{n\ge1}n\biggl(\frac{-4}n\biggr)\frac{nq^n}{1+q^{2n}}
=\frac{\eta_2^4\eta_8^4}{\eta_4^4},
\end{equation*}
а второй отвечает
\begin{align*}
\sum_{m,n\ge1}\frac{b(m)b(n)}n\,q^{mn}
&=\sum_{m,n\ge1}\frac{q^{mn}}n-\frac{q^{(2m)n}}n-\frac{q^{m(2n)}}{2n}+\frac{q^{(2m)(2n)}}{2n}
\\
&=\frac12\sum_{m,n\ge1}\frac{2q^{mn}-3q^{2mn}+q^{4mn}}n
\\
&=-\frac12\,\log\prod_{m\ge1}\frac{(1-q^m)^2(1-q^{4m})}{(1-q^{2m})^3}
=-\frac12\,\log\frac{\eta_1^2\eta_4}{\eta_2^3},
\end{align*}
так что
\begin{equation*}
L(E,2)
=\pi^2\int_0^\infty\frac{\eta_2^4\eta_8^4}{\eta_4^4}\bigg|_{\tau=iu}
\cdot\log\frac{\eta_1^2\eta_4}{\eta_2^3}\bigg|_{\tau=i/(32u)}\d u.
\end{equation*}
Применяя инволюцию \eqref{mod_inv} к эта-произведению под зн\'аком логарифма, мы получаем
\begin{equation*}
L(E,2)
=\pi^2\int_0^\infty\frac{\eta_2^4\eta_8^4}{\eta_4^4}
\,\log\frac{\sqrt2\eta_8\eta_{32}^2}{\eta_{16}^3}\bigg|_{\tau=iu}\d u.
\end{equation*}

На этом этапе в действие вступает ``модулярная магия'': выбирая модулярную функцию
$x(\tau)=\eta_2^4\eta_8^2/\eta_4^6$,
изменяющуюся от~0 до~1 при изменении $\tau$ от~0 до~$i\infty$,
несложно показать, что
$$
\frac1{2\pi i}\,\frac{x\,\d x}{2\sqrt{1-x^4}}
=-\frac{\eta_2^4\eta_8^4}{\eta_4^4}\,\d\tau
\quad\text{и}\quad
\biggl(\frac{\sqrt2\eta_8\eta_{32}^2}{\eta_{16}^3}\biggr)^2
=\frac{1-x}{1+x}.
$$

Подытоживая, мы получаем следующее интегральное представление.

\begin{proposition}
\label{ch7:th1}
Для эллиптической кривой $E$ кондуктора~$32$ имеет место формула
\begin{equation*}
L(E,2)
=\frac\pi8\int_0^1\frac{x}{\sqrt{1-x^4}}\,\log\frac{1+x}{1-x}\,\d x
=0.9170506353\dotsc.
\end{equation*}
\end{proposition}

\section{Общий случай}
\label{ch7:s2}

По-существу, проделанные нами манипуляции с $L(E,2)=L(f,2)$ в \S\,\ref{ch7:s1} сводятся к следующему:
сначала мы представляем параболическую форму $f(\tau)$ в виде произведения двух
рядов Эйзенштейна веса~1, а в завершении мы получаем произведение опять же двух рядов Эйзенштейна,
$g_2(\tau)$ и $g_0(\tau)$, но уже веса 2 и~0 соответственно, так что
$L(f,2)=c\pi L(g_2g_0,1)$ с некоторой алгебраической постоянной~$c$. Полученный объект обречен
быть периодом в смысле \cite{KZ}: действительно, в соответствии с теорией модулярных форм
$g_0(\tau)$ является логарифмом модулярной функции, а
$2\pi i\,g_2(\tau)\,\d\tau$ с точностью до умножения на модулярную функцию является дифференциалом модулярной функции;
наконец, любые две модулярные функции связаны алгебраическим соотношением над полем~$\overline{\mathbb Q}$ алгебраических чисел.

Описанный метод может быть формализован следующим образом.

Для двух \emph{ограниченных} последовательностей $a(m)$, $b(n)$
будем называть выражения вида
\begin{equation}
g_k(\tau)=a+\sum_{m,n\ge1}a(m)b(n)n^{k-1}q^{mn}
\label{k03}
\end{equation}
почти рядами Эйзенштейна веса $k$, при условии, что
$g_k(\tau)$ является модулярной формой какого-либо уровня, т.е.\ она преобразуется
в другой объект подобного же типа под действием инволюции
$\tau\mapsto-1/(N\tau)$ с некоторым выбором целого положительного~$N$.
Подобное выполняется автоматически в случае, когда $k>0$ и $g_k(\tau)$ является рядом Эйзенштейна
(например, с выбором $a(m)=1$ и $b(n)$ --- характера Дирихле по модулю~$N$
заданной четности, $b(-1)=(-1)^k$);
тогда $\widehat g_k(\tau)=g_k(-1/(N\tau))(\sqrt{-N}\tau)^{-k}$ опять же является рядом Эйзенштейна.
Приведенное определение почти рядов Эйзенштейна имеет замечательный смысл и для любого неотрицательного $k$.
Действительно, в случае $k=0$ мы получаем слабые модулярные формы $g_0(\tau)$ веса~0 --- логарифмы
модулярных функций, представимых в виде бесконечных произведений наподобие произведения для эта-функции Дедекинда.
Так же и в случае $k\le0$, примеры почти рядов Эйзенштейна даются так называемыми эйхлеровыми интегралами ---
первообразными порядка $1-k$ голоморфных рядов Эйзенштейна веса $2-k$; соответствующие правила
преобразования под действием модулярной группы следуют из знаменитой леммы Гекке~\cite[\S\,5]{We}.

Предположим теперь, что нас интересует $L$-значение $L(f,k_0)$ параболической формы $f(\tau)$ веса $k=k_1+k_2$,
которая может быть представлена в виде произведения (в общем случае в виде линейных комбинаций нескольких произведений)
двух почти рядов Эйзенштейна $g_{k_1}(\tau)$ и $\widehat g_{k_2}(\tau)$,
в котором первый ряд обращается в нуль на бесконечности (т.е.\ $a=g_{k_1}(i\infty)=0$ в~\eqref{k03}),
а второй --- в нуле (т.е.\ $\widehat g_{k_2}(i0)=0$). (Такое поведение на бесконечности и в нуле --- частные
случаи параболических точек --- может быть достигнуто, поскольку исходная модулярная форма $f(\tau)$ параболическая,
значит обращается в нуль во всех параболических точках.) В действительности нас интересует ряд не для самого $\widehat g_{k_2}(\tau)$,
а для его образа $g_{k_2}(\tau)=\widehat g_{k_2}(-1/(N\tau))(\sqrt{-N}\tau)^{-k_2}$ под действием инволюции:
$$
g_{k_1}(\tau)=\sum_{m,n\ge1}a_1(m)b_1(n)n^{k_1-1}q^{mn}
\quad\text{и}\quad
g_{k_2}(\tau)=\sum_{m,n\ge1}a_2(m)b_2(n)n^{k_2-1}q^{mn}.
$$

Имеем
\begin{align*}
L(f,k_0)
&=L(g_{k_1}\widehat g_{k_2},k_0)
=\frac1{(k_0-1)!}\int_0^1 g_{k_1}\widehat g_{k_2}\log^{k_0-1}q\,\frac{\d q}q
\\
&=\frac{(-1)^{k_0-1}(2\pi)^{k_0}}{(k_0-1)!}
\int_0^\infty g_{k_1}(it)\widehat g_{k_2}(it)t^{k_0-1}\,\d t
\displaybreak[2]\\
&=\frac{(-1)^{k_0-1}(2\pi)^{k_0}}{(k_0-1)!\,N^{k_2/2}}
\int_0^\infty g_{k_1}(it)g_{k_2}(i/(Nt))t^{k_0-k_2-1}\,\d t
\displaybreak[2]\\
&=\frac{(-1)^{k_0-1}(2\pi)^{k_0}}{(k_0-1)!\,N^{k_2/2}}
\int_{0}^\infty\sum_{m_1,n_1\ge1}a_1(m_1)b_1(n_1)n_1^{k_1-1}e^{-2\pi m_1n_1t}
\\ &\qquad\qquad\times
\sum_{m_2,n_2\ge1}a_2(m_2)b_2(n_2)n_2^{k_2-1}e^{-2\pi m_2n_2/(Nt)}t^{k_0-k_2-1}\d t
\displaybreak[2]\\
&=\frac{(-1)^{k_0-1}(2\pi)^{k_0}}{(k_0-1)!\,N^{k_2/2}}
\sum_{m_1,n_1,m_2,n_2\ge1}a_1(m_1)b_1(n_1)a_2(m_2)b_2(n_2)n_1^{k_1-1}n_2^{k_2-1}
\\ &\qquad\times
\int_{0}^\infty\exp\biggl(-2\pi\biggl(m_1n_1t+\frac{m_2n_2}{Nt}\biggr)\biggr)t^{k_0-k_2-1}\d t;
\end{align*}
замена интегрирования и суммирования является законной в виду экспоненциального
убывания подынтегрального выражения в концевых точках.
Осуществляя знакомую замену переменных $t=n_2u/n_1$
и меняя местами суммирование и интегрирование в обратном направлении, получаем
\begin{align*}
L(f,k_0)
&=\frac{(-1)^{k_0-1}(2\pi)^{k_0}}{(k_0-1)!\,N^{k_2/2}}
\sum_{m_1,n_1,m_2,n_2\ge1}a_1(m_1)b_1(n_1)a_2(m_2)b_2(n_2)n_1^{k_1+k_2-k_0-1}n_2^{k_0-1}
\\ &\qquad\times
\int_{0}^\infty\exp\biggl(-2\pi\biggl(m_1n_2u+\frac{m_2n_1}{Nu}\biggr)\biggr)u^{k_0-k_2-1}\d u
\displaybreak[2]\\
&=\frac{(-1)^{k_0-1}(2\pi)^{k_0}}{(k_0-1)!\,N^{k_2/2}}
\int_{0}^\infty\sum_{m_1,n_2\ge1}a_1(m_1)b_2(n_2)n_2^{k_0-1}e^{-2\pi m_1n_2u}
\\ &\qquad\qquad\times
\sum_{m_2,n_1\ge1}a_2(m_2)b_1(n_1)n_1^{k_1+k_2-k_0-1}e^{-2\pi m_2n_1/(Nu)}u^{k_0-k_2-1}\d u
\displaybreak[2]\\
&=\frac{(-1)^{k_0-1}(2\pi)^{k_0}}{(k_0-1)!\,N^{k_2/2}}
\int_0^\infty g_{k_0}(iu)g_{k_1+k_2-k_0}(i/(Nu))u^{k_0-k_2-1}\,\d u.
\end{align*}
В предположении наличия модулярного преобразования для почти ряда Эйзенштейна
$g_{k_1+k_2-k_0}(\tau)$ под действием инволюции $\tau\mapsto-1/(N\tau)$, мы можем записать
полученный интеграл в виде $c\pi^{k_0-k_1}L(g_{k_0}\widehat g_{k_1+k_2-k_0},k_1)$,
где $c$~--- некоторый алгебраический множитель
(плюс возможно дополнительные слагаемые, если $g_{k_1+k_2-k_0}(\tau)$ является эйхлеровым интегралом).
Альтернативно, если $g_{k_0}(\tau)$ удовлетворяет преобразованию при модулярной инволюции,
то после дополнительной замены переменной $v=1/(Nu)$ мы отождествляем исходное $L$-значение
с~$c\pi^{k_0-k_2}L(\widehat g_{k_0}g_{k_1+k_2-k_0},k_2)$. В обоих случаях мы получаем
нетривиальное равенство, связывающее исходное $L$-значение $L(f,k_0)$
с новым ``$L$-значением'', отвечающим (почти) модулярной форме того же веса.

Частный случай $k_1=k_2=1$, $k_0=2$, обсуждавшийся в \cite{RZ1}, \cite{RZ2} и в \S\,\ref{ch7:s1},
позволяет редуцировать соответствующие $L$-значения к периодам.
Как мы увидим в \S\,\ref{ch7:s3}, подобная редукция может достигаться
и в более общей ситуации с использованием неоднородных линейных уравнений, которым удовлетворяют
эйхлеровы интегралы.

\section({Формула Рамануджана и \$L(E,3)\$})%
{Формула Рамануджана и $L(E,3)$}
\label{ch7:s3}

Для преобразования $L(E,3)$, отвечающего эллиптической кривой кондуктора 32, мы опять воспользуемся
равенством $L(E,3)=L(f,3)$, где $f(\tau)=\eta_4^2\eta_8^2$, и запишем разложение~\eqref{f32} в следующем виде:
$$
f(it)
=\frac1{2t}\sum_{m_1,n_1\ge1}b(m_1)a(n_1)e^{-2\pi m_1n_1t}
\sum_{m_2,n_2\ge1}b(m_2)a(n_2)e^{-2\pi m_2n_2/(32t)}.
$$
Тогда
\begin{align*}
L(E,3)
&=L(f,3)
=\frac12\int_0^1f\,\log^2q\,\frac{\d q}q
=4\pi^3\int_0^\infty f(it)t^2\,\d t
\displaybreak[2]\\
&=2\pi^3\int_0^\infty\sum_{m_1,n_1,m_2,n_2\ge1}b(m_1)a(n_1)b(m_2)a(n_2)
\\ &\qquad\qquad\times
\exp\biggl(-2\pi\biggl(m_1n_1t+\frac{m_2n_2}{32t}\biggr)\biggr)t\,\d t
\displaybreak[2]\\
&=2\pi^3\sum_{m_1,n_1,m_2,n_2\ge1}b(m_1)a(n_1)b(m_2)a(n_2)
\\ &\qquad\times
\int_0^\infty\exp\biggl(-2\pi\biggl(m_1n_1t+\frac{m_2n_2}{32t}\biggr)\biggr)t\,\d t
\\ \intertext{(осуществляем замену переменной $t=n_2u/n_1$)}
&=2\pi^3\sum_{m_1,n_1,m_2,n_2\ge1}\frac{b(m_1)a(n_1)b(m_2)a(n_2)n_2^2}{n_1^2}
\\ &\qquad\times
\int_0^\infty\exp\biggl(-2\pi\biggl(m_1n_2u+\frac{m_2n_1}{32u}\biggr)\biggr)u\,\d u
\displaybreak[2]\\
&=2\pi^3\int_0^\infty
\sum_{m_1,n_2\ge1} b(m_1)a(n_2)n_2^2e^{-2\pi m_1n_2u}
\\ &\qquad\qquad\times
\sum_{m_2,n_1\ge1} \frac{b(m_2)a(n_1)}{n_1^2} e^{-2\pi m_2n_1/(32u)}\,u\,\d u.
\end{align*}

Имеют место следующие разложения в ряды Эйзенштейна:
\begin{align*}
\sum_{m,n\ge1}b(m)a(n)n^2q^{mn}
&=\sum_{\substack{m,n\ge1\\\text{$m$ неч.}}}\biggl(\frac{-4}n\biggr)n^2q^{mn}
=\frac{\eta_2^8\eta_8^4}{\eta_4^6},
\\
\sum_{m,n\ge1}b(m)a(n)m^2q^{mn}
&=\sum_{\substack{m,n\ge1\\\text{$m$ неч.}}}\biggl(\frac{-4}n\biggr)m^2q^{mn}
=\frac{\eta_4^{18}}{\eta_2^8\eta_8^4},
\end{align*}
откуда
\begin{equation}
r(\tau)=\sum_{m,n\ge1}\frac{b(m)a(n)}{n^2}\,q^{mn}
=\delta^{-2}\biggl(\frac{\eta_4^{18}}{\eta_2^8\eta_8^4}\biggr).
\label{rtau}
\end{equation}

Продолжая наши преобразования,
\begin{align*}
L(E,3)
&=2\pi^3\int_0^\infty \frac{\eta_2^8\eta_8^4}{\eta_4^6}\bigg|_{\tau=iu}
\cdot r(i/(32u))\,u\,\d u
\\ \intertext{(применяем инволюцию к эта-произведению)}
&=\frac{\pi^3}8\int_0^\infty\frac{\eta_4^4\eta_{16}^8}{\eta_8^6}\,r(\tau)\bigg|_{\tau=i/(32u)}
\frac{\d u}{u^2}
\\ \intertext{(осуществляем замену $u=1/(32v)$)}
&=4\pi^3\int_0^\infty\frac{\eta_4^4\eta_{16}^8}{\eta_8^6}\,r(\tau)\bigg|_{\tau=iv}\d v.
\end{align*}

На этом алгоритм из \S\,\ref{ch7:s2} завершается. Чтобы представить полученный интеграл
в виде периода, нам требуется еще один шаг. Как и в \S\,\ref{ch7:s1}, мы сделаем
модулярную подстановку, только на этот раз воспользуемся иной модулярной функцией $x(\tau)=4\eta_2^4\eta_8^8/\eta_4^{12}$,
которая изменяется от $0$ до $1$ при изменении $\tau$ от $i\infty$ до~$0$.
Тогда
$$
\delta x=\frac{4\eta_2^{12}\eta_8^8}{\eta_4^{16}},
\quad
(1-x^2)^{1/4}=\frac{\eta_2^4\eta_8^2}{\eta_4^6},
\quad
s(x)=\frac{(1-\sqrt{1-x^2})^2}{x(1-x^2)^{3/4}}
=\frac{16\eta_4^{10}\eta_{16}^8}{\eta_2^8\eta_8^{10}}.
$$
Кроме того, результат подстановки $z=x^2(\tau)$ в гипергеометрическую функцию
$$
F(z)={}_2F_1\biggl(\begin{matrix} \frac12, \, \frac12 \\ 1 \end{matrix}\biggm| z \biggr)
=\frac2\pi\int_0^1\frac{\d y}{\sqrt{(1-y^2)(1-zy^2)}}
$$
есть модулярная форма
$$
\varphi(\tau)
=F(x^2)
=\sum_{n=0}^\infty{\binom{2n}n}^2\Bigl(\frac x4\Bigr)^{2n}
=\frac{\eta_4^{10}}{\eta_2^4\eta_8^4}
$$
веса 1. Поскольку $F(z)$, а также $F(1-z)$, удовлетворяют гипергеометрическому дифференциальному уравнению
$$
z(1-z)\frac{\d^2F}{\d z^2}+(1-2z)\frac{\d F}{\d z}-\frac14F=0,
$$
несложно выписать и дифференциальный оператор второго порядка
$$
\mathcal L=x(1-x^2)\frac{\d^2}{\d x^2}+(1-3x^2)\frac{\d}{\d x}-x
$$
(в терминах $x$), для которого $\mathcal L\varphi=0$.

В этих обозначениях
\begin{align}
L(E,3)
&=\pi^3\int_0^\infty\frac{\eta_4^{10}\eta_{16}^8}{\eta_2^8\eta_8^{10}}\,
\varphi(\tau)\,r(\tau)\delta x\bigg|_{\tau=iv}\d v
\nonumber\\
&=\frac{\pi^3}{16}\int_0^\infty s(x(\tau))\,\varphi(\tau)\,r(\tau)\delta x\bigg|_{\tau=iv}\d v,
\label{LE3}
\end{align}
и на этом этапе заметим, что функция $h(\tau)=4\varphi(\tau)r(\tau)$
удовлетворяет неоднородному линейному уравнению
$$
\mathcal Lh=\frac1{1-x^2}\quad
\biggl(\text{которое есть не что иное \cite{ShVl}, \cite{Ya} как}\;
\frac{\delta^2r}{\delta x\cdot\varphi}=\frac{\eta_4^{24}}{4\eta_2^{16}\eta_8^8}\biggr),
$$
так что для этой функции можно применить метод вариации постоянных:
\begin{align*}
h
&=\frac\pi2\,\biggl(F(x^2)\int\frac{F(1-x^2)}{1-x^2}\,\d x
-F(1-x^2)\int\frac{F(x^2)}{1-x^2}\,\d x\biggr)
\\
&=\frac{\pi x}2\int_0^1\frac{F(x^2)F(1-x^2w^2)-F(1-x^2)F(x^2w^2)}{1-x^2w^2}\,\d w.
\\
&=x+\frac{5}9x^3+\frac{89}{225}x^5+\frac{381}{1225}x^7+\frac{25609}{99225}x^9+\frac{106405}{480249}x^{11}+\dotsb.
\end{align*}
В результате мы получаем выражение
$$
L(E,3)
=\frac{\pi^2}{128}\int_0^1s(x)\,h(x)\,\d x,
$$
которое очевидным образом записывается в виде (сложного) кратного действительного интеграла,
так что представление $L(E,3)$ в виде периода установлено.

Алгоритм реализации ряда Эйзенштейна отрицательного веса через решения неоднородных линейных уравнений
достаточно стандартен \cite{Ya} и успешно применим во многих аналогичных ситуациях.
Однако, возникший ряд Эйзенштейна \eqref{rtau} веса $-1$ допускает также другую реализацию,
благодаря формуле Рамануджана \cite[(2$\cdot$2)]{Duk}:
$$
r(\tau)=\sum_{\substack{m,n\ge1\\\text{$m$ неч.}}}\biggl(\frac{-4}n\biggr)\frac{q^{mn}}{n^2}
=\frac{\tilde x\,G(-{\tilde x}^2)}{4F(-{\tilde x}^2)},
$$
где $\tilde x(\tau)=4\eta_8^4/\eta_2^4$, $F(-{\tilde x}^2)=\eta_2^4/\eta_4^2$ и
$$
G(z)={}_3F_2\biggl(\begin{matrix} 1, \, 1, \, 1 \\ \frac32, \, \frac32 \end{matrix}\biggm| z \biggr)
=\frac14\int_0^1\!\!\!\int_0^1\frac{(1-x_1)^{-1/2}(1-x_2)^{-1/2}}{1-zx_1x_2}\,\d x_1\,\d x_2.
$$
Последний интеграл задает аналитическое продолжение гипергеометрического
$_3F_2$-ряда в область $\operatorname{Re}z<1$ (см.\ \cite[лемма~2]{Ne4});
замена $y=(1-x_1)^{1/2}$, $w=(1-x_2)^{1/2}$ приводит интеграл к виду
$$
G(z)=\int_0^1\!\!\!\int_0^1\frac{\d y\,\d w}{1-z(1-y^2)(1-w^2)}.
$$

Возвращаясь к модулярной функции $x(\tau)=4\eta_2^4\eta_8^8/\eta_4^{12}$ и замечая, что
$\tilde x=x/\sqrt{1-x^2}$, мы можем переписать \eqref{LE3} в виде
$$
L(E,3)
=\frac{\pi^3}{64}\int_0^\infty\frac{s(x(\tau))\,x(\tau)}{1-x(\tau)^2}\,G\biggl(-\frac{x(\tau)^2}{1-x(\tau)^2}\biggr)\delta x\bigg|_{\tau=iv}\d v.
$$
Осуществляя модулярную подстановку $x=x(\tau)$, мы получаем следующее интегральное представление.

\begin{proposition}
\label{ch7:th2}
Для эллиптической кривой $E$ кондуктора~$32$ имеет место формула
\begin{align*}
L(E,3)
&=\frac{\pi^2}{128}\int_0^1\frac{(1-\sqrt{1-x^2})^2}{(1-x^2)^{3/4}}\,\d x
\int_0^1\!\!\!\int_0^1\frac{\d y\,\d w}{1-x^2(1-(1-y^2)(1-w^2))}
\\
&=0.9826801478\dotsc.
\end{align*}
\end{proposition}

\section{Гипергеометрические представления}
\label{ch7:s4}

Важным достоинством~\cite{BoCr} интегралов в предложениях~\ref{ch7:th1} и~\ref{ch7:th2}
является возможность свести их к гипергеометрической форме.

\medskip
В интегральном представлении $L(E,2)$ в предложении~\ref{ch7:th1} запишем
$$
\log\frac{1+x}{1-x}
=\frac23x^3\,{}_2F_1\biggl(\begin{matrix} \frac34, \, 1 \\ \frac74 \end{matrix}\biggm| x^4 \biggr)
+2x\,{}_2F_1\biggl(\begin{matrix} \frac14, \, 1 \\ \frac54 \end{matrix}\biggm| x^4 \biggr)
$$
и произведем замену переменной $x^4=x_0$, чтобы получить
\begin{align*}
L(E,2)
&=\frac\pi{48}\int_0^1\biggl(x_0^{1/4}\,{}_2F_1\biggl(\begin{matrix} \frac34, \, 1 \\ \frac74 \end{matrix}\biggm| x_0 \biggr)
+3x_0^{-1/4}\,{}_2F_1\biggl(\begin{matrix} \frac14, \, 1 \\ \frac54 \end{matrix}\biggm| x_0 \biggr)\biggr)
(1-x_0)^{-1/2}\d x_0
\\
&=\frac{\pi^{3/2}}{48}\,\frac{\Gamma(\frac54)}{\Gamma(\frac74)}
\,{}_3F_2\biggl(\begin{matrix} \frac34, \, \frac54, \, 1 \\ \frac74, \, \frac74 \end{matrix}\biggm| 1 \biggr)
+\frac{\pi^{3/2}}{16}\,\frac{\Gamma(\frac34)}{\Gamma(\frac54)}
\,{}_3F_2\biggl(\begin{matrix} \frac14, \, \frac34, \, 1 \\ \frac54, \, \frac54 \end{matrix}\biggm| 1 \biggr);
\end{align*}
представление~\eqref{HLE2} теперь получается после применения преобразования Тома \cite[3.2.(1)]{Bai}
к каждому из двух $_3F_2$ рядов.

\medskip
В интеграле предложения~\ref{ch7:th2} положим $x_0=x^2$, $x_1=1-y^2$ и $x_2=1-w^2$:
\begin{align*}
L(E,3)
=\frac{\pi^2}{1024}
&\int_0^1x_0^{-1/2}\bigl((1-x_0)^{-3/4}-2(1-x_0)^{-1/4}+(1-x_0)^{1/4}\bigr)\,\d x_0
\\ &\qquad\times
\int_0^1\!\!\!\int_0^1\frac{(1-x_1)^{-1/2}(1-x_2)^{-1/2}}{1-x_0(1-x_1x_2)}\,\d x_1\,\d x_2.
\end{align*}

Сначала рассмотрим интеграл (см.\ \cite[1.5.(1) и 1.4.(1)]{Bai})
\begin{align*}
\int_0^1\frac{x_0^{-1/2}(1-x_0)^{a-1}}{1-x_0z}\,\d x_0
&=\frac{\Gamma(\frac12)\,\Gamma(a)}{\Gamma(a+\frac12)}\,
{}_2F_1\biggl(\begin{matrix} 1, \, \frac12 \\ a+\frac12 \end{matrix}\biggm| z \biggr)
\\
&=\frac{\Gamma(\frac12)\,\Gamma(a-1)}{\Gamma(a-\frac12)}\,
{}_2F_1\biggl(\begin{matrix} 1, \, \frac12 \\ 2-a \end{matrix}\biggm| 1-z \biggr)
\\ &\qquad
+\Gamma(a)\,\Gamma(1-a)\,(1-z)^{a-1}
{}_2F_1\biggl(\begin{matrix} a-\frac12, \, a \\ a \end{matrix}\biggm| 1-z \biggr)
\\
&=\frac{\sqrt\pi\,\Gamma(a-1)}{\Gamma(a-\frac12)}\,
{}_2F_1\biggl(\begin{matrix} 1, \, \frac12 \\ 2-a \end{matrix}\biggm| 1-z \biggr)
+\frac\pi{\sin\pi a}\,\frac{(1-z)^{a-1}}{z^{a-1/2}}
\end{align*}
при $z=1-x_1x_2$. Далее,
\begin{align*}
&
\int_0^1\!\!\!\int_0^1(1-x_1)^{-1/2}(1-x_2)^{-1/2}
{}_2F_1\biggl(\begin{matrix} 1, \, \frac12 \\ 2-a \end{matrix}\biggm| x_1x_2 \biggr)
\,\d x_1\,\d x_2
\\ &\quad
=4\,{}_4F_3\biggl(\begin{matrix} 1, \, 1, \, 1, \, \frac12 \\ 2-a, \, \frac32, \, \frac32 \end{matrix}\biggm| 1 \biggr)
\end{align*}
и
\begin{align*}
&
\int_0^1\!\!\!\int_0^1\frac{x_1^{a-1}(1-x_1)^{-1/2}x_2^{a-1}(1-x_2)^{-1/2}}{(1-x_1x_2)^{a-1/2}}\,\d x_1\,\d x_2
\\ &\quad
=\biggl(\frac{\Gamma(\frac12)\,\Gamma(a)}{\Gamma(a+\frac12)}\biggr)^2
{}_3F_2\biggl(\begin{matrix} a, \, a, \, a-\frac12 \\ a+\frac12, \, a+\frac12 \end{matrix}\biggm| 1 \biggr)
\\ &\quad
=2\pi^{1/2}\Gamma(a)\Gamma(\tfrac32-a)\cdot
{}_3F_2\biggl(\begin{matrix} \frac12, \, \frac12, \, \frac32-a \\ 1, \, \frac32 \end{matrix}\biggm| 1 \biggr),
\end{align*}
где на последнем этапе мы вновь воспользовались преобразованием Тома \cite[3.2.(1)]{Bai}.

Проведенное вычисление означает, что, полагая
\begin{align*}
I_1(a)&=\frac{\pi^{5/2}\Gamma(a-1)}{256\Gamma(a-\frac12)}\,
{}_4F_3\biggl(\begin{matrix} 1, \, 1, \, 1, \, \frac12 \\ 2-a, \, \frac32, \, \frac32 \end{matrix}\biggm| 1 \biggr),
\\
I_2(a)&=\frac{\pi^{7/2}\Gamma(a)\,\Gamma(\frac32-a)}{512\sin\pi a}\,
{}_3F_2\biggl(\begin{matrix} \frac12, \, \frac12, \, \frac32-a \\ 1, \, \frac32 \end{matrix}\biggm| 1 \biggr),
\end{align*}
находим
$$
L(E,3)
=\bigl(I_1(\tfrac14)-2I_1(\tfrac34)+I_1(\tfrac54)\bigr)
+\bigl(I_2(\tfrac14)-2I_2(\tfrac34)+I_2(\tfrac54)\bigr).
$$

Для $I_2(\frac34)$ суммирование Уотсона--Уиппла \cite[3.3.(1)]{Bai} приводит к
$$
{}_3F_2\biggl(\begin{matrix} \frac12, \, \frac12, \, \frac34 \\ 1, \, \frac32 \end{matrix}\biggm| 1 \biggr)
=\frac{\Gamma(\frac12)\,\Gamma(\frac54)}{\Gamma(\frac34)},
$$
так что $I_2(\tfrac34)=\pi^5/1024$. Кроме того,
\begin{align*}
I_2(\tfrac14)+I_2(\tfrac54)
&=\frac{\pi^{7/2}\Gamma(\frac14)\,\Gamma(\frac54)}{256\sqrt2}
\biggl({}_3F_2\biggl(\begin{matrix} \frac12, \, \frac12, \, \frac54 \\ 1, \, \frac32 \end{matrix}\biggm| 1 \biggr)
-{}_3F_2\biggl(\begin{matrix} \frac12, \, \frac12, \, \frac14 \\ 1, \, \frac32 \end{matrix}\biggm| 1 \biggr)\biggr)
\\
&=\frac{\pi^{7/2}\Gamma(\frac14)^2}{1024\sqrt2}
\sum_{n=0}^\infty\frac{(\frac12)_n^2\cdot\bigl((\frac54)_n-(\frac14)_n\bigr)}{(1)_n(\frac32)_nn!}
\displaybreak[2]\\
&=\frac{\pi^{7/2}\Gamma(\frac14)^2}{1024\sqrt2}
\sum_{n=1}^\infty\frac{(\frac12)_n^2(\frac54)_{n-1}}{(1)_n(\frac32)_n(n-1)!}
\displaybreak[2]\\
&=\frac{\pi^{7/2}\Gamma(\frac14)^2}{1024\sqrt2\cdot6}\,
{}_3F_2\biggl(\begin{matrix} \frac32, \, \frac32, \, \frac54 \\ 2, \, \frac52 \end{matrix}\biggm| 1 \biggr)
\\
&=\frac{\pi^{7/2}\Gamma(\frac14)^2}{1024\sqrt2\cdot6}\,\frac{\Gamma(\frac12)\,\Gamma(\frac14)\,\Gamma(\frac74)}{\Gamma(\frac54)^2}
=\frac{\pi^5}{512},
\end{align*}
где вновь было использовано суммирование Уотсона--Уиппла.

Подытоживая,
$L(E,3)=I_1(\tfrac14)-2I_1(\tfrac34)+I_1(\tfrac54)$, что в точности представляет из себя формулу~\eqref{HLE3}.

\medskip
Осуществленные вычисления приводят к удивительно схожим гипергеометрическим выражениям для $L(E,2)$ и $L(E,3)$.
В обозначении
$$
F_k(a)=\frac{\pi^{k-1/2}\Gamma(a)}{2^{3k-1}\Gamma(a+\frac12)}\,
\,{}_{k+1}F_k\biggl(\begin{matrix} \overbrace{1, \, \dots, \, 1}^{\text{$k$ раз}}, \, \frac12 \\
a+\frac12, \, \underbrace{\tfrac32, \, \dots, \, \tfrac32}_{\text{$k-1$ раз}} \end{matrix}\biggm| 1 \biggr)
$$
соотношения \eqref{HLE2} и \eqref{HLE3} могут быть переписаны в виде
\begin{equation}
L(E,2)=F_2(\tfrac54)+F_2(\tfrac34)
\quad\text{и}\quad
L(E,3)=F_3(\tfrac54)+2F_3(\tfrac34)+F_3(\tfrac14).
\label{LE23}
\end{equation}
Отметив еще формулу
$$
L(E,1)
=\frac{\pi^{-1/2}\Gamma(\frac14)^2}{8\sqrt2}
=\frac{\pi^{-1/2}\Gamma(\frac14)^2}{24\sqrt2}\,
\,{}_3F_2\biggl(\begin{matrix} 1, \, \frac12 \\ \frac74 \end{matrix}\biggm| 1 \biggr)
=2F_1(\tfrac54),
$$
можно заключить, что при $k=1$, 2 or~3 значение $L(E,k)$ записывается в виде (простой) $\mathbb Q$-линейной комбинации
$F_k(\frac74-\frac m2)$ с $m=1,\dots,k$. Экспериментальные вычисления показывают, что в случае $k>3$ подобного
простого выражения уже указать нельзя.

Формулы \eqref{LE23}, в свою очередь, могут быть переписаны как интегральные представления для $L(E,2)$
и $L(E,3)$, отличающиеся от приведенных в предложениях~\ref{ch7:th1} и~\ref{ch7:th2};
именно эти формулы указаны в выражениях \eqref{ILE2} и \eqref{ILE3} теоремы~\ref{th:7}.

\medskip
Отметим, что было бы интересно отождествить (линейные комбинации) гипергеометрических рядов
в~\eqref{HLE3} с линейными комбинациями логарифмических мер Малера многочленов, зависящих от трех переменных.
Результаты подобного типа обсуждаются в~\cite{Rog}, но там они приводятся в недостаточно гипергеометрической форме.

Принципиальная составляющая в выводе интегрального представления в теореме~\ref{ch7:th2}
--- уникальная гипергеометрическая форма для ряда Эйзенштейна отрицательного веса, доказанная в~\cite{Duk}.
Существуют ли другие подобные результаты?

\begingroup

\def\chaptername{}
\def\thechapter{}

\chapter{Заключительные замечания}
\label{final}

Основной целью \ifdisser диссертации \else монографии \fi является
изложение некоторых теорем автора, тесно связанных с теоремой Апери и
вопросами об арифметической природе значений дзета-функции Римана при целых
положительных значениях аргумента. В ``сердце'' всех этих результатов лежат
обобщенные гипергеометрические ряды, использование которых в теории чисел
и, в частности, в теории диофантовых приближений не является новым.
Новизна арифметико-гипергеометрического подхода --- полноценное применение гипергеометрии,
именно, различных формул суммирования и преобразования, интегральных
представлений и алгоритмов. Кроме того,
мы постарались собрать в \ifdisser диссертации \else монографии \fi
теоремы, обладающие ``особой эстетикой'': их доказательства несложные с технической точки зрения и достаточно короткие.

Автор не ставил целью полностью отразить современное состояние дел о применении гипергеометрических
рядов в теории диофантовых приближений даже касательно выбранной узкой тематики.
Важен уже тот факт, что интерес к ней не ослабевает --- за последнее десятилетие получено много новых результатов,
некоторые из которых дополняют (но не усиливают!) основные теоремы этой \ifdisser диссертации\else монографии\fi.
Мы также не обсуждаем ряд приложений гипергеометрической техники
в других задачах: полиномиальные рекурсии (типа Апери \eqref{eq:0.3}, \eqref{eq:0.12})
для математических постоянных \cite{CZ}, \cite{Z15}, \cite{Z23} и специальные дифференциальные уравнения
для периодов трехмерных многообразий Калаби--Яу~\cite{AZ}, \cite{ASZ}, формулы Рамануджана \cite{GuZ}, \cite{Z24}
для $1/\pi$ и их обобщения \cite{CWZ}, \cite{WaZ}, \cite{Z80}, кратные дзета-значения \cite{OZ}, \cite{Z13a},
и этот список далеко не полон. Обстоятельному изложению результатов в этих направлениях
посвящен написанный нами обзор ``Арифметические гипергеометрические ряды''~\cite{Z90}.

Возможности гипергеометрической методики далеко не исчерпаны, и мы искренне надеемся на появление ярких результатов,
синтезирующих гипергеометрию, арифметический метод и новые оригинальные идеи.

\bigskip
\begin{center}
{\it Feci quod potui}, {\it faciant meliora potentes}.
\end{center}

\endgroup

\newpage
\refstepcounter{chapter}

\def\bibname{Список литературы}


\begin{thebibliography}{PWZ}

\def\No{№}
\def\Bibitem[#1]#2{\bibitem{#2}}

\Bibitem[AZ]{AZ}
\textsc{G.~Almkvist, W.~Zudilin},
Differential equations, mirror maps and zeta values,
\emph{Mirror Symmetry V},
N.~Yui, S.-T.~Yau, and J.\,D.~Lewis (eds.),
AMS/IP Studies in Advanced Mathematics \textbf{38}
(International Press \& Amer. Math. Soc., Providence, RI 2007), 481--515.

\Bibitem[ASZ]{ASZ}
\textsc{G.~Almkvist, D.~van Straten, W.~Zudilin},
Generalizations of Clausen's formula and algebraic transformations of Calabi--Yau differential equations,
\emph{Proc. Edinburgh Math. Soc.} \textbf{54}:2 (2011), 273--295.

\Bibitem[An1]{An1}
\textsc{G.\,E.~Andrews},
Problems and prospects for basic hypergeometric functions,
\emph{Theory and application of special functions},
ed.~R.\,A.~Askey,
Proc. Advanced Sem., Math. Res. Center
(Univ. Wisconsin, Madison, Wis., 1975),
Math. Res. Center, Univ. Wisconsin, Publ. No.~35
(Academic Press, New York 1975), 191--224.

\Bibitem[An2]{An2}
\textsc{G.\,E.~Andrews},
The well-poised thread: An organized chronicle
of some amazing summations and their implications,
\emph{Ramanujan J.}
\textbf{1}:1 (1997), 7--23.

\Bibitem[AAR]{AAR}
\textsc{G.\,E.~Andrews, R.~Askey, R.~Roy},
\emph{Special functions},
Encyclopedia of Mathematics and its Applications
\textbf{71}
(Cambridge University Press, Cambridge 1999).

\Bibitem[AS]{AS}
\textsc{В.\,А.~Андросенко, В.\,Х.~Салихов},
Интеграл Марковеккио и мера иррациональности $\pi/\sqrt{3}$,
\emph{Вестник БГТУ} \textbf{34}:4 (2011), 129--132.

\Bibitem[Ap]{Ap}
\textsc{R.~Ap\'ery},
Irrationalit\'e de $\zeta(2)$ et $\zeta(3)$,
\emph{Ast\'erisque\/}
\textbf{61} (1979), 11--13.

\Bibitem[Ask]{Ask}
\textsc{R.~Askey},
The $q$-gamma and $q$-beta functions
\emph{Appl. Anal.} \textbf{8} (1978), 125--141.

\Bibitem[BBBC]{BBBC}
\textsc{D.\,H.~Bailey, D.~Borwein, J.\,M.~Borwein, R.\,E.~Crandall},
Hypergeometric forms for Ising-class integrals,
\emph{Experiment. Math.} \textbf{16}:3 (2007), 257--276.

\Bibitem[Bai0]{Bai0}
\textsc{W.\,N. Bailey},
Some transformations of generalized hypergeometric series, and contour-integrals of Barnes's type,
\emph{Quart. J. Math.} (\emph{Oxford}) \textbf{3}:1 (1932), 168--182.

\Bibitem[Bai]{Bai}
\textsc{W.\,N.~Bailey},
\emph{Generalized hypergeometric series},
Cambridge Math. Tracts
\textbf{32}
(Cambridge Univ. Press, Cambridge 1935);
2nd reprinted edition
(Stechert-Hafner, New York--London 1964).

\Bibitem[Bak]{Bak}
\textsc{A.~Baker},
The theory of linear forms in logarithms,
\emph{Transcendence theory: advances and applications},
Proc. Conf., Univ. Cambridge, Cambridge, 1976
(Academic Press, London 1977), 1--27.

\Bibitem[BC]{BC}
\textsc{A.~Baker, J.~Coates},
Fractional parts of powers of rationals,
\emph{Math. Proc. Cambridge Philos. Soc.}
\textbf{77} (1975), 269--279.

\Bibitem[BR]{BR}
\textsc{K.~Ball, T.~Rivoal},
Irrationalit\'e d'une infinit\'e de valeurs
de la fonction z\^eta aux entiers impairs,
\emph{Invent. Math.}
\textbf{146}:1 (2001), 193--207.

\Bibitem[Bar]{Bar}
\textsc{E.\,W.~Barnes},
A new development of the theory of
the hypergeometric functions,
\emph{Proc. London Math. Soc. II~Ser.}
\textbf{6} (1908), 141--177.

\Bibitem[Ben1]{Ben1}
\textsc{M.\,A.~Bennett},
Fractional parts of powers of rational numbers,
\emph{Math. Proc. Cambridge Philos. Soc.}
\textbf{114} (1993), 191--201.

\Bibitem[Ben2]{Ben2}
\textsc{M.\,A.~Bennett},
An ideal Waring problem with restricted summands,
\emph{Acta Arith.}
\textbf{66} (1994), 125--132.

\Bibitem[BBB]{BBB}
\textsc{J.\,L.~Berggren}, \textsc{J.\,M.~Borwein} and \textsc{P.~Borwein},
\emph{Pi: A source book},
3rd edn. (Springer, New York 2004).

\Bibitem[Be1]{Be1}
\textsc{F.~Beukers},
A note on the irrationality of~$\zeta(2)$ and~$\zeta(3)$,
\emph{Bull. London Math. Soc.}
\textbf{11}:3 (1979), 268--272.

\Bibitem[Beu]{Beu}
\textsc{F.~Beukers},
Fractional parts of powers of rationals,
\emph{Math. Proc. Cambridge Philos. Soc.}
\textbf{90} (1981), 13--20.

\Bibitem[Be2]{Be2}
\textsc{F.~Beukers},
Pad\'e approximations in number theory,
\emph{Lecture Notes in Math.}
\textbf{888} (Springer-Verlag, Berlin 1981), 90--99.

\Bibitem[Be3]{Be3}
\textsc{F.~Beukers},
Irrationality proofs using modular forms,
\emph{Ast\'erisque}
\textbf{147--148} (1987), 271--283.

\Bibitem[Bez]{Bez}
\textsc{J.-P.~B\'ezivin},
Ind\'ependence lin\'eaire des valeurs des solutions
transcendantes de certaines \'equations fonctionelles,
\emph{Manuscripta Math.}
\textbf{61} (1988), 103--129.

\Bibitem[BoCr]{BoCr}
\textsc{J.\,M.~Borwein} and \textsc{R.\,E.~Crandall},
Closed forms: what they are and why we care,
\emph{Notices Amer. Math. Soc.} \textbf{60} (2013), 50--65.

\Bibitem[Bo]{Bo}
\textsc{P.~Borwein},
On the irrationality of $\sum\frac1{q^n+r}$,
\emph{J. Number Theory}
\textbf{37} (1991), 253--259.

\Bibitem[Br]{Br}
\textsc{Н.\,Г.~де~Брёйн},
\emph{Асимптотические методы в анализе}
(Иностр. лит-ра, Москва 1961).

\Bibitem[BV]{BV}
\textsc{P.~Bundschuh, K.~V\"a\"an\"anen},
Arithmetical investigations of a certain infinite product,
\emph{Compositio Math.} \textbf{91} (1994), 175--199.

\Bibitem[BZ]{BZ}
\textsc{P.~Bundschuh, W.~Zudilin},
Irrationality measures for certain $q$-mathematical constants,
\emph{Math. Scand.} \textbf{101}:1 (2007), 104--122.

\Bibitem[CZ]{CZ}
\textsc{H.\,H.~Chan, W.~Zudilin},
New representations for Ap\'ery-like sequences,
\emph{Mathematika} \textbf{56}:1 (2010), 107--117.

\Bibitem[CWZ]{CWZ}
\textsc{H.\,H.~Chan, J.~Wan, W.~Zudilin},
Legendre polynomials and Ramanujan-type series for~$1/\pi$,
\emph{Israel J. Math.} \textbf{194}:1 (2013), 183--207.

\Bibitem[Ch]{Ch}
\textsc{G.\,V.~Chudnovsky},
Pad\'e approximations to the generalized hypergeometric functions. I,
\emph{J. Math. Pures Appl.} (9)
\textbf{58} (1979), 445--476.

\Bibitem[Ch2]{Ch2}
\textsc{G.\,V.~Chudnovsky},
On the method of Thue--Siegel,
\emph{Ann. of Math. II~Ser.}
\textbf{117}:2 (1983), 325--382.

\Bibitem[Da]{Da}
\textsc{Л.\,В.~Данилов},
Рациональные приближения некоторых функций в рациональных точках,
\emph{Матем. заметки}
\textbf{24}:4 (1978), 449--458.

\Bibitem[DD]{DD}
\textsc{F.~Delmer, J.-M.~Deshouillers},
The computation of $g(k)$ in Waring's problem,
\emph{Math. Comp.}
\textbf{54} (1990), 885--893.

\Bibitem[Du]{Du}
\textsc{А.\,К.~Дубицкас},
Оценка снизу величины $\|(3/2)^k\|$,
\emph{Успехи матем. наук}
\textbf{45}:1 (1990), 153--154.

\Bibitem[Duk]{Duk}
\textsc{W.~Duke},
Some entries in Ramanujan's notebooks,
\emph{Math. Proc. Camb. Phil. Soc.} \textbf{144} (2008), 255--266.

\Bibitem[Duv]{Duv}
\textsc{D.~Duverney},
Irrationalit\'e d'un $q$-analogue de $\zeta(2)$,
\emph{C.~R. Acad. Sci. Paris S\'er.~I Math.}
\textbf{321}:10 (1995), 1287--1289.

\Bibitem[DV]{DV}
\textsc{R.~Dvornicich, C.~Viola},
Some remarks on Beukers' integrals,
\emph{Colloq. Math. Soc. J\'anos Bolyai}
\textbf{51} (North-Holland, Amsterdam 1987), 637--657.

\Bibitem[Er]{Er}
\textsc{P.~Erd\H{o}s},
On arithmetical properties of Lambert series,
\emph{J. Indiana Math. Soc.}
\textbf{12} (1948), 63--66.

\Bibitem[Eu1]{Eu1}
\textsc{L.~Euler},
\emph{Comm. Acad. Sci. Imp. Petropol.}
\textbf{9} (1737), 160--188;
Reprint, 
\emph{Opera Omnia Ser.~I}
\textbf{15} (Teubner, Berlin 1927), 217--267.

\Bibitem[Eu2]{Eu2}
\textsc{L.~Euler},
Meditationes circa singulare serierum genus,
\emph{Novi Comm. Acad. Sci. Petropol.}
\textbf{20} (1775), 140--186;
Reprint,
\emph{Opera Omnia Ser.~I}
\textbf{15} (Teubner, Berlin 1927), 217--267.

\Bibitem[Ex]{Ex}
\textsc{H.~Exton},
\emph{$q$-Hypergeometric functions and applications},
Ellis Horwood Ser. Math. Appl.
(Ellis Horwood Ltd., Chichester 1983).

\Bibitem[Fin]{Fin}
\textsc{S.\,R.~Finch},
\emph{Mathematical constants},
Encyclopedia of Mathematics and its Applications
\textbf{94}
(Cambridge University Press, Cambridge 2003).

\Bibitem[Fi]{Fi}
\textsc{N.\,J.~Fine},
\emph{Basic hypergeometric series and applications},
Math. Surveys Monographs \textbf{27}
(Amer. Math. Soc., Providence, RI 1988).

\Bibitem[Fis]{Fis}
\textsc{S.~Fischler},
Irrationalit\'e de valeurs de z\^eta [d'apr\`es Ap\'ery, Rivoal, ...],
\emph{Ast\'erisque\/}
\textbf{294} (2004), 27--62.

\Bibitem[FZ]{FZ}
\textsc{S.~Fischler, W.~Zudilin},
A refinement of Nesterenko's linear independence criterion with applications to zeta values,
\emph{Math. Annalen} \textbf{347}:4 (2010), 739--763.

\Bibitem[Fr]{Fr}
\textsc{J.~Franel},
\emph{L'interm\'ediare des math\'ematiciens\/}
\textbf{1} (Gauthier-Villars, Paris 1894), 45--47;
Response 170,
\emph{L'interm\'ediare des math\'ematiciens\/}
\textbf{2} (Gauthier-Villars, Paris 1895), 33--35.

\Bibitem[GR]{GR}
\textsc{Г.~Гаспер, М.~Рахман},
\emph{Базисные гипергеометрические ряды}
(Мир, Москва 1993).

\Bibitem[Ge]{Ge}
\textsc{А.\,О.~Гельфонд},
\emph{Исчисление конечных разностей}, 3-е изд.
(Наука, Москва 1967).

\Bibitem[GKP]{GKP}
\textsc{R.\,L.~Graham, D.\,E.~Knuth, O.~Patashnik},
\emph{Concrete mathematics. A foundation for computer science},
2nd edition
(Addison-Wesley Publishing Company, Reading, MA 1994).

\Bibitem[Gui]{Gui}
\textsc{J.~Guillera},
WZ-proofs of ``divergent'' Ramanujan-type series,
\emph{Advances in Combinatorics, Waterloo Workshop in Computer Algebra, W80} (May 26--29, 2011),
I.~Kotsireas and E.\,V.~Zima (eds.) (Springer, New York 2013), 187--195.

\Bibitem[GuZ]{GuZ}
\textsc{J.~Guillera, W.~Zudilin},
Ramanujan-type formulae for $1/\pi$: The art of translation,
\emph{The Legacy of Srinivasa Ramanujan},
R.~Balasubramanian et al. (eds.),
Ramanujan Math. Soc. Lecture Notes Series \textbf{20} (Ramanujan Math. Soc., Mysore 2013), 181--195.

\Bibitem[Gu]{Gu}
\textsc{Л.\,А.~Гутник},
Об иррациональности некоторых величин, содержащих $\zeta(3)$,
\emph{Успехи матем. наук} \textbf{34}:3 (1979), 190;
\emph{Acta Arith.} \textbf{42}:3 (1983), 255--264.

\Bibitem[Ha]{Ha}
\textsc{L.~Habsieger},
Explicit lower bounds for $\|(3/2)^k\|$,
\emph{Acta Arith.}
\textbf{106} (2003), 299--309.

\Bibitem[HW]{HW}
\textsc{G.\,H.~Hardy, E.\,M.~Wright},
\emph{An introduction to the theory of numbers},
5th edition (Oxford Univ. Press, Oxford 1979).

\Bibitem[Ha1]{Ha1}
\textsc{M.~Hata},
Legendre type polynomials and irrationality measures,
\emph{J. Reine Angew. Math.}
\textbf{407}:1 (1990), 99--125.

\Bibitem[Ha1a]{Ha1a}
\textsc{M.~Hata},
Rational approximations to $\pi$ and some other numbers,
\emph{Acta Arith.} \textbf{63}:4 (1993), 335--349.

\Bibitem[Ha2]{Ha2}
\textsc{M.~Hata},
A note on Beukers' integral,
\emph{J. Austral. Math. Soc. Ser.~A}
\textbf{58}:2 (1995), 143--153.

\Bibitem[Ha3]{Ha3}
\textsc{M.~Hata},
A new irrationality measure for~$\zeta(3)$,
\emph{Acta Arith.} \textbf{92}:1 (2000), 47--57.


\Bibitem[He]{He}
\textsc{E.~Heine},
Untersuchungen \"uber die Reihe $\dots$,
\emph{J. Reine Angew. Math.} \textbf{34} (1847), 285--328.

\Bibitem[HP1]{HP1}
\textsc{Т.\,Г.~Хессами Пилеруд},
О линейной независимости векторов с полилогарифмическими координатами,
\emph{Вестник МГУ. Сер.~1. Матем., мех.}
\No~6 (1999), 54--56.

\Bibitem[HP2]{HP2}
\textsc{Т.\,Г.~Хессами Пилеруд},
\emph{Арифметические свойства значений гипергеометрических функций},
Дисс. канд. физ.-мат. наук
(Московский гос. ун-т, Москва 1999).

\Bibitem[JM]{JM}
\textsc{F.~Jouhet, E.~Mosaki},
Irrationalit\'e aux entiers impairs positifs d'un $q$-analogue de la fonction z\^eta de Riemann,
\emph{Intern. J. Number Theory} \textbf{6}:5 (2010), 959--988.

\Bibitem[KZ]{KZ}
\textsc{M.~Kontsevich, D.~Zagier},
Periods,
\emph{Mathematics unlimited\,---\,2001 and beyond} (Springer, Berlin 2001), 771--808.

\Bibitem[KR]{KR}
\textsc{C.~Krattenthaler, T.~Rivoal},
\emph{Hyperg\'eom\'etrie et fonction z\^eta de Riemann},
Mem. Amer. Math. Soc. \textbf{186}
(Amer. Math. Soc., Providence, RI 2007), no.~875.

\Bibitem[KRZ]{KRZ}
\textsc{C.~Krattenthaler, T.~Rivoal, W.~Zudilin},
S\'eries hyperg\'eom\'etriques basiques,
$q$-analogues des valeurs de la fonction z\^eta et formes modulaires,
\emph{Inst. Jussieu Math. J.}
\textbf{5}:1 (2006), 53--79.

\Bibitem[KW]{KW}
\textsc{J.~Kubina, M.~Wunderlich},
Extending Waring's conjecture up to $471600000$,
\emph{Math. Comp.}
\textbf{55} (1990), 815--820.

\Bibitem[Li]{Li}
\textsc{F.~Lindemann},
\"Uber die Zalh~$\pi$,
\emph{Math. Annalen}
\textbf{20} (1882), 213--225.

\Bibitem[Lu]{Lu}
\textsc{Ю.~Люк},
\emph{Специальные математические функции и их аппроксимации}
(Мир, Москва 1980).

\Bibitem[Ma]{Ma}
\textsc{K.~Mahler},
On the fractional parts of powers of real numbers,
\emph{Mathematika}
\textbf{4} (1957), 122--124.

\Bibitem[Mar]{Mar}
\textsc{R.~Marcovecchio},
The Rhin--Viola method for $\log2$,
\emph{Acta Arith.} \textbf{139}:2 (2009), 147--184.

\Bibitem[MV]{MV}
\textsc{T.~Matala-aho, K.~V\"a\"an\"anen},
On approximation measures of $q$-logarithms,
\emph{Bull. Austral. Math. Soc.} \textbf{58} (1998), 15--31.

\Bibitem[MVZ]{MVZ}
\textsc{T.~Matala-aho, K.~V\"a\"an\"anen, W.~Zudilin},
New irrationality measures for $q$-logarithms,
\emph{Math. of Comput.} \textbf{75} (2006), no.~254, 879--889.

\Bibitem[Me]{Me}
\textsc{F.~Mertens},
Ueber einige asymptotische Gesetze der Zahlentheorie,
\emph{J. Reine Angew. Math.}
\textbf{77}:4 (1874), 289--338.

\Bibitem[Ne1]{Ne1}
\textsc{Ю.\,В.~Нестеренко},
О линейной независимости чисел,
\emph{Вестник МГУ. Сер.~1. Матем., мех.}
\No~1 (1985), 46--54.

\Bibitem[Ne2]{Ne2}
\textsc{Ю.\,В.~Нестеренко},
Некоторые замечания о~$\zeta(3)$,
\emph{Матем. заметки}
\textbf{59}:6 (1996), 865--880.

\Bibitem[Ne3]{Ne3}
\textsc{Ю.\,В.~Нестеренко},
Модулярные функции и вопросы трансцендентности,
\emph{Матем. сб.} \textbf{187}:9 (1996), 65--96.

\Bibitem[Ne4]{Ne4}
\textsc{Yu.\,V.~Nesterenko},
Integral identities and constructions of approximations to zeta values,
\emph{J. Th\'eor. Nombres Bordeaux}
\textbf{15}:2 (2003), 535--550.

\Bibitem[Ne5]{Ne5}
\textsc{Ю.\,В.~Нестеренко},
О показателе иррациональности числа~$\ln2$,
\emph{Матем. заметки}
\textbf{88}:4 (2010), 549--564.

\Bibitem[Ni]{Ni}
\textsc{Е.\,М.~Никишин},
Об иррациональности значений функций $F(x,s)$,
\emph{Матем. сб.} \textbf{109}:3 (1979), 410--417.

\Bibitem[OZ]{OZ}
\textsc{Y.~Ohno, W.~Zudilin},
Zeta stars,
\emph{Commun. Number Theory Phys.} \textbf{2}:2 (2008), 325--347.

\Bibitem[PWZ]{PWZ}
\textsc{M.~Petkov\v sek, H.\,S.~Wilf, D.~Zeilberger},
\emph{$A=B$}
(A.\,K.~Peters, Ltd., Wellesley 1996).

\Bibitem[PS]{PS}
\textsc{Г.~Полиа, Г.~Сеге},
\emph{Задачи и теоремы из анализа}, Часть~II
(Наука, Москва 1978).

\Bibitem[Po]{Po}
\textsc{A.~van der Poorten},
A proof that Euler missed...
Ap\'ery's proof of the irrationality of~$\zeta(3)$,
\emph{Math. Intelligencer}
\textbf{1}:4 (1978/79), 195--203.

\Bibitem[Pr]{Pr}
\textsc{M.~Pr\'evost},
A new proof of the irrationality of $\zeta(3)$ using Pad\'e approximants,
\emph{J. Comput. Appl. Math.}
\textbf{67} (1996), 219--235.

\Bibitem[Pu1]{Pu1}
\textsc{Ю.\,А.~Пупырёв},
О линейной и алгебраической независимости $q$-дзета-значений,
\emph{Матем. заметки} \textbf{78}:4 (2005), 608--613.

\Bibitem[Pu2]{Pu2}
\textsc{Ю.\,А.~Пупырёв},
Эффективизация нижней оценки для $\|(4/3)^k\|$,
\emph{Матем. заметки} \textbf{85}:6 (2009), 927--935.

\Bibitem[RV1]{RV1}
\textsc{G.~Rhin, C.~Viola},
On the irrationality measure of~$\zeta(2)$,
\emph{Ann. Inst. Fourier} (\emph{Grenoble})
\textbf{43}:1 (1993), 85--109.

\Bibitem[RV2]{RV2}
\textsc{G.~Rhin, C.~Viola},
On a permutation group related to~$\zeta(2)$,
\emph{Acta Arith.}
\textbf{77}:1 (1996), 23--56.

\Bibitem[RV3]{RV3}
\textsc{G.~Rhin, C.~Viola},
The group structure for~$\zeta(3)$,
\emph{Acta Arith.}
\textbf{97}:3 (2001), 269--293.

\Bibitem[Rie]{Rie}
\textsc{Б.~Риман},
О числе простых чисел, не превышающих данной величины,
\emph{Сочинения} (ОГИЗ, Москва 1948), 216--224.

\Bibitem[Ri1]{Ri1}
\textsc{T.~Rivoal},
La fonction z\^eta de Riemann prend une infinit\'e
de valeurs irrationnelles aux entiers impairs,
\emph{C.~R. Acad. Sci. Paris S\'er.~I Math.}
\textbf{331}:4 (2000), 267--270.

\Bibitem[Ri2]{Ri2}
\textsc{T.~Rivoal},
\emph{Propri\'et\'es diophantiennes des valeurs
de la fonction z\^eta de Riemann aux entiers impairs},
Th\`ese de Doctorat
(Univ. de Caen, Caen 2001).

\Bibitem[Ri3]{Ri3}
\textsc{T.~Rivoal},
Irrationalit\'e d'au moins un des neuf nombres
$\zeta(5),\zeta(7),\dots,\zeta(21)$,
\emph{Acta Arith.} \textbf{103} (2002), 157--167

\Bibitem[RZ]{RZ}
\textsc{T.~Rivoal, W.~Zudilin},
Diophantine properties of numbers related to Catalan's constant,
\emph{Math. Annalen} \textbf{326}:4 (2003), 705--721.

\Bibitem[Rog]{Rog}
\textsc{M.\,D.~Rogers},
A study of inverse trigonometric integrals associated with three-variable Mahler measures, and some related identities,
\emph{J. Number Theory} \textbf{121} (2006), 265--304.

\Bibitem[RZ1]{RZ1}
\textsc{M.~Rogers, W.~Zudilin},
From $L$-series of elliptic curves to Mahler measures,
\emph{Compositio Math.} \textbf{148} (2012), 385--414.

\Bibitem[RZ2]{RZ2}
\textsc{M.~Rogers, W.~Zudilin},
On the Mahler measure of $1+X+1/X+Y+1/Y$,
\emph{Intern. Math. Research Notices} \textbf{2014} (2014), in press (\texttt{doi:\,10.1093/\allowbreak imrn/rns285}), 22~pp.

\Bibitem[Ru]{Ru}
\textsc{Е.\,А.~Рухадзе},
Оценка снизу приближения~$\ln2$ рациональными числами,
\emph{Вестник МГУ. Сер.~1. Матем., мех.}
\No~6 (1987), 25--29.

\Bibitem[Sa1]{Sa1}
\textsc{В.\,Х.~Салихов},
О мере иррациональности $\log3$,
\emph{Докл. РАН} \textbf{417}:6 (2007), 753--755

\Bibitem[Sa2]{Sa2}
\textsc{В.\,Х.~Салихов},
О мере иррациональности числа $\pi$,
\emph{Успехи матем. наук}
\textbf{63}:3 (2008), 163--164.

\Bibitem[SF]{SF}
\textsc{В.\,Х.~Салихов, А.\,И.~Фроловичев},
О кратных интегралах, представимых в виде линейной формы от $1,\zeta(3),\zeta(5),\dots,\zeta(2k-1)$,
\emph{Фундамент. и прикл. матем.} \textbf{11}:6 (2005), 143--178.

\Bibitem[Sc1]{Sc1}
\textsc{A.\,L.~Schmidt},
Generalized $q$-Legendre polynomials,
\emph{J. Comput. Appl. Math.} \textbf{49}:1--3 (1993), 243--249.

\Bibitem[Sc2]{Sc2}
\textsc{A.\,L.~Schmidt},
Legendre transforms and Ap\'ery's sequences,
\emph{J. Austral. Math. Soc. Ser.~A\/} \textbf{58}:3 (1995), 358--375.

\Bibitem[Sh]{Sh}
\textsc{А.\,Б.~Шидловский},
\emph{Трансцендентные числа}
(Наука, Москва 1987).

\Bibitem[ShVl]{ShVl}
\textsc{E.~Shinder, M.~Vlasenko},
Linear Mahler measures and double $L$-values of modular forms,
\emph{Preprint} \texttt{http://arxiv.org/abs/1206.1454} (2012), 22~pp.

\Bibitem[Sl]{Sl}
\textsc{L.\,J.~Slater},
\emph{Generalized hypergeometric functions}
(Cambridge Univ. Press, Cambridge 1966).

\Bibitem[SVA]{SVA}
\textsc{C.~Smet, W.~Van Assche},
Irrationality proof of a $q$-extension of $\zeta(2)$ using little $q$-Jacobi polynomials,
\emph{Acta Arith.} \textbf{138}:2 (2009), 165--178.

\Bibitem[So1]{So1}
\textsc{В.\,Н.~Сорокин},
Аппроксимации Эрмита--Паде для систем Никишина и иррациональность $\zeta(3)$,
\emph{Успехи матем. наук} \textbf{49}:2 (1994), 167--168.

\Bibitem[So2]{So2}
\textsc{В.\,Н.~Сорокин},
О мере трансцендентности числа~$\pi^2$,
\emph{Матем. сб.}
\textbf{187}:12 (1996), 87--120.

\Bibitem[So3]{So3}
\textsc{В.\,Н.~Сорокин},
Теорема Апери,
\emph{Вестник МГУ. Сер.~1. Матем., мех.}
\No~3 (1998), 48--52.

\Bibitem[So4]{So4}
\textsc{В.\,Н.~Сорокин},
Циклические графы и теорема Апери,
\emph{Успехи матем. наук} \textbf{57}:3 (2002), 99--134.

\Bibitem[Sta]{Sta}
\textsc{F.~Stan},
On recurrences for Ising integrals,
\emph{Adv. in Appl. Math.} \textbf{45} (2010), 334--345.

\Bibitem[St]{St}
\textsc{V.~Strehl},
Binomial identities---combinatorial and algorithmic aspects,
\emph{Discrete Math.} \textbf{136}:1--3 (1994), 309--346.

\Bibitem[As1]{As1}
\textsc{W.~Van Assche},
Approximation theory and analytic number theory,
\emph{Special functions and differential equations}
(Madras 1997), eds. K.~Srinivasa Rao et~al.
(Allied Publ., New Delhi 1998), 336--355.

\Bibitem[As2]{As2}
\textsc{W.~Van Assche},
Little $q$-Legendre polynomials and irrationality of certain Lambert series,
\emph{Ramanujan J.}
\textbf{5} (2001), 295--310.

\Bibitem[Vas]{Vas}
\textsc{О.\,Н.~Василенко},
Некоторые формулы для значения дзета-функции Римана в целых точках,
\emph{Теория чисел и ее приложения}
(Ташкент, 26--28~сентября 1990\,г.),
Тезисы докладов Республиканской научно-теоретической конференции
(Ташкент, Ташкентский гос. пед. институт 1990), 27.

\Bibitem[Va1]{Va1}
\textsc{Д.\,В.~Васильев},
Некоторые формулы для дзета-функции в целых точках,
\emph{Вестник МГУ. Сер.~1. Матем., мех.}
\No~1 (1996), 81--84.

\Bibitem[Va2]{Va2}
\textsc{D.\,V.~Vasilyev},
On small linear forms for the values of the Riemann
zeta-function at odd points,
Preprint \No~1\,(558)
(Nat. Acad. Sci. Belarus, Institute Math., Minsk 2001).

\Bibitem[Va]{Va}
\textsc{R.\,C.~Vaughan},
\emph{The Hardy--Littlewood method},
Cambridge Tracts in Mathematics
\textbf{125}
(Cambridge Univ. Press, Cambridge 1997).

\Bibitem[Vi]{Vi}
\textsc{C.~Viola},
Birational transformations and values of the Riemann zeta-function,
\emph{J. Th\'eor. Nombres Bordeaux} \textbf{15}:2 (2003), 561--592.

\Bibitem[VK]{VK}
\textsc{С.\,М.~Воронин, А.\,А.~Карацуба},
\emph{Дзета-функция Римана}
(Физматлит, Москва 1994).

\Bibitem[Wa]{Wa}
\textsc{Б.\,Л.~Ван дер Варден},
\emph{Алгебра}
(Наука, Москва 1976).

\Bibitem[WaZ]{WaZ}
\textsc{J.~Wan, W.~Zudilin},
Generating functions of Legendre polynomials: a tribute to Fred Brafman,
\emph{J. Approximation Theory} \textbf{164}:4 (2012), 488--503.

\Bibitem[We]{We}
\textsc{A.~Weil},
Remarks on Hecke's lemma and its use,
\emph{Algebraic number theory}, Kyoto Internat. Sympos., Res. Inst. Math. Sci., Univ. Kyoto, Kyoto 1976
(Japan Soc. Promotion Sci., Tokyo 1977), 267--274.

\Bibitem[Wh1]{Wh1}
\textsc{F.\,J.\,W.~Whipple},
A group of generalized hypergeometric series:
relations between 120 allied series of the type $F[a,b,c;d,e]$,
\emph{Proc. London Math. Soc. II~Ser.} \textbf{23} (1925), 104--114.

\Bibitem[Wh2]{Wh2}
\textsc{F.\,J.\,W.~Whipple},
On well-poised series, generalized hypergeometric series
having parameters in pairs, each pair with the same sum,
\emph{Proc. London Math. Soc. II~Ser.} \textbf{24} (1926), 247--263.

\Bibitem[Wi]{Wi}
\textsc{A.~Wiles},
Modular elliptic curves and Fermat's last theorem,
\emph{Annals of Math.} (2) \textbf{141}:3 (1995), 443--551.

\Bibitem[WZ]{WZ}
\textsc{H.\,S.~Wilf, D.~Zeilberger},
An algorithmic proof theory for hypergeometric
(ordinary and ``$q$'') multisum/integral identities,
\emph{Invent. Math.} \textbf{108}:3 (1992), 575--633.

\Bibitem[Ya]{Ya}
\textsc{Y.~Yang},
Ap\'ery limits and special values of $L$-functions,
\emph{J. Math. Anal. Appl.} \textbf{343} (2008), 492--513.

\Bibitem[Ze]{Ze}
\textsc{D.~Zeilberger},
Computerized deconstruction,
\emph{Adv. Appl. Math.} \textbf{31} (2003), 532--543.

\Bibitem[Zl1]{Zl1}
\textsc{С.\,А.~Злобин},
Интегралы, представляемые в виде линейных форм
от обобщенных полилогарифмов,
\emph{Матем. заметки}
\textbf{71}:5 (2002), 782--787.

\Bibitem[Zl2]{Zl2}
\textsc{С.\,А.~Злобин},
О некоторых интегральных тождествах,
\emph{Успехи матем. наук}
\textbf{57}:3 (2002), 153--154.

\Bibitem[Zl3]{Zl3}
\textsc{С.\,А.~Злобин},
Разложения кратных рядов в линейные формы,
\emph{Матем. заметки}
\textbf{77}:5 (2005), 683--706.

\Bibitem[Z1]{Z1}
\textsc{В.\,В.~Зудилин},
Разностные уравнения и мера иррациональности чисел,
\emph{Аналитическая теория чисел и приложения} (сб. статей),
Труды МИАН \textbf{218} (1997), 165--178.

\Bibitem[Z2]{Z2}
\textsc{В.\,В.~Зудилин},
Сокращение факториалов,
\emph{Матем. сб.}
\textbf{192}:8 (2001), 95--122.

\Bibitem[Z3]{Z3}
\textsc{В.\,В.~Зудилин},
Об иррациональности значений дзета-функции в нечетных точках,
\emph{Успехи матем. наук} \textbf{56}:2 (2001), 215--216.

\Bibitem[Z4]{Z4}
\textsc{В.\,В.~Зудилин},
Об иррациональности значений дзета-функции,
\emph{Современные исследования в математике и механике},
Материалы XXIII Конференции молодых ученых
мех.-мат. фак-та МГУ (9--14~апреля 2001\,г.), часть~2
(Изд-во мех.-мат\. фак-та МГУ, Москва 2001), 127--135.

\Bibitem[Z5]{Z5}
\textsc{В.\,В.~Зудилин},
Одно из восьми чисел $\zeta(5),\zeta(7),\dots,\zeta(17),\zeta(19)$
иррационально,
\emph{Матем. заметки} \textbf{70}:3 (2001), 472--476.

\Bibitem[Z6]{Z6}
\textsc{В.\,В.~Зудилин},
Одно из чисел $\zeta(5)$, $\zeta(7)$, $\zeta(9)$, $\zeta(11)$ иррационально,
\emph{Успехи матем. наук}
\textbf{56}:4 (2001), 149--150.

\Bibitem[Z7]{Z7}
\textsc{В.\,В.~Зудилин},
Об иррациональности $\zeta_q(2)$,
\emph{Успехи матем. наук}
\textbf{56}:6 (2001), 147--148.

\Bibitem[Z8]{Z8}
\textsc{В.\,В.~Зудилин},
Об иррациональности значений дзета-функции Римана,
\emph{Изв. РАН. Серия матем.}
\textbf{66}:3 (2002), 49--102.

\Bibitem[Z9]{Z9}
\textsc{W.~Zudilin},
Remarks on irrationality of $q$-harmonic series,
\emph{Manuscripta Math.} \textbf{107}:4 (2002), 463--477.

\Bibitem[Z10]{Z10}
\textsc{В.\,В.~Зудилин},
О мере иррациональности $q$-аналога $\zeta(2)$,
\emph{Матем. сб.} \textbf{193}:8 (2002), 49--70.

\Bibitem[Z11]{Z11}
\textsc{В.\,В.~Зудилин},
Совершенно уравновешенные гипергеометрические ряды и кратные интегралы,
\emph{Успехи матем. наук}
\textbf{57}:4 (2002), 177--178.

\Bibitem[Z12]{Z12}
\textsc{В.\,В.~Зудилин},
О диофантовых задачах для $q$-дзета-значений,
\emph{Матем. заметки} \textbf{72}:6 (2002), 936--940.

\bibitem{Z13a}
\textsc{В.\,В.~Зудилин},
Алгебраические соотношения для кратных дзета-значений,
\emph{Успехи матем. наук}
\textbf{58}:1 (2003), 3--32.

\Bibitem[Z14]{Z14}
\textsc{В.\,В.~Зудилин},
О функциональной трансцендентности $q$-дзета-значений,
\emph{Матем. заметки} \textbf{73}:4 (2003), 629--630.

\Bibitem[Z15]{Z15}
\textsc{W.~Zudilin},
Well-poised hypergeometric service for diophantine problems of zeta values,
\emph{J. Th\'eorie Nombres Bordeaux}
\textbf{15}:2 (2003), 593--626.

\Bibitem[Z16]{Z16}
\textsc{W.~Zudilin},
Heine's basic transform and a permutation group for $q$-harmonic series,
\emph{Acta Arith.} \textbf{111}:2 (2004), 153--164.

\Bibitem[Z17]{Z17}
\textsc{W.~Zudilin},
Arithmetic of linear forms involving odd zeta values,
\emph{J. Th\'eor. Nombres Bordeaux}
\textbf{16}:1 (2004), 251--291.

\Bibitem[Z18]{Z18}
\textsc{W.~Zudilin},
Well-poised hypergeometric transformations of Euler-type multiple integrals,
\emph{J. London Math. Soc.}
\textbf{70}:1 (2004), 215--230.

\Bibitem[Z19]{Z19}
\textsc{W.~Zudilin},
On a combinatorial problem of Asmus Schmidt,
\emph{Electron. J. Combin.} \textbf{11}:1 (2004), \#R22, 8~pages.

\Bibitem[Z19a]{Z19a}
\textsc{В.\,В.~Зудилин},
Эссе о мерах иррациональности $\pi$ и других логарифмах,
\emph{Чебышёвский сб.} (ТГПУ, Тула) \textbf{5}:2 (2004), 49--65.

\Bibitem[Z20]{Z20}
\textsc{В.\,В.~Зудилин},
Об обратном преобразовании Лежандра одного семейства последовательностей,
\emph{Матем. заметки} \textbf{76}:2 (2004), 300--303.

\Bibitem[Z20a]{Z20a}
\textsc{W.~Zudilin},
Approximations to $q$-logarithms and $q$-dilogarithms, with applications to $q$-zeta values,
\emph{Труды по теории чисел}, Записки научн. семинаров ПОМИ, СПб. \textbf{322} (2005), 107--124.

\Bibitem[Zu]{Zu}
\textsc{В.\,В.~Зудилин},
Формулы рамануджанова типа и меры иррациональности некоторых кратных числа~$\pi$,
\emph{Матем. сб.} \textbf{196}:7 (2005), 51--66.

\Bibitem[Z22]{Z22}
\textsc{W.~Zudilin},
A new lower bound for $\|(3/2)^k\|$,
\emph{J. Th\'eor. Nombres Bordeaux}
\textbf{19}:1 (2007), 313--325.

\Bibitem[Z23]{Z23}
\textsc{W.~Zudilin},
Approximations to -, di- and tri- logarithms,
\emph{J. Comput. Appl. Math.} \textbf{202}:2 (2007), 450--459.

\Bibitem[Z24]{Z24}
\textsc{W.~Zudilin},
Ramanujan-type formulae for $1/\pi$: A second wind?,
\emph{Modular Forms and String Duality}, N.~Yui, H.~Verrill, and C.\,F.~Doran (eds.),
Fields Inst. Commun. Ser. \textbf{54} (Amer. Math. Soc., Providence, RI 2008),  179--188.

\Bibitem[Z13]{Z13}
\textsc{W.~Zudilin},
Ap\'ery's theorem. Thirty years after,
\emph{Intern. J. Math. Computer Sci.} \textbf{4}:1 (2009), 9--19.

\Bibitem[Z80]{Z80}
\textsc{W.~Zudilin},
Ramanujan-type supercongruences,
\emph{J. Number Theory} \textbf{129}:8 (2009), 1848--1857.

\Bibitem[Z90]{Z90}
\textsc{В.\,В.~Зудилин},
Арифметические гипергеометрические ряды,
\emph{Успехи матем. наук} \textbf{66}:2 (2011), 163--216.

\Bibitem[Z96]{Z96}
\textsc{W.~Zudilin},
Period(d)ness of $L$-values,
\emph{Number Theory and Related Fields, In memory of Alf van der Poorten},
J.\,M.~Borwein et~al. (eds.),
Springer Proceedings in Math. \& Stat. \textbf{43} (Spinger, New York 2013), 381--395.

\Bibitem[Z100]{Z100}
\textsc{В.\,В.~Зудилин},
О мере иррациональности числа $\pi^2$,
\emph{Успехи матем. наук} \textbf{68}:6 (2013), 171--172.

\end{thebibliography}
\end{document}